\documentclass[12pt]{article}
\newcommand{\N}{\mbox{I$\!$N}}
\newcommand{\R}{\mbox{I$\!$R}}

\newcommand{\qed}{{\hfill {$\rlap{$\sqcap$}\sqcup$}}\\[0.2in]\hspace*{0.5in}}
\newcommand{\qedwh}{{\hfill {$\rlap{$\sqcap$}\sqcup$}}\\[0.2in]}
\newcommand{\bk}{\\[0.1in] \hspace*{0.5in} }
\newcommand{\btd}{\bigtriangledown}
\newcommand{\mfor}{\ \ \ \ {\mbox{for}} \ \ }

\begin{document}

\begin{center} {\LARGE   {\bf Construction of Blow-up Sequences for  the}}
\medskip \medskip \smallskip  \\ {\LARGE {\bf  Prescribed Scalar Curvature Equation}} \medskip \medskip \smallskip \\ {\LARGE {\bf  on $S^n$.\,   II. \,Annular Domains}}

\vspace{0.43in}

{\Large {Man Chun  {\LARGE L}}}{EUNG}\\[0.165in]

{\large {National University of Singapore}}${\,}^1$\\[0.1in]
{\tt matlmc@nus.edu.sg}\\
\end{center}
\vspace{0.3in}
\begin{abstract}
\vspace{0.33in}
\noindent
Using the Lyapunov\,-\,Schmidt reduction method,   we  describe how to use annular domains to construct (scalar curvature) functions on $\,S^n\,$ (\,$n \ge 6$\,)\,,\, so that each one of them enables the conformal scalar curvature equation to have  a blowing\,-\,up sequence of positive solutions.   The prescribed  scalar curvature function  is shown  to have  $\,C^{\,n - 1\,, \, \, \beta}$\, smoothness.

\end{abstract}

\vspace*{2.3in}

{\bf 2000 AMS MS Classification}\,: Primary 35J60\,;
\ \ Secondary 53C21. \\[0.01in]

$\overline{\ \ \ \ \ \ \ \ \ \ \ \ \ \ \ \ \ \ \ \ \  \ \ \ \ \ \ \ \ \ \ \ \ \ \ \ \ \ \ \ \ \ \ \ }$\\[0.05in]
\noindent{$\!\!{\!}^* {\small{\mbox{\, Department of Mathematics, National University of Singapore, 10, Lower Kent Ridge Rd.,}}}$}\\[0.03in]
 {\small{\mbox{\ \  \,Singapore 119076, Republic of Singapore\,.}}}

\newpage

\newpage

{\bf \large {\bf  1. \ \ Introduction.}}\\[0.1in]
In this article we apply the   Lyapunov\,-\,Schmidt reduction method to find non\,-\,constant functions  $\,{\cal K}\,$ so that the equation
$$
\Delta_1 \, u - {\tilde c}_n \,n \,(n - 1) \,u \ + \ ({\tilde c}_n \,{\cal K}) u^{{n + 2}\over {n - 2}} \ \,= \ 0 \ \ \ \ {\mbox{in}} \ \ S^n \leqno (1.1)
$$
has infinite number of positive solutions, which compose  a blow\,-\,up sequence of  solutions.  We refer to \S\,1\,c for the (rather standard) notations we use. \bk
As often with good insight, the Lyapunov\,-\,Schmidt reduction method works with  miraculous simplicity. It is used successfully by Ambrosetti, Berti and Malchiodi to construct many $C^k$\,-\,metrics $g$  on $S^n$ such that the Yamabe equation
$$
\Delta_g \, u \ - \ {\tilde c}_n R_g \,u \ + \ {\tilde c}_n \,n\,(n - 1) \,u^{{n + 2}\over {n - 2}} \ \,= \ 0 \ \ \ \ \ \ {\mbox{in}} \ \ \ S^n \leqno (1.2)
$$
has infinite number of positive solutions [\,here $\,g\,$ is non\,-\,conformal  to $\,g_1\,$ -- the standard metric on $\,S^n$\,;\, $\,R_g\,$ is the scalar curvature of $\,(S^n, \, g)$\,]\,.\,   See the renowned monograph \cite{Progress-Book}\,,\, in which we can find  their work on equation (1.1) as well. Initially, the degree of smoothness  on $g$ is  restricted to $\,2 \,\le \,k \,\le\, (n - 3)/4\,,\,$ but is improved to $\,C^\infty\,$ by Brendle in \cite{Brendle-1} for $\,n \,\ge\, 52\,$, and together with Marques \cite{Brendle-2} for the remaining cases $\,25 \,\le\, n \le \,51\,$. The result should be read with \cite{Compact}\,,\, in which compactness theorem for the Yamabe equation is obtained for $\,n\, \le \,24\,$ (see also \cite{Li-Zhang-Compactness},\, \cite{Li-Zhang-Compactness-3} and \cite{Marques} for earlier results). The unexpected critical dimension $\,n \,= \,24\,$ is both fascinating and intriguing.\bk
It is rewarding to view the prescribed scalar curvature equation (1.1) as a kind of ``dual" to the Yamabe equation (1.2). In contrast, blow\,-\,up sequences of positive solutions for (1.1) are widely studied as a mean to find solutions. However, to the knowledge of the author, no general method is known to construct  blow\,-\,up sequences of positive solutions for equation (1.1) [\,that is, in $S^n$\,,\, with fixed critical index $(n + 2)/(n - 2)$ and fixed ${\cal K}$\,]\,.\, Apparently the only previously known case is when ${\cal K}$ is equal to a positive constant. With the method presented in this article, we are able to  generate simple blow\,-\,ups (see  \cite{Compact} and \cite{Li-1}). The solutions are perturbations of the standard ones [\,cf. (1.4) below], and are adjusted to maintain the same scalar curvature function ${\cal K}\,.\,$ The degree of flatness of ${\cal K}$ at the blow\,-\,up point (the south pole) is up to $(n - 1)$\, [\,refer to \S\,3\,g \,and\, (3.30)\,]\,.\,\bk
Recalling the notations used in Part I \cite{I}\,,\, the reduced functional for equation (1.1) is defined by
$$\,
G_{|_{\bf Z}} ({\bf z}) \ = \ {\bar c}_{-1} \int_{\R^n} H (y) \left[ {\lambda\over {\lambda^2 + |\,y - \xi|^{\,2}}} \right]^n dy\, \ \ \ \ \mfor \ {\bf z} = V_{\lambda,\,\, \xi} \in {\bf Z}\,, \ \ \ \ \ \ \ \ {\mbox{where}} \leqno (1.3) $$
$$
{\bf Z} \ := \ \left\{ \   V_{\lambda,\  \xi} \,(y)\ := \ \left( {\lambda\over {\lambda^2 + |\,y - \xi|^{\,2}}} \right)^{\!\!{{n-2}\over 2}} \ \ \ \ \    \mbox{with} \ \  \ \ (\lambda, \ \xi)\,\in \,\R^+\times \R^n\, \right\}\,,  \leqno (1.4)
$$
$$K (y)\ :=\ {\cal K} \,({\dot{\cal P}}^{-1} (y)) \ = \ 4\,n\,(n - 1) \ + \ \varepsilon\, H (y)\ \ \mfor \ y \in \R^n\,.\,
\leqno (1.5)
$$

In (1.3)\,,\,
 ${\bar c}_{-1} \ =\ -\, [\,{\tilde c}_n \cdot (n - 2)]\,/(2n)\,.$\, For $\,\varepsilon \in \R\,$   small, the question on finding a positive solution of equation (1.1) is reduced   to finding a {\it stable\,} critical point of $G_{|_{\bf Z}}\,.\,$ See \cite{Progress-Book} or \cite{I}\,.\, It is shown in Part I \cite{I} (Proposition 7.10),  that the Hessian matrix of \,$G_{|_{\bf Z}}$\,  at a critical point is always trace\,-\,free, thus can never be positive or negative definite. Instead of looking for strict local    maximum or minimum (as in the Yamabe equation), we seek {\it saddle points\,.\,} Although it is a finite dimensional problem, practical examples show that it can be illusive to locate (local) critical points\,.\, Indeed, the existence of a critical point implies the fulfillment of the Kazdan\,-\,Warner condition for \,${\cal K}$\, (see Theorem 6.10 in \cite{I} for full details)\,. \bk
In \cite{Ambrosett-Scalar}, using the Lyapunov\,-\,Schimdt method and a degree counting method,  (finite number of) stable critical points of the $\,G_{|_{\bf Z}} $ are found.   The research is carried on in \cite{Ambrosetti-Malchiodi}\,,  where the authors explore  symmetries.  Each of these methods does not readily juxtapose   to produce a blow\,-\,up sequence of solutions for the scalar curvature equation (1.1). \bk
In this article, fixing an annular domain (see   \S\,1\,a below), we determine precisely the critical point of $\,G_{|_{\bf Z}}$\,,\, and show that the Hessian matrix at the critical point is non\,-\,degenerate (hence a stable critical point).
 By superimposing  concentric annular domains, and carefully estimating the gradient interference, we are able to find infinite number of stable critical points via degree theory for maps.  With the uniform cancelation property in the  Lyapunov\,-\,Schmidt  method\, (described in Part I \cite{I})\,,\, a blow\,-\,up sequence can be found. We obtain the claimed $\,C^{\,n - 1\,, \   \beta}\,$ regularity  in \S\,3\,g and \S\,3\,h  by choosing the ``strength" [\,see (1.6)] and the radii of the annular domains suitable. Regarding the degree of smoothness, cf.  also  Remark 7.11 in \cite{Leung-blow-up}\,.\, We observe that similar argument works for slightly offset, non\,-\,concentric annular domains.\bk
 The search for a construction on higher energy blow\,-\,ups (e.g. towering blow\,-\,ups, as well as the more complicated aggregated and clustered blow\,-\,ups)  is a challenging project for  research. See \cite{Leung-Supported} concerning a classification of blow\,-\,ups for (1.1)\, (cf. also \cite{Li-1} and \cite{Schoen-Notes}\,)\,.\, Refer to \cite{Leung-1st} \cite{Leung-2nd}  for non\,-\,compact  spaces.

 \vspace*{0.2in}


{\bf \S\,1\,a.   } {\it  Main Result.}\ \ \
Consider the annular domains
$$
B_o \,\left( {{1+ \eta}\over {\bf a}}\right)\,\bigg\backslash \  \overline{B_o \,\left( {{1- \eta}\over {\bf a}}\right)\!} \ , \ \  \cdot \ \cdot \ \cdot\ , \ \ \ B_o \,\left( {{1+ \eta}\over {{\bf a}^k}}\right)\,\bigg\backslash \   \overline{B_o \,\left( {{1- \eta}\over {{\bf a}^k}}\right)\!}\ , \  \ \cdot\  \cdot\  \cdot\ .
$$
We give explicit conditions on  the   numbers $\,{\bf a} \,>\, 1\,$ and $\,\eta \,\in\, (0, \ 1)\,$  in (1.9)\,.\,
Choose a small positive number $\,\sigma\,$  and fix a number
  $\tau\, \in\, (n- 1\,, \ n)\,.\,$ Let $\,H^S\,$ be given by
$$
(1.6) \ \ \ \ \ \  H^S\, (y) \ = \ \cases{ {\displaystyle{1\over {{\bf a}^{\tau k} }}}   \ \ & {\mbox if} \ \ \   $\displaystyle{ y \in B_o \,\left( {{1\,+ \, \eta  }\over {{\bf a}^k}}\right)\,\bigg\backslash \  B_o \,\left( {{1\,-  \, \eta  }\over {{\bf a}^k}}\right)\,, \ \ \ \ \ k = 1, \ 2, \cdot \cdot \cdot\,,}\ \ \ \ \ \ \ \ \ $\cr
\ & \ \cr
0 & {\mbox if } \ \ $\displaystyle{ y \,\not\in \  {{\bigcup}_{k = 1}^\infty}   \left\{ B_o \,\left( {{1\,+ \,(\eta  \,+\, \sigma)}\over {{\bf a}^k}}\right)\,\bigg\backslash \, B_o \,\left( {{1\,- \,(\eta  \,+\, \sigma)}\over {{\bf a}^k}}\right)\right\}}$\,.}
$$
In between
$\,
B_o \! \left( {{1\,+\, (\,\eta  \,+\, \sigma)}\over {{\bf a}^k}}\right) \!\bigg\backslash\! B_o \!\left( {{1\,+ \, \eta  }\over {{\bf a}^k}}\right)\ $\  and $\ B_o  \!\left( {{1\,- \, \eta  }\over {{\bf a}^k}}\right)\! \bigg\backslash \!B_o  \!\left( {{1\,- \, (\,\eta \, + \,\sigma)}\over {{\bf a}^k}}\right)\,,
$\,
we properly smooth\\[0.075in]
out
$\,H^S\,$ so that $\,H^S\,\in\, C^{\,n-1, \ \beta} \,(\R^n)\,$ (need not be rotationally symmetric\,; see \S\,3\,g    and \, \S\,3\,h)\,.\,
\vspace*{0.15in}

{\bf Main Theorem 1.7.} \ \   {\it For $\,n \ge 6\,$\,,\, let\,}  $H^S \in C^{\,n - 1\,, \ \beta} \,(\R^n)\,$ {\it be as described in }\,(1.6)\,,\, {\it  with the parameter $\,\tau\,$  and $\,\eta\,$ satisfying\,}

\vspace*{-0.25in}

$$
\tau \,\in\, (n - 1\,, \ n)\,, \ \ \ \  1 - A^2 \ >\ \eta \ >\  B^2 \ >\  0  \ \ \ \ {\it{and}} \ \ \ \
  {{1 + \eta}\over {1 - \eta}} \ \le \   {5\over 2}\,. \leqno (1.8)
  $$

  \vspace*{-0.05in}

{\it Here $\,A\,$ and $\,B\,$ are positive numbers\,.\, There exist positive constants \,$C$\,,\, $C_1$\,,\, $c$ \,and\, $\varepsilon_o$ so that if the parameters \,${\bf a}$\, and $\,\sigma$ satisfy}
$$
{\bf a}\  > \ C^2 \ \ \ \ {\it{and}} \ \ \ \ 0 \ < \ \sigma\  <\  c^2\,,
$$

\vspace*{-0.25in}

{\it then the equation }
$$
\Delta_o \,v + {\tilde c}_n  \left[\, 4n\,(n - 1) + \varepsilon \, H^S\,\right]  v^{{n + 2}\over {n - 2}} \ \,= \  0 \ \ \ \ \ \ \   {\it{in}} \ \ \  \,\R^n\leqno (1.9)$$
{\it has an infinite number of positive solutions\,} $\,\{ v_m \}_{m = 1}^\infty\, \subset\  C^{2\,,\ \bar\beta} \,(\R^n)$ {\it \,whenever\,} $|\,\varepsilon|\le \varepsilon_o\,.\,$ {\it Moreover,}

\vspace*{-0.25in}

 $$
 \Vert \,v_m - V_{\lambda_m\,, \ \xi_m} \Vert_\btd\ \le \ C_1\cdot\varepsilon \ \ \ \ \ \ {\it{for}} \ \ \ m = 1, \ 2, \cdot \cdot \cdot\,. \leqno (1.10)
 $$
{\it Here $\,\lambda_m \,\to\, 0\,$ and $\,|\,\xi_m| \,\to\, 0\,$ as $\,m \,\to \,\infty\,.\,$ As a result, \,$0\,$ is a blow\,-\,up point for $\,\{ v_m \}_{m = 1}^\infty$\ .}

\vspace*{0.2in}

\hspace*{0.5in}In (1.10), $\,\Vert \ \,\Vert_\btd\,$ represents the  $\,L^2$-\,norm on gradient for the Hilbert space ${\cal D}^{1\,, \ 2}$.\, See (2.7) in Part I \cite{I}\,.\, Via this $\,L^2$-\,norm on gradient and the Sobolev inequality, we discuss in \S\,3\,j\, how to transfer these solutions $\,\{ v_m \}_{m = 1}^\infty$\, back to $S^n$ as a blow\,-\,up sequence of solutions for equation (1.1).  [\,As for the number $5/2$ in (1.8)\,,\, it     appears naturally in some calculations. It is, however, not absolute or sharp.]

 \vspace*{0.2in}


{\bf \S\,1\,b. \  } {\it  Main features concerning the proof.}\ \  \
The key is to show that the   critical point with respect to a single annular domain is `stable' enough to withstand the interaction from other annular domains. Thus we are required to take \,${\bf a}$ to be large (big spacing between the annular domains \ $\Longrightarrow$ \ less interference\,)\,,\, and $\,\eta\,$ to be not too small (producing stronger effect on the gradient change for an individual annular domain\,)\,.\bk
In \S\,2, we  estimate the first derivatives  of \,$G_{|_{\bf Z}}$\, via the stereographic  projection back to $S^n$ (see \,\S\,2\,a and \,\S\,2\,b)\,.\, This geometric approach helps to visualize the location of the critical point (cf. Lemma 2.21), and the sharp changes in gradient after leaving the critical point (cf. Lemma 2.65).\bk
We remark that the method is stable under \,``\,small perturbation'' so that functions nearby $\,H^S$  can also be used in the construction. See \S\,3\,e for more detail.


\newpage

{\bf \S\,1\,c.  } {\it Conventions.}\ \
 Throughout this work,\, we  assume  the dimension $\,n \ge 3\,,\,$ except when otherwise is specifically mentioned, and let  \, ${\tilde c}_n = (n - 2)/[\,4\,(n - 1)]\,.\,$
We observe the practice of  using $\,C,\,$ possibly with sub-indices, to denote various positive constants,
which may be rendered {\it differently\,} from line to line according to the contents. {\it Whilst we use $\,{\bar c}\,$ and  $\,{\bar C},\,$ possibly with sub-index, to denote a\,}   fixed\,
{\it positive constant which always keeps the same value as it is first defined\,}.\, [\,The {\it negative\,} constant $\,{\bar c}_{-1}\,$ is defined in the sentence proceeding (1.5).\,]\\[0.075in]
$\bullet_1$ \ \, Denote by $\,B_y (r)\,$ the open ball in $\,\R^n$ with center at $\,y\,$ and radius $\,r > 0\,,\,$  and $\,\Vert \,  S^n \, \Vert\,$ the measure of  $\,S^n\,$ in $\,\R^{n + 1}$ with respect to the standard metric.\\[0.075in]
$\bullet_2$ \ \, $\Delta_g$ is the Laplace-Beltrami operator associated with the metric $g$ on $S^n$\,.\, Likewise, $\,\Delta_o$ is the Laplace-Beltrami operator associated with the Euclidean metric $g_o$ on $\R^n,\,$ and $\,\Delta_1$ is the
Laplace-Beltrami operator associated with the standard metric $g_1$ on $S^n$.\\[0.075in]
$\bullet_3$ \ \, Whereever there is no risk of misunderstanding, we suppress \,$dy$\, from the integral  expression.

\vspace*{0.2in}


{\large{\bf \S\, 2. \ \  Annular domains.}}\\[0.1in]
Let us start with a general
{\it smooth\,} and bounded domain $\,\Omega \subset \R^n$,\, and let  $H_\Omega$ be the characteristic function of  $\,\Omega$ (\,i.e., $H$   equals to $1$ within $\Omega$\,,\, and $0$ without). Expression (1.3) becomes

\vspace*{-0.27in}

$$\ \ \ \ \ \ \ \ \ \ \ \ \ \
G_{|_{\bf Z}} \, (\lambda\,,\, \xi\,; \ \Omega)  \ := \  \,{\bar c}_{-1} \int_{\Omega}   \left[ {\lambda\over {\lambda^2 + |\,y - \xi|^{\,2}}} \right]^n \ \ \ \ \ \ \ \ \ \ \ \ \ \ \ \ \ \ \ \ \ \ \ \     (H \,=\, H_\Omega)\,. \leqno (2.1)
$$
We remark that $G_{|_{\bf Z}}\, (\bullet\,,\, \bullet\,; \ \Omega)$ remains a smooth function on $\,\R^+ \times \R^n$.

\vspace*{0.15in}

%
%
{\bf \S\,2\,a.} \ \ {\it Rescaling and shifting.}\ \
Consider the change of variables $\,\bar y = \,\lambda^{-1}\, (y - \xi)\,$.\, We have

\vspace*{-0.27in}

$$
\int_{\Omega} \left(  {\lambda\over {\lambda^2 + |\,y - \xi|^{\,2}}}\right)^{\!\!n} \ = \
 \int_{  \left( {{ \Omega_{-\xi}}\over \lambda} \right) }  \left(  {1\over {1 + |\bar y|^{\,2}}}\right)^{\!\!n} \,d\bar y\,. \leqno (2.2)
$$
Here  we use a kind of
``{\it geometric notations}\,"\,:
$$
\Omega_{-\,\xi}\  := \ \{ \,\tilde y \in \R^n \ \ | \ \ \ \ \tilde y \ = \ y - \xi \mfor y \in \Omega\,\}\,, \leqno (2.3)
$$
 $$
 \left( \,{{ \Omega_{-\,\xi}}\over \lambda} \right)  \, := \, \{ \,\bar y \in \R^n \ |
\ \ \bar y = \lambda^{-1} \, \tilde y \mfor \tilde y \in
\Omega_{-\,\xi} \,\} = \left\{ \,\bar y \in \R^n \ | \ \ \bar y = {y\over \lambda} - {\xi \over \lambda} \!\!\!\mfor   y \in  \Omega  \right\}\!. $$

\vspace*{0.1in}

%
%
%
%
{\bf \S\,2\,b.} \ \ {\it Projection of annular regions.}\ \
 Let $\,\dot{\cal P}\,$ be the stereographic projection which sends the north pole $\,{\bf N}\,$ to $\,\infty\,.$\,
 The conformal factor is given by\,
 $\displaystyle{
g_1 \,(x) = \left[ \,{{4}\over {(1 + r^2)^{\,2}}} \right]\,g_o \,(y)\,\,}$ \,for   $\,y \,= \,\dot{\cal P}(x)\,.$\,   See, for instance, (2.1) and (2.3) in Part I \cite{I}\,.

We obtain

\vspace*{-0.27in}

$$
\int_{ \left({{ \Omega_{-\xi}}\over \lambda} \right)}  \left(  {1\over {1 + |\bar y|^{\,2}}}\right)^{\!\!n} \,d\bar y \ \,= \ {1\over {2^n}}
 \int_{ {\dot{\cal P} }^{-1} \left( {{ \Omega_{-\,\xi}}\over \lambda} \right)}\
  dS\,,\leqno (2.4)
  $$
where $\,dS\,$ refers to   the   standard   measure  on $\,S^n\,.$\,
Suppose
$
\,\overline{B_o (R_c)} \ \subset\, B_o (R_b)\,,\,
$
where $\,R_c\,$ and $\,R_b\,$ are positive numbers. Consider the annular domain
$$
 \Omega \ = \ B_o\,(R_b) \setminus  \overline{B_o (R_c)}\ \ \ \
 \Longrightarrow \ \  \ {{ \Omega_{-\,\xi}}\over \lambda} \ =
\  B_{\, {-\,\xi\over \lambda} } \left( \!\lambda^{-1} R_b \right) \, \big\backslash \,\,
 \overline{B_{ \,{-\,\xi\over \lambda} } \left( \lambda^{-1} R_c \right)} \ .
$$
Here we use (2.3)\,.\, We are interested in finding

\vspace*{-0.25in}

$$
\int_{ {{\dot{\cal P}}}^{-1}    \left(  B_{\! {{-\,\xi}\over \lambda} }\,
(\, \lambda^{-1} \,R_b ) \ \backslash \
 \overline{B_{ \!{{-\,\xi}\over \lambda} }\,(\lambda^{-1} \,R_c \,) }\,
\right) }   \ dS\,. \leqno (2.5)
$$
It is remarked that the above expression depends smoothly on $\,\lambda\,$ and $\,\xi\,.$\, \bk
%
%
%
%
The stereographic projection
$\dot{\cal P}$\,  sends an $(n - 1)\,-\,$sphere on $S^n \setminus \,\{ {\bf N} \}$
to an $(n - 1)\,-\,$sphere on $\R^n,\,$ and vis versa. See for example  \cite{Kulkarni}\,,\, pp. 7. It follows that the domain in the integral
(2.5) can be rendered as taking away two `caps' from $\,S^n$ .
 \vspace*{0.2in}
%

{\bf \S\,2\,c.} \   {\it  Volume of the cap in upper or lower hemisphere.}\ \
Consider a cap $\,{\cal C} \subset S^n$\,,\, bounded by the $(n - 1)$\,-\,sphere
  with radius  $\varrho$\,,\, {\bf measured in
the metric of\,} $\R^{n + 1}$.   {\it Here we only consider such a cap   that can be contained in a hemisphere\,.\,} It follows that $\varrho < 1\,.\,$  By symmetry, we view the cap as  upon the
southern hemisphere (resp. northern hemisphere)\,,\, ``\,centered\," at the respective pole\,:
\begin{eqnarray*}
{\mbox{Volume \,\, of \ the \ cap \ upon \ lower \ hemisphere  }}
& = & \Vert S^{n- 1} \Vert  \cdot\! \int_0^{\arcsin \varrho}  \ \,
[\,\sin\,\dot\phi\,]^{\,n - 1} \, d\,\dot\phi\,,\\
{\mbox{Volume\,\,\,of \, the \ cap \ upon \ upper \ hemisphere  }}
& = & \Vert S^{n- 1} \Vert\cdot\!
\int_0^{\arcsin \varrho} \ \,[\,\sin\,\phi\,]^{\,n - 1} \  d\,\phi\,.
\end{eqnarray*}
Here $\,\dot\phi\,$  is measured from the negative direction of the $\,x_{n + 1}\,-$axis,
and $\phi$ from the positive direction of $\,x_{n + 1}\,-$axis\,.\, We write the last expression as
$$
{\mbox{Vol\ Cap}} \,(\varrho) \ = \ \Vert \,S^{n-1} \Vert\cdot\!
\int_0^{\arcsin \, \varrho} \ [\,\sin \phi]^{\,n - 1} \,d\,\phi\,, \leqno (2.6)
$$
where we understand that $\varrho < 1$, and possibly after a rotation,  the cap is on the upper hemisphere. We can rewrite (2.5) as
$$
\int_{ {\cal P}^{-1}   \left(  B_{\! {{-\,\xi}\over \lambda} }
(\, \lambda^{-1} R_b ) \ \setminus \
 \overline{B_{ \!{{-\,\xi}\over \lambda} } (\lambda^{-1} R_c \,) }\
\right) }   \ dS \ = \ \Vert \,S^n \Vert  \ - \ {\mbox{Vol\ Cap}} \,(\varrho_b) \ - \ {\mbox{Vol\ Cap}} \,(\varrho_c)\,,\leqno (2.7)
$$
where $\varrho_b\,$ (respectively $\,\varrho_c$) is the radius of the sphere $\displaystyle{ {\cal P}^{-1} \,  \left( \partial  B_{\! {{-\,\xi}\over \lambda} }\,
(\, \lambda^{-1} R_b )
\right) }$ \ $\bigg[$\,{\mbox{respectively}}\\[0.055in]   $\displaystyle{  {\cal P}^{-1} \,  \left(\partial
  B_{ \!{{-\,\xi}\over \lambda} } \, (\lambda^{-1} R_c \,) \right) \, }\bigg]$\  \, (measured in the Euclidean metric of $\,\R^{n + 1}$\,)\,.

  \newpage

\hspace*{0.5in}For later application, let us consider  $\,\varrho\,$ depending smoothly  on a parameter $\delta$\,,\,  i.e., $\,\varrho \,= \, \varrho\, (\delta )\,< \,1$. We then have
\begin{eqnarray*}
(2.8) \ \ \ \  & \ &  {{d \,[\, {\mbox{Vol\ Cap}} \,(\varrho)]}\over {d\, \delta}} \ \,= \ \,
\Vert \,S^{n-1} \Vert\cdot [\,\sin \,\phi\,]^{\,n - 1}\bigg\vert_{\,\phi \,=\,\, \arcsin \, \varrho}\!
\times {1\over {\sqrt{1 - \varrho^2\,}}}
\,\cdot {{d\, \varrho (\delta) }\over {d\,\delta}}\ \ \ \ \ \ \ \ \ \ \ \ \ \ \ \ \ \ \ \ \ \ \ \ \\
& = &
\Vert \,S^{n-1} \Vert\, \cdot {{\varrho^{\,n - 1}}\over {\sqrt{1 - \varrho^2\,}}}
\,\cdot {{d\, \varrho(\delta) }\over {d\,\delta}} \ \,= \ \,
\Vert \,S^{n-1} \Vert\, \cdot
{{\varrho^{\,n - 1}}\over {\sqrt{1 - \varrho^2\,}}}
\,\cdot {1\over {2\,\varrho}} \,{{d\,[\,\varrho (\delta)]^2}\over {d\,\delta}}\\
& = &
{{\Vert \,S^{n-1} \Vert}\over 2}\, \cdot
{{\varrho^{\,n - 2}}\over {\sqrt{1 - \varrho^2\,}}}
  \cdot {{d\,[\,\varrho\,(\delta)]^2}\over {d\,\delta}} \ \ \ \ \ \ \ \ \ \  \ \ \ \ \  \ \ \ \ \ \ \ \ \ \ \ \ \ \ \ \ \ \ \ \ \ \ \ \ \ \ \ \ \ \ \  (\,\varrho < 1)\,.
\end{eqnarray*}
Note that the term
$$
{{\varrho^{\,n - 2}}\over {\sqrt{1 - \varrho^2\,}}} \ \ \ \
{\mbox{increases \ \ as }} \ \ \varrho \ \ {\mbox{increases  }} \ \ \ \ \ [\,{\mbox{here}} \ \ n \,\ge\, 3 \ \ {\mbox{and}} \ \ \varrho \,\in\, [\,0, \ 1)\,]\,. \leqno (2.9)
$$
 For example, when changing $\lambda$ in (2.7) while fixing $\xi$\, [\,$\varrho_b = \varrho_b \,(\lambda) < 1$ and $\varrho_c = \varrho_c \,(\lambda) < 1\,]\,,\,$ at the same time maintaining the cap
 $\, \displaystyle{ S^n \ \setminus \ {\cal P}^{-1}    \left(   B_{\! {{-\,\xi}\over \lambda} }\,
(\, \lambda^{-1} R_b )\,
\right) }$
to be on the northern hemisphere, and the cap  $\displaystyle{ \,{\cal P}^{-1}   \left(
  B_{ \!{{-\,\xi}\over \lambda} } \,(\lambda^{-1} R_c \,) \, \right) }$\,   on the southern hemisphere,  we then have

  \vspace*{-0.25in}

  \begin{eqnarray*}
(2.10) \     & \ &  {{\partial}\over {\partial \lambda}} \int_{ {\cal P}^{-1}   \left(  B_{\! {{-\,\xi}\over \lambda} }
(\, \lambda^{-1} R_b ) \ \setminus \
 \overline{B_{ \!{{-\,\xi}\over \lambda} } (\lambda^{-1} R_c \,) }\
\right) }   \ dS\\[0.075in]
&  \,&\ \ \ \ = \  {{\Vert \,S^{n-1} \Vert}\over 2}\, \cdot \left\{
  {{\varrho^{\,n - 2}_b}\over {\sqrt{1 - \varrho^2_b\,}}}
  \cdot \left[ -\, {{d\,[\,\varrho_b \,(\lambda)]^2}\over {d \lambda}} \right] \ \,+ \   \, {{\varrho^{\,n - 2}_c}\over {\sqrt{1 - \varrho^2_c\,}}}
  \cdot \left[ -\,{{d\,[\,\varrho_c \,(\lambda)]^2}\over {d \lambda}} \right] \right\}. \ \ \ \  \ \ \ \   \ \ \ \ \ \ \ \
  \end{eqnarray*}

  \vspace*{-0.1in}

 Bear in  mind that the bigger the upper cap, the smaller the value in (2.5). Likewise, the bigger the lower cap, the smaller the value in (2.5). Similar expression can be found on derivatives in (the components of) $\,\xi\,.$


\vspace*{0.1in}

{\bf \S\,2\,d.} \ \ {\it  Symmetric case and the geometric mean.}\ \
The geometric expressions  (2.5) and (2.6) provide light on the location of the
 critical point\,:
the boundary spheres
should have the same
radius \,(when seen in $\R^{n + 1}$)\, so as to cancel the increasing and
decreasing effect (when $\lambda$ is changed). In the following, we provide more details on this thought.\bk
In (2.5), take
$$
R_b \ = \ t \ + \ \Delta\,, \ \ \ \ \ \ \ \ R_c \ =\  t\  -\  \Delta\,, \ \ \ \ {\mbox{where \ \ }} t \ > \ \Delta \ > \ 0 \ \ {\mbox{are \ \ given}}\,. \leqno (2.11)
$$
 Observe that
$t$ equals  the {\it arithmetic  mean}\, of the outer and inner   radius. We can start  with (2.1).\, Putting $\xi \,= \, 0\,,\,$ we obtain
  \vspace*{-0.25in}

\begin{eqnarray*}
 (2.12) \!\!\!\!\!\!\!\!\!\!\!& \ & \ \ \ \ \int_{B_o\,(t + \Delta) \setminus \overline{B_o\,(t\,-\,\Delta)}}\    \left(  {\lambda\over {\lambda^2 + |\,y |^{\,2}}}\right)^n \ = \ \Vert \,S^{n- 1} \Vert \cdot \!\!\int_{t - \Delta}^{t + \Delta} \left(  {\lambda\over {\lambda^2 + r^2}}\right)^n r^{\,n - 1}\, dr\\
 & \ & \\
 & = & \! \Vert \,S^{n- 1} \Vert \cdot \!\!\int_{\arctan \,\left( {{t - \Delta}\over {\lambda}} \right)}^{\arctan \, \left({{t + \Delta}\over {\lambda}} \right)} \,\, {{ [\,\tan \theta\,]^{\,n - 1} \cdot [\,\sec \theta]^2}\over { [\sec \theta]^{2n} }} \  d\theta \ \ \ \ \ \ \ \ \ \ \ \ \ \ \ \ \ \  \ \ \ \ \ \ \ \   (\,r \ = \ \lambda \tan \theta\,)\ \ \ \ \ \ \ \ \ \ \ \ \\
  & = & \!\!\! \Vert \,S^{n- 1} \Vert \int_{\arctan \,
 \left( {{t - \Delta}\over {\lambda}} \right)}^{\arctan \,
 \left({{t + \Delta}\over {\lambda}} \right)} \,\,\, [\,\sin \theta\, \cos
 \theta]^{\,n - 1}  \, d\theta  = {{\Vert \,S^{n- 1} \Vert}\over {2^n}}\! \int_{\arctan \,
 \left( {{t - \Delta}\over {\lambda}} \right)}^{\arctan \,
 \left({{t + \Delta}\over {\lambda}} \right)}
  \, [\,\sin\, (2\theta)]^{\,n - 1}\,   \, d\,(2 \theta)\\
 & = & \!\!\!{{\Vert \,S^{n- 1} \Vert}\over {2^n}}
\int_{2\arctan \,\left( {{t - \Delta}\over
 {\lambda}} \right)}^{2\arctan \,
 \left({{t + \Delta}\over {\lambda}} \right)}
 \,\,\,[\,\sin \varphi]^{\,n - 1}\,   \, d\varphi \ \ \ \ \ \  \ \ \ \ \ \ \  \  \ \ \ \ \ \ \ \ \ \ \ \ \ \  \ \ \ \ \ \ \ \ \ \ \ \ \ \ \ \ \ \
 (\,\varphi = 2 \theta)\,.
\end{eqnarray*}
\begin{eqnarray*}
(2.13)\ \ \ \ & \cdot &   \!\!\!\!\!\!\cdot \cdot \cdot    \cdot \cdot \cdot \ \ \Longrightarrow \ \  {\partial\over {\partial \lambda}}
\int^{2\arctan \,\left( {{t\,+\, \Delta }\over {\lambda}}
\right)}_{2\arctan \,\left( {{t\,-\, \Delta }\over {\lambda}} \right)}
\,\,\,[\,\sin \varphi]^{\,n - 1}\,   \, d\varphi\\
 & \ &  \ \ \ \ \ \ \ \ \ \ \ = \ - \,{{2\,[\,t + \Delta]}\over {[\,{t + \Delta}]^2 + \lambda^2}}\,
(\,\sin \,{\overline{\varphi}}_+)^{\,n - 1}  \ + \,{{2\,[\,t - \Delta]}\over
 {[\,{t - \Delta}]^2 + \lambda^2}}\, (\,\sin \,{\overline{\varphi}}_-)^{\,n - 1}\,,\\
 & \ & \\
  {\mbox{where}} \ \ \ \ & \ & \,
 {\overline{\varphi}}_+ =  2\arctan  \left( {{t + \Delta}\over {\lambda}}
 \right)\,,    \ \ \ \ \ \ \ \  {\overline{\varphi}}_- = 2\arctan
 \left( {{t - \Delta}\over {\lambda}} \right)\,.
 \ \ \ \ \ \ \ \ \  \ \ \ \ \ \ \ \ \ \ \ \ \ \ \
  \ \ \ \ \   \ \ \ \ \ \ \ \ \ \   \ \ \ \ \ \ \ \ \ \ \ \ \ \ \ \ \ \ \ \ \ \ \ \ \ \ \
\end{eqnarray*}
Let us take
$$
 \lambda^2 \ =\  \lambda^2_M \ :=\  t^2 \ - \ \Delta^2 \ \ \ \ \ \  \Longrightarrow \ \ \ \lambda_M \ = \ \sqrt{(t \,+\, \Delta)(t \,-\, \Delta)\, }\,\ . \leqno (2.14)
$$
That is, $\,\lambda_M\,$ is the {\it geometric mean\,} of the inner and outer radius.
Consider   (2.13)\,:
\begin{eqnarray*}
 (2.15) \ \ \ \ \ \ \ \ \ \ \ \ \ \ \ \ \ \ \ \ \ \ \ \ \ \ \ \ \ \ \ \ \ \ \
 {{2\,[\,t \,+\, \Delta]}\over {[\,{t \,+\, \Delta}]^2
 + \lambda^2}} & = &  {{2}\over {(t \,+\, \Delta) \,+\, (t\, - \, \Delta) }}
 \ = \  {1\over t}\,,\ \ \ \ \ \ \ \ \ \ \ \ \ \ \ \ \ \ \ \ \ \ \ \ \ \ \ \ \ \  \\
(2.16) \ \ \ \ \ \ \ \ \ \   \ \ \ \ \ \ \ \ \ \ \ \ \ \ \ \ \ \ \ \ \ \ \ \
{{2\,[\,t \,- \,\Delta]}\over {[\,{t\, - \,\Delta}]^2 + \lambda^2_M}}& = &
{{2}\over {(t \,- \,\Delta)\, + \,(t \,+ \,\Delta) }} \ = \ {1\over t}\,, \ \ \ \ \ \ \ \ \ \
\end{eqnarray*}
\begin{eqnarray*}
 (2.17) \ \ \ \ \ \ \ \ \ \ \ \ \ \ \ \ \ \ \ \ \, {\bar \varphi}_{M_-} \ = \ 2\arctan   \left( {{t \,- \,\Delta}\over {\lambda_M}}
\right)  & = & 2\arctan \,\sqrt{ {{ t\, -\, \Delta}\over {t\, + \, \Delta}}\,}\,, \ \ \ \ \ \ \ \ \ \ \ \ \ \ \ \ \ \ \ \  \ \ \ \ \ \ \ \ \ \  \  \\
 (2.18) \ \ \ \ \ \ \ \ \ \ \ \ \ \ \ \ \ \ \ \ \,  {\bar \varphi}_{M_+}  \ = \ 2\arctan  \left( {{t \,+\, \Delta}\over {\lambda_M}}
\right)  & = & 2\arctan \,\sqrt{ {{ t \,+ \,\Delta}\over {t \,- \, \Delta}}\,}\,.
\end{eqnarray*}
We note that
\begin{eqnarray*}
\theta_- \!\!\!\!\ & = &\  \!\!\arctan \chi \ \ \ \ \,\Longleftrightarrow \ \ \ \ \tan \,
\theta_- \ =\  \chi\,,\\
\theta_+ \!\!\!\!\ & = & \!\!\ \arctan {1\over \chi} \ \ \ \ \Longleftrightarrow \ \ \ \
\tan \, \theta_+ \ =\  {1\over \chi}
\ \ \ \ \ \ \ \ \ \ \ \ \ (\,{\mbox{here}} \ \ \ 0 \,< \,\theta_-\,, \ \ \, \theta_+ \,< \,\pi/2\,)\\
& \Longrightarrow & [\tan \, \theta_-]\cdot [\tan \, \theta_+] \ =\  1\ \
 \Longrightarrow  \ \ {{ [\sin \, \theta_-]\cdot [\sin \, \theta_+]}\over {[\cos \, \theta_-]\cdot [\cos \, \theta_+] }} - 1 \ = \ 0\\
& \Longrightarrow & {{ [\sin \, \theta_-]\cdot [\sin \, \theta_+] - [\cos \, \theta_-]\cdot [\cos \, \theta_+] }\over {[\cos \, \theta_-]\cdot [\cos \, \theta_+] }}  \ =\  0 \ \
  \Longrightarrow  \ \  \cos \,(\,{\bar \varphi}_{M_+}    +  {\bar \varphi}_{M_-} ) \ =\ 0\\
& \Longrightarrow &  \!\! \ \theta_+  + \theta_-  = {\pi\over 2} \ \
\Longrightarrow \ \   \theta_+   \ = \ {\pi\over 2} -  \theta_-\ \ \ \ \ \ \ \ \ \ ({\mbox{as}} \ \ 0 \,<\, \theta_-\,, \ \, \theta_+ \,<\, \pi/2\,)\ \ \ \ \ \ \ \ \ \  \\
& \Longrightarrow &    \!\! {\bar \varphi}_{M_+}   \ = \ \pi -   {\bar \varphi}_{M_-}  \ \ \ \ (\Longrightarrow \ \ {\mbox{domain \ symmetric}}) \ \
  \Longrightarrow \     \sin\,{\bar \varphi}_{M_+}   \, = \, \sin\,{\bar \varphi}_{M_-} \ \ \ \ \ \\
(2.19)\cdot  \cdot \cdot  \!\! &\Longrightarrow &  (\,\sin\,{\bar \varphi}_{M_+} )^{\,n - 1}  \ = \ (\,\sin\,{\bar \varphi}_{M_-} )^{\,n - 1}.
\end{eqnarray*}
Together with (2.13), (2.15)\,--\,\,(2.19),  we arrive at
$$
\left[ {\partial\over {\partial \lambda}} \int^{2\arctan \,\left( {{t\,+ \,\Delta }\over {\lambda}} \right)}_{2\arctan \,\left( {{t\,-\, \Delta }\over {\lambda}} \right)} \,\,\ [\,\sin \varphi]^{\,n - 1}\,   \, d\varphi \right]_{ \lambda\, =\, \lambda_M\,} \  \, = \,0\,. \leqno (2.20)
$$

\vspace*{0.15in}

{\bf Lemma 2.21.} \ \ {\it Given any positive numbers \,$t$ and \,$\Delta$ with \,$t \,>\, \Delta$\,,\, let $\,H\,$ be equal to one in the annular domain $\Omega = {B_o\,(t \,+ \,\Delta) \setminus \overline{B_o\,(t\,-\,\Delta)}}\,,\,$ and zero outside, and $\lambda_M$  given by }\,

\vspace*{-0.3in}

$$
\lambda_M \ = \ \sqrt{(t \,+\, \Delta)(t \,-\, \Delta)\, }\,\,. \leqno (2.22)
$$
{\it Then}\, $\,(\lambda_M\,, \,\,{\vec{\,0}}\ )$ {\it  is a critical point for}\, $G_{|_{\bf Z}}\,(\bullet\,,\, \bullet\,;\, \Omega)$\, {\it given in\,} (2.1)\,.\\[0.2in]
{\bf Proof.} \ \ As $\,G_{|_{\bf Z}}\,(\bullet\,,\, \bullet\,;\, \Omega)$\, is differentiable, we first let $\,\xi = 0\,$ in the expression (2.1), and apply the calculations in (2.11)\,\,--\,\,(2.20) to conclude that
$\, \displaystyle{
{{\partial G_{|_{\bf Z}}}\over {\partial \lambda}} \bigg\vert_{(\lambda_M\,, \,\,{\vec{\,\,0}}\,)}\! = \ 0\,.\,
}$
In  light of Lemma 6.4 in Part I \cite{I}\,,\, we also have
$\ \displaystyle{
{{\partial G_{|_{\bf Z}}}\over {\partial \xi_j}} \bigg\vert_{(\lambda_M\,, \,\,{\vec{\,\,0}}\,)} \!= \ 0\,
}\,$
for $\,j = 1, \ 2, \cdot \cdot \cdot, \ n\,.$\,   \qed
With a closer look, we find that $\lambda_M$ is the only critical point [\,for the expressions in (2.12)]\,.\, As this point is not used in this article, we direct the interested readers to see \S\,A\,.12 in the e\,-\,Appendix.\, Next, we show that it is a non\,-\,degenerate critical point.\\[0.15in]
{\bf Lemma 2.23.} \ \ {\it Under the notations    and the conditions in Lemma\,} 2.21\,,\, {\it $(\lambda_M\,, \, {\vec{\,\,0}}\,)$ is a non-\,degenerate critical point for $\ G_{|_{\bf Z}}\,(\bullet\,,\, \bullet\,;\, \Omega)$\,.}

\vspace*{0.15in}

{\bf Proof.} \ \ We assert that

\vspace*{-0.3in}

$$
{{\partial^2 G_{|_{\bf Z}}}\over {\partial \lambda^2}} \, (\lambda_M\,, {\vec{\,\,0}}\,) \ \,< \ 0\,. \leqno (2.24)
$$
In view of   Lemma 7.7 and Proposition 7.10 in Part I \cite{I}\,,\, and the symmetry\,:
$$
{{\partial^2 G_{|_{\bf Z}}}\over {\partial \,\xi_\ell^2}} \, (\lambda_M\,,  {\vec{\,\,0}}\,) \ \,= \ \,{{\partial^2 G_{|_{\bf Z}}}\over {\partial\, \xi_j^2}} \, (\lambda_M\,, {\vec{\,\,0}}\,) \ \ \ \  \mfor \ \ \ 1 \,\le \,j\,,\ \, \ell \,\le\, n\,,
$$
(2.24) is enough to complete the proof. To do the task, we continue from (2.13),

\vspace*{-0.2in}

\begin{eqnarray*}
  & \ & {{\partial^2}\over {\partial \lambda^2}} \int^{2\arctan \,\left( {{t\,+\, \Delta }\over {\lambda}} \right)}_{2\arctan \,\left( {{t\,-\, \Delta }\over {\lambda}} \right)} \,\,[\,\sin \varphi]^{\,n - 1}\,   \, d\varphi\,\Bigg\vert_{\,\lambda \,=\, \lambda_M}\\[0.075in]
&\  &  \!\!\!\!\!\!\!\!\!\! =  \,{{4\,[\,t + \Delta] \,\lambda_M}\over {\{ [\,{t + \Delta}]^2 + \lambda^2_M\}^2}}\, [\,\sin {\bar \varphi}_{M_+}]^{\,n - 1} \,+\, {{4\,(t + \Delta)^2}\over {\{ [\,{t + \Delta}]^2 + \lambda^2_M\}^2}}\, (n - 1) \,[\,\sin {\bar \varphi}_{M_+}]^{\,n - 2} \cos {\bar \varphi}_{M_+}\\[0.075in]
& \  &  \!\!  - \,{{4\,[\,t - \Delta]\, \lambda_M}\over {\{ [\,{t - \Delta}]^2 + \lambda^2_M\}^2}}\, [\,\sin {\bar \varphi}_{M_-}]^{\,n - 1} - {{4\,(t - \Delta)^2}\over {\{ [\,{t - \Delta}]^2 + \lambda^2_M\}^2}}\, (n - 1)\, [\,\sin {\bar \varphi}_-]^{\,n - 2} \cos {\bar \varphi}_{M_-}\\
 & \ & \\[0.1in]
&\  &  \!\!\!\!\!\!\!\!\!\! =  \ {1\over {t^2}} \left( {{\lambda_M}\over {t + \Delta}} \right)
[\,\sin {\bar \varphi}_{M_+}]^{\,n - 1} \ \, - \ \,{1\over {t^2}} \left( {{\lambda_M}\over
{t - \Delta}} \right) [\,\sin {\bar \varphi}_{M_-}]^{\,n - 1}\\ & \ &  \hspace*{1.5in} \ \ \ \ \ \,
[\,  \uparrow \ \ {\mbox{the \ difference \ is \ negative\,, \ via}} \ \ (2.17\,-\,(2.19)\,]\\
& \ & \ \ \ \ \ \ \ \   +\ {{n - 1}\over {t^2}}\cdot [\,\sin \,\varphi_{M_\pm}]^{\,n - 2} \times [\, \cos {\bar \varphi}_{M_+} - \cos {\bar \varphi}_{M_-}] \ < \ 0\,.
 \end{eqnarray*}
Here we use   the notations introduced in (2.13), and (2.20)\,.\, In the last step above, we make use of (2.17) and (2.18) to  obtain
\begin{eqnarray*}
\ \pi \ > \  {\bar \varphi}_{M_+} \ > \ {\pi\over 2} \ \ \ \ \Longrightarrow \ \ \cos \, {\bar \varphi}_{M_+} < 0\,,\ \ \ \   \
 \ \ {\pi\over 2} >   {\bar \varphi}_{M_-} \ >  \ 0\ \ \ \ \Longrightarrow \ \ \cos \, {\bar \varphi}_{M_-}  \ > \ 0\,.
\end{eqnarray*}
With this, we arrive to (2.24)\,.\qed
It follows from (2.24), Lemma 7.7 and Proposition 7.10
in Part I \cite{I}\,,\,
that
$$
{\mbox{Deg}}\, (\,\btd\, G_{|_{\bf Z}}\,, \ B_{{\bf p_M}} \,(r)\,, \, \vec{\,\,0}\,) \ \,= \ - 1\,, \ \ \ \ \ \ {\mbox{where}} \ \ \ \ {\bf p_M} \ := \ (\lambda_M\,, \, \vec{\,\,0})\,,  \leqno (2.25)
$$
provided that  $\ \btd \,G_{|_{\bf Z}}\, (\lambda\,, \ \xi) \not= \vec{\,\,0}$\,\, for all  $\,(\lambda\,, \ \xi) \in \overline{B_{{\bf p_M}} (r)}\, \setminus \{ {\bf p_M} \}\,.\,$ See \cite{Fonseca-Gangbo} and also \S\,A.3 in the e\,-\,Appendix.\bk
In the process on superimposing annular domains, we must show that effects from other annular domains do not affect the stability of the one we focus on. For this, we require more information on the strength of the  first derivative [\,at points of the boundary of $\,B_{{\bf p_M}} \,(r)\,]\,.\,$ The following result is a prelude to this (cf. \S\,2\,h).

\vspace*{0.2in}

{\bf Lemma 2.26.} \ \ {\it Under the notations    and the conditions in Lemma\,} 2.21\,,\, {\it assume that}
$$
\sqrt{{5\over 2}\,} \ \ge  \ {{t \,+ \,\Delta}\over {\lambda_M}}   \  \ge  \ {{t \,- \,\Delta}\over {\lambda_M}}  \ \ge \   \sqrt{{2\over 5}\,}\,\,.\leqno (2.27)
$$
 {\it There exist positive constant $\,{\bar \varepsilon}_6\,$ and $\,{\bar \varepsilon}_7\,$ such that  }
$$
\bigg\vert \,{{\partial  \,G_{|_{\bf Z}}}\over {\partial \lambda }} \, (\lambda_M + s,\, {\vec{\,\,0}}\,) \bigg\vert \ \,> \ [\,{\bar C}_6  \cdot {\bar \varepsilon}_6\,] \cdot {{1}\over {\lambda_M}} \ \ \ \ \  \ \ \  \ \ \ {\it{for}} \ \ \ \ {\bar \varepsilon}_7 \ \ge \ {{|\,s|}\over {\lambda_M}} \ \ge \ {\bar \varepsilon}_6\,. \leqno (2.28)
$$
{\it The numbers $\,{\bar C}_6$,\, $\,{\bar \varepsilon}_6\,$ and $\,{\bar \varepsilon}_7\,$ can be chosen to be independent on $\,t\,$ and $\,\Delta\,$,\, as long as}\, (2.27) {\it is fulfilled.}\\[0.2in]
{\bf Proof.} \ \  From (2.14) \,[\,cf. (2.26)]\,,\, we obtain
\begin{eqnarray*}
 (2.29)  & \ &  {{\partial^2}\over {\partial \lambda^2}} \int^{2\arctan \,\left( {{t\,+\, \Delta }\over {\lambda}} \right)}_{2\arctan \,\left( {{t\,-\, \Delta }\over {\lambda}} \right)} \,\,[\,\sin \varphi]^{\,n - 1}\,   \, d\varphi\,\Bigg\vert_{\,\lambda \,=\, \lambda_M + s}\\
& \ & \!\!\!\!\!\!\!\!\!\!\!\!\!\!\!\!\!\!\!\!\!\!\!\!\!\!  = \ {{4\,[\,t + \Delta] \,[\,\lambda_M + s] }\over {\{ [\,{t + \Delta}]^2 + [\,\lambda_M + s]^2\}^2}}\cdot [\,\sin {\bar \varphi}_+]^{\,n - 1} \,+\ {{4\,(n - 1)\,(t + \Delta)^2 }\over {\{ [\,{t + \Delta}]^2 + [\,\lambda_M + s]^2\}^2}} \cdot [\,\sin {\bar \varphi}_+]^{\,n - 2} \cos {\bar \varphi}_+ \\
 & \ & \!\!\!\!\!\!\!\!\!\!\!\!\!\!\!\!\!\!\!\!\!\!\!   -\,    {{4\,[\,t - \Delta]\, [\,\lambda_M + s]}\over {\{ [\,{t - \Delta}]^2 + [\,\lambda_M + s]^2\}^2}}\cdot [\,\sin {\bar \varphi}_-]^{\,n - 1}   - {{4\,(n - 1)\,(t - \Delta)^2 }\over {\{ [\,{t - \Delta}]^2 + [\,\lambda_M + s]^2\}^2}}\cdot[\,\sin {\bar \varphi}_-]^{\,n - 2} \cos {\bar \varphi}_-\,\,.\ \ \ \ \ \ \ \ \ \ \ \ \ \ \ \
\end{eqnarray*}

As in (2.13),

\vspace*{-0.25in}

\begin{eqnarray*}
 \ \ \ {\bar \varphi}_- & = & 2\arctan   \left( {{t \,-\, \Delta}\over {\lambda_M + s}}
\right)  \ = \ 2\arctan   \left( {{t \,-\, \Delta}\over {\lambda_M  }} \cdot \left[ 1 \,+\, O \left({{s}\over {\lambda_M}} \right) \right]
\right)    \\
& = & 2\arctan   \left( {{t \,-\, \Delta}\over {\lambda_M  }}  \,+\, O \left({{s}\over {\lambda_M}}
\right) \right) \ \ \ \ \ \ \    [\,{\mbox{via}} \ \ (2.27)]\\
 &= & 2\arctan   \left( {{t \,-\, \Delta}\over {\lambda_M  }} \right) \ + \ O \left({{s}\over {\lambda_M}} \right) \ \ \  \ \ \ \ ({\mbox{first \   order \   approximation \    of \   arctan}})
\\
& = &  {\bar \varphi}_{M_-} \ + \  O \left({{s}\over {\lambda_M}}
\right)\ \ \ \ \ \ \ \ \ \ \ \ \ \ \ \ \ \ \ \ \   \mfor  \, \, \lambda_M^{-1}\cdot |\,s| \ \ {\mbox{small}} \ \ \ [\,{\mbox{cf.}} \ \ (2.17)]\,.
\end{eqnarray*}
Here we use the expansion
$$
{1\over {\lambda_M + s}} \ = \ {1\over {\lambda_M }}  \left( {1\over {1 + \lambda_M^{-1}\,s}} \right)\ =\  {1\over {\lambda_M }}  \left[ \, 1 \ + \ O \,\left(\lambda_M^{-1}\,s \right)\,\right] \ \ \  \ \ \ \ \ \  \ \mfor \lambda_M^{-1}\cdot |\,s| \ \ {\mbox{small}}\,,
$$

\vspace*{-0.1in}

and the bound on the derivative of $\arctan$. Similarly,
$\,\displaystyle{
 {\bar \varphi}_+   =   {\bar \varphi}_{M_+} +  O \left({{s}\over {\lambda_M}}
\right)}\,$ for $\,\lambda_M^{-1}\cdot |\,s|$\, small\,,\,
where we use (2.27) again.\, Moreover,
\begin{eqnarray*}
& \ & {{4\,[\,t \,+\, \Delta] \,[\,\lambda_M \,+\, s] }\over {\{ [\,{t \,+\, \Delta}]^2 + [\,\lambda_M \,+\, s]^2\}^2}}   \ = \   {1\over {\lambda_M^2}} \cdot {{4\,[(\lambda_M^{-1}\,(t \,+\, \Delta)] \,[\,1 \,+\,\lambda_M^{-1} \,s] }\over {\{ [\, \lambda_M^{-1}\,(t \,+\, \Delta)]^2 \ + \ [\,1 \,+\,\lambda_M^{-1}\, s]^2\}^2}}\\
& = &{1\over {\lambda_M^2}} \cdot \left[ {{4\,[(\lambda_M^{-1}\,(t \,+\, \Delta)]   }\over {\{ [\, \lambda_M^{-1}\,(t \,+\, \Delta)]^2 \,+\, 1\}^2}} \ +\  O \left({{s}\over {\lambda_M}}
\right)  \right] = {1\over {t^2}} \left( {{\lambda_M}\over {t \,+\, \Delta}} \right) \ + \ {1\over {\lambda_M^2}} \cdot  O \left({{s}\over {\lambda_M}}
\right)\\
& = & {1\over {t^2}} \left( {{\lambda_M}\over {t\, + \,\Delta}} \right) + {1\over {t^2}} \cdot  O \left({{s}\over {\lambda_M}}
\right)\,,\ \
  \ \ \ \ \ \ \ \ \ \ \ \ \   {\mbox{as \ \ \   (2.27)}} \ \ \   \Longrightarrow   \ \ \sqrt{{5\over 2}\,} \ \ge \  {{t  }\over {\lambda_M}}    \ \ge \ \sqrt{{2\over 5}\,}\,.
\end{eqnarray*}
We apply a similar method to estimate the other terms in (2.29), and obtain

\vspace*{-0.25in}

\begin{eqnarray*}
(2.30) \ \ \ \ \    & \ & {{\partial^2}\over {\partial \lambda^2}} \int^{2\arctan \,\left( {{t\,+\, \Delta }\over {\lambda}} \right)}_{2\arctan \,\left( {{t\,-\, \Delta }\over {\lambda}} \right)} \,\,[\,\sin \varphi]^{\,n - 1}\,   \, d\varphi\,\Bigg\vert_{\,\lambda \,= \,\lambda_M \,+\, s}\\[0.075in]
& \  & \!\!\!\!\!\!\!\!\!\!\!\!\!= \ \left[ \ {1\over {t^2}} \left( {{\lambda_M}\over {t \,+\, \Delta}} \right)
[\,\sin {\bar \varphi}_{M_+}]^{\,n - 1} \,-\ \, {1\over {t^2}} \left( {{\lambda_M}\over
{t \,-\, \Delta}} \right) [\,\sin {\bar \varphi}_{M_-}]^{\,n - 1} \,\right] \ \,+ \\
 & \ & \!\!\!\!\!\!\!\!\!\!\!\!\!\!\!\!\!\!\!  
[\, {\mbox{as \ \ in \ \ the \ \ proof \ \ of \ \ (2.24)}}\ \ \uparrow  {\mbox{  \ \ the \ difference \ is \ negative}}\,]\\
& \ & \\
& \ & \ \    +\ {{n - 1}\over {t^2}}\cdot  [\,\sin \,\varphi_{M_\pm}]^{\,n - 2} \cdot [\, \cos {\bar \varphi}_{M_+} \,- \,\cos {\bar \varphi}_{M_-}]  \ \,+ \ \, {1\over {t^2}} \cdot  O \left({{s}\over {\lambda_M}} \right)\,.\ \ \ \  \ \ \ \ \ \ \ \ \ \ \ \ \ \
\end{eqnarray*}
We make use of (2.17), (2.18), (2.27) to  deduce that there are positive  constants $\,c_1$,\, $\,c_2$,\, $\,c_3$\, and $\,c_4\,$ \, [\,independent on $\,t\,$ and $\,\Delta\,$\, as long as (2.27) is fulfilled\,] so that
\begin{eqnarray*}
  \pi - c^2_1 \ > \  {\bar \varphi}_+ \ > \ {\pi\over 2} + c^2_1 \ \ & \Longrightarrow &\ \ \cos \, {\bar \varphi}_+ \ < \ -c^2_2 \ \ \ \ {\mbox{and}} \ \ \    \sin \,\varphi_{M_+} \ge c^2_2\  >\  0\,,\\
  {\pi\over 2} - c^2_3 \ >  \  {\bar \varphi}_- \ > \  c^2_3\ > \  0\ \ &\Longrightarrow &\ \ \cos \, {\bar \varphi}_-  \ > \ c^2_4\ \ \ \ \ \  \,{\mbox{and}} \ \ \ \sin \,\varphi_{M_-} \ \ge\  c^2_4 \ > \ 0
\,.
\end{eqnarray*}
Hence we can find a positive number $\,{\bar \varepsilon}_7\,$ so that
$$
{{\partial^2}\over {\partial \lambda^2}} \int^{2\arctan \,\left( {{t\,+\, \Delta }\over {\lambda}} \right)}_{2\arctan \,\left( {{t\,-\, \Delta }\over {\lambda}} \right)} \,\,[\,\sin \varphi]^{\,n - 1}\,   \, d\varphi\,\Bigg\vert_{\,\lambda \,=\, \lambda_M + s}  \ \le\  -\, {{C_1}\over {t^2}} \ + \ {{C_2}\over {t^2}} \cdot  \left({{s}\over {\lambda_M}} \right) \ \,\le \  - \, {{C_3}\over {t^2}}
$$
for $\,\lambda_M^{-1}\cdot |\,s|\le {\bar \varepsilon}_7$\,.\, Consider the case $\,s > 0\,$ first.   Lemma 2.21 and an integration yield

\vspace*{-0.3in}

$$
{{\partial }\over {\partial \lambda }} \int^{2\arctan \,\left( {{t\,+\, \Delta }\over {\lambda}} \right)}_{2\arctan \,\left( {{t\,-\, \Delta }\over {\lambda}} \right)} \,\,[\,\sin \varphi]^{\,n - 1}\,   \, d\varphi\,\Bigg\vert_{\,\lambda \,= \, \lambda_M + s} \le  \ \, - {{C_3}\over {t^2}} \cdot s  \ \le\   - \,{{C_4}\over {\lambda_M }} \cdot{ s\over {\lambda_M}} \ \le \  - \,{{C_4}\over {\lambda_M }} \cdot {\bar \varepsilon}_6
$$
for $\ {\bar \varepsilon}_6 \,\le\, \lambda_M^{-1}\cdot  s  \,\le\, {\bar \varepsilon}_7\,.$\, Here we use (2.27) again. That is, we have (2.28) for $\,s > 0$\,.\, Similarly, we can handle the case $\,s \,< \,0\,.$ \qedwh
%
%
%
%
%
{\bf \S\, 2\,f.} \ \ {\it Perturbation in \,$\xi\,.$}\ \
In order  to apply Theorem 4.26 in Part I \cite{I}\,,\, we desire to provide lower bounds on $\,|\btd G_{|_{\bf Z}} \, (\lambda, \ \xi)|\,$\, off the critical point $(\,\lambda_M\,   {\vec{\,\,0}}\,)$. This means that we need to vary $\,\xi\,$ as well as $\,\lambda$\,.\, According to the geometric formula (2.5), $\,B_{\,\lambda^{-1} \, \xi} \,(\lambda^{-1} \,R)\,$
is the ball in $\,\R^n$ with center at $\,\lambda^{-1} \, \xi\,$
\,and radius \,$\lambda^{-1} \,R\,.\,$ In accordance,
$$
{\dot{\cal P}}^{-1}\, (\,\partial B_{\,\lambda^{-1} \, \xi} \, (\lambda^{-1} \,R ) ) \leqno (2.31)
$$
is a $\,(n - 1)$\,-\,sphere in $\,S^n \subset \R^{n + 1}.$\, [\, In (2.31)\,,\, $R\,$ is a positive number.]\, It is harder to visualize the effect when $\,\xi\,$ is changed.\, To proceed, we find the (Euclidean) radius $\,\varrho$\, of this $\,(n - 1)$\,-\,sphere\,.\,\bk
Via a rotation, we
assume that $\xi$ is along the {\it positive\,} $y_1$\,-\,axis. Granted this, let
$$
\xi \ = \ (\,\delta, \ 0, \cdot \cdot \cdot, \ 0)\,,
\leqno (2.32)
$$
where we consider  small perturbation so   that $\,0 \,\le\, \delta \,<\, R\,.\,$
The equation of
$\,\partial B_{\,\lambda^{-1} \, \xi} \,(\lambda^{-1} \,R )\,
$ is

\vspace*{-0.33in}

$$
(y_1 \,-\, \delta_\lambda)^2 \,+\, y_2^2 \,+\, \cdot \cdot \cdot \,+\, y_n^2 \ = \ R^{\,2}_\lambda\,, \leqno (2.33)
$$

\vspace*{-0.1in}

where

\vspace*{-0.3in}

$$
R_\lambda \ =\ {{R }\over \lambda} \ \ \ \ \ \ {\mbox{and}} \ \ \ \ \delta_\lambda \ =\  {{\delta} \over \lambda}\,.  \leqno (2.34)
$$
(2.33) and (2.34)  imply
$$
   r^2 \,- \,2y_1 \cdot \delta_\lambda\ = \ R^2_\lambda \,-\, \delta^2_\lambda \ \  \ \ \ \ \Longrightarrow   \ \ \ {{1\, +\, x_{n + 1}}\over {1 \,-\, x_{n - 1} }} \ - \
{{2\,\delta_\lambda\cdot x_1}\over {1 \,- \,x_{n + 1}}}\  =\  R^2_\lambda \,-\, \delta^2_\lambda\ .
$$
(As usual, $r^2 \,= \, y_1^2 \,+\, \cdot \cdot \cdot \,+ \, y_n^2$\,.)\,
See (2.1) and (2.2).\, The result is
$$
(-
2\,\delta_\lambda)\,x_1 \ +\   x_{n + 1}\,(1 \,+\, R^2_\lambda \,- \,\delta^2_\lambda)   \ = \
(R^2_\lambda \,-\, \delta^2_\lambda \,-\, 1)\,.\leqno (2.35)
$$
That is,
$${\dot{\cal P}}^{-1}\, (\,\partial B_{\,\lambda^{-1} \, \xi}\,
(\lambda^{-1} \,R) )=
S^n \ \cap \ \{ \,{\mbox{The \ hyperplane \ defined \ by \ equation}}  \ \
(2.35) \,\} \leqno (2.36)
$$

To find the Euclidean radius, we measure the
distance in $\,\R^{n + 1}$ of the inverse images  of the end points

\vspace*{-0.3in}

  $$-\,(R_\lambda \,-\, \delta_\lambda) \ \ \ \ \ {\mbox{and}} \ \ \ \ \ (R_\lambda \,+\, \delta_\lambda)\,, \leqno (2.37)
  $$
respectively (these two points lie on the $y_1$\,-\,axis)\,. In the plane defined by $\,x_1$\,-\,axis and $\,x_{n + 1}$\,-\,axis,  the inverse images (via ${\dot{\cal P}}^{-1}$) of these two points  have coordinates given by
$$
 \left(  -\,{{ 2 (R_\lambda \,-\, \delta_\lambda)}\over { (R_\lambda - \delta_\lambda)^2 \,+\, 1}} \,\,,
\  \
{{ (R_\lambda - \delta_\lambda)^2 \,-\, 1}\over { (R\, -\, \delta_\lambda)^2 \,+\, 1}}    \right)\,,
 \ \ \ \ \ \ \ \ \
 \left(   {{ 2 (R_\lambda\,+\,\delta_\lambda)}\over { (R \,+ \, \delta_\lambda)^2 \,+\, 1}} \,\,, \  \
  {{ (R_\lambda \,+\, \delta_\lambda)^2 \,- \,1}\over { (R_\lambda \,+ \,\delta_\lambda)^2 \,+ \,1}}    \right)\,,
  $$
respectively. Refer to (2.1) and (2.2).\, As $2\,\varrho$ is the distance between   these two points (measured in the standard
Euclidean metric in $\R^{n + 1}$\,),\, a calculation shows that (refer to \S\,A.13 in the e\,-\,Appendix)
$$
  (2 \varrho)^2 \ = \ {{ (4R_\lambda)^2\, [\, (1 \,+\, R_\lambda^{\,2} \,-\, \delta_\lambda^2)^2\,+\,  4\,\delta_\lambda^2\,]}\over
{[(R_\lambda \,+\, \delta_\lambda)^2 \,+\, 1]^2 \ [(R_\lambda \,-\, \delta_\lambda)^2 \,+ \,1]^2}}\,\,. \leqno (2.38)
$$
Compare  with the non-perturbed case [\,i.e., when $\,\xi \,= \, 0$\,]\, which can also be found directly by using (2.2) on the points $\,R_\lambda\,$ and $\,-R_\lambda\,$ along the $x_1$\,-\,axis\,:
$$
(2 \varrho)^2 = \left[ {{2 R_\lambda}\over {1 + R^2_\lambda}} \,-\,{{(-2 R_\lambda)}\over {1 + R^2_\lambda}}  \right]^2  =   {{(4R_\lambda)^2}\over {(1 + R_\lambda^{\,2})^2}}  \ \left(\! = (4R)^2 \cdot {{\lambda^2}\over {(R^2 + \lambda^2)^2}}\!\right)  \mfor \xi = 0\,.  \leqno (2.39)
$$
For later  reference in \S\,2\,h\,,\,  we differential both sides of (2.39) and obtain the following.
$$
{{\partial \,(2 \varrho)^2}\over {\partial \lambda}}  \ = \  {{ (2 \lambda)\,(4R )^2}\over {(R^2 \,+\, \lambda^2)^3}}\cdot (R^2 \,-\, \lambda^2) \ = \ {1\over \lambda} \,\left[\,{{  2 \cdot (4R_\lambda)^2}\over {(R^2_\lambda \,+ \,1)^3}} \  (R^2_\lambda - 1)\right] \ \  \,\mfor \,\xi \, = \, 0 \,.  \leqno (2.40)
$$

\vspace*{0.2in}

{\bf \S\, 2\,g.}  {\it Derivative in \ $\xi\,.\,$}\ \
Recall that we assume $\,\xi \, = \,(\delta\,, \ 0\,,\, \cdot \cdot \cdot\,,\ 0)$\,.\, When differentiating the denominator in (2.38), we make use of the following.
\begin{eqnarray*}
(2.41)  & \ & {{1}\over
{  [(R_\lambda+ \delta_\lambda)^2 \,+\, 1\,]  \cdot [(R_\lambda  - \delta_\lambda )^2 \,+\, 1\,]  }}\\
& = & {{1}\over
{  [\,(R_\lambda^{\,2} + \delta_\lambda^2 + 1) \,+\, 2 R_\lambda\cdot \delta_\lambda]\cdot [\,(R_\lambda^{\,2} + \delta_\lambda^2 + 1) - 2 R_\lambda\cdot \delta_\lambda\,] }}\\
& = & {{1}\over
{   (R_\lambda^{\,2} + \delta_\lambda^2 + 1)^2  \,-\, 4 R_\lambda^{\,2}\cdot \delta_\lambda^2   }}\ \
 \ \ (\,{\mbox{cancelation \ \ of \ \ first \ \ order \ \ terms \ \ in }}\ \  \delta_\lambda) \ \ \ \ \\
& =& {{1}\over
{   (1 +  R_\lambda^{\,2}\,)^2\, + \,\delta_\lambda^2 \,( 2 + \delta_\lambda^2 - 2 R_\lambda^{\,2}\,)  }} \ \ \ \ \ \ \ \ \ \ \ \ \ \ \ \ \ \ \ \ \ \ \ \ \ \ \ \ \ \ \  \left({\mbox{recall}} \ \  \delta_\lambda \,= \,{\delta\over \lambda} \right)\,.
\end{eqnarray*}

One finds that

\newpage

\begin{eqnarray*}
(2.42) & \ &  {{\partial \,(2\varrho)^2}\over {\partial \delta}}
\ = \ -\ {\delta_\lambda\over \lambda}\cdot
{{  (4R_\lambda)^2}\over {[(R_\lambda + \delta_\lambda)^2 \,+ \, 1\,]^2
\,[(R_\lambda - \delta_\lambda)^2 \,+\, 1\,]^2}} \,\times \\
& \ &  \times
\left\{ 4   \,[\,(R_\lambda^{\,2} - 1) -\delta_\lambda^2\,] - \,{{8  \,[(1 + R_\lambda^{\,2} - \delta_\lambda^2)^2 \,+\,
4\,\delta_\lambda^2\,]}\over
{[(R_\lambda + \delta_\lambda)^2 \,+\, 1\,] \,[(R_\lambda - \delta_\lambda)^2 \,+\, 1\,]}} \cdot
\left[\,(R_\lambda^{\,2} - 1) \,\,-\delta_\lambda^2\,\right]    \right\} \ \ \ \ \ \ \ \ \\[0.1in]
&=&  -\ {\delta_\lambda\over \lambda}\cdot[\,(R_\lambda^{\,2} - 1) -\delta_\lambda^2\,] \,\times
{{  (4R_\lambda)^2}\over {[(R_\lambda + \delta_\lambda)^2 \,+ \, 1]^2
\,[(R_\lambda - \delta_\lambda)^2 \,+\,1\,]^2}} \times \ \ \ \ \ \ \\
& \ &  \ \ \ \ \  \ \ \ \ \ \ \ \   \ \ \ \ \ \  \ \ \ \ \ \ \ \  \ \ \ \ \ \ \ \  \ \ \ \ \ \ \times
\left\{ 4    \  - \,{{8  \,[(1 \,+\, R_\lambda^{\,2} - \delta_\lambda^2)^2 +
4\,\delta_\lambda^2\,]}\over
{[(R_\lambda + \delta_\lambda)^2 \,+ \,1\,] \,[(R_\lambda - \delta_\lambda)^2 \,+\, 1\,]}} \, \,
   \right\}\,.
\end{eqnarray*}

Continuing from the last line in (2.41), one can approximate
\begin{eqnarray*}
(2.43) \ \ \ \ \ \ \ & \ & {{1}\over
{  [(R_\lambda+ \delta_\lambda)^2 \,+\,1\,]  \cdot [(R_\lambda  - \delta_\lambda )^2 \,+\,1\,]  }}\\
& =& {{1}\over
{   (1 +  R_\lambda^{\,2}\,)^2 }}\cdot {1\over {1 + \delta_\lambda^2 \cdot  {{ ( \,2 \,+\, \delta_\lambda^2 \,-\, 2 R_\lambda^{2}) }\over {(1 \,+\,  R_\lambda^{2})^2}} }}  \ = \ {{1}\over
{   (1 +  R_\lambda^{\,2}\,)^2}} \ +\   O \,(\delta_\lambda^2)\ \ \ \ \ \ \ \ \ \ \ \ \ \ \ \ \ \ \ \ \ \ \ \ \ \ \ \
\end{eqnarray*}
for $\,\delta_\lambda > 0$\, small\,,\, and  \,$R_\lambda \le \sqrt{ {5\over 2} \,}\,.$\,\,
It follows that

\begin{eqnarray*}
(2.44) \ \ \ \ \ \   {{\partial \,(2\varrho)^2}\over {\partial \delta}} & = & -{\delta_\lambda\over \lambda}\cdot[\,(R_\lambda^{\,2} - 1) -\delta_\lambda^2\,] \cdot \left[ {{(4R_\lambda)^2 }\over {(1 + R_\lambda^{\,2})^4  }} + O \,(\delta_\lambda^2)\right]\cdot [\, 4 - 8 +  O \,(\delta_\lambda^2)]\ \ \ \ \ \  \ \ \ \ \ \  \\
  & = & {\delta_\lambda\over \lambda}\cdot \left[  (R_\lambda^{\,2} - 1)  \cdot   {{4\cdot(4R_\lambda)^2 }\over {(1 + R_\lambda^{\,2})^4  }} + O \,(\delta_\lambda^2) \right]
\end{eqnarray*}
for $\,\delta_\lambda > 0\,$ small\,,\, and $\, R_\lambda \le \sqrt{ {5\over 2} \,}\,.\ $  Formula (2.44) indicates that
$$
{{\partial \,(2\varrho)^2}\over {\partial \delta}} > 0 \ \   {\mbox{when}} \ \ R_\lambda > 1\,; \ \ \ \  {\mbox{whereas}} \ \ {{\partial \,(2\varrho)^2}\over {\partial \delta}} < 0 \ \  {\mbox{when}} \ \ R_\lambda < 1
$$
for $\,\delta_\lambda$ small enough (relative to $R_\lambda$)\,.

\vspace*{0.15in}

{\it Expansion for small variation of \,$\lambda_M$\,.\,} \ \  Take
\begin{eqnarray*}
(2.45)    \!\!\!\!\!\!& \ & \ \ \ \lambda = \lambda_M + s \ \ {\mbox{and}} \ \ R = t + \Delta\\[0.075in]
\Longrightarrow  & \ & \!\!\!\!\!\!R_\lambda^{\,2} = \left(\, {{t + \Delta}\over {\lambda }} \right)^2 \,=\, \left( {{t + \Delta}\over {\lambda_M + s}} \right)^2 \,=\, \left( {{t + \Delta}\over {\lambda_M  }} \right)^2 + O \left( {s\over {\lambda_M}} \right) \,=\, {{t + \Delta}\over {t - \Delta}} + O \left( {s\over {\lambda_M}} \right)\\
& \ & \ \ \ \ \ \ \  \ \ \ \ \ \  \ \ \ \ \ \ \  \ \ \ \ \ \  \ \ \ \ \ \ \  \ \ \ \ \ \  \ \ \ \ \ \ \  \ \ \ \ \ \  \ \  \ \   \left( \ \uparrow \ \ {\mbox{cf. \ \ (2.43)}} \right)\\
\Longrightarrow & \ & \!\!\!\!\!\!\!\!\!\!(R_\lambda^{\,2} - 1)  \cdot   {{4\cdot(4R_\lambda)^2 }\over {(1 + R_\lambda^{\,2})^4  }} \, = \,    {{ {{128 \Delta}\over {t - \Delta}} \cdot{{t + \Delta}\over {t - \Delta}}  }\over {  \left( {{2 t}\over {t - \Delta}} \right)^4 }}  + O \left(  \lambda_M^{-1} s \right)  =     {{  128 \Delta  \,(t + \Delta)(t - \Delta)^2 }\over {  \left(  2 t  \right)^4 }} + O \left(  \lambda_M^{-1} s \right)\\
\Longrightarrow  & \ & {{\partial \,(2\varrho)^2}\over {\partial \delta}} \bigg\vert_{\,R \,=\, t \,+\, \Delta} \ = \   {\delta_\lambda\over \lambda}\cdot \left[  4^3  \cdot  {{  2 \Delta  \,(t + \Delta)(t - \Delta)^2 }\over {  \left(  2 t  \right)^4 }} \ + \ O \left( {s\over {\lambda_M}} \right) \ + \ O \,(\delta_\lambda^2) \right]
\end{eqnarray*}
for $\,\delta_\lambda$\, and \,$\lambda_M^{-1}\, |\,s|$\, small\,,\, and $ \,R_\lambda\, \le\, \sqrt{ {5\over 2} \,}\,.$\, Likewise, take $\,\lambda\, = \,\lambda_M \,+\, s$\, and $\,R \,=\, t\, - \,\Delta\,,\,$ we obtain
$$
 {{\partial \,(2\varrho)^2}\over {\partial \delta}} \bigg\vert_{\,R \,= \,t\, - \,\Delta} \ = \ -\  {\delta_\lambda\over \lambda}\cdot \left[  4^3  \cdot  {{  2 \Delta  \,(t - \Delta)(t + \Delta)^2 }\over {  \left(  2 t  \right)^4 }}\ + \ O \left( {s\over {\lambda_M}} \right)\  + \ O \,(\delta_\lambda^2) \right]
\leqno (2.46)
$$
for $\,\delta_\lambda$\, and $\,\lambda_M^{-1}\, |\,s|\,$ small, \,and   $\,R_\lambda \,\le \,\sqrt{ {5\over 2} \,}\ .$\
In the following discussion, we restrict ourselves to the case
\vspace*{-0.3in}

$$
1 \,-\, A^2  \ \ge\  {{ \Delta}\over {t  }} \ \ge\  B^2 \ > \ 0  \leqno (2.47)
 $$
for some (fixed) positive numbers \,$A$ and $B$\,,\,  and
 $$ \sqrt{5\over 2\,} \ \ge \ {{t \,+\, \Delta}\over {\lambda_M}} \  > \ {{t \,- \,\Delta}\over {\lambda_M}}   \ \ge \ \sqrt{2\over 5\,}\,\,. \leqno (2.48)
$$

\vspace*{0.15in}

{\bf Lemma 2.49.} \ \ {\it Under the notations    and conditions in Lemma\,} 2.21\,,\, {\it assume also  \,} (2.47) \,{\it and} \, (2.48)\,.\,
{\it Let}
$$
  \lambda \ = \ \lambda_M \,+ \,s \ \ \ \ {\it{and}} \ \ \ \ \delta \ = \ |\,\xi|\,.
$$
{\it There exist  positive constants \,${\bar \varepsilon}_8$\,,\, $\,{\bar \varepsilon}_9$\,,\, ${\bar \varepsilon}_{10}$  and \,${\bar C}_7$\,}  {\it such that if }

\vspace*{-0.25in}

$$ {\bar \varepsilon}_8 \ \ge \ \lambda^{-1}_M \,  \,\delta  \ \ge\  {\bar \varepsilon}_9 \ >\  0 \ \ \ \ \ \  {\it{and}} \ \ \ \   \lambda_M^{-1}\,|\,s|   \ \le\  {\bar \varepsilon}_{10}\,,  \leqno (2.50)
$$
{\it   then we  have}

\vspace*{-0.35in}

$$
\Vert  \btd_\xi \, G_{|_{\bf Z}} \, (\lambda_M + s\,, \ \xi)\Vert \ = \ \ \sqrt{ \ \sum_{\ell = 1}^n \,\bigg\vert\, {{\partial G_{|_{\bf Z}} }\over {\partial \xi_\ell}} \  (\lambda\,, \, \xi) \bigg\vert^{\,2} \,} \ \ \, \ge   \   [\,{\bar C}_7\ \cdot {\bar \varepsilon}_{9}\,]\cdot {1 \over {\lambda_M}}\ . \leqno (2.51)
$$
{\it In\,} (2.50) {\it and\,} (2.51), {\it the numbers \, ${\bar \varepsilon}_8$\,,\, ${\bar \varepsilon}_9$\,,\, ${\bar \varepsilon}_{10}$  and \,${\bar C}_7$\,    are independent on $t$ and $\,\Delta$   as long as\,} (2.47)  {\it and\,} (2.48)  {\it are fulfilled}\,.\\[0.2in]
{\bf Proof.} \ \ As the situation is rotationally symmetric, we may assume that
 $$
 \xi \ =  \ (\,\delta\,,\, 0\,,\, \cdot \cdot \cdot\,,\, 0)\,, \ \ \ \ {\mbox{where}} \ \ \delta \,>\, 0\,.
 $$

 \vspace*{-0.25in}

Moreover,

\vspace*{-0.3in}

$$
\lambda \ = \ \lambda_M \,+\, s \ \ \ {\mbox{and}} \ \ \ \lambda_M^{-1} \, |\,s| \le {\bar\varepsilon}_{10} \ \ \ \ \Longrightarrow \ \ (1 \,-\, {\bar\varepsilon}_{10}) \, \lambda_M \ \le \ \lambda \ \le\  (1 \,+ \,{\bar\varepsilon}_{10}) \, \lambda_M\,, \leqno (2.52)
$$
In formula (2.38)\,,\, we let
 $$
 \,\varrho_b\, \ \ \ {\mbox{when}} \ \ \   \,R \ =\  R_b = t \,+\, \Delta\,;\, \ \ \ \ \ {\mbox{and \ \ by}} \ \  \,\varrho_c\, \ \ \ {\mbox{when}} \ \  \,  R \ = \ R_c \ = \ t \,- \, \Delta\,.\,
 $$
Clearly \,$\varrho_b$\, and \,$\varrho_c$\, depend on $\,\lambda\,$ and $\,\xi\,$.\,
 Using (2.43) and a argument similar to that in   (2.45), we obtain

 \newpage

\begin{eqnarray*}
 (2.53) \ \ \ \ \ & \ &      (2 \,\varrho)^2 \ \, = \  {{ (4 R_{\lambda_M})^2}\over { (1 + R_{\lambda_M}^2)^2}} \ +  \ O \left( {s\over {\lambda_M}} \right) \ + \ O \,(\delta_\lambda^2)\\[0.075in]
\Longrightarrow \ \ \   \varrho_b &=& {{\lambda_M}\over t} +  O \left( {s\over {\lambda_M}} \right) \ + \ O \,(\delta_\lambda^2) \ \ \ \ {\mbox{and}} \ \ \  \,\varrho_c  \ = \  {{\lambda_M}\over t} \ + \  O \left( {s\over {\lambda_M}} \right) \ + \ O \,(\delta_\lambda^2) \ \ \ \\
  & \ & \ \ \ \ \ \ \ \ \ \ \ \ \     \   \mfor  \ \delta_\lambda   \ \ {\mbox{and}}\, \ \ \lambda_M^{-1}\, |\,s| \ \ \ {\mbox{small}}\,,\, \ \ {\mbox{and}} \ \  R_{\lambda_M} \le \sqrt{ {5\over 2} \,} \ \ \ \ \ \ \ \ \ \   \\
\Longrightarrow \ \  \varrho_b  &= &  \varrho_c  \ + \  O \left( {s\over {\lambda_M}} \right)\  + \  O \,(\delta_\lambda^2\,)\,.
\end{eqnarray*}
Applying (2.48) and (2.49), we have
\begin{eqnarray*}
  1 - {\Delta\over t} \ \ge \ a^2\ \ \ \
 \Longrightarrow   \  \   {{R_b}\over {\lambda}} &= & {{R_b}\over {\lambda_M \,+\, s}} \ = \ {{R_b}\over {\lambda_M  }}  \ + \ O \left( {|\,s|\over {\lambda_M}} \right)\ \
   \ \ \ \ \ \ \ \  \ \ \ \ \ \ \ \     \left[\,{\mbox{cf.}} \ \ (2.52)\,\right] \\
& = & \sqrt{{{t + \Delta}\over {t - \Delta}}\,  } \ + \ O\! \left( {|\,s|\over {\lambda_M}} \right) \ = \ \sqrt{1 + {{2\, (\Delta/t)}\over {1 - (\Delta/t)}}  } \ \ + \ O\! \left( {|\,s|\over {\lambda_M}} \right)\\
 \Longrightarrow    \   \ {{R_b}\over {\lambda}}   & \ge & \sqrt{1 + {{2\, b^2}\over {a^2}}  \,}\ \, +  \ \, O \!\left( {|\,s|\over {\lambda_M}} \right)\,.
\end{eqnarray*}
Observe that
$\displaystyle{\
{{R_b}\over {\lambda_M}} \ =\  \sqrt{ {{t \,+\, \Delta}\over {t \,-\, \Delta}}\,} \  \le \   \sqrt{ 5\over 2\,}\,.\ \,}
$
Thus if we choose $\,{\bar \varepsilon}_8\,$ and $\,{\bar \varepsilon}_{10}\,$ small enough (they depend on $\,{\bf a}\,$ and \,$b$\,)\,,\, we can find a positive number $\,c\,$ such that

\vspace*{-0.05in}

$$
{{R_b}\over {\lambda}} \ \ge\  1 \,+\, c^2 \ \ \ \ \Longrightarrow \ \ {\displaystyle{B_o (1) \ \subset\      B_{  {{-\,\xi}\over \lambda} }\,
(\, \lambda^{-1} R_b\, )}}\,.
$$

\vspace*{-0.02in}

Hence the boundary sphere
$ \,
\displaystyle{  \partial \,{\dot {\cal P}}^{-1} \,  \left(   B_{\! {{-\,\xi}\over \lambda} }\,
(\, \lambda^{-1} R_b )
\right) }\,
$
lies inside the northern hemisphere. Let us denote its radius (measured in the Euclidean metric in $\R^{n + 1}$) by $\,{\varrho_{b_\lambda}}\,.\,$\, Together with the upper bound in (2.48), we can find positive constants $\,c_1\,$ and $\,c_2\,$ so that
$$
1 \,-\, c_1^2  \,\,\ge\,\, r_{b_\lambda} \,\,\ge\,\, c_2^2  \,\,> \,\,0\,. \leqno (2.54)
$$

\vspace*{-0.1in}

Likewise, the boundary sphere

\vspace*{-0.25in}

$$
\displaystyle{  \partial \,{\dot {\cal P}}^{-1} \,  \left(   B_{ {{-\,\xi}\over \lambda} }\,
(\, \lambda^{-1} R_c )
\right) }\leqno (2.55)
$$

\vspace*{-0.1in}

lies inside the southern hemisphere when we suitably choose $\,{\bar \varepsilon}_8\,$ and $\,{\bar \varepsilon}_{10}$\,.\, Denote by $\,{\varrho_{c_\lambda}}\,$ the radius (again measured in the Euclidean metric in $\R^{n + 1}$\,) of the sphere in (2.55).\, Similarly, via (2.47) and (2.48), there are positive numbers $\,c_3\,$ and $\,c_4$\, such that
$$
1 \,- \,c_3^2  \,\,\ge\,\, r_{b_\lambda}\,\, \ge\,\, c_4^2 \,\, > \,\,0\,. \leqno (2.56)
$$
Thus we are justified to use (2.10), which, together with  (2.45) and (2.46),\,  yield

\newpage

\begin{eqnarray*}
 \!\!& \ & {{\partial G_{|_{\bf Z}} }\over {\partial \xi_1}} \,(\lambda\,,\,\xi)     \ = \    {{\Vert\,S^{n-1}\Vert }\over 2}\, \left\{
{{\varrho_{b_\lambda}^{\,n - 2}}\over {\sqrt{1 - \varrho_{b_\lambda}^2\,}}}
  \cdot \left[  -\, {{d\,[\,\varrho\, (\delta)]^2}\over {d\,\delta}} \bigg\vert_{\,\varrho \,=\, \varrho_{b_\lambda}} \right]\right.\\
  & \ & \ \ \ \ \  \ \ \ \ \   \ \ \ \ \ \ \ \ \ \ \ \ \ \ \ \  \ \ \ \ \ \ \ \ \    \ \ \ \ \ \ \ \    \ \ \ \ \ \ \ \  \left. +\
{{\varrho_{c_\lambda}^{\,n - 2}}\over {\sqrt{1 - \varrho_{c_\lambda}^2\,}}}
  \cdot \!\left[ -\, {{d\,[\,\varrho \,(\delta)]^2}\over {d\,\delta}} \bigg\vert_{\,\varrho \,= \,\varrho_{c_\lambda}} \right] \right\}\\
  & \ & \\[0.1in]
& \  & \!\!\!\!\!\!\!\!\!\!\!=\ -\ {{\Vert\,S^{n-1}\Vert }\over 2} {{\varrho_{b_\lambda}^{\,n - 2}}\over {\sqrt{1 - \varrho_{b_\lambda}^2\,}}} \left\{\! {{d\,[\,\varrho (\delta)]^2}\over {d\,\delta}} \bigg\vert_{\,\varrho \,= \, \varrho_{b_\lambda}} \!\!\!\!\!\!+ {{d\,[\,\varrho (\delta)]^2}\over {d\,\delta}} \bigg\vert_{\,\varrho \,= \, \varrho_{c_\lambda}}\!\right\} + O \left( {{\delta_\lambda}\over {\lambda}} \right) \left[ O \left( {s\over {\lambda_M}} \right) + O \,(\delta_\lambda^2) \right] \ \ \ \ \ \ \\
  & \ & \\[0.1in]
& \  & \!\!\!\!\!\!\!\!\!\!\!=\ {{4^4 \Vert\,S^{n-1}\Vert }\over 2} \cdot {{\delta_\lambda}\over {\lambda}} \cdot {{\varrho_{b_\lambda}^{\,n - 2}}\over {\sqrt{1 - \varrho_{b_\lambda}^2\,}}} \cdot {{\Delta^2\,(t + \Delta) (t - \Delta)}\over {4 \,t^4}}    \  +\  O \!\left( {{\delta_\lambda}\over {\lambda}} \right) \left[ O\! \left( {s\over {\lambda_M}} \right) + O \,(\delta_\lambda^2) \right]\,. \ \ \ \ \ \ \ \ \ \ \ \ \
\end{eqnarray*}
It follows from (2.47), (2.48), (2.54) and (2.56) that
$$
{{r_{b_\lambda}^{\,n - 2}}\over {\sqrt{1 - r_{b_\lambda}^2\,}}} \cdot {{\Delta^2\,(t \,+\, \Delta) (t - \Delta)}\over {4 \,t^4}}\ \ \ge \ {{c_2^{\,2(n - 2)} }\over { \sqrt{1 - (1 - c_1)^2\,} }} \cdot {{ b^4}\over 4} \cdot \left[\,1 - (1 - a^2)^2 \right] \ \ge\  c^2_5\ > \ 0
$$
Hence we have
$$
\bigg\vert\, {{\partial \,G_{|_{\bf Z}} }\over {\partial \,\xi_1}} \,(\lambda\,,\,\xi) \bigg\vert\  \ \ge \  c_6^2 \cdot\delta_\lambda \cdot  {{1}\over {\lambda_M}}  \ \ \ \  \mfor  \ \  \lambda_M^{-1}\, \delta   \ \ \ {\mbox{and}} \ \ \ \lambda_M^{-1}\, |\,s| \ \ \ {\mbox{small}} \ \ \  \ \left( \delta_\lambda  = {\delta\over \lambda}\right)\,.
$$
Combining    with (2.34) and (2.50)   we obtain (2.51). \qedwh
{\bf \S\, 2\,h.}  {\it Derivative in $\lambda$\,.}\ \
Let us express (2.38) in terms of $\,\lambda:\,$
\begin{eqnarray*}
(2.57)  \ \  \ \  (2\,\varrho)^2 & = &   {{ (4R_\lambda)^2\, [\, (1 \,+\, R_\lambda^{\,2} \,-\, \delta_\lambda^2)^2\,+\,  \delta_\lambda^2]}\over
{[(R_\lambda \,+\, \delta_\lambda)^2 \,+\,1\,]^2 \,[(R_\lambda \,-\, \delta_\lambda)^2 \,+\,1\,]^2}} \ \ \ \ \ \ \  \  \left( R_\lambda \ = \ {R\over \lambda}\,, \ \  \delta_\lambda \ =\ {\delta\over \lambda } \right)\ \ \ \ \ \ \ \ \ \ \ \ \ \ \ \ \ \ \ \ \ \ \ \ \ \\[0.1in]
& \  &  \!\!\!\!\!\!\!\!\!\!\!\!\!\!\!\!\!\!\!\!\!\!\!\!\!\!\!\!\!\!\!\!\!\!\! = \  {{\  \lambda^2\,(4R )^2\, [\, (\lambda^2 + R^2  - \delta^2 )^2
+  \delta^2 \,\lambda^2\,]\,}\over
{[(R  + \delta )^2 + \lambda^2]^2 \,[(R  - \delta )^2
+ \lambda^2]^2}}  \ = \ {{ \ \lambda^2\,(4R )^2\, [\, (\lambda^2 + R^2  \,-\, \delta^2 )^2
+  \delta^2 \,\lambda^2\,]\,}\over
{[\, (R^2 \,+\, \delta^2 \,+\, \lambda^2)^2 \,-\, 4R^2 \delta^2\,]^{\,2}}}\,.\   \ \ \ \ \ \ \ \ \ \ \ \
\end{eqnarray*}
Cf. the second last line in (2.41).
The expression suggests that $\,\varrho \to 0\,$ as $\,\lambda \to 0^+\,$ (or $\,\lambda \to \infty\,),\,$ corresponding to pushing the boundary sphere toward the north (resp. south) pole.
 (2.57) guides us to
\begin{eqnarray*}
(2.58) \ \ \ \ \ {{\partial \,(2\,r)^2}\over {\partial \lambda}}  & = &
{{  2\lambda\,(4R)^2\,  }\over
{\{ [(R + \delta)^2 + \lambda^2]  \cdot [(R - \delta)^2
+ \lambda^2]\,\}^2}}\,\times \\
&  & \\
& \ & \!\!\!\!\!\!\!\!\!\!\!\!\!\!\!\!\!\!\!\!\!\!\!\!\!\!\!\!\!
\!\!\!\!\!\!\!\!\!\!\!\!\!\!\!\!\!\!\!\!\!\!\!\!\!\! \times \left[
(\lambda^2 + R^2 - \delta^2)^2
+ 2 \lambda^2 (\lambda^2 + R^2 )\,- \,{{[\,4\,\lambda^2]
[\, (\lambda^2 + R^2 - \delta^2 )^2 + \delta^2\, \lambda^2]\cdot
[R^2 + \delta^2  + \lambda^2] }\over
{  [(R + \delta)^2 + \lambda^2]  \cdot [(R - \delta)^2
+ \lambda^2] }}\right]. \ \ \ \ \ \
\end{eqnarray*}
We note that when $\,\delta  \ =\  0\,$\, (i.e.,\, $\,\xi \,\,= \,\, 0\,$)\, and
$$
\ \ \ \ \ \ \ \ \ \ \ \ \ \ \ \ \  \lambda \,\,=\,\, \lambda_M \ \ \ \ \ \ \ \ \ \  ({\mbox{\,recall \ \ that}} \ \ \lambda_M \,\,=\,\,  \sqrt{ (t\,+\, \Delta)(t\,-\, \Delta)\,}\,, \ \ \ \  t\,\, > \,\, \Delta)\,,
$$
(2.56) \,  provides  the following information.
When $\,R \ =\ R_b\  := \ t \,+ \,\Delta\,: $
\begin{eqnarray*}
{{\partial \,[(2\,\varrho)^2]}\over {\partial \lambda}} \, \bigg\vert_{\lambda \,  = \,  \lambda_M}\!\!
 &  = & {{ (2 \lambda_M)\cdot  (4 R_b)^2}\over { (R_b^2 + \lambda_M^2)^4}} \times \left[ \,(R_b^2 + \lambda_M^2)^2 + 2 \lambda_M^2 (R_b^2 + \lambda_M^2) - 4\lambda_M^2 (R_b^2 + \lambda_M^2) \right] \ \ \ \ \ \ \ \ \ \
\\[0.075in] & = &
    {{(2 \lambda_M) \cdot (4R_b)^2    }\over {(R^2_b + \lambda^2_M)^3}}\, (R_b^2 - \lambda_M^2) \ =\   {{8\, \Delta \, \lambda_M}\over {t^3}}\,.
    \end{eqnarray*}
When $R  \ =\ R_c \ := \ t \,-\, \Delta\,: $
$$
{{\partial \,[(2\,r)^2]}\over {\partial \lambda }} \, \bigg\vert_{\lambda = \lambda_M} = {{(2 \lambda_M) \cdot (4R_c)^2   }\over {(R^2_c + \lambda^2_M)^3}}\ (R_c^2 - \lambda_M^2) \,=\, -\, {{8\, \Delta\, \lambda_M}\over {t^3}} \,.  \leqno (2.59)
 $$
Cf. (2.40).\, Combining with the information on the radius of the caps obtained in (2.39),  together with (2.10), the changes at the boundaries of the two caps cancel each other. This provides another way to see (2.20) and Lemma 2.21.\bk
%
%
%
%
%
As the situation is rotationally symmetric, we continue to  assume that we arrange
$
\,\xi \ = \ (\delta\,, \ 0\,,\,\cdot \cdot \cdot, \ 0),\,
$
where $\,\delta \,\ge\, 0\,.\,$
With the domain $\,\Omega\,$ as in (2.1) being an annular domain, clearly,
$\,
G_{|_{\bf Z}} (\lambda, \ \xi)\, = \ G_{|_{\bf Z}} (\lambda, \  -\,\xi)\,.
$\,
It follows that
$$
G_{|_{\bf Z}} (\lambda, \ \xi)  \ = \ G_{|_{\bf Z}} (\lambda, \ 0)  + O \,(|\,\xi|^2) \ \ \ \ \ \ \ \ \ {\mbox{when}} \ \ \ |\,\xi| \ = \ |\,\delta| \ \ {\mbox{is \ \ small}}\,.
$$
In what follows, we provide a more formal argument.

  \vspace*{0.15in}

{\it Expansion for small \,$\delta\,.$}\ \ \
Rewrite \,(2.58)\, in terms of $\,R_\lambda\ \, (\,= \,R \, \lambda^{-1})\,\,$ and $\,\delta_\lambda \     (\,=\,\delta  \,\lambda^{-1})$\,:
\begin{eqnarray*}
(2.60) \ \ \ \ \ {{\partial \,(2\,r)^2}\over {\partial \lambda}}  & = & {1\over \lambda} \cdot
{{ 2 \,(4R_\lambda)^2\,  }\over
{\{ [(R_\lambda + \delta_\lambda )^2 \,+\,1\,]  \cdot [(R_\lambda- \delta_\lambda )^2
+ 1\,]\,\}^2}}\,\times \\
&  & \\
& \ & \!\!\!\!\!\!\!\!\!\!\!\!\!
\!\!\!\!\!\!\!\!\!\!\!\!\!\!\!\!\!\!\!\!\!\!\!\!\!\!\!\!\!\! \times\ \left[
(1+ R_\lambda^{\,2} - \delta_\lambda^2)^2
\ + \ 2\,  (1+ R_\lambda^{\,2} )\,- \,{{4\,
[\, (1 + R_\lambda^{\,2} - \delta_\lambda^2 \,)^2 \,+\, \delta_\lambda^2  \,]\times
[R_\lambda^{\,2} + \delta_\lambda^2 \,+\,1\,] }\over
{  [\,(R_\lambda+ \delta_\lambda)^2 \,+\,1\,]  \cdot [(R_\lambda - \delta_\lambda )^2
\,+\,1\,] }}\right]. \ \ \ \ \ \ \ \
\end{eqnarray*}
Using  the expansion in (2.43)  we derive   the following [\,cf. also (2.40)].
\begin{eqnarray*} (2.61) \ \ \ \ \ \
{{\partial \,(2\,r)^2}\over {\partial \lambda}}  & = & {2\over \lambda} \cdot {{ (4 R_\lambda)^2}\over { (1 + R^2_\lambda)^3}} \cdot [\,(R_\lambda^{\,2} - 1) + O \,(\delta^2_\lambda)\,]\\
& \ & \\
 & \ &  \!\!\!\!\!\!\!\!\!\!\!\!\!  \!\!\!\!\! \!\!\!\!\!\!\!\!\!\!\!\!\!  \!\!\!\!\!\!\!\!\!\!\!\!\!  \!\!\!\!\!\!\!\!\!\!\!\!\! \!\!\!\!\!    \left(\!   =     {{ (2 \lambda)\cdot(4 R)^2}\over { (\lambda^2 + R^2 )^3}} \cdot\! \left[\,(R^2 - \lambda^2) + (\lambda^2)\cdot O \,(\delta^2_\lambda)\,\right] \!\right)\!\!\!
\mfor   \delta_\lambda \ge 0 \ \ {\mbox{small}}\,,\, \  {\mbox{here}} \    R_\lambda \,\le\, \sqrt{ {5\over 2} \,} \,+\, c^2.\ \ \ \ \ \ \ \
\end{eqnarray*}
In the above formula, $\,c\,$ is a fixed and small positive number.
Comparing with (2.40) and using the discussion in \S\,2\,f\,,\, we recognize that
the leading order term in (2.61) is   the derivative in the unperturbed case. Combining with (2.10), we find that (\,for $\,\lambda\,$ close to $\,\lambda_M\,$)
$$
{{\partial G_{|_{\bf Z}}}\over {\partial \lambda}}\, (\lambda\,, \, \ \xi) \ = \ {{\partial G_{|_{\bf Z}}}\over {\partial \lambda}}\, (\lambda\,, \, \ 0) \,+\, {1\over {\lambda}} \,  \, O \,(\delta^2_\lambda) \mfor   \delta_\lambda \ge 0 \ \ {\mbox{small}}\,,\,  \ \  R_\lambda \,\le\, \sqrt{ {5\over 2} \,} + c^2\,. \leqno (2.62)
$$
Applying Lemma 2.26 to (2.62), we obtain the following result.


\newpage

{\bf Lemma 2.63.} \ \ {\it Under the notations    and the conditions in Lemma\,} 2.21\,,\, {\it let \,$\delta \,=\, |\,\xi|$\, and $\,\lambda \,=\, \lambda_M \,+ \,s\,,\,$ and assume\,} (2.48)\,,\,
{\it there exist  positive constants \,${\bar \varepsilon_{11}}\,,\,$ ${\bar \varepsilon_{12}}\,,\,$ ${\bar \varepsilon_{13}}\,,\,$ and ${\bar C}_8$\,,\,  such that if }

\vspace*{-0.3in}

$$    \lambda^{-1}_M \, |\,\delta| \,\,\le\,\, {\bar \varepsilon_{11}}\,, \ \ \ \ \  {\bar \varepsilon_{12}} \,\,\ge\,\, \lambda^{-1}_M \, |\,s|  \,\,\ge \,\,  {\bar \varepsilon_{13}} \,\,>\,\, 0\,, \leqno (2.64)
$$
{\it   then we  have}
$$
\bigg\vert\, {{\partial G_{|_{\bf Z}} }\over {\partial \lambda}} \, (\lambda_M + s\,, \, \xi) \bigg\vert \ \,\ge \ \, {\bar C}_8\cdot {\bar \varepsilon_{13}} \cdot  {1\over {\lambda_M}}\ .
$$
{\it Here the positive constant $\,{\bar C}_8\,$ is independent on $\,t\,$,\, $\,\Delta\,$ and $\,\lambda\,$ as long as\,} (2.48)\, {\it and\,} (2.64) {\it are fulfilled.}
\vspace*{0.2in}

%
%
%
%
%
%
%
{\bf \S\,2\,i.  } {\it Stability under perturbation.} \ \
Recall that we denote the critical point  by
$$ {\bf p_M} \ = \ (\lambda_M\,, \ {\vec{\,0}}\,)\,,\, \ \ \ \ \ \ \ {\mbox{where \ \  (as usual) \ \ }} \lambda_M \ = \ \sqrt{(t \,+\, \Delta)(t \,-\, \Delta)\,}\ .$$

\vspace*{0.15in}

{\bf Lemma 2.65.} \ \ {\it Under the notations    and the conditions in Lemma\,} 2.21\,,\, {\it assume also conditions} (2.47) {\it and\,} (2.48)\,.\,
{\it   We can find positive numbers  $\,\gamma \,\in\, (0, \ 1)$\, and\, $\,{\bar C}_9\,$ such that}
$$
\min_{\partial B_{\,{\bf p_M}} (\gamma\, \lambda_M)}   \Vert \btd G_{ |_{\bf Z}}(\lambda\,, \ \xi\,; \ \Omega) \Vert
\ \,\ge\ \, {\bar C}_9 \cdot \gamma\cdot {{1}\over {\lambda_M}}\ ,\, \ \ \ \ \ {\it{where}}
$$
 $$
 \partial B_{{\bf p_M}} (\gamma\, \lambda_M\,) \,= \,\{\, (\lambda\,, \ \xi\,) \,\in\, \R^+ \times \R^n \ \ \big\vert  \ \ \ |\,s|^2 + |\,\xi|^2 = (\gamma \cdot \lambda_M)^2\,, \ \ {\it{with}} \ \ s = \lambda - \lambda_M\,\}\,.
 $$
 {\it  In addition, the constants $\,\gamma\,$ and \,${\bar C}_9\,$ do not depend on \,$t\,$ and $\,\Delta\,$ as long as\,} (2.47)\, {\it and}\, (2.48)\, {\it  are fulfilled.}\\[0.2in]
%
%
%
{\bf Proof.} \ \    We take $\,\gamma\,$ to be small  so that
$$(\lambda, \ \xi) \ \in \  \partial B_{{\bf p_M}} (\gamma\, \lambda_M\,)  \ \ \ \  \Longrightarrow  \ \  \sqrt{{5\over 2}\,} \,\,+\, c^2 \ \ge \  {{t \,+\, \Delta}\over {\lambda}}\  > \ {{t -  \Delta}\over {\lambda}} \ \ge \  \sqrt{{2\over 5}\,} - c^2 \ > \ 0 \,.
$$
For a point $(\lambda, \ \xi)\ \in \ \,\partial B_{\,{\bf p_M}} (\gamma\, \lambda_M)\,$,\, with $\ s \,=\, \lambda  \,-\, \lambda_M$\, and $\,|\,\delta| \,=\, |\,\xi|$\,,\, we have
$$
 |\,s| \ > \ {{ \,\gamma}\over {\sqrt{2}}}\,\,\lambda_M \ \ \ \ {\mbox{and}} \ \ \ \ |\,\delta| \ \le \  {{ \,\gamma}\over {\sqrt{2}}}\,\,\lambda_M \ \ \ \
  \Longrightarrow
\ \ \bigg\vert \,{{\partial \,G_{ |_{\bf Z}}(\lambda\,, \ \xi\,; \ \Omega)}\over {\partial \lambda}}
  \,\bigg\vert \ \ \ge \ \,{{C\,\gamma}\over {\lambda_M}}\ ,
  \leqno {\mbox{either}}
  $$
where we use Lemma  2.63\,;
$$
{\mbox{or}}  \ \ \ \ \ \ \ \ \ \   |\,\delta| \ > \ {{ \,\gamma}\over {\sqrt{2}}}\,\,\lambda_M \ \ \ \ {\mbox{and}} \ \ \  |\,s|\ \le \ {{ \gamma}\over {\sqrt{2}}}\,\,\lambda_M   \ \ \ \ \  \  \Longrightarrow
\ \ \Vert \,\btd_\xi\, G_{ |_{\bf Z}}(\lambda\,, \ \xi\,; \ \Omega) \,\Vert
\ \ge \ {{C'\,\gamma}\over {\lambda_M}} \ \ \ \ \ \ \ \  \ \ \ \ \
$$
(via Lemma 2.49). Thus we have the desired estimate.\qedwh

\newpage

{\large \bf \S 3. \ \ Infinite number of solutions.}\\[0.15in]
We begin by understanding the effect caused by  other annular domain $\,B_o \,(R) \setminus \overline{B_o \,(\rho)}\,\,$  at a  point $\,(\lambda_c\,, \ {\vec{\,0}}\,)$ (\,can be thought of as a critical point for another annular domain).\, As before, $\,R \,>\, \rho\,$ are given positive numbers. \\[0.2in]
{\bf \S 3\,a.  }{\it  $C^o$--\,effect}\,:\, {\it the case $\,\lambda_c^{-1} R \,\,\gg\,\, 1\,$ and $\,\lambda_c^{-1} \,\rho \,\,\gg\,\, 1$\,.}\ \
Let
$\displaystyle{\,
{R\over {\lambda_c}} \,\,>\,\, {\rho\over {\lambda_c}}   \,\,\gg \,\, 1\,.}\,
$
In this\\[0.07in]
case  the annular domain, when pulled back to \,$S^n$\, via $\,\dot{\cal P}\,$,\,\,   shrinks around the {\it north\,} pole. To be more precise, we consider
the following computation.
%
%
 \begin{eqnarray*}
  {\mbox{Let}} \ \ \ \   \tan\, \theta_R  \ = \ {R\over {\lambda_c}} \
 & \Longleftrightarrow & \ \  \theta_R \ :=\  \arctan \left( {R\over {\lambda_c}} \right)\\
 &  \Longrightarrow & \ \  \sin \,
\left({\pi\over 2} \,-\,  \theta_R\right)
\ = \ \left( 1 \,+ \,{{R^2}\over {{\lambda_c}^2}}\right)^{\!\!- {1\over 2}}
 \ = \ O \left( {{\lambda_c}\over R} \right)\\
& \Longrightarrow & \sin \,
\left({\pi\over 2} \,-\,  \theta_R\right) \ = \  {{\lambda_c}\over R}
 \left( 1 \,+\, {{\lambda_c}\over R}\right)^{\!\!-{1\over 2}} \ = \  {{\lambda_c}\over R}
 \left[ \,1 \,+\, O \left( {{\lambda_c}\over R} \right)\, \right]\!,\\
 & \ & \\
(3.1) \ \   \varphi_R \ :=\  {\pi\over 2} -  \theta_R  & \Longrightarrow &
 \ \ \sin \, \varphi_R \ = \ \varphi_R \,[\,1 \,+\, O \,(\varphi_R)] \ = \ \varphi_R \,\left[ \,1 \,+\, O \left( {{\lambda_c}\over R} \right) \right]\\
 & \Longrightarrow & \ \ \varphi_R \,\,= \,\,{{\lambda_c}\over R}
 \left[ 1 \,+\, O \left( {{\lambda_c}\over R} \right) \right] \ \ \ \  \ \ \   \ \ \ \ \ \  \mfor \, \ \ \ {R\over {\lambda_c}}  \ \ \ \ {\mbox{large}}\,. \ \ \ \ \ \ \
 \end{eqnarray*}
 It follows that
 $$
 0\  < \ \sin \,2\,\theta_R \ = \ \sin \,(\pi - 2\theta_R) \ = \ \sin \,2
\left({\pi\over 2}\, -\,  \theta_R\right) \ \le\  2 \sin
\left({\pi\over 2} \,-\,  \theta_R\right) \, = \ O\! \left(  {{\lambda_c} \over R} \right)\,.\leqno (3.2)
$$
Likewise, we define \,$\theta_\rho$\, is a similar fashion, and obtain
$$
0 \ < \  \sin \,2\theta_\rho \ = \ \sin \,(\pi - 2\theta_\rho)\ = \  \sin \,2
\left({\pi\over 2} \,- \, \theta_\rho\right) \le 2 \sin
\left({\pi\over 2} \,-\,  \theta_\rho\right) \, =\,\, O \!\left(  {{\lambda_c} \over \rho} \right) \leqno (3.3)
$$
for $\,\lambda_c^{-1} \, \rho \,\,\gg \,\,1\,.\,$
It follows that
\begin{eqnarray*}
(3.4) \ \ \ \ \ & \ & \int_{B_o (R) \setminus B_o (\rho)} \left( {{\lambda_c}\over
{{\lambda_c}^2 + |\,y|^{\,2}}}\right)^n \, \\
& = &   {{ \Vert S^n \Vert}\over {2^n}} \,
\int^{2\arctan \,{R\over {\lambda_c}}}_{2\arctan
\,{\rho\over {\lambda_c}}}   \
 [\,\sin \,\theta'\,]^{\,n - 1}\,d\theta' \ \ \ \  \ \ \ \ \  \ \ \ \ \   (\,{\mbox{here}} \ \ |\,y| = \tan\, \theta\,, \ \ \theta' = 2 \, \theta\,) \ \ \ \ \ \ \ \ \\
 & \le & C  \left({{\lambda_c} \over \rho} \right)^{n-1} \cdot \left[ \, \arctan \left({R\over {\lambda_c}} \right)
 \ - \ \, \arctan \left({\rho\over {\lambda_c}} \right)\right] \ \ \ \ \ \ \ \ \    [\,{\mbox{via}} \ \ (3.2)\, \ \& \ \,(3.3)\,] \\
 & \le & C  \left({{\lambda_c} \over r} \right)^{n-1}\cdot
 \left\{ \,\left[ \,{\pi\over 2} \ - \ \arctan \left({\rho\over {\lambda_c}} \right) \,\right] \
 -\  \left[{\pi\over 2} \ - \ \arctan \left({R\over {\lambda_c}} \right) \,\right]\,\right\} \\
 & \le & C_1 \left({{\lambda_c} \over \rho} \right)^{n-1}\cdot
 \left[  \, \left({{\lambda_c} \over \rho} \right)\
 - \  \left( {{\lambda_c} \over R} \right)\  + \ O \left( {{{\lambda_c}^2}\over {\rho^2}}\right)\,\right] \ \ \ \ \ \ [\,{\mbox{using}} \ \ (3.2)\, \ \& \ \,(3.3)\,] \\
 & \le & C_2 \left({{\lambda_c} \over \rho} \right)^{n-1} \cdot  {{{\lambda_c}\,(R -\rho)}\over {R\,\rho}}  \ \ \ \   \ \ \ \ \ \ \ \ \   \ \ \ \ \  \ \ \    \ \ \ \ \   {\mbox{when \ \ }} \lambda_c^{-1} \, R \ > \ \lambda_c^{-1} \, \rho \ \gg \  1\,.
\end{eqnarray*}

 \vspace*{0.15in}


%
%
%
%
{\bf \S 3\,b.   }{\it  $C^o$--\,effect}\,:\,  {\it the case \,$\lambda_c^{-1} \,R \,\,\approx\,\, 0$\, and \,$\lambda_c^{-1} \,\rho \,\,\approx\,\, 0$\,.}\ \
Let
$$
\ \ \varepsilon^2 \,\,\ge\,\, {R\over {\lambda_c}}\ > \ {\rho\over {\lambda_c}}   \ > \ 0\,, \ \ \ \ \ \ \ \  {\mbox{where}} \ \ \, \varepsilon \ \ {\mbox{is \ \ a \ \ (small) \ \ number}}\,.\leqno (3.5)
$$
The situation is upside down (in $S^n$) of the consideration in \S3\,a\,.
\begin{eqnarray*}
(3.6) \!\!\!\!\!\!\!\!\!\!& \ & \int_{B_o (R) \setminus B_o (\rho)} \left( {{\lambda_c}\over
{{\lambda_c}^2 + |\,y|^{\,2}}}\right)^n \,
  \le   \  C\, {{ \Vert S^n \Vert}\over {2^n}} \ \cdot
\int^{2\arctan \,{R\over {\lambda_c}}}_{2\arctan
\,{\rho\over {\lambda_c}}}   \
 [\,\theta'\,]^{\,n - 1}\,d\theta'  \ \ \ \ \  (\,\theta' = 2 \, \theta)\\
 & \le & C' \left[ \left({R\over {\lambda_c}} \right)^n
 - \left({\rho\over {\lambda_c}} \right)^n \right] \
 \ \ \ \ \ \ \ \ \ \ \ \ \ \ ({\mbox{domain \ \ gathers \ \ around \ \ the \ \ south \ \ pole}}) \ \ \ \ \
 \\
 & \le & C'' \left({R\over {\lambda_c}} \right)^{n-1} \left[
 \left({R\over {\lambda_c}} \right)
 - \left({\rho\over {\lambda_c}} \right)  \right]={{C'' \cdot  R^{\,n - 1} \, (R- \rho)}\over {{\lambda_c}^{\!\!n}}}  \ \    \mfor \varepsilon^2 \,\,\ge\,\, {R\over {\lambda_c}}\ > \ {\rho\over {\lambda_c}}   \ > \ 0\,.
\end{eqnarray*}

 \vspace*{0.15in}

{\bf \S 3\,c.  }{\it  $C^1$--\,effect\,.}\ \
Based on (6.6) in Part I \cite{I}, we have
\begin{eqnarray*}
\ \ \bigg\vert\, {{\partial G_{|_{\bf Z} } }\over {\partial \xi_j}} \,(\lambda_c\,, \ \vec{\,0}) \bigg\vert  &\le & C_1\,\int_{B_o (R) \setminus B_o (\rho)}  \ \  {{{\lambda_c}^n\cdot |\,y|}\over
({{\lambda_c}^2 + |\,y|^{\,2})^{n + 1}}}\  \, \ = \  C_2 \int^R_\rho   {{{\lambda_c}^n\cdot r^n \,dr} \over
{({\lambda_c}^2 + r^2)^{n + 1} }}
\\
& = &  {{C_2}\over {\lambda_c}}    \cdot \int^{\arctan \,{R\over {\lambda_c}}}_{\arctan
\,{r\over {\lambda_c}}}   \
[\,\cos \,\theta]^n\,[\,\sin \,\theta]^n\,d\theta =   {{C_3}\over {{\lambda_c}}}   \cdot
\int^{2\arctan \,{R\over {\lambda_c}}}_{2\arctan
\,{r\over {\lambda_c}}}   \
 [\,\sin \,\theta'\,]^n\,d\theta' \ \ \ \ \ \ \ \ \ \  \ \ \ \ \
\end{eqnarray*}
for $j = 1\,, \ 2\,, \ \cdot \cdot \cdot\,, \ n\,.\,$ In the above, $\,r = \tan \theta\,$ and $\,\theta' = 2\, \theta\,.\,$
We can carry out similar argument found in \S\,3\,a and \S\,3\,b.
As in (3.4), we have
$$
 \int_{B_o (R) \setminus B_o (\rho)}   {{{\lambda_c}^n\cdot r}\over
({{\lambda_c}^{\!\!\!2} + |\,y|^{\,2})^{n + 1}}}\   \,  \le\  {{ \bar C}\over {{\lambda_c}}}
\cdot   \left({{\lambda_c} \over \rho} \right)^n
{{{\lambda_c}\,(R -\rho)}\over {R\,\rho}} \mfor \,\lambda_c^{-1} \,R \gg 1\,, \ \ \lambda_c^{-1} \,\rho \gg 1\,. \leqno (3.7)
$$
 On the other hand,
\begin{eqnarray*}
(3.8) \  \int_{B_o (R) \setminus B_o (\rho)}   {{{\lambda_c}^n\cdot r}\over
({{\lambda_c}^2 + |\,y|^{\,2})^{n + 1}}}\  \,  & \le & {{C''}\over {\lambda_c}}  \left({R\over {\lambda_c}} \right)^n \left[
 \left({R\over {\lambda_c}} \right)
 - \left({\rho\over {\lambda_c}} \right)  \right]\\[0.1in]
 & = & {{C''}\over {\lambda_c}} \cdot {{R^n \, (R- \rho)}\over {{\lambda_c}^{n+1}}}  \mfor   \lambda_c^{-1}\, R  \ \ {\mbox{and}} \ \ \lambda_c^{-1} \,\rho  \ \ {\mbox{small}}\,. \ \ \ \  \ \ \ \ \ \ \  \ \ \ \ \ \ \ \ \
\end{eqnarray*}
Likewise, using (6.1) in Part I \cite{I},  we provide similar estimates for
$ \displaystyle{\,
\bigg\vert \,{{\partial G_{|_{\bf Z} } }\over {\partial \lambda}} \,(\lambda_c\,, \ \vec{\,0}) \,\bigg\vert\,.\,}
$
Observe that the first term inside the bracket in (6.1) is similar to the integral in (3.4) and (3.6). Knowing that this becomes the dominating term,  we have the following results.
\begin{eqnarray*}
(3.9) \ \ \ \ \bigg\Vert \btd  \left[\int_{B_o (R) \setminus B_o (\rho)}     \left( {\lambda\over {\lambda^2 + |\,y - \xi|^2}} \right)^n\,\right]_{(\lambda_c\,,\ \vec{\,0})}
\bigg\Vert & \le &   {C\over {\lambda_c}}\, \left({{\lambda_c} \over \rho} \right)^{n-1} {{1}\over R}\cdot
{{ R -\rho }\over { \rho}} \cdot \lambda_c\\
& \ & \\
& \ &    \!\!\!  {\mbox{for}} \ \
 \lambda_c^{-1}\, R \gg 1 \ \ {\mbox{and}} \ \ \lambda_c^{-1} \,\rho \gg 1\,.  \\
 & \ & \\
(3.10) \ \ \bigg\Vert \btd  \left[\int_{B_o (R) \setminus B_o (\rho)}     \left( {\lambda\over {\lambda^2 + |\,y - \xi|^2}} \right)^n\,\right]_{(\lambda_c\,,\ \vec{\,0})}
\bigg\Vert & \le &  {C\over {\lambda_c}} \cdot \left( {R\over {\lambda_c}}\right)^{\,n - 1}\!\!\cdot (R -\rho ) \cdot {1\over {\lambda_c}}\\
& \ & \\
& \ & \!\!\!\!\!      \mfor   \lambda_c^{-1}\, R  \ \ {\mbox{and}} \ \ \lambda_c^{-1} \,\rho  \ \ {\mbox{small}}\,. \ \ \ \  \ \ \ \ \ \ \  \ \ \ \ \ \ \ \ \
\end{eqnarray*}

\vspace*{0.15in}


{\bf \S\,3\,d.} {\it  Effect from $\,\,\xi$.\,}\ \
Consider the expression\,:
 $$
 \left( {{\lambda_c}\over {\lambda_c^2 + |\,y - \xi|^{\,2}}}  \right)^n
=\ \left[ \, {\lambda_c\over {\lambda_c^2 + |\,y |^{\,2} - 2 \,y \cdot \xi +   |\, \xi|^{\,2}}}  \right]^n = \  \left( {\lambda_c\over {\lambda_c^2 + |\,y |^{\,2}}}  \right)^n \cdot
 \left[ {1\over{1 +   {{|\, \xi|^{\,2} - 2\,y \cdot \,\xi }\over
{ |\,y |^{\,2} + \lambda_c^2  }}     }}\right]^n\!\!.
$$
{\bf Case one}\,: there exist positive numbers $\,c_1\,$ and $\,C\,$ such that
\begin{eqnarray*}
& \ & \lambda_c^{-1} \,|\,\xi| \le c_1 \ \ \ \ {\mbox{and}} \ \ \ \ \lambda_c^{-1} \,\rho \ge  C \ \  \Longrightarrow  \ \ \Bigg\vert {{|\, \xi|^{\,2} - 2\,y \cdot \,\xi }\over
{ |\,y |^{\,2} + \lambda_c^2  }}  \Bigg\vert \ \le \ {1\over 2}\\
 \Longrightarrow  & \ &   \left({{\lambda_c}\over {\lambda_c^2 + |\,y - \xi|^{\,2}}}  \right)^n \le \ 2^n \!\cdot\! \left( {{\lambda_c}\over {\lambda_c^2 + |\,y |^{\,2}}}  \right)^n   \ \ \ \   \ \ \ \   \ \ \ \  \ \ \ \ \ \ \mfor \rho \le |\,y| \le R.
\end{eqnarray*}
{\bf Case two}\,: there exist positive numbers $\,c_1\,$ and \,$\varepsilon\,$ such that
\begin{eqnarray*}
 \lambda_c^{-1} \,|\,\xi| \le c_1 & {\mbox{and}} & \!\! \lambda_c^{-1} \,R \le \varepsilon\ \
 \Longrightarrow       \bigg\vert {{|\, \xi|^{\,2} - 2\,y \cdot \,\xi }\over
{ |\,y |^{\,2} + \lambda_c^2  }}  \bigg\vert  = \bigg\vert \, {{\lambda_c^{-2} |\,\xi|^2 - 2 \, (\lambda_c^{-1} y)\cdot (\lambda_c^{-1} \xi)}\over {\lambda_c^{-2}|\,y|^2 + 1}}  \bigg\vert \ \le \ {1\over 2}\\
 & \  &  \!\!\!\!\!\!\!\!\!\!\!\!\!\!\!\!\!\!\!\!\!\!\!\!\!\!\!\!\!\!\!\!\!\!\!\!\!\!\!\!\!\!\!
 \!\!\!\!\!\ \  \  \Longrightarrow   \ \ \left( {{\lambda_c}\over {\lambda_c^2 + |\,y - \xi|^{\,2}}}  \right)^n \  \le \ 2^n \!\cdot\! \left( {{\lambda_c}\over {\lambda_c^2 + |\,y |^{\,2}}}  \right)^n\ \
 \  \ \ \ \ \ \ \ \ \ \ \ \ \ \ \ \ \ \ \ \    \mfor \rho \le |\,y| \le R.
\end{eqnarray*}
Under the condition $\,\lambda_c^{-1} \,|\,\xi| \le c_1$\,,\, the argument leading to (3.9) and (3.10)  can be applied to the situation $\,\xi \,\not=\, 0\,$ after multiplying the estimates by
a fixed constant. Moreover, so far we assume that $\,|\,H|\le 1\,$.\, The results can be rescaled correspondingly to $\,|\,H|\,\le\, {\cal A}\,$ for any positive constant $\,{\cal A}\,$\,.\\[0.2in]
{\bf Lemma 3.11.} \ \ {\it Given any positive numbers \,$R\,$ and \,$\rho\,$ with \,$R > \rho$\,,\, in\,} (1.3), {\it assume $\,|\,H|\, \le\, {\cal A}\,$ in the annular domain $\,{B_o\,(R) \setminus \overline{B_o\,(\rho)\!}}\,,\,$ and that $\,H\,$ is equal to zero outside. Then we can find  positive constants $\,c\,$, $\,\varepsilon\,$\,, and $\,C\,$ such that for }

\newpage

\begin{eqnarray*}
(3.12)\ \ \,\,   & \ &  \lambda_c^{-1} \,|\,\xi| \,\le \,c\,, \ \ \ \
 \lambda_c^{-1}\, R \,\ge \,C \ \ \ \ {\it{and}} \ \ \ \ \lambda_c^{-1} \,\rho \,\ge \,C \ \ \ \ \ \ \ \ \ \ \ \ \ \ \ \ \ \ \ \ \ \ \ \ \ \ \ \ \ \ \ \  \ \ \ \ \ \ \ \ \ \ \ \ \ \ \ \  \\
   & \ & \\
\Longrightarrow   \!\bigg\Vert\, \btd &\!&\!\!\!\!\!\!\!\!\!\!\!  \left(\int_{B_o (R)\, \setminus B_o (\rho)} H (y)\,   \left( {\lambda\over {\lambda^2 \,+\, |\,y \,-\, \xi|^2}} \right)^n\,\right)\bigg\vert_{(\lambda_c\,,\ \xi)}
\bigg\Vert   \  \le        {{C_1}\over {\lambda_c}}\, \left({{\lambda_c} \over \rho} \right)^{n-1} {{{\cal A} }\over R}\cdot
{{ R \,-\,\rho }\over { \rho}}\cdot  \lambda_c\,. \ \ \ \ \ \ \ \ \ \ \\
& \ & \\
(3.13)\ \ \,\,   & \ &    \lambda_c^{-1} \,|\,\xi| \,\le\, c\,, \ \ \
\varepsilon \,\ge\, \lambda_c^{-1}\, R\, >\, 0 \ \ \ \ {\it{and}} \ \ \ \  \varepsilon \,\ge\, \lambda_c^{-1} \,\rho\, >\, 0\\
   & \ & \\
\!\!\!\! \Longrightarrow\!
\bigg\Vert \,\btd &\!&\!\!\!\!\!\!\!\!\!\!\!  \left(\int_{B_o (R) \, \setminus B_o (\rho)}    H (y) \left( {\lambda\over {\lambda^2 \,+\, |\,y \,-\, \xi|^2}} \right)^n\,\right)\!\!\bigg\vert_{(\lambda_c\,,\ \xi)}
\bigg\Vert    \le   {{C_2}\over {\lambda_c}} \cdot \!\left( {R\over {\lambda_c}}\right)^{\,n - 1}\!\!\!\cdot (R -\rho )\cdot {\cal A}\cdot {1\over {\lambda_c}}\,.
\end{eqnarray*}
{\it Moreover, the positive constants $\,C_1\,$ and $\,C_2\,$ do not depend on $\,\lambda_c\,$,\, $\,\xi\,$,\, $R\,$ and $\,\rho\,$ as long as the conditions in\,} (3.12) {\it and\,} (3.13) {\it are fulfilled\,} ({\it respectively\,}).

\vspace*{0.2in}


{\bf \S\,3\,e.} \  {\it  Effect from changing $H$.\,}\ \
Suppose that $H$ and $H_1$ are bounded   and $L^2$-integrable on $\R^n$,\, and

\vspace*{-0.35in}

$$
\int_{\R^n} |\,H (y) - H_1 (y)|^{\,2} \  \ \le \ \ \varepsilon^n \! \cdot
\left( \,\sup_{\R^n} \,|\,H|\,\right)^2. \leqno (3.14)
$$

\vspace*{-0.05in}

To distinguish  the reduced functionals, we denote by
$$
G_{|_{\bf Z}} (H_1) \,({\bf z}) \ =\ {\bar c}_{-1} \int_{\R^n} H_1 (y) \left[ {\lambda\over {\lambda^2 + |\,y - \xi|^{\,2}}} \right]^n \, \ \ \ \mfor {\bf z} \,=\, V_{\lambda,\,\, \xi} \,\in\, {\bf Z}\,. \leqno (3.15)$$
Likewise, we define $G_{|_{\bf Z}} (H)\,.\,$

\vspace*{0.1in}

{\it $C^o$--\,estimate.} \ \ Applying the H\"older  inequality, we have
\begin{eqnarray*}
(3.16) \   & \ & |\, G_{|_{\bf Z}} (H) - G_{|_{\bf Z}} (H_1) | \\
& = & \Bigg\vert\  {\bar c}_{-1} \int_{\R^n} H (y)\left(  {\lambda\over {\lambda^2 + |\,y - \xi|^{\,2}}}\right)^n\  - \ \, {\bar c}_{-1} \int_{\R^n} H_1 (y)\left(  {\lambda\over {\lambda^2 + |\,y - \xi|^{\,2}}}\right)^n\,\Bigg\vert \ \ \ \ \ \ \ \ \ \ \ \  \ \ \ \ \ \ \ \ \ \ \  \ \ \ \   \\
& \le & |\,{\bar c}_{-1}|\cdot  \!\left[\int_{\R^n} |\,H (y) - H_1 (y)|^{\,2} \,
\right]^{1\over 2} \cdot \left[  \int_{\R^n}  \left(  {\lambda\over {\lambda^2 + |\,y - \xi|^{\,2}}}\right)^{2n}
 \, \right]^{1\over 2}\\
& \le & {\bar C}_{o,\  n}\!\cdot \left( {\varepsilon\over \lambda}
\right)^{n\over 2}\cdot \left( \,\sup_{\R^n} \,|\,H|\,\right)\,.
\end{eqnarray*}
Precisely,
$ \displaystyle{\
{\bar C}_{o,\  n} \ = \ |\,{\bar c}_{-1}|\cdot \! \left(\Vert \,S^{n-1} \Vert \, \int_0^{\pi\over 2}   \,[ \sin \theta]^{\,n-1} \,\,  [\cos \theta]^{\,3\,(n - 1)} \ d \,\theta \right)^{\!\!{1\over{\, 2}}}\!.}
$


\newpage

{\it $C^1$--\,estimate.} \ \ Consider the first derivative in $\,\xi$\,.\, Refer to (6.3)  in Part I \cite{I}\,.\,  As in (3.16),  we have

\vspace*{-0.2in}

\begin{eqnarray*}
(3.17)\!\!\!\!\!\!\!& \ & \bigg\vert\, {{\partial G_{|_{\bf Z}} (H)}\over {\partial \xi_j}} \ - \ {{\partial G_{|_{\bf Z}} (H_1)}\over {\partial \xi_j}}  \bigg\vert\\
& \le  & |\,{\bar c}_{-1}|\cdot  \left[\int_{\R^n} |\,H (y) \,-\, H_1 (y)|^{\,2} \,
\right]^{1\over 2}\left[  \int_{\R^n}    {{\lambda^{2n}\, |\,\xi_j \,- \,y_j|^{\,2}}
\over {(\lambda^2 \,+\, |\,y \,-\, \xi|^{\,2})^{2n + 2}}} \,\,  \right]^{1\over 2} \ \ \ \ \ \ \ \ \ \ \ \ \  \ \ \ \ \ \  \\
& \le &   C\,  \varepsilon^{n\over 2}\cdot \left( \,\sup_{\R^n} \,|\,H|\,\right) \left[      \int_{\R^n}
 {{\lambda^{2n}\, |\,\xi  \,-\, y |^{\,2}}\over
 {(\lambda^2 \,+ \,|\,y \,-\, \xi|^{\,2})^{2n + 2}}} \,\, \right]^{1\over 2}\\
& =& C'\,  \varepsilon^{n\over 2} \cdot \left( \,\sup_{\R^n} \,|\,H|\,\right) \left[   \int_0^\infty {{\lambda^{2n} r^{n+1} \,dr}
\over {(\lambda^2 + r^2)^{2n + 2}}} \right]^{1\over 2}\\
&=& C'\,  \varepsilon^{n\over 2} \cdot \!\left( \,\sup_{\R^n} \,|\,H|\,\right)\left[
\int_0^{\pi\over 2} {{\lambda^{3n  + 2} [\,\tan
\theta\,]^{\,n+1} \,\sec^2 \theta \, d \theta}\over
{\lambda^{4(n + 1)} \,[\,\sec^2 \theta]^{2n + 2} }} \right]^{1\over 2}  \,=   \,{\bar C}_{1,\  n} \left( {\varepsilon\over \lambda}
\right)^{n\over 2} {1\over {\lambda}}\cdot  \sup_{\R^n} \,|\,H|\,. \ \ \ \ \ \
\end{eqnarray*}
Here $\,j \,=\,1,\, \cdot \cdot \cdot, \ n$\,,\, and $\ {\bar C}_{1,\  n}\,$ has a similar expression as $\,{\bar C}_{o,\  n}\,$ [\,in (3.16)\,]. From here, we know that by counting the power of $\,\lambda\,$ and $\,r\,$,\, the other expressions in the first derivative of $\,\lambda\,$ can be estimated accordingly, and it is also   of order $\displaystyle{ \,O \left( \left[\, {\varepsilon\over \lambda} \right]^{n\over 2} \!\cdot  {1\over {\lambda}}\,\right)}\,$.\, We summarize the discussion of this   section in the following.\\[0.15in]
%
%
%
{\bf Lemma 3.18.} \ \ {\it Suppose that $H$ and $H_1$ are bounded   and $L^2$ functions on}\, $\R^n$.\, {\it Assume that  } (3.14) {\it holds. Then}
$$
\bigg\Vert \,\btd \, [\,G_{|_{\bf Z}} (H) - G_{|_{\bf Z}} (H_1)\,] \,\bigg\vert_{\,(\lambda\,, \ \xi)} \,\bigg\Vert \ \ \le \ \  C  \left( {\varepsilon\over \lambda}
\right)^{n\over 2}\!\cdot \left( \,\sup_{\R^n} \,|\,H|\,\right) \cdot {{  1 }\over \lambda}  \leqno (3.19)
$$
{\it for all} $\ (\lambda\,,\ \xi) \,\in\, \R^+ \times \R^n.$

\vspace*{0.25in}

{\bf \S\,3\,f.}  {\it Spacing.}\ \
For a number $\,{\bf a} > 1$\,,\, let the annular domains be set up at
$$
B_o \,\left( {{1\,+\, \eta}\over {\bf a}}\right)\,\bigg\backslash\  \overline{B_o \,\left( {{1\,- \, \eta}\over {\bf a}}\right)}\ , \ \ \ \cdot \cdot \cdot\,, \ \
B_o \,\left( {{1\,+ \,\eta}\over {{\bf a}^k}}\right)\,\bigg\backslash \  \overline{B_o \,\left( {{1\,- \, \eta}\over {{\bf a}^k}}\right)}\ , \  \cdot \cdot \cdot\,. \leqno (3.20)
$$
Given positive numbers $\,\tau\,$ and $\,\eta$\,  (to be specified later), define
%
%
%
\begin{eqnarray*}
(3.21) \ \
H (y) & = & {1\over {{\bf a}^{\tau \, k}}} \ \ \  \mfor \ \ y \,\in \,B_o \,\left( {{1\,+\, \eta}\over {{\bf a}^k}}\right)\,\bigg\backslash \  \overline{B_o \,\left( {{1\,-\, \eta}\over {{\bf a}^k}}\right)}\,, \ \ \ \ k \,= \,1, \ 2,  \cdot \cdot \cdot\,; \ \ \ \ \ \ \ \  \ \ \ \ \\
& \ & \\
H (y) & = & 0 \ \ \ \ \ \ \ \mfor \ \ y \,\not\in \,  \bigcup_{k = 1}^\infty \ \left\{ B_o \,\left( {{1\,+ \,\eta}\over {{\bf a}^k}}\right)\,\setminus \, \overline{B_o \,\left( {{1\,-\, \eta}\over {{\bf a}^k}}\right)} \ \right\}\,.
\end{eqnarray*}
 We first
observe that if we   choose \,$\eta \in (0, \ 1)$ so that
$$
(1 \,+\, \eta) \ < \  (1 \,-\, \eta) \, {\bf a} \ \ \ \ \Longrightarrow \ \ \
  {{ 1 \,+\, \eta}\over {{\bf a}^{k+1}}}
\ < \   {{1 \,-\,\eta}\over {{\bf a}^k}} \ \ \ \ {\mbox{and}} \ \ \ \
  {{1\,+\,\eta}\over {{\bf a}^{k}}}
\ < \  {{1\,-\, \eta}\over {{\bf a}^{k-1}}}\,.\leqno (3.22)
$$
(That is, the next annular domain is surrounded by the present one.)

\newpage

\hspace*{0.5in}For each individual annular domain indexed by a positive integer $m$\,,\, the critical point appears in $\,(\lambda_{M_m}\,, \ \vec{\,0})$\,,\, where
 $$
\lambda_{M_m} \,\,:=\,\, \sqrt{{{1\,+\, \eta}\over {{\bf a}^m}} \cdot {{1\,-\, \eta}\over {{\bf a}^m}}\ } \ = \ {{ \sqrt{\,1 - \eta^2\,} }\over {{\bf a}^m}} \ \ \ \ \ \ \ \ \  \ \ \ \ \ \ \ \ \  (\,{\mbox{via \ \ Lemma 2.21}}\,)\,. \leqno (3.23)
$$
For the reader's convenience, we describe the effects of other annular domains first, before smoothing out the ``\,edges".

\vspace*{0.2in}

{\it Effects from outside annular domains.} \ \  Here we fix $m \ge 2$\,.\, We prepare to apply (3.12) in Lemma 3.11 to estimate the $C^1$-\,effect. For a point $(\lambda\,, \ \xi) \in \R^+ \times \R^n$ with
$$
|\, \lambda \,- \,\lambda_{M_m}|\,\le \,c \, \lambda_{M_m}\ \   \left[\,{\mbox{i.e.}} \ \  (1 \,-\, c) \,\lambda_{M_m} \ \le\  \lambda \,\le\,  (1 \,+\, c) \,\lambda_{M_m}\,\right]\,, \ \ \  {\mbox{and}} \ \ \ \  \lambda_{M_m}^{-1}\, |\,\xi| \,\le\, c \,,\,
$$
where $c \in (0\,, \ 1)$ is a  constant, we obtain
\begin{eqnarray*}
& \ & \lambda^{-1} \,\times \,({\mbox{inner  \ \ radius \ \ of \ \ the \ \ nearest \ \ outside \ \ annular \ \ domain}})\\
 & \ge  & {1\over {1 \,+\, c}} \cdot{{ 1\, - \,\eta}\over {{\bf a}^{\,m - 1} }} \cdot {{ 1}\over {\lambda_{M_m}}}  \ = \ {{\bf a}\over {1 \,- \,c}}   \cdot \sqrt{\,{{ 1 - \eta} \over {1 \,+\, \eta}}\,}\,.
\end{eqnarray*}
When $\,{\bf a}\,$ is large enough, $\,c\,$ is small enough, and $\,\eta \,<\, 1/2\,$,\, we can apply (3.12) in Lemma 3.11 to obtain

\vspace*{-0.3in}

\begin{eqnarray*}
(3.24)\!\!\!\!\!\!\!\!\!\!\!\!\!\!& \ & \ \ \ \ \ \  \sum_{k = 1}^{m - 1}\   \bigg\Vert \btd  \left[\,\int_{B_o \,\left( {{1\,+\, \eta}\over {{\bf a}^k}}\right)\,\big\backslash \  \overline{B_o \,\left( {{1\,-\, \eta}\over {{\bf a}^k}}\right)}} \  \, H (y)\cdot    \left( {\lambda\over {\lambda^2 + |\,y \,- \,\xi|^2}} \right)^n\,\right]\bigg\vert_{(\lambda\,,\ \xi)}\,
\bigg\Vert  \ \ \ \ \ \ \ \ \ \ \\[0.075in]
  &\le & {{  C  }\over {\lambda }}\  \sum_{k  = 1}^{m - 1}\, \left\{ \left[ {\lambda\over {  (1 \,- \,\eta)\, {\bf a}^{-k} }}\right]^{\,n - 1} \cdot   \left[ {\lambda\over {  (1 \,+\, \eta)\, {\bf a}^{-k} }}\right] \left( {{2 \eta}\over {1\, - \, \eta}} \right)\cdot \left[{1\over {{\bf a}^{\tau\, k} }} \right]  \right\}\\
  & \le & {{  C_{\eta\,,\,c} }\over {\lambda_{M_m}  }} \  \sum_{k  = 1}^{m - 1}   \,\left\{
\left({1\over {{\bf a}^m}} \right)^n \cdot (\,{\bf a}^k )^n \cdot \left[{1\over {{\bf a}^{\tau \,k} }} \right] \right\} \\
& = & {{ C_{\eta\,,\,c}}\over {\lambda_{M_m}  }} \cdot   {1\over {{\bf a}^{m\,n}}}\,\sum_{k = 1}^{m - 1} \,(\,{\bf a}^k \,)^{\,n   -\tau} \ \ = \ \   {{  C_{\eta\,,\, c} }\over {\lambda_{M_m}  }} \cdot {1\over {{\bf a}^{m\,n}}}\,\sum_{k = 1}^{m - 1} [\,{\bf a}^{n   -\tau} \,]^k\\[0.1in]
& = &  {{  C_{\eta\,,\,c}}\over {\lambda_{M_m}  }} \cdot  {1\over {{\bf a}^{m\,n}}}\cdot {{  [\,{\bf a}^{n  - \tau} \,]^m - {\bf a}^{n  -\tau} }\over
{ {\bf a}^{n - \tau}  - 1}}
  \le    {{  C_{\eta\,,\,c} }\over {\lambda_{M_m}  }}  \cdot{1\over {{\bf a}^{m\,n}}} \cdot [\,{\bf a}^{n  -\tau} \,]^{\,m - 1} \ \ \   (\,{\mbox{provided}} \ \ {\bf a}^{n -\tau} > 2) \ \ \ \ \ \ \ \ \ \ \   \\[0.1in]
& = & \left( {{  1 }\over {\lambda_{M_m}  }}    \cdot  {1\over {{\bf a}^{\tau m} }}  \right) \cdot  {{ C_{\eta\,,\,c}}\over {{\bf a}^{n  - \tau}}} \!\!\mfor \,0 \,<\, \eta \,<\, {1\over 2}\,,\,\   |\, \lambda \,-\, \lambda_{M_m}|\,\le \,c \, \lambda_{M_m}\,  \ \    {\mbox{\&}} \   \ \,\lambda_{M_m}^{-1}\, |\,\xi| \,\le \,c\,.\, \ \ \ \ \
\end{eqnarray*}

\vspace*{0.01in}

Roughly speaking, the first term in the last expression above is the strength of the gradient from the contribution of the \,$m$\,-\,th annular domain (cf. Lemma 2.65)\,.\, Once we take $\,\tau > n\,$\,,\, the last  term in the above, that is, $\, C_{\eta\,,\ c}\,/ \,{\bf a}^{n  - \tau}\,$,\, can be made small when we choose $\,{\bf a}\,$ large enough (maintaining $\,0 \,<\, \eta\, <\, 1/2\,$)\,.\\[0.1in]
{\it Effects from inner annular domains.} \ \  Similarly,
\begin{eqnarray*}
& \ & \lambda^{-1} \,\cdot \,(\,{\mbox{outer  \ \ radius \ \ of \ \ the \ \ nearest \ \ inside \ \ annular \ \ domain}})\\
 & \le  & {1\over {1 \,-\, c}} \cdot{{ 1 \,+\, \eta}\over {{\bf a}^{m + 1} }} \cdot {{ 1}\over {\lambda_{M_m}}}  \ \,=\, \  {1\over {1 \,-\, c}} \cdot {1\over {\bf {\bf a}}} \cdot \sqrt{\,{{ 1\, + \, \eta} \over {1\, -\, \eta}}\,}\,.
\end{eqnarray*}
When $\,{\bf a}$\, is large enough, $c$ is small enough, and $\eta < 1/2$,
we apply (3.13) in Lemma 3.11 to form the following estimate.
\begin{eqnarray*}
(3.25)  & \ & \sum_{k = m + 1}^{\infty}\   \bigg\Vert \btd  \left(\int_{B_o \,\left( {{1+ \eta}\over {{\bf a}^k}}\right)\,\big\backslash\  \overline{B_o \,\left( {{1- \eta}\over {{\bf a}^k}}\right)}} \  H (y)\cdot    \left( {\lambda\over {\lambda^2 + |\,y - \xi|^{\,2}}} \right)^n\,\right)\bigg\vert_{\,(\lambda\,,\ \xi)}\
\bigg\Vert  \ \ \ \ \ \ \ \ \ \ \ \ \ \ \ \ \ \ \ \ \\
& \ & \\
  &\ & \!\!\!\!\!\!\!\!\!\!\!\!\!\!\!\!\!\!\!\!\!\!\! \le \ {{  C'  }\over {\lambda }}\  \sum_{k  = m + 1}^{\infty} \left\{ \left[ { {  (1 + \eta)\, {\bf a}^{-k}  }\over \lambda} \right]^{\,n - 1} \cdot    \left( {{2 \,\eta\, {\bf a}^{-k}}\over {\lambda}} \right)\cdot \left[{1\over {{\bf a}^{\tau\, k} }} \right]  \right\} \ \  \le\ \
  {{  C_{\eta\,, \ c}' }\over {\lambda_{M_m}}} \cdot {\bf a}^{m\,n} \cdot \!\! \sum_{k  = m + 1}^{\infty}  \left[ {1\over {{\bf a}^{n  + \tau} }} \right]^k\\[0.075in]
  &\ & \!\!\!\!\!\!\!\!\!\!\!\!\!\!\!\!\!\!\!\!\!\!\! \le \
  {{ C_{\eta\,, \ c}'}\over {\lambda_{M_m}}} \cdot {\bf a}^{m\,n} \, \, \left[ {1\over {{\bf a}^{n  + \tau} }} \right]^{m + 1}\!
\cdot {1\over {1 - {1\over {{\bf a}^{n  +\tau}}} }}\ \  \le \ \    \left( {{ 1}\over {\lambda_{M_m}}}  \cdot  {1\over {{\bf a}^{\tau m} }} \right) \cdot {{ C_{\eta\,,\ c}'  }\over {{\bf a}^{n  +\tau + 1}}}\\
\end{eqnarray*}
for $\,0 \,< \,\eta \,<\, {1\over 2}$\,,\, $|\, \lambda \,-\, \lambda_{M_m}|\,\le\, c\, \lambda_{M_m}\,\,$ and $\,\, \lambda_{M_m}^{-1}\, |\,\xi| \,\le\, c\,.$
Likewise, the last term $\,C_{\eta\,,\, c}' \,/ \,{\bf a}^{\,n + \tau + 1}\,$ can be made small when we choose $\,{\bf a}\,$ large.

 \vspace*{0.2in}

{\bf \S\,3\,g.}  {\it Smoothing out the edges.}\ \
We thicken the annular domain by bringing in the adjustment factor  $\,\sigma> 0$\,:

\vspace*{-0.35in}

$$
B_o \,\left( {{1\,+\, (\eta \,+\,  \sigma) }\over {{\bf a}^k}}\right)\,\bigg\backslash \ \overline{B_o \,\left( {{1\,-\, (\eta \,+ \,\sigma)}\over {{\bf a}^k}}\right)\!}\ \,. \leqno (3.26)
$$
Here $1\, -\, (\eta \,+\, \sigma)\, >\, 0\,.$\,
Recall that $\,H\,$ is defined in (3.21)\,.\,
We make use of the extra  space available in  (3.26) to smooth out   the inner and outer `edges'. To do that, fix a $\,C^\infty$\,-\,function \,${\cal F}$ \,on \,$\R^n$ with the following properties.
\begin{eqnarray*}
(3.27)  \ \  \ \ \  {\cal F} \,(y) & = & 1 \ \ \mfor 1\, -\, \eta \ \le \ |\,y| \ \le\  1 \,+ \,\eta\,,\\
     {\cal F} (y) &= &0  \ \ \mfor   |\,y| \ \ge \ 1 \,+\, (\eta \, + \,\sigma) \ \ \ {\mbox{or}} \ \ |\,y|\ \le\  1 \,- \, (\eta\, + \,\sigma)\,,\\
  & \ &   \!\!\!\!\!\!\!\!\!\!\!\!\!\!\!\!\!\!\!\!\!\!\!\!\!\!\!\!
  \!\!\!\!\!\!\!\!\!\!\!\!\!\!\!\!\!\!\!\!\!\!\!\! 1\  \ge \ {\cal F} (y) \ \ge\    0   \mfor 1 \,-\,\left(  \eta \,+ \,  \sigma \right) \,\,< \, |\,y| \, < \,\,1 \,-\, \eta \ \  {\mbox{or}} \ \, \,1 \,+ \,\left( \eta \,+\,  \sigma \right)\,  \ge\,  |\,y| \, \ge\,  1\, + \, \eta\,. \ \ \ \ \ \  \ \ \ \ \ \  \ \ \ \ \ \
\end{eqnarray*}
For an integer $\,m \ge 1\,$,\, \,let\, $\varsigma  \,=\, {\bf a}^{-m}$\,.\,
Then consider the rescaling
\begin{eqnarray*}
(3.28) \ \ \  & \ & \ \ \ {\cal S} \,(y) \,\,=\,\,{1\over{{\bf a}^{\tau \,m} }}\cdot {\cal F} \,(\,\varsigma^{-1} \,y\, )\\
  \Longrightarrow & \ &\
\ \ \sup_{\R^n} \Vert \btd^{(\,h)} {\cal S} \Vert \ \,\le \  \, {1\over{{\bf a}^{\tau \,m} }}\cdot  {{C \,(\sigma)}\over {\varsigma^k}}\ \, = \ \, {{C \,(\sigma)}\over {{\bf a}^{m\,(\tau - h)}} } \\[0.1in]
  \Longrightarrow  & \ & \ \ \ \sup_{\R^n} \Vert \btd^{(\,h)} {\cal S} \Vert \  \to \ 0  \   \mfor h \,\le\, n - 1 \ \ \ \   {\mbox{as}} \ \  \ m\to \infty \ \    \ \ (n - 1 \,<\, \tau\, < \,n)\,. \ \ \ \
\end{eqnarray*}
Denote the smoothened function by $H^S$\,,\,  which has the following property.
\begin{eqnarray*}
(3.29) \ \ \ \ \ \ \ \ \ \
 & \ & \!\!\!\!\!\!\!\!\!\!\!\!\!\!\!\!\!\!\!\!\!H^S\, (y)\ =\ {1\over {{\bf a}^{\tau \, k}}} \ \,\mfor y \in B_o \,\left( {{1+ \eta}\over {{\bf a}^k}}\right)\,\bigg\backslash \  \overline{B_o \,\left( {{1- \eta}\over {{\bf a}^k}}\right)}\,, \ \ \ \ k = 1, \ 2,  \cdot \cdot \cdot\,; \ \ \ \ \ \ \ \ \\
& \ & \\
0 \   \le  \  H^S& \ &\!\!\!\!\!\!\!\!\!\!\!\!\!(y)   \le  {1\over {{\bf a}^{\tau \, k}}} \mfor y \in B_o \,\left( {{1+ (\eta + \sigma)}\over {{\bf a}^k}}\right)\,\bigg\backslash \  \overline{B_o \,\left( {{1- (\eta + \sigma)}\over {{\bf a}^k}}\right)}\,,\\
H^S \,(y) & = & 0   \mfor y \,\not\in  \,\bigcup_{k = 1}^\infty \ \left\{ B_o \,\left( {{1+ (\eta + \sigma)}\over {{\bf a}^k}}\right)\,\bigg\backslash \  \overline{B_o \,\left( {{1- (\eta + \sigma)}\over {{\bf a}^k}}\right)} \  \right\}.
\end{eqnarray*}
$$
 \Vert \btd^{(\,h)}  H^S \,(y) \Vert   \to   0  \ \ \ \ \ \mfor \ \  h \ \le \ n - 1 \ \ \ {\mbox{as}} \ \ |\,y| \to 0^+. \leqno (3.30)
$$

\vspace*{0.2in}

{\bf \S\,3\,h.} {\it Regularity of $\,H^S$ at $\,0$\,.}\ \
It follows from (3.30) that   $\,H^S$\, has $\,(n - 1)$\,-\,th order of flatness at $0$\,.\, As
$$
n > \tau > (n - 1) \ \ \Longrightarrow \ \ \tau - (n - 1) > 0\,, \leqno (3.31)
$$
we can gain a bit of H\"older regularity at the origin, arriving to
$H^S\in C^{\,n -1, \ \beta} \ (\R^n)\,.$\, Currently, this appears to be the best balance between strength [\,cf. (3.29)\,]  and spacing [\,cf. \S\,3\,f \,and \,(2.48)\,]\,.

\vspace*{0.2in}


{\bf \S\,3\,i.} \ \  {\it Existence of an infinite number of stable critical points and the proof of Main  Theorem\,} 1.7.\ \
Let us now fix   \,$\tau\, \in\, (n - 1, \ n)\,,\,$ and make it precise that
\begin{eqnarray*}
(3.32) \ \ \ \ \ \ \ \ \ \ \ \ \ \ \ \Omega_m &:=& B_o \,\left( {{1+ \eta}\over {{\bf a}^m}}\right)\,\bigg\backslash \ \overline{B_o \,\left( {{1- \eta}\over {{\bf a}^m}}\right)}\,,\\
  H_m\, (y) \!\!& := &\!\! {1\over {{\bf a}^{\tau \, m}}} \mfor y \in B_o \,\left( {{1+ \eta}\over {{\bf a}^m}}\right)\,\bigg\backslash \ \overline{B_o \,\left( {{1- \eta}\over {{\bf a}^m}}\right)}\,,\\
& \ & \\
H_m \,(y) & \!\!:= &\!\! 0 \ \ \ \ \mfor y \not\in   \,B_o \,\left( {{1+ \eta }\over {{\bf a}^m}}\right)\,\bigg\backslash \  \overline{B_o \,\left( {{1-  \eta  }\over {{\bf a}^m}}\right)}\ . \ \ \ \ \ \ \ \ \ \ \ \ \ \ \ \ \ \ \ \ \ \ \ \ \ \ \ \ \ \ \ \ \
\end{eqnarray*}
Here $m \in {\N} \,   \setminus \{ 0 \}\,.\,$ With definition (2.1), and
via   Lemma 2.21 and Lemma 2.23\,,\, $G_{|_{\bf Z}} (\bullet\,, \ \bullet\,; \ \Omega_m)$ has a non\,-\,degenerate critical point at

\vspace*{-0.3in}

$$
{\bf p_{M_m}} : = (\lambda_{M_m}\,, \ \vec{\,0})\,, \ \ \ \ \ \ \ \ \  {\mbox{where}} \ \ \ \ \  \lambda_{M_m} \ = \  {{\sqrt{(1 \,+\, \eta) (1 \,-\, \eta)\,}}\over {{\bf a}^m}}\ . \leqno (3.33)
$$
It follows that
$$
{{ (1 \,+\, \eta)\,{\bf a}^{\,-m} }\over {\lambda_{M_m}}} \ = \  \sqrt{{{1 \,+\, \eta}\over {1 \,-\, \eta}}\,}\,, \ \ \ \ \ \ {{ (1 \,- \,\eta)\,{\bf a}^{\,-m} }\over {\lambda_{M_m}}}  \ = \   \sqrt{{{1 \,-\, \eta}\over {1 \,+ \,\eta}}\,}\,, \ \ \  \ \ \ {{   \eta\cdot{\bf a}^{\,-m} }\over {{\bf a}^{-m}}}  \,\,= \,\, \eta\,.\leqno (3.34)
$$
$$
\ \ \ \ \ \ \  {{t + \Delta}\over {\lambda_M}}\,, \ \ \ \ \ \ \ \ {{t - \Delta}\over {\lambda_M}}\,, \ \ \ \ \ \ \ \ \   {{ \Delta}\over {t}} \ \ \ \ \ \ \ \ \ \   {\mbox{in \ \ Lemma \ \ 2.65\,.}} \leqno {\mbox{Cf.}}
$$
We now find and fix  \,$\eta \in (0, \ 1)$ so that
$$
1\, -\, A^2 \  > \ \eta \,\,>\,\, B^2 \,\,>\,\, 0   \ \ \ \
 {\mbox{and}}
 \ \ \ \ \
 {5\over 2}\  \ge\
  {{1 \,+ \,\eta}\over {1 \,- \,\eta}} \,\,> \,\,  {{1\,- \, \eta}\over {1 \,+\, \eta}} \,\,\ge \,\,{2\over 5}\,.\leqno (3.35)
$$
One can verify that with (3.34)\,--\,(3.36)\,,\, conditions (2.47) and (2.48)  are fulfilled.\, (We emphasize that the choice of $\,\eta\,$ is once and for all $\,m= 1, \ 2, \cdot \cdot \cdot\,.$\,)\bk
After a rescaling (\,$ H \,= \, 1 \,\,\to \,\,H \,= \,{\bf a}^{-\tau \, m}$)\,,\,  the conclusion in Lemma 2.65 implies that

\vspace*{-0.3in}

$$
\min_{\partial B_{\,{\bf p_{M_m}}} (\gamma\, \lambda_{M_m})}   \Vert \btd\, G_{ |_{\bf Z}} (\lambda\,, \ \xi\,; \ \Omega_m) \,(\lambda\,, \ \xi) \Vert
\,\,\ \ge\,\, \    {{{\bar C}_9\,\,\gamma} \over {\lambda_{M_m}}}  \cdot {1\over { {\bf a}^{\,\tau\, m} }}\ . \leqno (3.36)
$$
Here the positive constants $\,{\bar C}_9\,$ and $\,\gamma\,$ do not depend on $\,m\,.\,$ \bk
To keep the notation neat, we introduce
$  G_{ |_{\bf Z}} (H_m) $
and   $G_{ |_{\bf Z}} (H^S)$ as in (3.15)\,.\, {\it We claim that}\\[0.05in]
(3.37)  \centerline{``\,{\it $G_{|_{\bf Z}} (H^S)$ \,has a stable critical point inside  $\,B_{\,{\bf p_{M_m}}} (\gamma\, \lambda_{M_m})\,.$"}}\\[0.1in]

Let us begin with
\begin{eqnarray*}
(3.38)\!\!\!\!\!\!& \ & \ \ \bigg\Vert \btd \,\{ G_{ |_{\bf Z}} (H^S) \,- G_{ |_{\bf Z}} (H_m) \}  \bigg\vert_{\,(\lambda\,, \ \xi)} \,\bigg\Vert\\
\le &  \ & \!\!\!\!\!\!\!\!\!C_1  \left\{   \   \bigg\Vert \btd  \left(\int_{B_o \,\left( {{1+ (\eta + \sigma)}\over {{\bf a}^m}}\right)\,\setminus \, \overline{B_o \,\left( {{1- (\eta + \sigma)}\over {{\bf a}^m}}\right)}} \ \   |\, H^S  (y) - H \,(y)| \cdot    \left( {\lambda\over {\lambda^2 + |\,y - \xi|^2}} \right)^n\,\right)\bigg\vert_{(\lambda\,,\ \xi)} \,\bigg\Vert \right.\\[0.075in]
& \ &  \     + \,\sum_{k = 1}^{m-1}\   \bigg\Vert \btd  \left(\int_{B_o \,\left( {{1+ (\eta + \sigma)}\over {{\bf a}^k}}\right)\,\setminus \, \overline{B_o \,\left( {{1- (\eta + \sigma)}\over {{\bf a}^k}}\right)}} \ \,   H^S  (y)\cdot    \left( {\lambda\over {\lambda^2 + |\,y - \xi|^2}} \right)^n\,\right)\bigg\vert_{(\lambda\,,\ \xi)}
\bigg\Vert\\[0.075in]
& \ & \ \     \ \ \ \  +  \!\!\left.  \sum_{k = m + 1}^{\infty}    \bigg\Vert \btd  \left(\int_{B_o \,\left( {{1+ (\eta + \sigma)}\over {{\bf a}^k}}\right)\,\setminus \, \overline{B_o \,\left( {{1- (\eta + \sigma)}\over {{\bf a}^k}}\right)}} \ \,  H^S   (y)\cdot    \left( {\lambda\over {\lambda^2 + |\,y - \xi|^2}} \right)^n\,\right)\bigg\vert_{(\lambda\,,\ \xi)}\,
\bigg\Vert\,
  \right\}. \ \ \ \
\end{eqnarray*}
Starting with the first group in the right hand side of  (3.38), we proceed   with
\begin{eqnarray*}
& \ & \int_{B_o \,\left( {{1+ (\eta + \sigma)}\over {{\bf a}^m}}\right)\,\setminus \, \overline{B_o \,\left( {{1- (\eta + \sigma)}\over {{\bf a}^m}}\right)}} \ \ \ |\, H^S \,(y) \ - \ H_m \,(y)|^2\\[0.1in]
& \le & \left( {2\over {{\bf a}^{\tau m}}}  \right)^2 \cdot \int_{B_o \,\left( {{1+ (\eta + \sigma)}\over {{\bf a}^m}}\right)\,\setminus \, \overline{B_o \,\left( {{1- (\eta + \sigma)}\over {{\bf a}^k}}\right)}} \ \ \  dy  \ \ \ \ \  \ \ \ \ \ \  \ \ \ \ \  \  \ \ \ \  [\,{\mbox{using}} \ \ (3.29) \ \ \& \ \ (3.30)\,]\\[0.075in]
& \le &  4\, \Vert \,S^{\,n - 1} \Vert  \cdot \left[\, \int_{{{1+  \eta }\over {{\bf a}^m}}}^{{{1+  (\eta + \sigma)}\over {{\bf a}^m}}} r^{\,n - 1} \, dr \ \ + \ \int^{{{1-  \eta }\over {{\bf a}^m}}}_{{{1-  (\eta + \sigma)}\over {{\bf a}^m}}} \ r^{\,n - 1} \, dr\right] \cdot \left( {1\over {{\bf a}^{\,\tau\, m}}}  \right)^2\\
& \le & C_2 \left( {1\over {{\bf a}^m}} \right)^n \cdot \left\{ \,[\,1+  (\eta + \sigma)\,]^{\,n} - [\,1+    \eta  ]^{\,n} \ + \ [\,1- \eta ]^n - [\,1- (\eta + \sigma) \, ]^{\,n}\,\right\} \cdot \left( {1\over {{\bf a}^{\,\tau\, m}}}  \right)^2 \\
& \le & C_3 \left( {1\over {{\bf a}^m}} \right)^{\!n} \cdot \sigma \times \left( {1\over {{\bf a}^{\,\tau\, m}}}  \right)^{\!2}\ \ \ \ \ \ \ \ \ \ \ \ \ \ \ \ \ \ \ \ \ \ \ \ \ \ \ \ \ \ \ \ \ \ \ \ \ \ (\,{\mbox{for}} \ \  \sigma > 0 \ \ {\mbox{ small}} \,)\\[0.1in]
& \le & \left( {{C_4 \  \sigma^{1\over n} }\over {{\bf a}^m}} \right)^{\!n}\cdot \left(\, \sup_{\R^n} |\, H_m| \right)^2 \ \ \le \ \   \left(  C_5 \  \sigma^{1\over n} \cdot \lambda_{M_m} \right)^n\cdot \left( \,\sup_{\R^n} |\, H_m| \right)^{\!2} \ \ \ \ \ \ [\,\mbox{via} \ \ (3.33)]\,.
\end{eqnarray*}
 The positive constant $C_5$ does not depend on $m$\,.\, Cf.  condition (3.14), which is required in Lemma 3.18.\,
The argument leading to Lemma 3.18 shows that
\begin{eqnarray*}
(3.39)\!\!\!\!\!\!\!\!& \ &  \bigg\Vert \btd  \left(\int_{B_o \,\left( {{1+ (\eta + \sigma)}\over {{\bf a}^m}}\right)\,\setminus \, \overline{B_o \,\left( {{1- (\eta + \sigma)}\over {{\bf a}^m}}\right)}} \ \ \ |\, H^S \,(y) - H \,(y)| \cdot    \left( {\lambda\over {\lambda^2 + |\,y - \xi|^2}} \right)^n\,\right)\bigg\vert_{(\lambda\,,\ \xi)} \bigg\Vert \\
& \ & \\
& \le & C_6 \, \sqrt{\sigma\,} \cdot {1\over {\lambda_{M_m} }} \cdot \left\{ {1\over {{\bf a}^{\,\tau\, m}}}  \right\} \ \ \ \ \  \ \ \  \ \ \ \ \  \ \ \   \ \ \ \ \  \ \ \  \ \  [\,\lambda_{M_m} \ \ {\mbox{is \ \ expressed \ \ in \ \ }} (3.33)\,]
\ \ \ \ \ \end{eqnarray*}
\hspace*{0.9in}for $\,(\lambda\,, \ \xi) \,\in\, \overline{B_{\,{\bf p_{c_m}}} (\gamma\, \lambda_{M_m})}\!$ \ \ \  \ [\,\,$\Longrightarrow \ \ ( 1 \,-\, \gamma) \,\lambda_{M_m}\,\, \le\,\, \lambda \,\,\le\,\, (1 \,+\, \gamma) \,\lambda_{M_m}\,].$\bk
%
%
%
 %
%
%
The second and third groups in the right hand side of (3.38) can be estimated as in (3.24) and (3.25) by changing $\,(1 \,+\, \eta)\, \to \,[\,1 \,+\, (\eta \,+\, \sigma)\,]\,$ and $\,(1 \,-\,\eta) \,\to\, [\,1 \,-\, (\eta \,+\, \sigma)\,]\,$ [\,when $\sigma > 0\,$ is small,  the arguments in (3.24) and (3.25) are not affected\,]\,.\, This leads to \\[0.05in]
(3.40)

\vspace*{-0.25in}

$$
\bigg\Vert \btd \,\{ G_{ |_{\bf Z}} (H^S) \,- G_{ |_{\bf Z}} (H_m) \}  \bigg\vert_{\,(\lambda\,, \ \xi)} \,\bigg\Vert \ \ \le \  \ \left[\, {{1}\over {\lambda_{M_m} }} \cdot   {1\over {{\bf a}^{\,\tau\, m}}} \right]  \cdot \left( C_6 \,\sqrt{\sigma\,} \,  + \,{{C_7}\over{{\bf a}^{n - \tau}}} + {{C_8}\over {{\bf a}^{n + \tau + 1} }} \right).
$$

\vspace*{-0.05in}

The positive constants do not depend on $m$,\, ${\bf a}$,\, $\eta$\, and \,$\eta$\,,\, as long as the conditions in (3.35) are fulfilled, and $\eta$ is small enough. Recall that $n - 1 < \tau < n$\,.\, Hence we can  choose `${\bf a}$' to be large enough, and $\sigma$ to be small enough, so that (3.40) yields
 $$
 \bigg\Vert \btd \,\{ G_{ |_{\bf Z}} (H^S) \,- G_{ |_{\bf Z}} (H_m) \}  \bigg\vert_{\,(\lambda\,, \ \xi)} \,\bigg\Vert \
 \ \le \  \  {1\over 2} \cdot {{{\bar C}_9\,\gamma}\over {\lambda_{M_m} }}     \cdot   {1\over {{\bf a}^{\,\tau\, m}}}  \leqno (3.41)
  $$
for  $\,(\lambda\,, \ \xi) \in \overline{ B_{\,{\bf p_{c_m}}} (\gamma\, \lambda_{M_m})}\,.\ $
  (3.36) and (3.41) imply that the gradient estimate
 $$
 \bigg\Vert \btd \,\{ G_{ |_{\bf Z}} (H^S)    \bigg\vert_{\,(\lambda\,, \ \xi)} \,\bigg\Vert
 \ \ \ge \   \ {1\over 2} \cdot {{{\bar C}_9\,\gamma}\over {\lambda_{M_m} }}     \cdot   {1\over {{\bf a}^{\,\tau\, m}}}  \ \ \ \mfor \ \  (\lambda\,, \ \xi) \in \partial  B_{\,{\bf p_{c_m}}} (\gamma\, \lambda_{M_m})\,. \leqno (3.42)
  $$
Together with (2.25), and $\,C^o\,$ property of degree  [\,refer to \cite{Fonseca-Gangbo}\,;\, \,cf. also (A.3.5), (A.3.6), (A.3.7) and (A.3.8) in the e\,-Appendix\,], we conclude that
$$
{\mbox{Deg}} \, (\,\btd\,G_{ |_{\bf Z}} (H^S)\,, \ B_{\,{\bf p_{M_m}}} (\gamma\, \lambda_{M_m})\,, \ {\vec{\,0}}\,)  \,=\, {\mbox{Deg}} \, (\,\btd\,G_{ |_{\bf Z}} ( H_m)\,, \ B_{\,{\bf p_{M_m}}} (\gamma\, \lambda_{M_m})\,, \ {\vec{\,0}}\,) = -\,1\,.
$$
(For more information on the degree, see, for example,  \S\,A.3.3 in the e\,-Appendix\,.\,)
Hence $\,G_{|_{\bf Z}} (H^S)\,$ has a stable critical point  [\,denoted by $\,(\lambda_m\,, \ \xi_m)$\,]\, inside  $\,B_{\,{\bf p_{M_m}}} (\gamma\, \lambda_{M_m})$\,.\, \{\,See \cite{Fonseca-Gangbo}\,;\, cf. also (A.3.4) and \S\,A.3.4 in the e\,-Appendix.\}\bk
As the above argument works for all $\,m \,\in\, {\N}\,,\,$ we apply (3.42) in above, together with Theorem 4.26  and   Lemma 5.4  in Part I \cite{I}, to deduce the existence part in Main Theorem 1.7. Moreover, \,(1.10) in the Main Theorem is a consequence of Theorem 4.17 in Part I \cite{I} and the fact that the stable critical point    $\,(\lambda_m\,, \ \xi_m) \,\in\, B_{\,{\bf p_{M_m}}} (\gamma\, \lambda_{M_m}) \ \Longrightarrow \ \lambda_m \to 0^+\,$ \,and \,$\xi_m \,\to\, \vec{\,0}\,$ as $\,m \to \infty$\,.\, It also follows that $\,|\,\xi_m| \,\le \,C\, \lambda_m\,$ for $\,m \,\ge\, 1\,.$\\[0.15in]
{\bf \S\,3\,j\,.}   {\it Back to $S^n$.}\ \
When we  transfer the solutions $\,\{ v_m \}_{m = 1}^\infty$ back to $\,S^n\,$ via the stereographic projection as solutions $\,\{ u_m \}_{m = 1}^\infty\,$ of equation (1.1), we need to show that the north pole \,${\bf N}$\, is a removable singularity.
This follows from the fact that the solutions $\,\{ v_m \}_{m = 1}^\infty\,$  appeared in Main Theorem 1.7, where we obtain in \S\,3\,h, have bounded  $\,\Vert \ \,\Vert_\btd\,$ -\,norms (\,recall (2.8) in Part I \cite{I})\,.\,   In fact, using   the Sobolev inequality [\,(2.9) in Part I \cite{I}\,], the Kelvin transform as used in the proof of Theorem 3.16 in \cite{Leung-2nd}, together with a result of Brezis and Kato \cite{Brezis-Kato},  we show that \,${\bf N}$\, is a removable singularity.       See \S\,A.14 in the e\,-Appendix. From the form of $\dot{\cal P}$\,  [\,(2.1) in Part I \cite{I}\,]\,,\, it can be seen that the south pole ${\bf S}$ is the blow\,-\,up point for the sequence  $\,\{ u_m \}_{m = 1}^\infty$\,.

 \vspace*{0.2in}

{\large \bf e\,-Appendix}\, \ \ can be found in
$$
{\mbox{www.math.nus.edu.sg/}}\, \tilde{\ }\, {\mbox{matlmc/e\,-Appendix.pdf}}
$$

\end{document}


\centerline{\bf \LARGE {\bf Appendix to the Articles  }}

\medskip \smallskip

\begin{center} {\LARGE   {\bf Construction of Blow-up Sequences for   }}
\medskip \medskip \smallskip \smallskip \\ {\LARGE {\bf  the Prescribed Scalar Curvature }} \medskip \medskip \smallskip \smallskip\\ {\LARGE {\bf Equation  on $S^n$.\,   I \,and\, II. }}

\vspace{0.43in}

{\Large {Man Chun  {\LARGE L}}}{EUNG}\\[0.165in]

{\large {National University of Singapore}} \\[0.1in]
{\tt matlmc@nus.edu.sg}\\
\end{center}
\vspace{0.33in}


\vspace*{0.2in}

{\it In this appendix we follow the notations, conventions, equation numbers, section numbers, lemma, proposition and theorem numbers as  used in the main articles\,} (Part\, I \cite{I} and Part  II \cite{II}\,), {\it unless otherwise is specifically mentioned} ({\it \, for instances, those equation numbers starting with}\, A\,). {\it The references for the appendix are listed at the last three pages.}\\[0.4in]
{\bf \S\,A.\,1 \ \ A proof of the integration by parts formula (2.11) in Part I.}\\[0.1in]
The Hilbert space $\,{\cal D}^{1,\,2}$ coincides with
the completion of $C^\infty_o\,(\R^n)$ with respect to
the $L^{\,2}$\,-\,norm of the gradient \cite{Ambrosett-Scalar}. That is, for any $f\in {\cal D}^{1, \, 2}$, there exists a sequence of compactly supported smooth functions $\{f_{o, \, i} \}$ so that
$$
\Vert \,f - f_{o, \, i} \Vert_\btd^{\,2} \ \ = \ \int_{\R^n} | \btd  (f - f_{o, \, i})|^{\,2} \ \to \ 0 \ \ \ \ \ {\mbox{as}} \ \ i \to \infty\,.\leqno (A.1.1)
$$
Using the   Green identity (compact support $\Longrightarrow$ no boundary term), we find
$$
\int_{\R^n} \langle \btd  \,f_{o, \, i}\,,\, \btd  \,h\rangle  \ = \ - \int_{\R^n}  f_{o, \, i}\cdot \Delta_o \,h    \ \ \mfor \ i \,=\, 1, \ 2, \cdot \cdot \cdot\,.
$$
It follows that\\[0.05in]
\begin{eqnarray*}
&\ & \!\!\!\bigg\vert\int_{\R^n} \langle \btd  \,f\,,\, \btd  \,h\rangle   \ +\   \int_{\R^n} ( f\, \Delta_o \,h ) \, \bigg\vert\\
& \ & \\
& \le & \!\!\! \bigg\vert\int_{\R^n}  \left\{ \,\langle \btd  \,(f - f_{o, \, i})\,,\, \btd  \,h\rangle\,  +  \langle \btd  \,f_{o, \, i}\,,\, \btd  \,h\rangle\,    + \,   ( f_{o, \, i}\, \Delta \,h )\,  +  \, (  \, [\,f - f_{o, \, i} \, ]\, \Delta_o \,h ) \, \right\} \bigg\vert\\
& \ & \\
& \le &\!\!\!\int_{\R^n}   |\btd  \,(f - f_{o, \, i})|\cdot| \btd  \,h | \    +  \  \int_{\R^n}   |\,f - f_{o, \, i}|\cdot | \,\Delta_o \,h \,| \\
& \ & \\
& \le &\!\!\!\left[\, \int_{\R^n}   |\btd  \,(f  \, - \,  f_{o, \, i})|^{\,2} \right]^{1\over 2} \cdot\, \left[\,\int_{\R^n}    | \btd  \,h |^{\,2}  \right]^{1\over 2} \ \ \ \ \ \ \ \ \ \ \ \ \ \ \ \ \ \ \
  \left(\ \to  0\ \ \ \ {\mbox{as}} \ \ i  \, \to \,  \infty \right) \\[0.1in]
 & \ & \ \ \ \ \  \    +\, \left[\, \int_{\R^n}   |\,f  \, - \,  f_{o, \, i}|^{{2n}\over {n - 2}} \right]^{{n - 2}\over {2n}} \cdot \, \left[\, \int_{\R^n}   | \,\Delta_o \,h | ^{{2n}\over {n + 2}} \right]^{{n + 2}\over {2n}}\ \ \ \
  \left(\ \to  0\ \ \ \ {\mbox{as}} \ \ i  \, \to \,  \infty \right).
\end{eqnarray*}
Here we use (A.1.1),   H\"older's  inequality,
the  Sobolev inequality (2.9)\,,\, and (2.10), both in Part I\,. \qedwh
 \vspace*{0.01in}


{\bf    \S\,A.\,2.   \ Riesz representation theorem. \  }\\[0.1in]
 For every bounded linear functional\,  $F\,: {\cal H} \to \R$ on a Hilbert space ${\cal H}$ over the scalar field   $\R$, equipped with inner product $\,\langle \  \,, \ \rangle_{\cal H}\,,\,$   there is a unique element $v \in {\cal H}$ so that
$$
F\, (u) \ = \ \langle u\,, \  v \rangle_{\cal H} \ \ \ \ \ \ \ \ \ \  {\mbox{for \ \ all}}  \ \  \,u \,\in \,{\cal H}\,. \leqno (A.2.1)
$$
Moreover, we may take
 $$
v \ = \ {{F (w)}\over {\Vert\,w \Vert^{\,2}_{\cal H}}}\,\cdot w \ \ \ \  \ \  {\mbox{for \ \ all}}  \ \   w \,\perp\,{\bf N}\,:=\, \{ \,F (\bullet)\, = \,0\}\,,\leqno (A.2.2)
$$
 and $\ \Vert\, F \Vert^{\,2}_{\cal H} = \langle v\,, \  v\,\rangle_{\cal H}\,.\,$ See, for example, \cite{McOwen}\,.

\vspace*{0.2in}

{\bf Corollary A.2.3.} \ \ {\it Under the conditions and notations described in}\, (A.2.1) {\it and}\, (A.2.2), {\it suppose $\,w \,\in \,{\cal H}\,$ with $\,F (w) \,\not=\, 0\,$,\, then}
$$
\Vert \,F \Vert^{\,2} \ \ge\  {{|\,F (w)|^{\,2}}\over {\Vert\,w \Vert^{\,2}_{\cal H}}}\,.  \leqno (A.2.4)
$$
({\it Here $w$ may not be perpendicular to the null space of\,} $F$\,.\,)

\newpage

{\bf Proof.} \ \ We have
$$
w \ = \ w_\perp \,+\, w_o \ \ \ \ {\mbox{with}} \ \ w_\perp  \perp {\bf N} \ \ {\mbox{and}} \ \ w_o \in {\bf N} \ \ \Longrightarrow \ \ \Vert\,w \Vert^{\,2}_{\cal H} \ = \ \Vert\,w_\perp \Vert^{\,2}_{\cal H}  \, + \,  \Vert\,w_o \Vert^{\,2}_{\cal H}.
$$
Hence
$$
\Vert F \Vert^{\,2} \ = \ {{|\,F (w_\perp)|^{\,2}}\over {\Vert\,w_\perp \Vert^{\,2}_{\cal H}}} \ = \ {{|\,F (w)|^{\,2}}\over {\Vert\,w_\perp \Vert^{\,2}_{\cal H}}} \ \ge \  {{|\,F (w )|^{\,2}}\over {\Vert\,w  \Vert^{\,2}_{\cal H}}}\,.
$$
The result follows. \qed
Regarding the functional $I_o$\, in (2.13)\,,\, Part I,  consider the following. For ${\bf z} \in {\bf Z}\,,\,$
\begin{eqnarray*}
(A.2.5) \  \  \ \  v \in {\mbox{Ker}}\, I''_o\,({\bf z})\ \ &\Longrightarrow &\ \ (I_o'' ({\bf z})\,[v] \,f)  \ = \ 0 \ \ \ \ \mfor {\mbox{all}} \ \ f \in  {\cal D}^{1,\, 2}\\[0.1in] &\Longrightarrow &\ \ (I_o'' ({\bf z})\,[f] \,v)  \ = \ 0 \ \ \ \ \ \ \ \ [\,{\mbox{using  \ \ (2.21)\,, \ Part \ I\,}}]\\[0.1in]
&\Longrightarrow &\ \ v \in {\mbox{Ker}}\,  I_o'' ({\bf z})\,[f] \\[0.1in]
&\Longrightarrow &\ \ \langle h, \ v \rangle_\btd  \ = \ 0\\
 & \ & \ \ \ \ \   \ \ \ \ \ \ \ \ \ \ \ \ \   \ \  [\,{\mbox{for}} \ \ h\,,\, \ {\mbox{refer  \ \ to \ \ (2.23)\,, \ Part \ I}}\,] \ \ \ \ \ \ \ \ \ \ \\[0.1in]
&\Longrightarrow &\ \ h \,\perp \,{\mbox{Ker}}\, I''_o({\bf z})\,.
\end{eqnarray*}

Moreover,
\begin{eqnarray*}
(A.2.6) \  \ \ \ \  u \in {\mbox{Coker}}\, I''_o({\bf z}) \ \ &\Longleftrightarrow &\ \  \bigg\langle {{ (I_o'' ({\bf z})\,[f] \,h) }\over {\Vert h \Vert_\btd^{\,2}}}\,\cdot h\,, \ u \bigg\rangle_{\btd} = \ 0\\[0.1in]
\ \ &\Longleftrightarrow &\ \  (I_o'' ({\bf z})\,[f] \,u) \ = \ 0 \ \ \  \mfor {\mbox{all}} \ \ f \,\in\, {\cal D}^{1, \, 2}\\
\ \ &\Longleftrightarrow &\ \  (I_o'' ({\bf z})\,[u] \,f) \ = \ 0 \ \ \  \mfor {\mbox{all}} \ \ f \,\in\, {\cal D}^{1, \, 2} \ \ \ \\[0.1in]
& \ & \  \ \ \ \ \ \ \ \ \ \  \ \ \ \  \  \ \ \ \ \ \ \ \ \  \ \ \ \ \ \  \ \    [\,{\mbox{via  \ \ (2.21)\,, \ Part \ I\,}}]  \  \ \ \ \ \ \ \ \ \ \  \ \ \ \ \ \ \ \  \ \ \ \ \\[0.1in]
\ \ &\Longleftrightarrow &\ \  u \,\in \,{\mbox{Ker}}\, I''_o({\bf z})\,.
\end{eqnarray*}
This symmetry allows us to show that $I''_o\,(u)$, treated in the above way, is a Fredholm map with index zero once we know that the kernel has finite dimension. See Lemma 4.1 in \cite{Progress-Book},\, pp. 47\,--\,50\,.

 \newpage


{\bf \S\,A.\,3 \ \ Non-degenerate critical points and fundamental notions }\\[0.05in]
\hspace*{0.59in}{\bf   \,in
the degree theory.} \\[0.2in]
{\bf \S\,A.3.\,a.} \ \  {\it Jacobian and the degree of a regular value for a $C^1$\,-\,map\,.}\ \
See \cite{Fonseca-Gangbo}, pp. 6\,.\, Here we restrict ourselves to an open and bounded domain  $\Omega$  in $\R^n$ with smooth boundary $\partial \,\Omega\,.\,$  Let
\begin{eqnarray*}
\Phi\,: \overline{\,\Omega} & \to &  \R^n\\
\cap   & \ & \\
\R^n\!\!\!   & \ &
\end{eqnarray*}
be a $C^1$-map (up to the boundary $\partial \,\Omega)\,.\,$ Suppose $\vec{\,\,0}\,\in\, \Phi (\Omega) \setminus  \Phi (\partial \,\Omega)$ is  a regular  value (pre\,-\,images having non\,-\,zero  Jacobians), then
$$
{\mbox{Deg}} \, (\Phi\,, \ \Omega\,, \ {\vec{\,\,0}}) \ := \sum_{y \, \in \,\,\Phi^{-1} \,( \!{\vec{\,\,0}})} {\mbox{sign}}\,  J_\Phi (y) \,. \leqno (A.3.1)
$$
Here $J_\Phi$ is the Jacobian (determinant) of the map $\Phi$\,.\, \bk
Intuitively, the degree counts the `left/right' cutting of the map at $\!{\vec{\,\,0}}$\,.\, In case $\!{\vec{\,\,0}}$ is a not a regular value, the degree is defined by using a ``near-by" regular value. See \cite{Fonseca-Gangbo}.\\[0.2in]
{\bf \S\,A.3.\,b.} \ \  {\it Degree for a continuous function.} \ \ Let
$$
\Phi_o\,: \overline{\,\Omega} \ \to \ \R^n \leqno (A.3.2)
$$
be $C^o$ (up to the boundary) and
$$
{\vec{\,\,0}}\, \not\in \,\Phi_o\, (\partial \,\Omega)\,.\ \ \ \ {\mbox{Indeed\,, \ \ assume \ \ that }}  \ \ \delta \, :=\, {\mbox{dist}} \  ({\vec{\,\,0}}\,, \ \Phi_o \,(\partial \,\Omega)) \,>\, 0\,.\, \leqno (A.3.3)
$$
Define
$$
{\mbox{Deg}} \, (\Phi_o\,, \ \Omega\,, \ {\vec{\,\,0}}) \ = \ {\mbox{Deg}} \, (\Psi\,, \ \Omega\,, \ {\vec{\,\,0}})
$$
for any $\Psi \in C^1 (\overline{\Omega})$ such that
$$
|\,\Phi_o\, (y) - \Psi (y)| \ < \ \delta \mfor {\mbox{all}} \ \ y \,\in\, \overline\Omega \ = \ \Omega \,\cup \,\partial \,\Omega\,.
$$
See Definition 1.18 in \cite{Fonseca-Gangbo}\,,\, pp. 17.

\vspace*{0.2in}

{\bf \S\,A.3.\,c.} \ \  {\it $C^o$ property of degree and existence of solution.}\ \
Let $\Phi_o$ be as in (A.3.2) satisfying  (A.3.3). Then we have the following properties. (See Theorems 2.1 and 2.3 in \cite{Fonseca-Gangbo}\,,\, pp. 30.)
$$
{\mbox{Deg}} \, (\,\Phi_o\,, \ \Omega\,, \ {\vec{\,\,0}})\  \not= \ 0 \ \ \Longrightarrow \ \ \Phi_o\, (y_o)\ = \ {\vec{\,\,0}} \mfor {\mbox{a}} \ \ y_o \,\in\, \Omega\,.  \leqno (A.3.4)
$$
That is, the equation $\Phi_o\, (y)\  =\  {\vec{\,\,0}}$ \,has a solution $\,y_o \in \Omega\,.\,$
\begin{eqnarray*}
 (A.3.5) \ \ \ \ \ \ \ \ \  \ & \ & {\mbox{For \   \ any \  \ }}   \Psi_o \in C^o (\overline{\,\Omega}) \  \ \ \ {\mbox{with}} \ \  \ \Vert\, \Phi_o - \Psi_o \Vert_{C^o ( \Omega \, \cup \,\partial \, \Omega\,)} \ < \ \delta\ \ \ \ \ \ \ \ \ \ \ \ \ \ \ \ \  \ \ \ \ \ \ \ \  \\
   & \ & \\
   \Longrightarrow  &  \ & {\mbox{Deg}} \, (\,\Phi_o\,, \ \Omega\,, \ {\vec{\,\,0}})  \ =\ {\mbox{Deg}} \, (\,\Psi_o\,, \ \Omega\,, \ {\vec{\,\,0}})\,.
\end{eqnarray*}
 Here $\delta $ is given in (A.3.3)\,.\, \\[0.2in]
{\bf \S\,A.3.\,d.} \ \  {\it Hessian and stable critical point of a $C^{\,2}$-\,function.}\ \
Let ${\cal O}$ be a non-empty open set in $\R^n$ and
$$
f \,: {\cal O}\, \to \,\R
$$
be a $C^{\,2}$-function\,.\, Suppose $y_c \in {\cal O}$ a critical point for $f$\, (i.e., $\btd \, f (y_c) = \vec{\,0}\,)\,.\,$  $y_c$ is called a {\it stable critical point}\, if there is an open bounded domain $\Omega \subset\!\!\subset {\cal O}$ such that $y_c \in \Omega$ and
$$
{\mbox{Deg}} \, (\btd \, f\,, \ \Omega\,, \ \vec{\,0}\,) \,\,\not= \,\,0\,.
$$
(We assume that $\partial \,\Omega$ is smooth\,.)\,
Here we render
\begin{eqnarray*}
 \ \ \ \ \ \  \ \ \  \ \ \ \ \ \  \ \ \ \ \ \ \ \ \  \ \ \  \btd \, f \,: {\cal O}  \subset \R^n &\,\to\,& \R^n \\[0.075in]
y & \mapsto & \btd \, f (y)\,.\ \ \ \ \ \  \ \ \ \ \ \ \ \ \  \ \ \ \ \ \ \  \ \ \ \ \ \ \ \ \  \ \ \ \ \ \ \ \ \ \ \  \ \ \
\end{eqnarray*}
Refer to \cite{Progress-Book}, pp. 25.\,
The Jacobian of the map $\btd\, f$ is then the determinant of the Hessian matrix of the function $f$. See, for example, \cite{Pressley}, pp. 93. \bk
A critical point of a $C^{\,2}$-function $f$ is called non-degenerate if the Hessian is a {\it non-degenerate} bilinear form, or equivalently, the determinant of the Hessian is non-zero. As a consequence, the critical point is isolated. That is, there is a ball $B_{y_c} (r)$ such that there is no other critical point in $B_{y_c} (r)$ besides\, $y_c$\,.\, Take $\Phi = \btd \, f$ in (A.3.1), we obtain\\[0.05in]
(A.3.6)
$$
{\mbox{Deg}} \, (\btd \, f\,, \ B_{y_c} (r)\,, \ \vec{\,0}\,) \, = \,\,{\mbox{sign}} \  J_{\btd f} \,(y_c) \, =\, {{ {\mbox{det\ Hessian \ at}} \ \,y_c }\over { |\,{\mbox{det\ Hessian \   at \   }} y_c\, |}} \, \not= \,0 \ \ \Longrightarrow \ \, ``{\mbox{stable}}".
$$
Moreover, let
$$
\bar \delta \,\,:= \,\,\min_{y \in \partial B_{y_c} (r)} |\, \btd f (y)|\,.\leqno (A.3.7)
$$
For a $C^1$-function $h$ defined in a neighborhood of\, $\overline{B_{y_c} (r)}\,,\,$ via (A.3.5), we have
\begin{eqnarray*}
(A.3.8) \ \ \ \ \ \   \ \ \ \ \ \ & \ & ``\,|\btd f (y)  \,-\, \btd \,h (y)| \,\,< \,\,\bar \delta \ \  \ \  \mfor {\mbox{all}} \ \ y \,\in\, \overline{B_{y_c} (r)}\ \, "\\
& \ & \\
\Longrightarrow  \  & \ &  ``\, {\mbox{Deg}} \, (\btd \, f\,, \ B_{y_c} (r)\,, \ \vec{\,0}\,)\  = \ {\mbox{Deg}} \, (\btd \, g\,, \ B_{y_c} (r)\,, \ \vec{\,0}\,)\ "\,.\ \ \ \ \ \   \ \ \ \ \ \ \ \ \ \ \ \   \ \ \ \ \ \
\end{eqnarray*}

\vspace*{0.3in}

{\bf   \S\,A.\,4. \ \ Curvature, gradient, Hessian, cancelation \& matching\,.} \\[0.1in]
{\bf  \S\,A.4.a.} \ \ {\it Calculation of the  sectional curvature.}\ \
Refer to \S\,3\,b in Part I.\,  On $\,\R^+ \times \R^n\,$ with coordinates
\begin{eqnarray*}
& \ & (\lambda\,, \ \xi_1\,, \cdot \cdot \cdot\,, \ \xi_n)\,,\\
& \ & \!\!\!\!\!\!\!\!\!\!\!\!\!\!\!\!\!\!\!\!\!\!\!\!\!\!\!\!\!\!\!\!\!\!{\mbox{indexed zero}}   \ \,\uparrow \ \
\end{eqnarray*}
we denote
$$
\partial_o \  = \ {\partial \over {\partial \lambda}}\,, \ \ \partial_1 \  = \ {\partial \over {\partial \xi_1}}\,,\, \cdot \cdot \cdot\,, \ \ \partial_n \  = \ {\partial \over {\partial \xi_n}}\,.
$$
Consider a metric of the form
$$
g_{{\bf z}_{ij}} \  = \ {{ c_i \,\delta_{ij}}\over {\lambda^{\,2}}}\,, \ \ \ \ 0\, \le i\,, \ j \le \,n\,, \ \    {\mbox{with}} \ \  c_o \  = \ {\bar c}_o  > 0 \ \ {\mbox{and}} \ \ c_1 = c_2 = \cdot  \cdot \cdot = c_n = {\bar c}_1 > 0\,.
$$
In \cite{Progress-Book}\,, pp. 73\,-\,74, we find a list of standard formulas for Riemannian curvature tensor. Let us start with the
  Christoffel symbol, which  is given by
\begin{eqnarray*}
& \ & \Gamma^w_{ij} \ = \ {1\over 2} \left[\,D_i \,h_{wj} + D_j \,h_{wi} - D_w \,h_{ij}\,\right]\cdot {g_{\bf z}}^{ww}\\
& \Longrightarrow & \Gamma^o_{jj} \ = \ {{c_o^{-1}\, c}\over {\lambda}}\,, \ \ \ \ \Gamma^o_{oo} = \Gamma^w_{ow} = - {1\over \lambda}\,, \ \ \Gamma^w_{oo} = \Gamma^o_{io} = {\mbox{others}} = 0 \mfor w \not= 0, \ \ i \not= 0\,.
\end{eqnarray*}
Recall  that
$$
R\,(\partial_u, \ \partial_v)\,\partial_w = \btd_{\partial_u} \btd_{\partial_v} \partial_w - \btd_{\partial_v} \btd_{\partial_u} \partial_w - \btd_{[\partial_u\,,\ \partial_v]} \,\partial_w\,.
$$
In local coordinates $x^\mu$ the Riemann curvature tensor is given by
\begin{eqnarray*}
& \ & R^w_{\ jkl}  \ = \   \,dy^w (R(\partial_k, \ \partial_l)\,\partial_j)\,,\,\\[0.075in]
& \ & R^w_{\ jkl}  \ = \  D_k\, \Gamma^w_{lj} - D_l\, \Gamma^w_{kj} +  \Gamma^w_{km}  \,\Gamma^m_{lj} - \Gamma^w_{lm} \, \Gamma^m_{kj} \\[0.075in]
& \Longrightarrow &   R^w_{\ jij}  \ = \  D_i\, \Gamma^w_{jj} - D_j\, \Gamma^w_{ij} +  \Gamma^w_{im} \, \Gamma^m_{jj} - \Gamma^w_{jm}  \, \Gamma^m_{ij} = \Gamma^w_{io} \, \Gamma^o_{jj} - \Gamma^w_{jo}\,  \Gamma^o_{ij}\\[0.075in]
 & \ &\ \ \ \ \ \ \   \ = \  \left({{c_o^{-1}\, c}\over {\lambda}} \right) \left( - {1\over \lambda} \right)  \ \ \   \ \ \ \ \ \ \  \ \ \ \ \ \ \ \   (\,i \,\not= 0\,, \ \ j \,\not=\, 0\,, \ \ i \,\not= \,j\,, \ \ w \,=\, i)\,.
\end{eqnarray*}
The related covariant tensor is obtained via
\begin{eqnarray*}
& \ & R_{ijkl}  \ = \ {g_{\bf z}}_{iw}\,R^w_{\ jkl} \\
& \Longrightarrow &  R_{ijij}  \ = \  {g_{\bf z}}_{iw}\,R^w_{\ jij}  \ = \  {c\over {\lambda^{\,2}}} \, \left({{c_o^{-1}\, c}\over {\lambda}} \right) \left( - {1\over \lambda} \right)  \\
& \ & \ \ \ \ \ \ \  \ \ \ \ \ \ \   \ \ \ \ \ \ \   \ \ \ \ \ \ \  \ \ \ \ \ \ \   \ \ \ \ \ \ \   \ \ \ \ \ \  (\,i \,\not= \, 0\,, \ \ j\, \not= \,0\,, \ \ i \,\not= \,j\,)\,.
\end{eqnarray*}
The sectional curvature \cite{Lee} is given by
\begin{eqnarray*}
K (i, \, j) & = & {{ \langle\, R(\partial_i\,, \ \partial_j) \, \partial_j\,, \ \partial_i \rangle_{g_{\bf z}}   }\over {\langle\, \partial_i\,,\ \partial_i\rangle_{g_{\bf z}} \langle \, \partial_j\,,\ \partial_j\rangle_{g_{\bf z}} - \langle \, \partial_i\,,\ \partial_j\rangle^{\,2}_{g_{\bf z}} }} = {{{g_{\bf z}}_{iw}\, R^w_{\ jij} } \over {{g_{\bf z}}_{ii}\, {g_{\bf z}}_{jj} - {g_{\bf z}}_{ij}^{\,2} }}\\
\Longrightarrow \ \ K (i, \, j) & = & {{R_{ijij}}\over {{g_{\bf z}}_{ii} \,{g_{\bf z}}_{jj} - {g_{\bf z}}_{ij}^{\,2} }} =
-\,{1\over {\,c_o}}  \mfor  \,i \not= 0\,, \ \ j \not= 0\,, \ \ i \not= j\,.
\end{eqnarray*}
\begin{eqnarray*}
R_{ojoj} &= & {g_{\bf z}}_{oo}\, R^o_{\ joj}\,,\\
R^o_{\ joj} & = & {{\partial \,\Gamma^o_{jj}}\over {\partial \lambda}} - 0 \,+\, \Gamma^o_{oo}\,\Gamma^o_{jj}   - \Gamma^o_{jj}\,\Gamma^j_{oj}\\
 &= & - \left({{c_o^{-1}\, c}\over {\lambda^{\,2}}} \right) - \left({{c_o^{-1}\, c}\over {\lambda^{\,2}}} \right) + \left({{c_o^{-1}\, c}\over {\lambda^{\,2}}} \right) = - \left({{c_o^{-1}\, c}\over {\lambda^{\,2}}} \right)\\
 \Longrightarrow R_{ojoj} & = & -\, {{c_o}\over {\lambda^{\,2}}}   \left({{c_o^{-1}\, c}\over {\lambda^{\,2}}} \right)\\[0.075in]
K (0,\, j ) & = &  -\,{1\over {\,c_o}}  \ \ \ \ \mfor  \ \ j \,\not= \,0\,.
\end{eqnarray*}

\vspace*{0.2in}

{\bf  \S\,A.4.b.} \ \ {\it Gradient and Hessian.}\ \
It follows from the  standard expression for gradient \cite{Book} \cite{Lee} (pp. 28) that
\begin{eqnarray*}
 &\,&\btd_{g_{\bf z}} \,  {\mathcal{F}} \,\, =\,\,  g_1^{ij} \,\left[\, \partial_i \,  {\mathcal{F}}  \right] \, \partial_j \,\,= \,\,{\lambda^{\,2}}   \left( {1\over {c_o}} {{\partial \,  {\mathcal{F}} }\over {\partial \lambda }} \,\partial_o \ +  \    {1\over {c_1}} {{\partial \,  {\mathcal{F}} }\over {\partial \xi_1 }} \,\partial_1\ + \cdot \cdot \cdot  + \ \ {1\over {c_1}} {{\partial \,  {\mathcal{F}} }\over {\partial \xi_n }}\,\partial_n\right)\\[0.1in]
 \Longrightarrow  & \ & \btd_{g_{\bf z}} \,  {\mathcal{F}} (\bar p) \,\,=\,\, 0 \ \ \ \Longleftrightarrow \ \ {{\partial \,  {\mathcal{F}} }\over {\partial \lambda }} (\bar\lambda\,, \ \bar\xi) = {{\partial \,  {\mathcal{F}} }\over {\partial \xi_1 }} (\bar\lambda\,, \ \bar\xi) = \cdot \cdot \cdot  = 0\,.
\end{eqnarray*}
Here,
$$
{\mbox{the \ image \ of }} \ \ (\bar\lambda\,, \ \bar\xi)\ \ \ {\mbox{ is \ denoted \ by}} \ \ \bar p\ \ {\mbox{ via the parametrization.}} \leqno (A.4.1)
$$
The Hessian is given by \cite{Lee}
\begin{eqnarray*}
{\mbox{Hess}}_{\,{g_{\bf z}}}\, {\mathcal{F}}   & = & \btd^{\,2}_{g_1} \,{\mathcal{F}}   \ = \ (\btd_{\,{g_{\bf z}}}  \, {\mathcal{F}} )^i_{\, ;\,j} \ \partial_i \otimes dy^j\\[0.075in]
{\mbox{with}} \ \ \ \ \ \  \ \ \ \ \ \  \ \ \ \ \ \ \ \ \  (\btd_{\,{g_{\bf z}}} \, {\mathcal{F}} )^i_{\, ;\,j} & = & \partial_j \, (\btd_{\,{g_{\bf z}}} \, {\mathcal{F}} )^i \ + \ (\btd_{\,{g_{\bf z}}} \, {\mathcal{F}} )^k \,\Gamma_{jk}^i\,, \ \ \ \ \ \  \ \ \ \ \ \  \ \ \ \ \ \  \ \ \ \ \ \ \ \ \ \ \ \
\end{eqnarray*}
where $\Gamma_{jk}^i$ is the Christoffel symbol \cite{Lee}.
In particular,

\newpage

\begin{eqnarray*}
 \btd_{g_{\bf z}} \,  {\mathcal{F}} (\bar p) & = & 0  \\[0.075in]
 \Longrightarrow  \      (\btd_{\,{g_{\bf z}}} \, {\mathcal{F}} )^i_{\, ;\,j}    \,(\bar p)   &=&   \left[\, {{\partial^{\,2} {\mathcal{F}} \, (\bar\lambda\,, \ \bar\xi)}\over {\partial \xi_i\  \partial \xi_j}} \right] \, \partial_i \otimes dy^j \mfor i \not= 0 \ \ {\mbox{and}} \ \ j \not= 0\,, \\[0.075in] (\btd_{\,{g_{\bf z}}} \, {\mathcal{F}} )^i_{\, ;\,j}    \,(\bar p)   &=&    \left[\, {{\partial^{\,2} {\mathcal{F}} \, (\bar\lambda\,, \ \bar\xi)}\over {\partial \xi_j\  \partial \lambda}} \right] \, \partial_o \otimes dy^j \mfor i  \,=\, 0\,, \ \ j \,\ge\, 1\,,\\ (\btd_{\,{g_{\bf z}}} \, {\mathcal{F}} )^0_{\, ;\,0}    \,(\bar p)   &=&  \left[\, {{\partial^{\,2} {\mathcal{F}} \, (\bar\lambda\,, \ \bar\xi)}\over { \partial \lambda^{\,2}}} \right] \, \partial_o \otimes d\lambda \mfor i  \,=\, 0\,, \ \ j \,\ge\, 1\,, \\[0.1in]
 \Longrightarrow \  {\mbox{Hess}}_{\,{g_{\bf z}}} \   {\mathcal{F}} (\bar p)  &{\mbox{is}}& {\mbox{ non-degenerate}} \ \  \Longleftrightarrow \ \ {\mbox{Hess}}_{\,(\lambda\,, \ \xi)} \ {\mathcal{F}} \  (\bar\lambda\,, \ \bar\xi)  \  {\mbox{is \  non-degenerate\,.}} \ \ \ \ \ \ \ \ \
\end{eqnarray*}


\vspace*{0.2in}

{\bf  \S\,A.4.\,c.} \ \
Refer to \S\,3\,d,\, Part I.  First we find\\[0.1in]
$$
\Delta_y \,\left[ \,{{ \lambda^{\,2} - r^{\,2} }\over { (\lambda^{\,2} + r^{\,2})^{n \over 2} }} \right] \ \ \ \ \ \ \ \ \ \ \ \ \ \  (\,r^{\,2} = y_1^{\,2} + \cdot \cdot \cdot + y_n^{\,2}\,)\,.
$$
$$
{{\partial }\over {\partial r}} \left[ {{\lambda^{\,2} - r^{\,2}}\over {(\lambda^{\,2} + r^{\,2})^{n \over 2} }}\right] \ = \  {{  - 2r }\over {(\lambda^{\,2} + r^{\,2})^{n \over 2} }} - n {{r\,(\lambda^{\,2} - r^{\,2})}\over {(\lambda^{\,2} + r^{\,2})^{{n \over 2}+1} }}
$$

$$
 {{\partial^{\,2} }\over {\partial r^{\,2}}} \left[ {{\lambda^{\,2} - r^{\,2}}\over {(\lambda^{\,2} + r^{\,2})^{n \over 2} }}\right] \ = \ {{  - 2}\over {(\lambda^{\,2} + r^{\,2})^{n \over 2} }} -  {{ n\,\lambda^{\,2}  }\over {(\lambda^{\,2} + r^{\,2})^{{n \over 2}+1} }} + {{ 3n\, r^{\,2} }\over {(\lambda^{\,2} + r^{\,2})^{{n \over 2}+1} }}
 $$

$$ \ \ \ \ \ \ \ \ \ \ \ \  \ \ \ \ \ \ \ \ \ \ \ \  \ \ \ \ \ \ \ \ \ \ \ \  \ \ \ \ \ \ \ \ \ \ \ \ \ \   +\  {{   2n\, r^{\,2} }\over {(\lambda^{\,2} + r^{\,2})^{{n \over 2} + 1} }} \ + \  {{ n\,(n +2) \,r^{\,2}\,(\lambda^{\,2} - r^{\,2}) }\over {(\lambda^{\,2} + r^{\,2})^{{n \over 2}+2} }}\ .
$$

It follows that
\begin{eqnarray*}
& \ & \Delta_y \,\left[ {{\lambda^{\,2} - r^{\,2}}\over {(\lambda^{\,2} + r^{\,2})^{n \over 2} }}\right]\,\, =\,\, \left[ {{\partial^{\,2} }\over {\partial r^{\,2}}} + {{n - 1}\over r} {{\partial  }\over {\partial r }} \right]\,\left[\, {{\lambda^{\,2} - r^{\,2}}\over {(\lambda^{\,2} + r^{\,2})^{n \over 2} }}\,\right]\\[0.1in]
 & = & {1\over {(\lambda^{\,2} + r^{\,2})^{{n \over 2}+2} }} \left[ - 2n (\lambda^{\,2} + r^{\,2})^{\,2} + n \,(n + 4)\, r^{\,2}(\lambda^{\,2} + r^{\,2})\right.\\
  & \ & \ \ \ \ \ \ \ \  \ \ \ \ \ \ \ \ \  \  \ \ \ \ \ \ \ \  \ \ \ \ \ \ \ \ \  \  \left. - \,n^{\,2} \lambda^{\,2} (\lambda^{\,2} + r^{\,2}) +  n\,(n + 2) \,r^{\,2}\,(\lambda^{\,2} - r^{\,2}) \,\right]\\[0.1in]
 & = & {1\over {(\lambda^{\,2} + r^{\,2})^{{n \over 2}+2} }} \left[ \, n\,(n + 2) \lambda^{\,2} r^{\,2} - n \,(n + 2) \lambda^4 \,\right] \ \ \ \ \ \ \ \ \  (\,{\mbox{cancelation + matching}}\,)\\[0.1in]
 & = & {{-\,n\,(n + 2) \,\lambda^{\,2}\, (\lambda^{\,2} - r^{\,2})}\over {(\lambda^{\,2} + r^{\,2})^{{n \over 2}+2} }}\,\,.
\end{eqnarray*}
Hence
$$
\Delta_y \,\phi_o\,\, = \,\,-\,n\,(n + 2) \lambda^{\,2}\cdot{{ (\lambda^{\,2} - |\,y - \xi|^{\,2})}\over {(\lambda^{\,2} + |\,y - \xi|^{\,2})^{{n \over 2}+2} }}
$$
As for
\begin{eqnarray*}
\Delta_y \,\phi_1 & = & \Delta_y \,\left[ {{  (\xi_1 - y_1)}\over {(\lambda^{\,2} + |\, y - \xi|^{\,2})^{{n }\over 2} }} \right]\\
& = & 0 + 2\, \btd_y (\xi_1 - y_1) \cdot \btd_y \,\left[{1\over {(\lambda^{\,2} + |\, y - \xi|^{\,2})^{{n }\over 2} }}\right]\\
 & \ & \ \ \ \ +\, (\xi_1 - y_1) \cdot \Delta_y\, \left[{1\over {(\lambda^{\,2} + |\, y - \xi|^{\,2})^{{n }\over 2} }}\right]\,.\\[0.1in]
\Delta_y \, \left[{1\over {(\lambda^{\,2} + r^{\,2})^{{n }\over 2} }}\right] & = & {{2n\,r^{\,2} - n^{\,2} \lambda^{\,2}}\over {(\lambda^{\,2} + r^{\,2})^{{n \over 2}+2} }}\\[0.075in]
&  & \\
\Longrightarrow \ \ \Delta_y \,\phi_1 & = &  {{2n\,(y_1 - \xi_1)}\over {(\lambda^{\,2} + r^{\,2})^{{n \over 2}+1} }} + (\xi_1 - y_1) \,{{2n |\, y - \xi|^{\,2} - n^{\,2} \lambda^{\,2}}\over {(\lambda^{\,2} + |\, y - \xi|^{\,2})^{{n \over 2}+2} }}\\[0.075in]
& = &  \,{{ - n^{\,2} \lambda^{\,2}(\,\xi_1 - y_1)}\over {(\lambda^{\,2} + |\, y - \xi|^{\,2})^{{n \over 2}+2} }}\ .
\end{eqnarray*}
Likewise, we can find other terms in (3.8).

\vspace*{0.2in}

{\bf  \S\,A.4.\,d.} \ \ Proof of Corollary 3.13, Part I.  Recall that \
$$ \displaystyle{
\Vert w \Vert_{\,\btd}^2\ = \ \int_{\R^n}|\btd  w|^{\,2} \,\,\,, \ \ \mbox{and} \ \ \ \ \ \Delta\,V_{\lambda\,, \, \xi}+ n (n - 2)\, V_{\lambda\,, \,\xi}^{{n + 2}\over {n - 2}}\ = \ 0} \ \ \ \ {\mbox{in}} \ \ \R^n\,.$$    Via the integration by parts formula (2.11) in Part I, we obtain
\begin{eqnarray*}
  \ \ \   \int_{\R^n} V_{\lambda\,, \ \xi}^{{n + 2}\over {n - 2}} \, w &   =  &  -\,{1\over {n \,(n - 2)}} \, \int_{\R^n} (\Delta V_{\lambda\,, \  \xi}) \, w\,\ = \ {1\over {n \,(n - 2)}} \, \int_{\R^n} \langle \btd \,V_{\lambda\,, \  \xi}\,,\, \btd \,w\, \rangle\,\ \ \\
&\ = \ & {1\over {n \,(n - 2)}} \  \langle V_{\lambda\,, \  \xi}\,,\, w \rangle_{\,\btd} \ = \   {1\over {n \,(n - 2)}} \  \langle \,{\bf z}\,,\, w \rangle_{\,\btd}\,.
\end{eqnarray*}
Using   Proposition 3.11, Part I,  and the fact that
$
 \,w \perp ( \,{\mbox{span}} \, \{{\bf z} \} \oplus T_{\bf z}\, {\bf Z}\,) \    \Longrightarrow \   \langle \,{\bf z}\,,\, w \rangle_{\,\btd}\ = \ 0\,,
$
we obtain the inequality (3.14).  Cf. Remark 4.2 in \cite{Progress-Book}, pp. 50, and the proof of Lemma 2.10 (\,loc. cit. pp. 21).


\newpage

{\bf   \S\,A.\,5. \ \    Uniform approximation.}\\[0.1in]
{\bf  \S\,A.5.\,a.} \ \ Proof of Proposition 4.6,  Part I. \ \ For $f \in (T_{\bf z} \,{\bf Z})^{\,\perp}$\,,\,  write
$$f =a\, V_{\lambda\,,\  \xi} + {\bar f}\ = \  a\, {\bf z} \,+\, {\bar f}\,,\,$$
where $\bar f \perp ( \,{\mbox{span}} \, \{ {\bf z} \} \oplus T_{\bf z}\, {\bf Z}\,)$.\, Here $a \in \R\,.\,$ As
$$
I''_o \,({\bf z})\, [f]\ = \ a\,I''_o \,({\bf z})\, [{\bf z}] \ + \ I''_o \,({\bf z})\, [{\bar f}\,]\,,
$$
 using (2.21), Part I, the Rieze Representation Theorem [\,cf. (2.23)\,,\,Part I\,], and calculation similar to that in (3.15), Part I,  we have
$$
 I''_o \,({\bf z})\, [\,a\,{\bf z}]\ = \ - {{4 }\over {n - 2}} \,\, (a\,{\bf z})\,.
$$
Thus (4.7), Part I, holds when $\bar f\ = \ 0$\,.\, Next, we assume that $\bar f \not= 0\,.\,$
Write
$$
{\bar f}\ = \ {\bar f}_p \,+\, {\bar f}_o\,, \ \ \ \ {\mbox{where}} \ \ {\bar f}_p \ \perp \ {\mbox{Ker}}\, I''_o \,({\bf z}) \,[{\bar f}\,] \ \ \ \ {\mbox{and}} \ \ f_o \,\in\, {\mbox{Ker}}\, I''_o \,({\bf z}) \,[\,{\bar f}\,]\,. \leqno (A.5.1)
$$
Observe that
\begin{eqnarray*}
{\mbox{Corollary \  3.13\,,\, Part I}} \ \ &\Longrightarrow & \ \, (I''_o \,({\bf z}) \, [{\bar f}\,]\,{\bar f}\,) \ \ge \  2 \,{\bar c}_2\, \Vert \bar f \Vert_{\,\btd}^2\  > \  0 \\[0.1in] & \Longrightarrow & \ \, \bar f \  \not\in \ {\mbox{Ker}}\, I''_o \,({\bf z}) \,[{\bar f}\,] \ \ \Longrightarrow \ \, {\bar f}_p \ \not= \ 0\,.
\end{eqnarray*}
Using (2.23), Part I, we obtain
$$
I''_o \,({\bf z}) \, [\,{\bar f}\,]\ = \ {{ (I_o'' ({\bf z})\,[{\bar f}\,] \,{\bar f}_p) }\over {\Vert {\bar f}_p \Vert_{\,\btd}^2}}\cdot   {\bar f}_p \ = \ {{ (I_o'' ({\bf z})\,[{\bar f}\,] \,{\bar f} ) }\over {\Vert {\bar f}_p \Vert_{\,\btd}^2}}\cdot {\bar f}_p \ \ \ \ \ \ \  \ \ \ \    [\,{\mbox{as}} \ \ (I_o'' ({\bf z})\,[{\bar f}\,] \,{\bar f}_o \, )\ = \ 0\,]\,.
$$
Via Corollary 3.13, Part I, we have
$$
\bigg\Vert \left(\, {{ (I_o'' ({\bf z})\,[{\bar f}\,] \,{\bar f} ) }\over {\Vert {\bar f}_p \Vert_{\,\btd}^2}}\,\right) {\bar f}_p \, \bigg\Vert_{\,\btd} \ \ge \ {{ 2 \,{\bar c}_2 \, \Vert \bar f\Vert_{\,\btd}^2}\over {\Vert {\bar f}_p \Vert_{\,\btd}^2}} \, \cdot \Vert {\bar f}_p \Vert_{\,\btd} \ \ge \ 2 \,{\bar c}_2\, \Vert {\bar f}  \Vert_{\,\btd}\ ,
$$
as $\,\Vert \,{\bar f} \Vert_{\,\btd} \ \ge \ \Vert \,{\bar f}_p \Vert_{\,\btd}\,.\,$
Together with (2.21), integration by parts formula (2.11) (both are found in  Part I) and the fact that $\bar f \,\perp \,  {\bf z}$\,,\, we infer that
$$
\ I''_o \,({\bf z}) [{\bar f}\,] \,({\bf z})\ = \ 0 \ \ \Longrightarrow \ \ {\bf z} \in {\mbox{Ker}}\, I''_o \,({\bf z}) \,[{\bar f}\,] \ \ \Longrightarrow \ \ {\bf z} \, \perp \, {\bar f}_p \ \ \ \ \ \ \ \  \ \  \ [\,{\mbox{via}} \ \ (A.\,5.1)\,]\,.
$$
It follows that
\begin{eqnarray*}
& \ & I''_o \,({\bf z})\, [f] \ = \ - \left(\, {{4a}\over {n - 2}} \right) {\bf z} \ + \   \left(\, {{ (I_o'' ({\bf z})\,[{\bar f}\,] \,{\bar f} ) }\over {\Vert {\bar f}_p \Vert_{\,\btd}^2}}\,\right) {\bar f}_p\\[0.05in]
& \Longrightarrow & \Vert I''_o \,({\bf z})\, [f]\,  \Vert_{\,\btd}^2 \,\ge\,  \left[ \left(\, {{4 }\over {n - 2}} \right)^2  \Vert\, a \,{\bf z} \Vert_{\,\btd}^2 + 4 \,{\bar c}^2_2\, \Vert\,  {\bar  f} \Vert_{\,\btd}^2 \right]   \ \ge\ {\bar c}_3^2 \,\Vert \,   f  \Vert_{\,\btd}^2\,,\\
 & \ & \ \ \ \ \  \ \ \ \ \  \ \ \ \ \  \ \ \ \ \ \  \ \ \ \ \  \ \ \ \ \  \   \ \ \  \ \ \ \ \  \ \ \ \ \  \   \ \ \ \ \  \left(\,\,{\mbox{as}} \ \  \Vert   \,  f  \Vert_{\,\btd}^2\ = \   \Vert \,  a\, {\bf z}  \Vert_{\,\btd}^2 + \Vert    {\bar f}\,  \Vert_{\,\btd}^2\,\right)\ .
\end{eqnarray*}
Here we can take
$\,\displaystyle{\
{\bar c}_3\ = \ \min \, \left\{{{4 }\over {n - 2}}\,, \ \  2 \,{\bar c}_2  \right\}\,.\
}$

\vspace*{0.3in}

{\bf  \S\,A.5.\,b.} \ \  Refer to  \S\,4\,g,\, Part I. See also \cite{Progress-Book} \ (pp. 53), and \cite{Primer}.
As in Lemma 2.21 in \cite{Progress-Book} (pp. 27\,)\,,\, we seek  to show that the remainder
$$
   R_{\bf z}\,(w) \  := \ I'_o\,({\bf z}+ w) \,-\, I''_o\,({\bf z})\, [w]
 $$
has the uniform property
 $$
 \Vert R_{\bf z}(w) \Vert \ = \ o\,(\Vert  w \Vert_{ \btd} ) \ \ \ \  \ \ {\mbox{when}} \ \ \ \Vert  \, w \Vert_{ \btd} \, \to \, 0 \ \ \ (\,{\mbox{uniform \ \ in \ \ }} {\bf z} \, \in\,  {\bf Z})\,. \leqno (A.5.2)
 $$
Recall that, as operators acting on $[\bullet]$, we have
\begin{eqnarray*}
I'_o\,({\bf z}+ w) & = & \int_{\R^n} \left[\,\langle \btd \,({\bf z}+ w)\,,  \btd\, [\bullet]\rangle_\btd  \ - \ n \, (n - 2) ({\bf z}+w)^{{n + 2}\over {n - 2}}\cdot [\bullet]\,\right]\,,\ \ \ \ \ \  \ \ \ \\[0.1in]
{\mbox{and}} \ \ \ \ \ \  \ \ \  I''_o\,({\bf z})[w] & = & \int_{\R^n} \left[\,\langle \btd \,w\,,  \btd\, [\bullet] \rangle_\btd  \ - \ n \, (n + 2)\, {\bf z}^{4\over {n - 2}} \,w \cdot [\bullet]\,\right]\,.
\end{eqnarray*}

It follows that\\[0.1in]
(A.5.3)
\begin{eqnarray*}
R_{\bf z}(w) & = & \int_{\R^n} \left[\,\langle \btd \,z\,,  \btd\, [\bullet] \,\rangle_\btd \, -\, n \, (n - 2) ({\bf z}+ w)_+^{{n + 2}\over {n - 2}} \cdot [\bullet] \, + \, n\,  (n + 2)\, {\bf z}^{4\over {n - 2}}\, w \cdot [\bullet]\,\right] \\[0.1in]
& = & \int_{\R^n} \left[\,-\,\Delta \,{\bf z} -n (n - 2) ({\bf z}+ w)_+^{{n + 2}\over {n - 2}}   +  n (n + 2)\, {\bf z}^{4\over {n - 2}} \,w \,\right]\cdot[\bullet] \\[0.1in]
& = & - n \, (n - 2) \int_{\R^n} \left[     ({\bf z}+ w)_+^{{n + 2}\over {n - 2}} \, -\, {\bf z}^{{n + 2}\over {n - 2}}  \,  - \, {{n + 2}\over {n - 2}} \, {\bf z}^{4\over {n - 2}}\, w \,\right]\cdot[\bullet] \\[0.1in]
& \ & \ \ \ \ \ \ \ \ \ \  \ \ \ \ \ \ \ \ \    (\, {\mbox{applying \ \ the \  \ equation \ \ }} \Delta\, {\bf z}\, + \, n \, (n - 2) \,{\bf z}^{{n + 2}\over {n - 2}} \, = \, 0\,)\,.
\end{eqnarray*}
One can proceed  to the remark made in the proof of (ii) in Lemma 4.7 in \cite{Progress-Book} \, (\,pp. 53\,, that is for the subcritical case). The proof for the critical case is similar.  We provide the following argument.\bk
Using the Taylor expansion we observe that
\begin{eqnarray*}
& \ & \!\!\!\!\!\!\!\!\!(1 + a)^{{n + 2}\over {n - 2}} \ = \ 1 \,+ \,{{n + 2}\over {n - 2}} \cdot a \,+ \,O (1) \,a^{\,2} \ \ \ \  \mfor \ \ |\,a|\,\le \,{1\over 2}\\[0.1in]
& \Longrightarrow & \ \ \bigg\vert \,(1 + a)^{{n + 2}\over {n - 2}} - 1 -  {{n + 2}\over {n - 2}} \cdot a  \bigg\vert \le C a^{\,2} \ \ \ \ \mfor \ \ |\,a|\,\le\, {1\over 2}\,,\\[0.1in]
& \  & \!\!\!\!\!\!\!\!\!\!\!\!\!\!\!\!(b + c)^{{n + 2}\over {n - 2}} = b^{{n + 2}\over {n - 2}} \left(\, 1 + {c\over b}\right)^{{n + 2}\over {n - 2}} =  b^{{n + 2}\over {n - 2}} \left[ 1 + {{n + 2}\over {n - 2}} \cdot {c\over b} + O (1) \cdot \left(\, {c\over b}\right)^{\,2} \right] \mfor \bigg\vert \, {c\over b} \bigg\vert \le {1\over 2} \\[0.1in]
& \Longrightarrow & \ \ \bigg\vert\,  (b + c)^{{n + 2}\over {n - 2}} \,-\, b^{{n + 2}\over {n - 2}} \,- \,{{n + 2}\over {n - 2}}\,\, b^{4\over {n - 2}} \,c \bigg\vert \ \le\  C \left(\, {c\over b}\right) b^{4\over {n - 2}} \cdot c \ \ \mfor \ \ \bigg\vert \, {c\over b} \bigg\vert \,\le\, {1\over 2}\,.
\end{eqnarray*}
Here $a$, $b$ and $c$ are numbers\,,\, whereas $C$ is a positive constant depending on the dimension $n\,.\,$
 Given a (small) number $\gamma$,\, as long as
\begin{eqnarray*}
 & \ & \!\!\!\!\!\bigg\vert \, {{w(y)}\over {{\bf z}(y)}} \bigg\vert \ < \  \gamma \ \ \ \ \ \  \ \ \ \ \ \   \ \ \ \ \ \ (\,{\mbox{recall \ \ that \ \ }} {\bf z} \ \ {\mbox{is \ \ a \ \ positive \ \ function}})\\[0.1in]
  \Longrightarrow & \ &  \!\!\!\!\!\bigg\vert \,    [\,{\bf z}(y)+ w(y)]_+^{{n + 2}\over {n - 2}} -[\,{\bf z}(y)]^{{n + 2}\over {n - 2}}   - {{n + 2}\over {n - 2}} \, [\,{\bf z}(y)]^{4\over {n - 2}}\, w (y) \,\bigg\vert \le C \cdot\gamma \,[\,{\bf z}(y)]^{4\over {n - 2}}\, w(y)\,.
\end{eqnarray*}
 Let
 $$
 \Omega_\gamma \ :=\  \left\{ \,y \,\in\, \R^n \ \ \big\vert \ \ \ \bigg\vert \, {{w(y)}\over {{\bf z}(y)}} \bigg\vert \,< \, \gamma \right\}\,.
 $$

 Hence
 \begin{eqnarray*}
& \ &  \int_{\Omega_\gamma}  \bigg\vert \,    ({\bf z}+ w)_+^{{n + 2}\over {n - 2}} -{\bf z}^{{n + 2}\over {n - 2}}   - {{n + 2}\over {n - 2}} \, {\bf z}^{4\over {n - 2}}\, w \,\bigg\vert \cdot |\, [\bullet]| \\[0.075in]
& \le & C \cdot \gamma \int_{\Omega_\gamma} {\bf z}^{4\over {n - 2}}\, |\, w| \cdot |\, [\bullet]|  \le   C \,\gamma \int_{\R^n} {\bf z}^{4\over {n - 2}}\, |\, w|\cdot |\, [\bullet]|  \\[0.075in]
 & \le & C \cdot \gamma\left(\, \int_{\R^n}  | \,{\bf z}^{4\over {n - 2}}\, w |^{{2n}\over {n + 2}}    \right)^{{n + 2}\over {2n}} \left(\, \int_{\R^n}  [\bullet]^{{2n}\over {n - 2}}   \right)^{{n - 2}\over {2n}}\\[0.075in]
& \le & C_1 \cdot \gamma\left(\, \int_{\R^n}  | \,z|^{ {{2n}\over {n - 2}} \cdot {4\over {n + 2}}} \, |\,w |^{{2n}\over {n + 2}}    \right)^{{n + 2}\over {2n}} \Vert\, [\bullet] \Vert_\btd \\[0.075in]
& \le & C_1 \cdot \gamma \left[ \left(\, \int_{\R^n}  | \,z|^{{2n}\over {n - 2}} \right)^{4\over {n + 2}} \, \left(\,\int_{\R^n}  |\,w |^{ {{2n}\over {n + 2}}  \cdot {{n + 2}\over {n - 2}}}   \right)^{{n - 2}\over {n + 2}} \right]^{{n + 2}\over {2n}} \Vert\, [\bullet] \Vert_{ {\cal D}^{1\,\, 2}}\\[0.075in]
& \le & C_2 \cdot \gamma \,\Vert \,w \Vert_\btd\cdot \Vert\, [\bullet] \Vert_\btd \ \ \ \ \ \ \ \ \ \ \ \ \ \ \ \ \ \  \ \ \ \ \ \ \   ({\mbox{applying \ the \ Sobolev \ inequality}})\,.\\[0.02in]
\end{eqnarray*}
On the other hand,
\begin{eqnarray*}
& \ &  \int_{\R^n \setminus \,\Omega_\gamma}  \bigg\vert \,    ({\bf z}+ w)_+^{{n + 2}\over {n - 2}} \,-\,{\bf z}^{{n + 2}\over {n - 2}}   \,- \,{{n + 2}\over {n - 2}} \, {\bf z}^{4\over {n - 2}}\, w \,\bigg\vert \cdot |\, [\bullet]| \\
& \le & \int_{\R^n \setminus \,\Omega_\gamma}  \left[ \,    ({\bf z}+ |\, w|)_+^{{n + 2}\over {n - 2}} + {\bf z}^{{n + 2}\over {n - 2}}  \,+ \,{{n + 2}\over {n - 2}} \, {\bf z}^{4\over {n - 2}}\, |\,w | \,\right] \cdot |\, [\bullet]|\\
& \le & \int_{\R^n \setminus \,\Omega_\gamma}  \left[ \,    (1 + \gamma^{-1} )^{{{n + 2}\over {n - 2}}} |\, w|^{{n + 2}\over {n - 2}}\, + \,\gamma^{-{{n + 2}\over {n - 2}}} |\, w|^{{n + 2}\over {n - 2}}  \,+ \,{{n + 2}\over {n - 2}} \, \gamma^{-{4\over {n - 2}}} |\, w|^{{n + 2}\over {n - 2}}  \,\right] \cdot |\, [\bullet]|\\
& \ & \ \ \ \ \ \ \ \ \ \ \ \ \ \  \ \ \  [\,|\,w (y)| \ge \gamma \cdot {\bf z} \,(y) \ \ \ {\mbox{for}} \ \ \ y \in \R^n \setminus \Omega_\gamma \ \ \Longleftrightarrow \ \ \gamma^{-1} |\,w (y)| \ge   {\bf z}(y)\,]\\
& \le & C_1  \,\gamma^{-{{n + 2}\over {n - 2}}} \int_{\R^n} |\, w|^{{n + 2}\over {n - 2}}   \cdot |\, [\bullet]|\\
& \le & C_1  \,\gamma^{-{{n + 2}\over {n - 2}}} \left(\, \int_{\R^n} |\, w|^{{2n}\over {n - 2}}  \right)^{{n + 2}\over {2n}}  \cdot \left(\, \int_{\R^n} |\, [\bullet]|^{{2n}\over {n - 2}}  \right)^{{n - 2}\over {2n}}\\
& \le & C_2 \left[\, \gamma^{-{{n + 2}\over {n - 2}}}  \cdot \Vert\,w \Vert_{ \btd}^{4\over {n - 2}} \right]\cdot \Vert\,w \Vert_{ \btd} \cdot \Vert [\bullet] \Vert_{ \btd}\,.
 \end{eqnarray*}
Thus as long as
$$
\gamma^{-{{n + 2}\over {n - 2}}} \cdot \Vert \,w \Vert_{ \btd}^{4\over {n - 2}} \ \le\  \gamma \ \ \Longleftrightarrow \ \ \Vert\,w \Vert_{ \btd}\ \le \ \gamma^{n\over 2}\ , \leqno (A.5.4)
$$
we have
$$
\int_{\R^n \setminus \,\Omega_\gamma}  \bigg\vert \,    ({\bf z}+ w)_+^{{n + 2}\over {n - 2}} -{\bf z}^{{n + 2}\over {n - 2}}   - {{n + 2}\over {n - 2}} \, {\bf z}^{4\over {n - 2}}\, w \,\bigg\vert \cdot \Vert\, [\bullet]\Vert_\btd \le C_3 \cdot \gamma \,\Vert \,w \Vert_{ \btd} \cdot \Vert \,[\bullet] \Vert_{ \btd} \leqno (A.5.5)
$$

\vspace*{-0.05in}

$$
\Longrightarrow \ \ \int_{\R^n }  \bigg\vert \,    ({\bf z}+ w)_+^{{n + 2}\over {n - 2}} -{\bf z}^{{n + 2}\over {n - 2}}   - {{n + 2}\over {n - 2}} \, {\bf z}^{4\over {n - 2}}\, w \,\bigg\vert \cdot \Vert\, [\bullet]\Vert_\btd \le C_4 \cdot \gamma \,\Vert\, w \Vert_{ \btd} \cdot \Vert \,[\bullet] \Vert_{ \btd}
$$

\vspace*{-0.05in}

$$
\Longrightarrow \ \ \Vert \,R_{\bf z}(w)  \Vert = \sup_{\Vert\, [\bullet] \Vert_\btd \not=0 } {{ |\,R_{\bf z}(w) [\bullet] \,|}\over { \Vert\, [\bullet] \Vert_{ \btd}   }}  = C_5 \cdot \gamma \,\Vert\, w \Vert_\btd\,.
$$

\vspace*{0.1in}

As $\gamma$ can be chosen to be small [\,cf.  (A.5.4)] when $\Vert\,w \Vert_{ \btd} \,\to \,0$, we obtain
$$
  \Vert \,R_{\bf z}(w)  \Vert \ = \ o \,(1) \cdot \Vert \,w \Vert_\btd  \ =\ o \,(\,\Vert \,w \Vert_{ \btd}) \ \ \ \ \ \ {\mbox{when}} \ \ \ \Vert\,w \Vert_{ \btd} \,\to\, 0\,.
$$

\vspace*{0.4in}

{\bf   \S\, A.\,6. \ \ Bounded projection.}\\[0.1in]
In Part I,  refer to (4.11) and condition {\bf (iii)} in Lemma 4.11\,.\, From (2.14),
$$
G \,(f) \ = \  {\bar c}_{-1}   \int_{\R^n} H f_+^{{2n}\over {n - 2}} \, \ \ \mfor \ f \in {\cal D}^{1\,,\, 2} \,.
$$
The Fr\'echet derivative of $G$ is given by
$$
G' ({\bf z} + w) ([\bullet]) \ = \ -\,{\bar c}_n    \int_{\R^n} H \cdot ({\bf z} + w)_+^{{n+2}\over {n - 2}} \cdot [\bullet] \ \ \ \ \mfor [\bullet] \in {\cal D}^{1\,,\, 2} \,. \leqno (A.6.1)
$$
It follows that
\begin{eqnarray*}
 (A.6.2) \ \ \ \ |\,G' ({\bf z} + w) ([\bullet])| &\le &C_1 \int_{\R^n} |\,{\bf z} + w|^{{n+2}\over {n - 2}} \,|\,[\bullet]|\ \ \ \ \ \ \ \ \ \ \ \ (\,|\,H| \ \ {\mbox{is \ \ bounded\,}}) \\
& \le & C_1 \left(\, \int_{\R^n}  | \, {\bf z} + w|^{{{n + 2}\over {n - 2}} \times {{2n}\over {n + 2}} }   \right)^{{n + 2}\over {2n}} \,\cdot\!\left(\, \int_{\R^n}  |\,[\bullet]|^{{2n}\over {n - 2}}  \right)^{{n - 2}\over {2n}}\\
& \le & C_2 \left(\, \int_{\R^n}  ( \, {\bf z} + |\, w| )^{ {{2n}\over {n-2}} }  \right)^{{n + 2}\over {2n}} \cdot \Vert \, [\bullet] \Vert_{ \btd}\,.
\end{eqnarray*}
Using the inequality
$$
( \, {\bf z} + |\, w|\,)^{{2n}\over {n - 2}} \le C(n) \left[\,{\bf z}^{{2n}\over {n - 2}}  + |\, w|^{{2n}\over {n - 2}}\, \right] \ \  \ \ \  \ \ ({\mbox{recall \ \ that \ \ }} {\bf z} \in {\bf Z} \ \ {\mbox{is \ \ positive}}),
$$
we have
\begin{eqnarray*}
|\,G' ({\bf z} + w)\, ([\bullet])| & \le & C_2 \left(\, \int_{\R^n}  | \, {\bf z}|^{ {{2n}\over {n-2}} }  +  | \, w|^{ {{2n}\over {n-2}} } \, \right)^{{n + 2}\over {2n}} \cdot \Vert \, [\bullet] \Vert_{ \btd}\\
& \le & C_2 \cdot [\, C(n) + C_3]^{{n + 2}\over {2n}} \cdot \Vert [\bullet] \Vert_{ \btd} \ \ \ \ \mfor \Vert \,w \Vert_{ \btd} \le 1\,.
\end{eqnarray*}
Here we use the Sobolev inequality.
Therefore
$$
\Vert \,G' ({\bf z} + w) \Vert \ \le \  C \mfor {\mbox{all}} \ \ {\bf z} \in {\bf Z} \ \ \ {\mbox{and}} \ \ w \ \ {\mbox{with}} \ \  \Vert\, w \Vert \le 1\,. \leqno (A.6.3)
$$

\vspace*{0.3in}


{\bf   \S\, A.7. \ \ More uniform bounds.}\\[0.1in]
Refer to \S\,4\,d, Part I.  Consider the following uniform bounds.
\begin{eqnarray*}
   (A.7.1) \ \ \ \ \ \   \ \ \ \ \ \     I_o\,({\bf z} + w_\varepsilon ({\bf z})) & = & I_o\,({\bf z}) + I_o'({\bf z})[w_\varepsilon ({\bf z})] + o\,(\Vert  w_\varepsilon ({\bf z}) \Vert) \ \ \ \ \ \  \ \ \ \ \ \ \ \ \ \ \ \  \ \ \ \ \ \   \\
        & = &     {\bar c}_4 +   o\,(\Vert   w_\varepsilon ({\bf z}) \Vert)\,; \\
        & \ & \\
       (A.7.2) \ \ \ \ \ \   \ \ \ \ \ \ \,\, G ({\bf z} + w_\varepsilon ({\bf z}))     &= &G ({\bf z})
        + [ G' ({\bf z}) \,| \,w_\varepsilon ({\bf z})] + o\,(\Vert\,w_\varepsilon ({\bf z}) \Vert)\,; \\
             & \ & \\
          (A.7.3) \ \ \ \ \ \   \ \ \ \ \ \ \, I_o' ({\bf z} + w_\varepsilon ({\bf z})) & = &    [\,I_o''({\bf z})\,| \,w_\varepsilon ({\bf z})]
           + o\,(\Vert  w_\varepsilon ({\bf z}) \Vert)\,;   \\
                & \ & \\
          (A.7.4) \ \ \ \ \ \   \ \ \ \ \ \   G' ({\bf z} + w_\varepsilon ({\bf z}))     &= &G' ({\bf z})
           + [ G'' ({\bf z}) \,| \,w_\varepsilon ({\bf z})] + o\,(\Vert\,w_\varepsilon
           ({\bf z}) \Vert)\,.\\
         \end{eqnarray*}

         \vspace*{-0.2in}

See \S\,A\,.5\, for \,(A.7.3).\, Here we demonstrate the arguments leading to (A.7.1) and (A.7.4)\,.
\begin{eqnarray*}
 {\bf I.} \ \ & \ & I_o\,({\bf z} + w)\\
  & = & \int_{\R^n} \left(\, \ {1\over 2} \,\langle\, \btd\, ({\bf z} + w)\,, \  \btd\, ({\bf z} + w)\, \rangle -{{n - 2}\over {2n}} \cdot n (n - 2) ({\bf z} + w)^{{2n}\over {n - 2}} \right) \\
 & = & {1\over 2} \int_{\R^n} \langle\, \btd \,{\bf z}\,, \  \btd \,{\bf z}\, \rangle  +  \int_{\R^n} \langle\, \btd\, w\,, \  \btd \,{\bf z}\, \rangle  + {1\over 2} \int_{\R^n} \langle\, \btd\, w\,, \  \btd\, w\, \rangle \\
 &  \ & \ \ \ \ \ \ \ \  \ \ \ \ -\, {{n - 2}\over {2n}} \cdot n (n - 2)  \int_{\R^n}({\bf z} + w)^{{2n}\over {n - 2}} \\
 & \ & \\
 & = & I_o\,({\bf z}) -   \, {{n - 2}\over {2n}} \cdot n (n - 2)  \int_{\R^n}\left[ \,({\bf z} + w)^{{2n}\over {n - 2}} -{\bf z}^{{2n}\over {n - 2}} - {{2n}\over {n - 2}} \cdot {\bf z}^{{n + 2}\over {n - 2}}\cdot w \right]  \\
 & \ & \ \ \ \ \ \ \ \ \ \ \ \ \ \ \ \ \ \ \ \ \ \ \ \ \ \ \ \ \ \ \ \ \ \ \ \ \ \ \ \  \ \ \ \ \ \ \   (\,{\mbox{cf. \ \ the \ \ argument \ \ in}} \ \ \S\,A\,.5 \,)\\
 & = & I_o\,({\bf z})     + o \,\left(\,\int_{\R^n} {\bf z}^{{n + 2}\over {n - 2}} \cdot w \right) \ \ \ \ \ \  \ \ \ \ \ \   \ \ \ \ \ \   \ \ \ \ \ \  \left(\, \, p ={{2n}\over {n + 2}}\,, \ \ q = {{2n}\over {n - 2}}\, \right)\ \ \ \ \ \ \\
 & \ & \\
 & = & I_o\,({\bf z})  + o \,(\Vert\,w \Vert_{ \btd })  \ \  \ \ \ \ \ \ \ \ \ \ \ \ \ \ \ \ \ \ \ \ \ \ \ \ \ \   (\,{\mbox{uniformly, \ \ \ \ as}} \ \ \Vert \,w \Vert_{\btd} \to 0)\,.
         \end{eqnarray*}
Here
 \begin{eqnarray*}
I_o\,({\bf z}) & = &  {1\over 2} \int_{\R^n} \left[\,  \langle \btd_o \,{\bf z}\,,\, \btd_o\,{\bf z} \rangle \,-\,   (n - 2)^{\,2}   \,{\bf z}^{{2n}\over {n - 2}} \,\right]\\[0.075in]
 &= &{1\over 2} \int_{\R^n} \left[\,  [\,-\,\Delta_o \,{\bf z}]\cdot {\bf z} \, -\,   (n - 2)^{\,2}   \,{\bf z}^{{2n}\over {n - 2}} \,\right]\\[0.075in]
&= & {1\over 2} \int_{\R^n} \left[\,  n\,(n -2)\, {\bf z}^{{2n}\over {n - 2}}  \,-\,   (n - 2)^{\,2}   \,{\bf z}^{{2n}\over {n - 2}} \,\right]\\[0.075in]
 &= & (n - 2) \int_{\R^n} \left(\,{1\over {1 + |\,y|^{\,2}}} \right)^n \ = \ {\bar c}_4\,,
 \end{eqnarray*}
 via a translation and rescaling.
Likewise,
\begin{eqnarray*}
 {\bf II.} \ \ \ \ \ & \ &  (\,G' ({\bf z} + w) - G'({\bf z}) - G''({\bf z})[\,w])\,[\bullet] \\[0.075in] & = & C\, \int_{\R^n} \left(\,({\bf z} + w)^{{n + 2}\over {n - 2}}\ - \ {\bf z}^{{n + 2}\over {n - 2}} \ - \  {{n + 2}\over {n - 2}}\  {\bf z}^{4\over {n - 2}}\cdot w \right)\,[\bullet]\,\ \ \ \ \ \ \ \ \ \ \ \ \ \ \ \ \ \ \ \ \ \ \ \ \\[0.075in]
 & = &  o \,(\Vert \,w \Vert_{ \btd }) \ \ \ \ \ \ \ \ \ \ \ \  \ \ \ \ \ \ \ \ \ \ \ \ \ \ \ \ \ \ \ \   (\,{\mbox{uniformly, \ \ as}} \ \ \Vert\, w \Vert_{ \btd } \,\to \, 0)\,.\ \ \ \ \ \ \ \ \ \ \ \ \ \ \ \ \ \ \ \ \ \ \ \ \
         \end{eqnarray*}
Here indeed the formulas allow $w \in {\cal D}^{1, \ 2}$ with suitable bounds on $\Vert\, w \Vert_{ \btd }\,.$

\vspace*{0.4in}

{\bf   \S\,A.\,8. \ \ Argument toward Proposition 4.15 in Part I.}\\[0.1in]
Let $\,\{ q_\ell^{\varepsilon\,, \, i} \}\,$ be an orthonormal basis of $\,T_{z_{\varepsilon, \,i}}  {\bf Z}\,.$\, See (3.5). Here
$$
{\bf z}_{\varepsilon\,, \,i} (y) = \left(\, {{\lambda_{\varepsilon\,, \,i} }\over { \lambda_{\varepsilon\,, \,i}^{\,2} + | \,y - \xi_ {\varepsilon\,, \,i}  | }} \right)^{{n - 2}\over 2}.
$$
Via the  Lemma 4.1, given any $\delta > 0$, there exists a positive number $\,\bar \varepsilon\,$ such that
$$
\Vert\,w_\varepsilon ({\bf z}) \Vert_{\btd} \le  \delta \mfor {\mbox{all}} \ \  z \in {\bf Z} \ \ {\mbox{and}} \ \ |\,\varepsilon|\le \bar \varepsilon\,. \leqno  (A.8.1)
$$
We follow closely the argument in the proof of Theorem 2.12, as well as the notations used in \cite{Progress-Book}. As the argument is the same at every critical point $z_{\varepsilon\,, \,i}$, we simplify the notation by taking away the suffix  and superfix   ${\varepsilon\,, \,i}$. That is, in what follows
\begin{eqnarray*}
q_\ell^{\varepsilon, \,i} &\to& {\tilde q}_\ell\\
{\bf z}_{\varepsilon, \,i} &\to& {\tilde z}\\
w_{\varepsilon, \,i} &\to& {\tilde w} \\
\lambda_{\varepsilon, \,i} &\to& {\tilde \lambda} \\
& {\mbox{etc.}}&
\end{eqnarray*}
Set
$$
{\cal B}_{j, \, \ell}   \ :=\ {\tilde \lambda} \cdot \langle D_j \, w \,,\, \ {\tilde q}_\ell  \, \rangle \ \ \ \ \ \ \ \ \ \ \  \mfor j\,, \ \,\ell \ =\  0, \ 1, \ 2, \cdot \cdot \cdot, \ n\,,
$$
and
$$
D_o    \,=\, {{\partial  }\over {\partial \lambda}}\,, \ \ \ D_k\, \,=\, {{\partial }\over {\partial \xi_k}} \ \ \mfor \ \  k \,=\,  1, \ 2, \cdot \cdot \cdot, \ n\,.
$$
As $w_\varepsilon \,({\bf z})$, the solution to the auxiliary equation, is  perpendicular to $ T_{\bf z} \,{\bf Z}$\,,\,
one has
\begin{eqnarray*}
(A.8.2)\!\!\!\!\!\!& \ & \langle \,\tilde w \,(\tilde {\bf z})\,, \ {\tilde q}_\ell\,\rangle_{\btd} = 0 \ \ \ \ \ \ \mfor \ \  \ell = 0, \ 1, \ 2, \cdot \cdot \cdot, \ n\\[0.075in]
& \Longrightarrow & \ \ \langle \,D_j \,{\tilde w} \, (  \tilde {\bf z})\,, \ {\tilde q}_\ell\,\rangle_{\btd}\, +\, \langle \,{\tilde w}  (\tilde z)\,, \ D_j\,{\tilde q}_\ell\,\rangle_{\btd} \ =\ 0\\[0.075in]
& \Longrightarrow & \ \   |\,\langle \,D_j \,\tilde w \,(\tilde {\bf z})  \,, \ {\tilde q}_\ell\,\rangle_{\btd}\,|\le    \Vert \tilde w \,(\tilde {\bf z}) \Vert_{\btd} \cdot \Vert D_j\,\,{\tilde q}_\ell\,  \Vert_{\btd}  \le {{C(n)}\over {\tilde \lambda}} \cdot  \Vert\,  \tilde w \,(\tilde {\bf z}) \Vert_{\btd}\\[0.075in]
& \ & \ \ \ \ \ \ \ \ \  \ \ \ \ \ \ \ \ \   \ \ \ \ \ \ \ \ \ \ \ \ \ \ \ \ \ \   \ \ \ \ \ \ \ \ \   \ \ \ \ \ \ \ \ \   \ \   \ \ \ \ \ \ \  ({\mbox{using \ \ (3.6))}}\   \ \ \ \ \ \ \ \ \   \ \  \   \ \ \ \ \ \ \ \ \   \ \   \\
& \Longrightarrow & \ \ |\,{\cal B}_{j,\,\ell} | \ \le \   C (n) \cdot  \Vert  \tilde w \,(\tilde {\bf z})    \Vert_{\btd} \ \le \ C(n)   \cdot \delta\,.
\end{eqnarray*}
Here we apply (A.8.1). As $\tilde{\bf z} $ is a critical point of $\Phi_\varepsilon$ (see the chart in \S\,2\,e), we have
\begin{eqnarray*}
(A.8.3)\!\!\!  & \ & \langle \,I'_\varepsilon (\tilde{\bf z}+ \tilde w\, (\tilde{\bf z})) \,, \ \ D_j\,\tilde{\bf z}  + D_j \,\tilde w\, (\tilde{\bf z}) \, \rangle_{ \btd} = 0 \mfor j = 0, \ 1, \cdot \cdot \cdot, \ n\\[0.1in]
& \Longrightarrow & \bigg\langle \,I'_\varepsilon (\tilde{\bf z}+ \tilde w\, (\tilde{\bf z})) \,, \ \ {{ D_j\,\tilde{\bf z}  + D_j \,\tilde w\, (\tilde{\bf z}) }\over {  \Vert D_j\,\tilde{\bf z}  \Vert_{\btd} }} \, \bigg\rangle_{ \btd}\ =\  0\\[0.075in]
& \Longrightarrow &  \langle \,I'_\varepsilon (\tilde{\bf z}+ \tilde w\, (\tilde{\bf z})) \,, \ \ {\tilde q}_j   \rangle_{ \btd}
+  {\tilde C}_j \,\tilde \lambda \ \bigg\langle \,I'_\varepsilon (\tilde{\bf z}+ \tilde w\, (\tilde{\bf z})) \,, \ \      D_j \,\tilde w\, (\tilde{\bf z})    \bigg\rangle_{ \btd}
\ =\ 0\,. \ \ \ \ \ \ \ \ \ \ \ \ \ \ \ \ \ \ \
\end{eqnarray*}
Here we use (3.2) and (3.3). Write
$$
I'_\varepsilon \,(\tilde{\bf z}+ \tilde w\, (\tilde{\bf z})) \ = \ \sum_{\ell = 0}^n A_\ell \,\,{\tilde q}_\ell\,.
$$
$$
``\,(A.8.3)\," \Longrightarrow \ \ A_j  \ + \ \sum_\ell \, A_\ell\, [\,{\cal B}_{j,\, \ell} \cdot {\tilde C}_j ] \ = \ 0\,.
$$
Other part of the proof follows as in the proof of Theorem 2.12 (pp. 22 in \cite{Progress-Book})\,.\qedwh
{\bf Remark.} \ \ The smallness of $\,{\cal B}_{j,\,\ell}$\,
basically said that the ``tangent spaces" at\,
${\bf z}$ and at \,${\bf z} + w({\bf z})$ are  almost parallel\,.

\newpage

{\bf   \S\,A.\,9. \ \    $G'$.}\\[0.1in]
Refer to the proof of Lemma 4.16 in Part I.
Here
$$
G \,(u) \ = \ -\,{\bar c}_n {{n - 2}\over {2n}}  \int_{\R^n} H u_+^{{2n}\over {n - 2}} \,  \,.
$$
It follows that
$$
G' ({\bf z} + w) ([\bullet]) \ = \ -\,{\bar c}_n    \int_{\R^n} H \cdot ({\bf z} + w)_+^{{n+2}\over {n - 2}} \cdot [\bullet] \mfor {\bf z} \in {\bf Z} \ \ {\mbox{and}} \ \ w \in {\cal D}^{1, \ 2}\,. \leqno (A.9.1)
$$
Hence
\begin{eqnarray*}
|\,G' ({\bf z} + w) \,([\bullet])| &\le &C_1 \int_{\R^n} |\,{\bf z} + w|^{{n+2}\over {n - 2}} \,|\,[\bullet]|\ \ \ \ \ \ \ \ \ (\,{\mbox{as}} \ \ H \ \ {\mbox{is \ \ bounded\,}}) \\[0.1in]
& \le & C_1 \left(\, \int_{\R^n}  | \, z + w|^{{{n + 2}\over {n - 2}} \cdot {{2n}\over {n + 2}} }   \right)^{{n + 2}\over {2n}} \left(\, \int_{\R^n}  |\,[\bullet]|^{{2n}\over {n - 2}}  \right)^{{n - 2}\over {2n}}\\[0.1in]
& \le & C_2 \left(\, \int_{\R^n}  ( \, z + |\, w| \,)^{ {{2n}\over {n-2}} }  \right)^{{n + 2}\over {2n}} \cdot \Vert \,[\bullet] \Vert_\btd\,.
\end{eqnarray*}
Using the inequality
$$
( \, {\bf z} + |\, w|\,)^{{2n}\over {n - 2}}\  \le\  C_3 \left[\,{\bf z}^{{2n}\over {n - 2}}  \,+\, |\, w|^{{2n}\over {n - 2}}\, \right]\,,
$$
we obtain
\begin{eqnarray*}
(A.9.2)\ \ \  \ \ |\,G' \,(\,{\bf z} + w) ([\bullet])| & \le & C_4 \left(\, \int_{\R^n} {\bf z}^{ {{2n}\over {n-2}} }  \,+\,  | \, w|^{ {{2n}\over {n-2}} } \,\right)^{{n + 2}\over {2n}} \cdot \Vert \,[\bullet] \Vert_\btd\\[0.1in]
& \le & C_4\, [\, C_5 \,+\, C_6]^{{n + 2}\over {2n}} \cdot \Vert \, [\bullet] \Vert_\btd \ \ \ \ \mfor \Vert \,w \Vert_\btd \le 1\,. \ \ \ \ \ \ \
\end{eqnarray*}
Here we use the Sobolev inequality.
(A.9.2) leads to
$$
\Vert \,G' \,(z + w) \Vert \ \le\  C \ \ \mfor {\mbox{all}} \ \ {\bf z} \in {\bf Z} \ \ {\mbox{and}} \ \ w \ \ {\mbox{with}} \ \  \Vert \,w \Vert_\btd\ \le \ 1\,. \leqno (A.9.3)
$$

\newpage

{\bf \S\,A.\,10. \ \ Kazdan\,-\,Warner condition.}\\[0.1in]
For equation (1.1),
 a guiding relation is revealed when we differentiate $\,{\cal K}\,$ with respect to a {\it conformal Killing vector field\,} $\,X\,$ -- one that generates a family of conformal transformations. In this way we obtain the renowned Kazdan-Warner formula \cite{BourguignonEzin}
$$
\int_{S^n}  X ({\cal K}) \, u^{{2n}\over {n - 2}} \, dV_{g_1}\ =\  0\,, \ \ \ \ \ \ {\mbox{where}} \ \ X ({\cal K})\, =\,  \langle\,X\,, \ \btd_{g_1} \,{\cal K}\rangle_{g_1}\,.
$$
Simple and elegant, the Kazdan-Warner formula encapsulates a central character of the equation, namely, the {\it balance}\, between $\,\btd_{g_1} \,{\cal K}\,$ and $\,u\,.$ \\[0.2in]
{\bf Definition A.10.1.} \ \ {\it A function\,} $\hat{\cal K} \in C^1 (\R^n)\,$ {\it is said to satisfy the K\,-W condition if there exists a\,} ({\it single}) {\it positive  function\,} $f \in C^o (S^n)$ {\it such that}\\[0.1in]
(A.10.2)
$$
\int_{S^n}  X (\hat{\cal K}) \, f^{{2n}\over {n - 2}} \  dV_{g_1} \ = \ 0  \ \ \ \ \ \ {\it{for \ \ all \ \ conformal \ \ Killing \ \ vector \ \ fields \ \ }} \ \ X\,.
$$

\vspace*{0.15in}

\hspace*{0.5in}The collection of all conformal Killing vector field on $S^n$ can be regarded as a linear space of dimension $(n + 1)(n + 2)/2$, with a basis formed by the generators  of the dilations ($n + 1$ dimension),\,
denoted by $X_1,\cdot \cdot \cdot, \, X_{n + 1}\,,\,$
and of the rotations [\,$n (n + 1)/2$ dimension],\,
denoted by $\,X_{n + 2}\,,\cdot \cdot \cdot, \, X_{{(n + 1)(n + 2)}\over 2}\,$.
Refer to \cite{Han-Li}.\\[0.2in]
{\it Dilations.} \ \ By the homogeneity of $S^n$, these are generated by
$$\btd_{g_1} \,x_\ell \mfor \ell \ = \ 1, \ 2, \cdot \cdot \cdot, \ n + 1\,.\,$$

\vspace*{0.1in}

{\it Rotations.} \ \ Let  $\Theta$ be a rotation of $S^n$, which is an isometry. Consider the integral
$$
\int_{S^n}  [ \,\hat{\cal K} \circ \Theta\,] \, f^{{2n}\over {n - 2}} \, dV_{g_1} \ = \ \int_{S^n}   \hat{\cal K} \circ   \, [\,f \circ \Theta^{-1}]^{{2n}\over {n - 2}} \  dV_{g_1}
$$
Here we apply the change of variables formula, noticing that the rotation of the standard metric $g_1$ on $S^n$ is isometric to itself. In particular, if
$$
f (x) \ = \ \left(\,{{\lambda}\over{\lambda^{\,2} + |\,y - \xi|^{\,2}}} \right)^{{n-2}\over {2}}\,,
$$
then, via the Obata theorem \cite{Obata}, we can express
$$
f (\Theta^{-1} (x))\ = \ \left(\,{{\lambda'}\over{{\lambda'}^{\,2} + |\,y - \xi'|^{\,2}}} \right)^{{n-2}\over {2}} \,. \leqno (A.10.3)
$$

\vspace*{0.15in}

{\bf \S\,A.10.1.} \ {\it Pohozaev identity.} \ \  Expressing (A.10.1) in $\R^n$ via  the stereographic projection, we obtain
$$
\Delta_o \,v \,+ \,{\tilde c}_n  \, K  v^{{n + 2}\over {n - 2}} \ = \  0\ \ \ \ {\mbox{in}} \ \ {\R}^n,
\leqno (A.10.4)
$$
where
$$
K  (y)\ := \ {\cal K} \,({\dot{\cal P}}^{-1} (y)) \ \ \ \ \ \ {\mbox{and}} \ \ \ \ \ v \,(y)\ =\ u \,({\dot{\cal P}}^{-1} (y)) \,\left(\, {2\over{1 + |y|^{\,2}}} \right)^{{n-2}\over 2}.  \leqno (A.10.5)
$$
[\,Cf. (1.5) and (2.4)\,.]\,
Corresponding to $X_{n + 1} = \btd_{g_1} \,x_{n + 1}$\,,\, note that
$$
\btd_{g_1}   \,x_{n + 1} \ = \ g_1^{ij} \,\left[\, {{\partial \,[ \,x_{n + 1}  \circ {\dot{\cal P}}^{-1} (y)]}\over {\partial y_i}}\right] \, \partial_j\ = \ {1\over 4} \,(1 + r^{\,2})^{\,2} \left[ {{4 y_i}\over {(1 + r^{\,2})^{\,2}}} \right]\partial_i \ = \ \sum_{i = 1}^n y_i \,\partial_i\,,
$$
as
$$
x_{n + 1} \ = \ {{r^{\,2} \,-\, 1}\over {r^{\,2} \,+\, 1}}\,.
$$
Likewise,
$$
\btd_{g_1}  {\cal K} \ = \  g_1^{ij} \,\left[\, {{\partial \,{\cal K}  \circ {\dot{\cal P}}^{-1} (y)]}\over {\partial y_i}}\right] \, \partial_j\ = \  {1\over 4} \,(1 + r^{\,2})^{\,2} \left[\, {{\partial \,{\cal K}  \circ {\dot{\cal P}}^{-1} (y)]}\over {\partial y_i}}\right] \partial_i
$$
Hence
$$
X ({\cal K} ) \ = \  \sum_{i = 1}^n  {4\over {1 + r^{\,2}}} \,\left[y_i \,\partial_i  \right] \cdot \left[ {1\over 4} \,(1 + r^{\,2})^{\,2} \left[\, {{\partial \,{\cal K}  \circ {\dot{\cal P}}^{-1} (y)}\over {\partial y_i}}\right] \partial_i \right] = \sum_{i = 1}^n  y_i \,{{\partial \,{\cal K}  \circ {\dot{\cal P}}^{-1} (y)}\over {\partial y_i}}\,.
$$
Using the stereographic transformation, (A.10.1), (A.10.4) and (A.10.5) lead to  the following radial Pohozaev identities
$$
\int_{\R^n} r {{\partial   K}\over {\partial r}}\, v^{{2n}\over {n - 2}}  \ = \ \sum_{i = 1}^n \, \int_{\R^n}   y_i  {{\partial  K}\over {\partial y_i }} v^{{2n}\over {n - 2}}  \ = \  0 \ \ \ \ \ [\,  K (y) \ = \ {\cal K} \circ {\dot{\cal P}}^{-1} (y) ]\,. \leqno (A.10.6)
$$
This corresponds to the rescalings by a positive number $\sigma$\,:
$$
(r, \ \vartheta) \,\mapsto\, (r \,\sigma, \ \vartheta)\,, \ \ \ \ {\mbox{where}} \ \ r \,\ge\, 0 \ \ {\mbox{and}} \ \ \vartheta \,\in\, S^{n - 1}.\leqno (A.10.7)
$$
As for
$X_i \ = \ \btd_{g_1} \,x_i$\, when \,$i \ = \  1, \ 2, \cdot \cdot \cdot, \ n\,,$\, we have (for example)
\begin{eqnarray*}
x_1 &= & {{2y_1}\over {1 + r^{\,2}}}\,, \\
{\mbox{and}} \ \ \ \ \ \btd_{g_1}   x_1 &= &g_1^{ij} \,\left[\, {{\partial \,[ \,x_1  \circ {\dot{\cal P}}^{-1} (y)]}\over {\partial y_i}}\right] \, \partial_j\ = \  {1\over 4} \,(1 + r^{\,2})^{\,2} \left[ {{\partial}\over {\partial y_i}} \left(\, {{2y_1  }\over {1 + r^{\,2} }} \right) \right]\partial_i\ \ \ \ \ \ \ \ \ \ \\
 &= &  \left[ \,{1\over 2} \,(1 + r^{\,2}) \ - \ \  y_1^{\,2}\right]\,\partial_1 \ -\, \sum_{i = 2}^n y_1\,y_i \,\partial_i\,.
\end{eqnarray*}
From (A.10.1), (A.10.4) and (A.10.5), we obtain
\begin{eqnarray*}
(A.10.8) \ \ \ \ \ \ \ \ {1\over 2} \int_{\R^n}  {{\partial K}\over {\partial y_1}} \, v^{{2n}\over {n - 2}} & + &
{1\over 2} \int_{\R^n} r^{\,2} {\,{\partial K}\over {\partial y_1}} \, v^{{2n}\over {n - 2}} \,
- \int_{\R^n} y_1^{\,2} \,{{\partial K}\over {\partial y_1}} \, v^{{2n}\over {n - 2}} \, \ \ \ \ \ \   \ \ \ \ \ \  \ \ \ \ \ \  \\[0.1in]
&  \  & \ \ \ \ \ \   \ \ \ \ \ \  \ \ \ \ \ \   -\ 
\sum_{i = 2}^n \, \int_{\R^n} y_1\,  y_i  {{\partial K}\over {\partial y_i }} \,v^{{2n}\over {n - 2}} \, \ =\ 0\,.
\end{eqnarray*}
It follows that if we have the translational   Pohozaev   identity\,:
$$
\ \int_{\R^n}  {{\partial K}\over {\partial y_1}} \, v^{{2n}\over {n - 2}} \, = 0\,, \leqno (A.10.9)
$$
\begin{eqnarray*}
{\mbox{then}} \ \ \ \ \ \ \ \ \ (A.10.8)
  & \Longrightarrow  &
{1\over 2} \int_{\R^n} r^{\,2} {\,{\partial K}\over {\partial y_1}} \, v^{{2n}\over {n - 2}} \,
 -
\sum_{i = 1}^n \, \int_{\R^n} y_1\,  y_i  {{\partial \hat   K}\over {\partial y_i }} \,v^{{2n}\over {n - 2}} \, \ =\  0 \ \ \ \ \ \ \ \ \\[0.1in]
  & \Longrightarrow  &
  \int_{\R^n} r^{\,2} {\,{\partial K}\over {\partial y_1}} \, v^{{2n}\over {n - 2}} \,
\ =\ 2
 \, \int_{\R^n} y_1\, \left[ \sum_{i = 1}^n  y_i  {{\partial K}\over {\partial y_i }}\right] \,v^{{2n}\over {n - 2}} \,\\[0.1in]
(A.10.10)\,\cdot \cdot \cdot \cdot \cdot \cdot  \ \ \ \    & \Longrightarrow  &
  \int_{\R^n} r^{\,2} {\,{\partial K}\over {\partial y_1}} \, v^{{2n}\over {n - 2}} \,
\ =\ 2
 \, \int_{\R^n} y_1\, \left[ r  {{\partial K}\over {\partial r }}\right] \,v^{{2n}\over {n - 2}} \,\,.\ \ \ \ \ \ \ \ \ \ \ \ \ \ \ \ \ \ \ \ \ \ \ \ \ \ \
\end{eqnarray*}

\vspace*{0.2in}

{\bf Theorem A.10.11.} \ \ {\it For a given ${\cal H} \in C^1 (S^n)$\,,\, let $H (y) = {\cal H} \, ({\dot{\cal P}}^{-1} (y))$\,.\, Suppose $G_{|_{\bf Z}}$\,} [\,{\it corresponding to}\, $H$] {\it has a critical point. Then }
$$
{\hat{\cal K}} \ = \ 1 \,+\, \varepsilon \,{\cal H}
$$
{\it fulfills the  K\,-W  condition\,}\, (A.10.\,2)
{\it for any\,}  $\varepsilon \in \R$\,. \\[0.2in]
%
%
%
{\bf Proof.} \ \ Let $(\lambda_c\,, \ \xi_c)$ be a critical point of $G_{|_{\bf Z}}$\,.\,  It follows from  formulas (6.1) and (6.3) that
\begin{eqnarray*}
 (A.10.12) \ \ \, {{\partial G_{|_{\bf Z}} }\over {\partial \lambda }} (\lambda_c\,, \ \xi_c) & = & 0 \ \ \Longrightarrow \ \  n\,\int_{\R^n}   H (y) \,  {{ \lambda_c^{n-1}\,  }\over {\left(\,\lambda_c^{\,2}\, +\, |\,y - {\xi_c}|^{\,2}\right)^n}}\,\, \\
& \  & \ \ \ \ \ \ \ \ \ \ \  - 2n\,  \int_{\R^n}   H (y)\,   {{\lambda_c^{n+1}\, }\over {\left(\,\lambda_c^{\,2} + |\,y - {\xi_c}|^{\,2}\right)^{n+1}}}\ =\   0\,,\ \ \ \ \ \ \ \ \ \ \ \ \ \ \  \ \\
 (A.10.13) \ \  {{\partial G_{|_{\bf Z}}}\over {\partial \xi_i}}\, (\lambda_c\,, \ \xi_c) & = & 0 \ \ \Longrightarrow \ \   \int_{\R^n}  H (y) \,  {{\lambda_c^n\, ({\xi_c}_i - y_i)}\over {\left(\,\lambda_c^{\,2} + |\,y - {\xi_c}|^{\,2}\right)^{n+1}}} \ =\   0\,.
\end{eqnarray*}
Here $H (y) = {\cal H} \, ({\dot{\cal P}}^{-1} (y))$\,.\, \\[0.2in]
{\it Translational Case.} \ \ From  (A.10.13) we have
$$
 \int_{\R^n} H (y) \, {\partial\over {\partial y_i}} \left[ \left(\, {{\lambda_c}\over {\lambda_c^{\,2} \,+\, |\,y - \xi_c|^{\,2}}} \right)^n\right] \ =\  0\,.\leqno (A.10.14)
$$
We have
$$
|\,H| \ \le \ C_1 \ \ \ {\mbox{and}} \ \ \ | \btd \,  H (y)|     \le {{C_2}\over {|\,y|^{\,2}}} \mfor |\,y|\gg 1
$$
(see Lemma 9.37 in \cite{Leung-Supported})\,.\, Also,
$$
\lim_{|\,y_1| \to \infty} \left[ H \,(y_1\,, \cdot \cdot \cdot, \ y_n) \cdot \left(\, {{\lambda_c}\over {\lambda_c^{\,2} + |\,y - \xi_c|^{\,2}}} \right)^n \right] \ =\  0 \mfor (y_2\,, \cdot \cdot \cdot, \ y_n) \in \R^{\,n - 1}.
$$
Via  Fubini's theorem and integration by parts, we obtain  from (A.10.13)
the translational   Pohozaev   identity\,:
$$
\ \int_{\R^n}  {{\partial H}\over {\partial y_1}} \, w^{{2n}\over {n - 2}} \ =\  0 \leqno (A.10.15)
$$
with $$
w (y) \ =\   \left(\, {{\lambda_c}\over {\lambda_c^{\,2} + |\,y - \xi_c|^{\,2}}} \right)^{{n-2}\over 2}.\leqno (A.10.16)
$$
{\it Radial Case.} \ \ Consider the integral\\[0.1in]
(A.10.17)
\begin{eqnarray*}
& \ & \sum_{i = 1}^n \int_{\R^n} H (y)\, {\partial\over {\partial y_i}} \left[(y_i - {\xi_c}_i) \cdot \left(\, {{\lambda_c}\over {\lambda_c^{\,2} + |\,y - \xi_c|^{\,2}}} \right)^n\right]  \,\\[0.1in]
& = & \sum_{i = 1}^n \int_{\R^n} H (y)\,   \left[ \, \left(\, {{\lambda_c}\over {\lambda_c^{\,2} + |\,y - \xi_c|^{\,2}}} \right)^n - {{2n\, (y_i - {\xi_c}_i)^{\,2}\lambda_c^n}\over {[\,\lambda_c^{\,2} + |\,y - \xi_c|^{\,2}]^{n + 1}}} \right]  \,\\[0.1in]
& = & n\, \int_{\R^n}  H (y)\,     \left(\, {{\lambda_c}\over {\lambda_c^{\,2} + |\,y - \xi_c|^{\,2}}} \right)^n  \ - \ 2n\, \int_{\R^n}  H (y)\,  {{  |\,y - {\xi_c}|^{\,2}\lambda_c^n}\over {[\,\lambda_c^{\,2} + |\,y - \xi_c|^{\,2}]^{n + 1}}}  \  \\[0.1in]
& = & n\, \int_{\R^n}  H (y)\,     \left(\, {{\lambda_c}\over {\lambda_c^{\,2} + |\,y - \xi_c|^{\,2}}} \right)^n \ -\  2n\, \int_{\R^n} H (y)\,  {{ [|\,y - {\xi_c}|^{\,2} + \lambda_c^{\,2} - \lambda_c^{\,2}]\,\lambda_c^n}\over {[\,\lambda_c^{\,2} + |\,y - \xi_c|^{\,2}]^{n + 1}}}  \  \\[0.1in]
& = & - n\, \int_{\R^n}   H (y)\,     \left(\, {{\lambda_c}\over {\lambda_c^{\,2} + |\,y - \xi_c|^{\,2}}} \right)^n \ +\ 2n\, \int_{\R^n}  H (y)\,  {{ \lambda_c^{n+2} }\over {[\,\lambda_c^{\,2} + |\,y - \xi_c|^{\,2}]^{n + 1}}} \  \\[0.1in]
& = & 0 \ \ \ \   \ \ \ \ \ \  \ \ \ \   \ \ \ \ \ \ \ \ \ \   \ \ \ \ \ \  \ \ \ \   \ \ \ \ \ \  \ \ \ \   \ \ \ \ \ \ \ \ \ \   \ \ \ \ \ \ \ \ \ \ \ \ \ \ \ \ \  [\,{\mbox{using \ \ (A.10.12)}}]\,.
\end{eqnarray*}

\newpage

On the other hand, via  Fubini's theorem and integration by parts, we obtain
\begin{eqnarray*}
& \ & \sum_{i = 1}^n \int_{\R^n} H (y)\, {\partial\over {\partial y_i}} \left[(y_i - {\xi_c}_i) \cdot \left(\, {{\lambda_c}\over {\lambda_c^{\,2} + |\,y - \xi_c|^{\,2}}} \right)^n\right]  \,\\
& = & - \sum_{i = 1}^n \int_{\R^n}   {{\partial H}\over {\partial y_i}}\, \left[(y_i - {\xi_c}_i) \cdot \left(\, {{\lambda_c}\over {\lambda_c^{\,2} + |\,y - \xi_c|^{\,2}}} \right)^n\right]  \,\\
& = & - \sum_{i = 1}^n \int_{\R^n}   y_i \, {{\partial H}\over {\partial y_i}}\,   \cdot \left(\, {{\lambda_c}\over {\lambda_c^{\,2} + |\,y - \xi_c|^{\,2}}} \right)^n   \, + \sum {\xi_c}_i \int_{\R^n}   {{\partial H}\over {\partial y_i}}\,   \cdot \left(\, {{\lambda_c}\over {\lambda_c^{\,2} + |\,y - \xi_c|^{\,2}}} \right)^n   \, \\
& = & -  \int_{\R^n} \left(\, \sum_{i = 1}^n y_i \, {{\partial H}\over {\partial y_i}}   \right)   \cdot \left(\, {{\lambda_c}\over {\lambda_c^{\,2} + |\,y - \xi_c|^{\,2}}} \right)^n   \,   \ \ \ \ \ \ \ \  \ \ \ \ \ \ \ \  \ ({\mbox{using   \      (A.10.15)}}\\
   & = & -   \int_{\R^n}   r\, {{\partial H}\over {\partial r}}\,   \cdot \left(\, {{\lambda_c}\over {\lambda_c^{\,2} + |\,y - \xi_c|^{\,2}}} \right)^n   \, = 0 \ \ \ \ \ \ \ \ \ \ \ \ \ \ \ \ \ \  \ \ \ ({\mbox{via}} \ \ (A.10.17))\,.
\end{eqnarray*}
Cf. the radial Pohozaev identity (A.10.6). Running the argument backward as in (A.10.6), we have\\[0.1in]
(A.10.18)
$$
\int_{S^n} X_{n + 1} \,({\cal H})\, f^{{2n}\over{n - 2}}\,dV_{g_1} \ =\   0\,, \ \ \ \ {\mbox{where}} \ \ f (x) \ =\    \left(\, {{\lambda_c}\over {\lambda_c^{\,2} + |\,y - \xi_c|^{\,2}}} \right)^{{n-2}\over 2} \cdot \left(\,{{1 + |\,y|^2}\over 2 } \right)^{{n - 2}\over 2}.
$$
Cf. also (A.10.5)\,.\\[0.2in]
{\it $X_1\,, \cdot \cdot \cdot\, \ X_n$ \,directions.} \ \
From (A.10.7) , (A.10.8), (A.10.9) and (A.10.10), we need to show that
$$
\int_{\R^n} r^{\,2} {\,{\partial \hat H}\over {\partial y_1}} \, \left(\,{\lambda_c\over {\lambda_c^{\,2} + |\,y - \xi_c|^{\,2}}} \right)^n  \ =\  2
 \, \int_{\R^n} y_1\, \left[ r  {{\partial \hat H}\over {\partial r }}\right] \, \left(\,{\lambda_c\over {\lambda_c^{\,2} + |\,y- \xi_c|^{\,2}}} \right)^n \,  \leqno (A.10.19)
 $$
 in order to obtain \\[0.1in]
(A.10.20)
 $$
\int_{S^n} X_{i} \,({\cal H})\, f^{{2n}\over{n - 2}}\,dV_{g_1} \ =\   0\,, \ \ \ \ {\mbox{for}} \ \ i \ =\   1, \ 2, \ \cdot \cdot \cdot, \ n\,, \ \ {\mbox{where}} \ \ f (x) \ \ {\mbox{is \ \ given \ \ in \ \ (A.10.18)}}\,.
$$
Using Fubini's theorem and integration by parts, we obtain
\begin{eqnarray*}
(A.10.21)\!\!\!\!\!\!& \ & \int_{\R^n} r^{\,2} {\,{\partial   H}\over {\partial y_1}} \, \left(\,{\lambda_c\over {\lambda_c^{\,2} + |\,y - \xi_c|^{\,2}}} \right)^n  \,\\[0.1in]
& = & \int_{\R^n}  H {\,{\partial  }\over {\partial y_1}} \left[\ r^{\,2}  \, \left(\,{\lambda_c\over {\lambda_c^{\,2} + |\,y - \xi_c|^{\,2}}} \right)^n \right] \,\\[0.1in]
& = & -2 \int_{\R^n} y_1  \, H \left(\,{\lambda_c\over {\lambda_c^{\,2} + |\,y - \xi_c|^{\,2}}} \right)^n dy\ + \ 2n \! \int_{\R^n} r^{\,2}   H \cdot {{(y_1 - {\xi_c}_1) \,\lambda_c^n}\over {(\lambda_c^{\,2} + |\,y - \xi_c|^{\,2})^{n + 1}}} \  dy\,.
\end{eqnarray*}

\newpage

On the other side,

\begin{eqnarray*}
(A.10.22) \ \    & \ & 2
 \, \int_{\R^n} y_1\, \left[ r  {{\partial H}\over {\partial r }}\right] \, \left(\,{\lambda_c\over {\lambda_c^{\,2} + |\,y- \xi_c|^{\,2}}} \right)^n \,\\[0.15in]
 & = & 2
 \, \int_{S^{n - 1}} \int_0^\infty   {{\partial H}\over {\partial r }}  \, \left[ r^n y_1\,  \left(\,{\lambda_c\over {\lambda_c^{\,2} + |\,y- \xi_c|^{\,2}}} \right)^n \right]\, dr \, dS_\vartheta\\[0.15in]
 & = & -2n
 \, \int_{S^{n - 1}} \int_0^\infty     {H } \, y_1\,  r^{n-1}   \left(\,{\lambda_c\over {\lambda_c^{\,2} + |\,y- \xi_c|^{\,2}}} \right)^n \, dr \, dS_\vartheta\\
 & \   & \ \ \ \  -\, 2
 \, \int_{S^{n - 1}} \int_0^\infty     {H } \,  y_1\, r^{n-1}   \left(\,{\lambda_c\over {\lambda_c^{\,2} + |\,y- \xi_c|^{\,2}}} \right)^n \, dr \, dS_\vartheta\\
 & \   & \ \ \ \   \ \ \ \  +\, 2n
 \, \int_{S^{n - 1}} \int_0^\infty     {H } \,  y_1\, r^n   \left(\,{{\lambda_c^n \left(\, 2 r - 2 \sum {\xi_c}_i \!\cdot \! {{y_i}\over r} \right)} \over {[\,\lambda_c^{\,2} + r^{\,2} - 2 \sum {\xi_c}_i \!\cdot \! y_i + |\,{\xi_c}|^{\,2}\,]^{n + 1} }} \right)  \, dr \, dS_\vartheta\,,
\end{eqnarray*}
where we apply integration by parts formula.
\begin{eqnarray*}
& \ & {\mbox{The \ \ last \ \ term \ \ in \ \ (A.10.22)}}\\
& \ &\\
 &=  &   2n
 \, \int_{S^{n - 1}} \int_0^\infty     {H } \, y_1\,   r^{n-1}  \left(\,{{\lambda_c^n \left(\, 2 r^{\,2} - 2 \sum {\xi_c}_i \!\cdot \! y_i \right)} \over {(\,\lambda_c^{\,2} + |\,y - \xi_c|^{\,2} )^{n + 1} }} \right)  \, dr \, dS_\vartheta\\[0.15in]
 &=  &   2n
 \, \int_{S^{n - 1}} \int_0^\infty     {H } \, y_1\,  r^{n-1}   \left[\,{{\lambda_c^n\cdot \left(\,\lambda_c^{\,2}  +  r^{\,2} - 2 \sum {\xi_c}_i \!\cdot \! y_i +|\,\xi_c|^{\,2} + r^{\,2} - |\,\xi_c|^{\,2} - \lambda_c^{\,2} \right)  } \over {(\,\lambda_c^{\,2} + |\,y - \xi_c|^{\,2} )^{n + 1} }} \right] \, dr \, dS_\vartheta\\[0.1in]
& = & 2n
 \, \int_{S^{n - 1}} \int_0^\infty     {H } \,  y_1\,  r^{n-1}  \left(\,{\lambda_c\over {\lambda_c^{\,2} + |\,y- \xi_c|^{\,2}}} \right)^n \, dr \, dS_\vartheta\\
 & \   & \ \ \ \  + \,2n
 \, \int_{S^{n - 1}} \int_0^\infty     {H } \, y_1\,  r^{\,2} \cdot r^{n-1}  \left(\,{\lambda_c^n\over {(\lambda_c^{\,2} + |\,y- \xi_c|^{\,2})^{n + 1}}} \right)  \, dr \, dS_\vartheta\\
 & \ & \ \ \ \ \ \ \ \  - 2n \,(|\,\xi_c|^{\,2} + \lambda_c^{\,2}) \, \int_{S^{n - 1}} \int_0^\infty     {H } \,  y_1\,   \cdot r^{n-1} \left(\,{\lambda_c^n\over {(\lambda_c^{\,2} + |\,y- \xi_c|^{\,2})^{n + 1}}} \right)  \, dr \, dS_\vartheta\,.
\end{eqnarray*}

Hence
\begin{eqnarray*}
(A.10.23)  & \ & 2
 \, \int_{\R^n} y_1\, \left[ r  {{\partial H}\over {\partial r }}\right] \, \left(\,{\lambda_c\over {\lambda_c^{\,2} + |\,y- \xi_c|^{\,2}}} \right)^n \,\\[0.15in]
 & = & -\, 2
 \, \int_{S^{n - 1}} \int_0^\infty     {H } \,  y_1\, r^{n-1}   \left(\,{\lambda_c\over {\lambda_c^{\,2} + |\,y- \xi_c|^{\,2}}} \right)^n \, dr \, dS_\vartheta\\
 & \   & \ \ \ \  + 2n
 \, \int_{S^{n - 1}} \int_0^\infty     {H } \,  r^{\,2} \cdot r^{n-1} y_1\,  \left(\,{\lambda_c^n\over {(\lambda_c^{\,2} + |\,y- \xi_c|^{\,2})^{\,n + 1}}} \right)  \, dr \, dS_\vartheta\\
 & \ & \ \ \ \ \ \ \ \  - 2n\, (|\,\xi_c|^{\,2} + \lambda_c^{\,2}) \, \int_{S^{n - 1}} \int_0^\infty     {H } \,    y_1\,  r^{n-1} \left(\,{\lambda_c^n\over {(\lambda_c^{\,2} + |\,y- \xi_c|^{\,2})^{\,n + 1}}} \right)  \, dr \, dS_\vartheta\\[0.15in]
 & = & -\, 2
 \, \int_{\R^n}    {H } \,   y_1\,  \left(\,{\lambda_c\over {\lambda_c^{\,2} + |\,y- \xi_c|^{\,2}}} \right)^n \,\\
 & \   & \ \ \ \  + 2n
 \, \int_{\R^n}     {H } \,  r^{\,2}  ( y_1 - {\xi_c}_1)\,  \left(\,{\lambda_c^n\over {[\,\lambda_c^{\,2} + |\,y- \xi_c|^{\,2}]^{\,n + 1}}} \right)  \,\\
 & \ & \ \ \ \ \ \ \ \  - 2n \,(|\,\xi_c|^{\,2} + \lambda_c^{\,2}) \, \int_{\R^n}    {H } \cdot    (y_1 -{\xi_c}_1)\,  \left(\,{\lambda_c^n\over {[\,\lambda_c^{\,2} + |\,y- \xi_c|^{\,2}]^{\,n + 1}}} \right)  \,\\
 & \   & \ \ \ \  + 2n
 \, \int_{\R^n}     {H } \,  r^{\,2}  \,  {\xi_c}_1 \,  \left(\,{\lambda_c^n\over {[\,\lambda_c^{\,2} + |\,y- \xi_c|^{\,2}]^{n + 1}}} \right)  \,\\
 & \ & \ \ \ \ \ \ \ \  - 2n \,(|\,\xi_c|^{\,2} + \lambda_c^{\,2}) \, \int_{\R^n}    {H } \,    {\xi_c}_1\,  \left(\,{\lambda_c^n\over {[\,\lambda_c^{\,2} + |\,y- \xi_c|^{\,2}]^{\,n + 1}}} \right)  \,\,.
\end{eqnarray*}
The third term above is equal to zero via (A.10.13). The first two terms cancel with the two terms in (A.10.21). Thus we are required to show the sum of the last two terms in (A.10.23) is equal to zero.\\[0.1in]
(A.10.24)
\begin{eqnarray*}
& \ & {\mbox{The \ \ \ \ sum \ \   of \ \ the \ \ last \ \ two \ \ terms \ \ in  \ \ (A.10.23)}}\\[0.1in]
 &=  & 2n\,{\xi_c}_1\,  \int_{\R^n}    {H } \, \left[\,r^{\,2} - |\,\xi_c|^{\,2} - \lambda_c^{\,2}\right] \,   \left(\,{\lambda_c^n\over {(\lambda_c^{\,2} + |\,y- \xi_c|^{\,2})^{n + 1}}} \right)  \,\\[0.1in]
 &=  & 2n\,{\xi_c}_1\,  \int_{\R^n}    {H } \, \left[\,(\lambda_c^{\,2} + |\, y - \xi_c|^{\,2}) - 2\,\lambda_c^{\,2}   - 2\,|\,\xi_c|^{\,2} + 2\sum y_i \cdot {\xi_c}_i \right] \,   \left(\,{\lambda_c^n\over {(\lambda_c^{\,2} + |\,y- \xi_c|^{\,2})^{n + 1}}} \right)  \,\\[0.1in]
 & = & 2n\,{\xi_c}_1\,  \lambda_c \left[\int_{\R^n}    {H } \, {\lambda_c^{n-1}\over {(\lambda_c^{\,2} + |\,y- \xi_c|^{\,2})^{n}}}  \  - \ 2 \int_{\R^n}    {H } \, {\lambda_c^{n+1}\over {(\lambda_c^{\,2} + |\,y- \xi_c|^{\,2})^{n +1}}}  \  \right]\\[0.075in]
 & \ & \ \ \ +\,2n\,{\xi_c}_1\,   \left[ \int_{\R^n}    {H } \, { {\lambda_c^n [\,2\,\sum y_i \cdot {\xi_c}_i] }\over {(\lambda_c^{\,2} + |\,y- \xi_c|^{\,2})^{n+1}}}  \  - \ \int_{\R^n}    {H } \, {{\lambda_c^n\, [\, 2 |\,\xi_c|^{\,2}]}\over {(\lambda_c^{\,2} + |\,y- \xi_c|^{\,2})^{n +1}}}  \   \right]\\[0.1in]
 & = & 2n\,{\xi_c}_1\,   \int_{\R^n}    {H } \, { {\lambda_c^n [\,2\,\sum (y_i - {\xi_c}_i)\cdot {\xi_c}_i] }\over {(\lambda_c^{\,2} + |\,y- \xi_c|^{\,2})^{n+1}}}  \   = \ 2n\,{\xi_c}_1\,  \sum {\xi_c}_i \,\int_{\R^n}    {H } \, { {\lambda_c^n [\, y_i - {\xi_c}_i] }\over {(\lambda_c^{\,2} + |\,y- \xi_c|^{\,2})^{n +1}}}  \  \\[0.075in]
 & = & 0\,.
\end{eqnarray*}
Here we use (A.10.13).\, Combining (A.10.21)\,--\,(A.10.24), we obtain (A.10.20)\,.\\[0.2in]
{\it Rotations.} \ \
Let $\Theta_t$ be a family of rotations so that
$$
{{d \Theta_t}\over {dt}} \bigg\vert_{\,t \,= \,0 } \ = \  X_k \ \ \ \ \mfor \  k \in {\N} \ \ {\mbox{with}}  \ = \  n + 2 \,\le\, k \,\le\,  (n + 1) (n + 2)/2\,.
$$
As in (A.10.4)\,,\, we have
\begin{eqnarray*}
\int_{S^n} X_k (K) \,f^{{n + 2}\over {n - 2}} dV_{g_1} & = & {{d }\over {dt}} \int_{S^n}  {\cal K} (\Theta_t) \,f^{{n + 2}\over {n - 2}} \, dV_{g_1} \bigg\vert_{t = 0 } \ = \  {{d }\over {dt}} \int_{S^n}  K \,[\,f(\Theta_t^{-1}) ]^{{n + 2}\over {n - 2}} \,dV_{g_1} \bigg\vert_{\,t = 0 }\\
& = & {{d }\over {dt}}\int_{\R^n} {H} \left(\,{{\lambda_t}\over {\lambda_t^{\,2} + |y - |\,\xi_t|^{\,2}}} \right)^n\,\,\bigg\vert_{\,t = 0 } \ = \  0\\
 & \ & \ \ \ \ \ \ \ \ \ \ \ \   \ \ \ \ \ \ \ \ \ \ \ \ \ \ \ \ \   \ \ \ \ \  \ \ \ \ \ \ \ \ \ \ \ \  (\lambda_o \ = \  \lambda_c\,, \ \ \xi_o \ = \  \xi_c)
\end{eqnarray*}
via the chain rule and (A.10.12) and (A.10.13).       \qed
%
%
%
%
Here we point out  that, in the special case
  $(\,\lambda_c,\,\, \xi_c)\ = \ (1, \ {\vec{\,0}}\,)\,,\,$
\begin{eqnarray*}
& \ & ``\,(6.5) \ \ \& \  \,(6.6)\ \   ({\mbox{Part I}})\,"   \Longrightarrow   \ \ \int_{S^n} {\cal H} (x) \, x_\ell \ dS_x\ = \ 0\\[0.1in]
& \ & \ \ \ \ \Longrightarrow \ \ \int_{S^n} {\cal H} (x) \, \Delta_1 \,(x_\ell) \,dS_x\ = \ 0
 \ \ \Longrightarrow  \ \int_{S^n} \langle \btd \,{\cal H} (x)\,,\,\,  \btd_{g_1}  (x_\ell) \rangle\,dS_x\ = \ 0\\[0.1in]
 & \ & \ \ \ \ \Longrightarrow \ \,  \int_{S^n} X_\ell\, ({\cal H})  \,dS_x\ = \ 0
\end{eqnarray*}
for $\,\ell\ = \ 0, \ 1, \cdot \cdot \cdot\,, \ n$\,.\, Here $\,X_\ell\ = \ \btd_{g_1} \,x_\ell\,$ is the conformal Killing vector field on $S^n$ generating the dilations (\,see \cite{Han-Li})\,.\, Observe that when  $\,(\,\lambda_c\,, \  \xi_c)\ = \ (1\,, \ {\vec{\,0}}\,)\,,\,$ $\,\psi \,\equiv\, 1\,$ as described in Theorem 6.10 (Part I).\bk
We conclude this section with the following (well-known) remark.\\[0.2in]
{\bf Proposition A.10.25.} \ \ {\it Suppose\,}  ${\cal K} \in C^1 (S^n)\,$  {\it satisfy the K\,-W condition and $\Phi\,: S^n \to S^n$ is a conformal transformation. then ${\cal K} \circ \Phi$ also satisfies the K\,-W condition\,.}\\[0.2in]
{\bf Proof.} \ \ We first note that $\Phi$ induces an isomorphism on the collection of all conformal Killing vector field. Suppose $\phi_t \,: S^n \to S^n$ is a parameter of conformal transformations which generate $X$.\, That is,
$$
{{d\,\phi_t}\over {dt}} \bigg\vert_{\,t \,=\, 0} \ = \  X\,.
$$
This results
$$
\int_{S^n}  X ({\cal K}) \, f^{{2n}\over {n - 2}} \, dV_{g_1}  = \left[ {{d}\over {dt}} \int_{S^n}  [ {\cal K}\circ \phi_t\,] \, f^{{2n}\over {n - 2}} \, dV_{g_1} \right]_{\,t \,=\, 0}\,.
$$
Then
\begin{eqnarray*}
\int_{S^n}  X ({\cal K} \circ \,\Phi) \, f^{{2n}\over {n - 2}} \, dV_{g_1}  & = &  \left[\, {{d}\over {dt}} \int_{S^n}  [ {\cal K} \circ \Phi \circ \phi_t\,] \, f^{{2n}\over {n - 2}} \, dV_{g_1} \right]_{\,t\, =\, 0}\\[0.1in]
& =&  \int_{S^n}  \tilde X ({\cal K} \circ\, \Phi) \, f^{{2n}\over {n - 2}} \, dV_{g_1} = 0\,,
\end{eqnarray*}
where
$$
\tilde X \ = \ {{d\,[\Phi \circ \phi_t\,]}\over {dt}} \bigg\vert_{\,t \,=\, 0}\,.
$$

\vspace*{0.1in}

Hence ${\cal K} \circ \,\Phi$ also satisfies the K\,-W condition.\qed
%
%
%
%
%
 \vspace*{0.1in}


{\bf   \S\,A.\,11. \ \  Deriving (7.14) and (7.15)\,,\, Part I.} \\[0.1in]
Following (7.3) and (7.4) in Part I\,,\, a direct calculation shows that
\begin{eqnarray*}
 \ \ \ \ \ \ \ {{\partial^2 {G_{|_{\bf Z}} }}\over {\partial \lambda^2 }} \,(1, \ {\vec{\,0}}\,) & = &   - 2\,n \, (n + 2) \,{\bar c}_{-1} \int_{\R^n}  {{ H (\,y) }\over {\left(\,1 + |\,y   |^{\,2}\right)^{n+1}}}\, \,\,\\[0.075in]
 & \ & \ \ \ \ \ \ \ +  \,4n\,(n + 1)  \,{\bar c}_{-1} \int_{\R^n}
 \,  {{ H (\,y) }\over {\left(\,1 + |\,y  |^{\,2}\right)^{n+2}}}\,
 \,\,\\[0.15in]
 & = & -\,{{n\,(n + 2)}\over {2^n}}\,{\bar c}_{-1}\,\int_{S^n}
 (1 - x_{n + 1}) \,H(x) \, dS_x\\[0.075in]
 & \ & \ \ \ \ \ \ \ \ \ +\,\,{{n\,(n + 1)}\over {2^n}}\,{\bar c}_{-1}\,\int_{S^n}
 (1 - x_{n + 1})^2 \,H(x) \, dS_x\\[0.15in]
 & = & {{n\,(n + 1)}\over {2^n}}\,{\bar c}_{-1} \left[\,   \int_{S^n}
  x_{n + 1}^2 \,H(x) \, dS_x - {1\over {n + 1}} \int_{S^n}
   \,H(x) \, dS_x\right]\\
   & \ & \\
   & \ & \ \ \ \ \ \ \ \ \ \ \ \ \ \ \ \ \ \ \ \  \ \ \ \  \ \ \ \  \ \ \
   \left(\,{\mbox{using}} \int_{S^n}
  x_{n + 1}  \,H(x) \, dS_x \,=\, 0\right)\,,
\end{eqnarray*}
In addition,
\begin{eqnarray*}
 \ \ \ \ \ \ \ \ {{\partial^2 {G_{|_{\bf Z}} }}\over {\partial \xi_\ell^2}}\,
(1, \ {\vec{\,0}}\,) & = & -\,2 n\,{\bar c}_{-1} \int_{\R^n}  {{H (\,y)   }\over {\left(\,1+ |\,y |^{\,2}\right)^{n+1}}}\  \,\\
& \  & \ \ \ \ \ +\,4 n\,(n + 1)\,{\bar c}_{-1}
\int_{\R^n} H (\,y)   {{  ( y_\ell)^2}\over
{\left(\,1 + |\,y |^{\,2}\right)^{n+2}}}\  \,\\[0.075in]
  =   -\,{{n}\over {2^n}}\,{\bar c}_{-1}\,\int_{S^n} (1\!\!\! &-&\!\!\! x_{n + 1})
\,H(x) \, dS_x \ + \ {{n\,(n + 1)}\over {2^n}}\,{\bar c}_{-1}\,\int_{S^n} (x_\ell )^2  \,H(x) \, dS_x \ \ \ \ \ \ \ \ \ \ \ \ \ \ \ \
\\[0.075in]
  =   -\,{{n}\over {2^n}}\,{\bar c}_{-1}\,\int_{S^n}& \ &
\!\!\!\!\!\!\!\!\!\!\!\!\!\!H(x) \, dS_x\  + \  {{n\,(n + 1)}\over {2^n}}\ {\bar c}_{-1}\,\int_{S^n} (x_\ell )^2
\,H(x) \, dS_x\,.
\end{eqnarray*}
Here we use
$$
x_{n + 1} \ = \ {{r^2 - 1}\over {r^2 + 1}} \ \ \ \ \Longrightarrow
\ \ {1\over {1 + r^2}}\  = \ {{1 - x_{n + 1}}\over 2}\,.
$$

\newpage

{\bf   \S\,A.\,12. \ \  $\lambda_M =
\sqrt{(t + \Delta)(t - \Delta) \,}\ \  $ is the only critical point in expression }\\[0.05in]
\hspace*{0.7in}{\bf  \,(2.12), Part II.}\\[0.1in]
Refer to paragraph proceeding the proof of Lemma 2.21, Part II.  The argument is based on the following.\\[0.1in]
{\bf (i)} \ \ The function $\displaystyle{\, {a \over {a^{\,2} + \lambda^{\,2}}}\, }$ is decreasing on $\lambda\,;$ it is decreasing in $a$ when $a \,>\, \lambda\,.$\\[0.1in]
{\bf (ii)}  \ \,$\displaystyle{\, {{t + \Delta} \over \lambda } \ > \ {{t -\Delta } \over \lambda } \ \Longrightarrow \ {\bar \phi}_+ \ > \ {\bar \phi}_-\,. }$\,\, See (2.13) in Part II\,. \\[0.1in]
We divide the argument into three parts.\\[0.1in]
{\bf (a)} \ \ When $\displaystyle{\, \lambda < t - \Delta\,.\, }$ We have
$$
{{t + \Delta} \over \lambda } \ > \ {{t -\Delta } \over \lambda } > 1 \ \Longrightarrow \ \pi \ > \ {\bar \phi}_+ \ > \ {\bar \phi}_- > {\pi\over 2} \ \Longrightarrow \ \sin {\bar \phi}_+ \ <\  \sin {\bar \phi}_-\,.
$$
Thus we have
\begin{eqnarray*}
\lambda \ < \ t - \Delta \ \ \Longrightarrow   &\ & -\, {{2\,[\,t + \Delta]}\over {[\,{t + \Delta}]^{\,2} + \lambda^{\,2}}}\  [\,\sin {\bar \phi}_+]^{n - 1} + {{2\,[\,t - \Delta]}\over {[\,{t - \Delta}]^{\,2} + \lambda^{\,2}}}\ [\,\sin {\bar \phi}_-]^{n - 1} > 0\\
& \ &  \ \ \ \ \ \uparrow \ \ ({\mbox{smaller}}) \ \ \ \uparrow \ \ ({\mbox{smaller}})\\
&\ &\\
\Longrightarrow   &\ &  {\partial\over {\partial \lambda}}
\int^{2\arctan \,\left(\, {{t+ \Delta }\over {\lambda}}
\right)}_{2\arctan \,\left(\, {{t- \Delta }\over {\lambda}} \right)}
\,\,\,[\,\sin \varphi]^{n - 1}\,   \, d\varphi\ >\ 0  \ \ \ \ \  [\,{\mbox{via}} \ \ (2.13)\,, \ {\mbox{Part \ II}}\,]\,.
\end{eqnarray*}

{\bf (b)} \ \ Likewise,
\begin{eqnarray*}
\lambda > t + \Delta   \Longrightarrow   &\ & \!\!\!\!\!-\, {{2\,[\,t + \Delta]}\over {[\,{t + \Delta}]^{\,2} + \lambda^{\,2}}}\ [\,\sin {\bar \phi}_+]^{\,n - 1} + {{2\,[\,t - \Delta]}\over {[\,{t - \Delta}]^{\,2} + \lambda^{\,2}}}\ [\,\sin {\bar \phi}_-]^{\,n - 1} \ < \ 0\\
& \ &  \ \ \ \ \ \uparrow \ \ ({\mbox{bigger}}) \ \ \ \uparrow \ \ ({\mbox{bigger}})\\
& \ & \\
\Longrightarrow   &\ &  {\partial\over {\partial \lambda}}
\int^{2\arctan \,\left(\, {{t+ \Delta }\over {\lambda}}
\right)}_{2\arctan \,\left(\, {{t- \Delta }\over {\lambda}} \right)}
\,\,\,[\,\sin \varphi]^{\,n - 1}\,   \, d\varphi \ < \ 0 \ \ \ \ \  [\,{\mbox{via}} \ \ (2.13)\,, \ {\mbox{Part \ II}}\,]\,.
\end{eqnarray*}

{\bf (c)} \ \ {\it Cross over}\,. \ \ When $\displaystyle{\,t + \Delta >  \lambda > t - \Delta\,.\, }$ We have
$$
0 \ < \ {\bar \phi}_- \ < \ {\pi\over 2} \ \ \ \ {\mbox{and}} \ \ \ \ {\pi\over 2}\  <\  {\bar \phi}_+ \ < \ \pi\,.
$$
Thus
$$
\sin {\bar \phi}_- \ \downarrow \  {\mbox{from \ }} 1\ \ \ \ {\mbox{and}} \ \ \ \  \sin {\bar \phi}_+ \ \uparrow \  {\mbox{to \ }} 1 \ \ \ \ {\mbox{when}} \ \ \lambda \ \ {\mbox{increases \ in \ this \ range}}\,.
$$

\newpage

While
\begin{eqnarray*}
& \ & {{[\,t + \Delta]}\over {[\,{t + \Delta}]^{\,2} + \lambda^{\,2}}} \ < \ {{[\,t - \Delta]}\over {[\,{t - \Delta}]^{\,2} + \lambda^{\,2}}} \ \ \ \ {\mbox{when}} \ \ \lambda\ = \ t - \Delta\\
\Longrightarrow   &\ &  {\partial\over {\partial \lambda}}
\int^{2\arctan \,\left(\, {{t+ \Delta }\over {\lambda}}
\right)}_{2\arctan \,\left(\, {{t- \Delta }\over {\lambda}} \right)}
\,\,\,[\,\sin \varphi]^{n - 1}\,   \, d\varphi\ <\  0\\
 & \ & \ \ \ \ \ \ \ \ \ \ \ \ \  \ \ \ \ \ \ \ \  \ \ \  \ \ \ \ \ \ \ \ \ \ \ \ \  \ \ \ \ \ \ \  [\,{\mbox{via}} \ \ (2.13)\,, \ {\mbox{Part \ II}}\,]\,. \\
& \ & \\
& \ & {{[\,t + \Delta]}\over {[\,{t + \Delta}]^{\,2} \ + \ \lambda^{\,2}}} \ > \ {{[\,t - \Delta]}\over {[\,{t - \Delta}]^{\,2} + \lambda^{\,2}}} \ \ \ \ {\mbox{when}} \ \ \lambda \ =\  t + \Delta\\
\Longrightarrow   &\ &  {\partial\over {\partial \lambda}}
\int^{2\arctan \,\left(\, {{t+ \Delta }\over {\lambda}}
\right)}_{2\arctan \,\left(\, {{t- \Delta }\over {\lambda}} \right)}
\,\,\,[\,\sin \varphi]^{\,n - 1}\,   \, d\varphi \ > \ 0 \ \ \ \ \ \  [\,{\mbox{via}} \ \ (2.13)\,, \ {\mbox{Part \ II}}\,]\,.\\
\end{eqnarray*}
$$
{\mbox{Moreover,}} \ \ \ \ \ \left[\,{{[\,t\, +\, \Delta]}\over {[\,{t \,+\, \Delta}]^{\,2}
+ \lambda^{\,2}}} \ -\  {{[\,t \,-\, \Delta]}\over {[\,{t \,- \, \Delta}]^{\,2} \,+ \,\lambda^{\,2}}}
\right] \ \  \uparrow \ \ {\mbox{in}} \ \ \lambda \ \ \uparrow  \  (\,> 0)\,. \ \ \ \ \  \ \ \ \ \ \ \ \ \ \  \ \ \ \ \ \ \ \ \ \  \ \ \ \ \ \ \ \ \ \  \ \ \ \ \ \ \ \ \ \  \ \ \ \
$$

\vspace*{0.2in}

Thus we see that there can only be one critical point in $\,(\,t \,-\, \Delta\,, \ \,t \,+\, \Delta)\,.$

\vspace*{0.4in}

{\bf   \S\,A.\,13. \ \  Deriving (2.39) in Part II.}\\[0.1in]
 From (2.38), Part II, we have
\begin{eqnarray*}
  (2 \varrho)^2
 & = &   \left[ {{ 2 (R_\lambda+\delta_\lambda)}\over { (R_\lambda + \delta_\lambda)^2 + 1}} \,\,-\,\,
\left(\, -\,{{ 2 (R_\lambda - \delta_\lambda)}\over { (R_\lambda \,-\, \delta_\lambda)^2 + 1}} \right)\right]^2\\
 & \ & \ \ \ \ \ \ \ \  \ \ \ \ \ \ \ \ \   +\  \left[  {{ (R_\lambda + \delta_\lambda)^2 - 1}\over { (R_\lambda + \delta_\lambda)^2 + 1}}\,\, -\,\,  {{ (R_\lambda - \delta_\lambda)^2 - 1}\over { (R_\lambda - \delta_\lambda)^2 + 1}}    \right]^2\\
& = & {1\over {[(R_\lambda\, + \,\delta_\lambda)^2 \,+ \,1]^2 \,[(R_\lambda \,-\, \delta_\lambda)^2 \,+\, 1]^2}} \times\\
& \ &   \!\!\!\!\!\!\times\left(\, \left\{ 2 (R_\lambda + \delta_\lambda) [(R_\lambda \,- \,\delta_\lambda)^2\, + \,1] \,\,+ \,\,2 (\,R_\lambda \,-\, \delta_\lambda) [\,(R_\lambda\, + \,\delta_\lambda)^2 \,\,+\,\, 1]\right\}^2
\right. \ \ \ \ \ \\
& \ &  \left. \ \ \ \ \ \  \ \ \ \ \ \ \  \  \ \ \ \ \ \ \ \ \    +\ [\,(\,{\mathcal{A}}  - 1)(\,{\mathcal{B}}  + 1) - (\,{\mathcal{B}} - 1) (\,{\mathcal{A}}  + 1)]^2\right)\\
& = & {{  4\,[ \,2R_\lambda + (R_\lambda + \delta_\lambda)(R_\lambda - \delta_\lambda)(2R_\lambda)]^2\,\,+ \,\, [\,2(\,{\mathcal{A}}  - \,{\mathcal{B}})]^2        }\over {[(R_\lambda + \delta_\lambda)^2 + 1]^2 \,[(R_\lambda - \delta_\lambda)^2 + 1]^2}}   \\
& = & {{4 \cdot
\left(\,\, [ 2R_\lambda (1 + R_\lambda^2 - \delta_\lambda^2)]^{\,2}\,\,+\,\,
[(R_\lambda + \delta_\lambda)^2 - (R_\lambda - \delta_\lambda)^2]^{\,2}\right)      } \over {[(R_\lambda + \delta_\lambda)^2 + 1]^2 \,[(R_\lambda - \delta_\lambda)^2 + 1]^2}}\,.
\end{eqnarray*}
In the above,
$$
{\mathcal{A}} \,\,= \,\,(R_\lambda \, +\, \delta_\lambda)^2 \ \ \ \ {\mbox{and}} \ \ \ \ {\mathcal{B}} \,\, = \,\,(R_\lambda  \,-\, \delta_\lambda)^2.
$$
It follows that
$$
  (2 \varrho)^2 \ = \ {{ (4R_\lambda)^2\, [\, (1 + R_\lambda^2 - \delta_\lambda^2)^2\,\,+\,\,  4\,\delta_\lambda^2\,]}\over
{[(R_\lambda + \delta_\lambda)^2 \, +\,  1]^2 \,\,[(R_\lambda - \delta_\lambda)^2 \, + \, 1]^2}}\,\,.
$$

\vspace*{0.4in}

{\bf   \S A.\,14.} \ \ {\bf Back to $S^n$ \,(refer to \S\,3\,j in Part II)\,.}\\[0.2in]
{\bf Lemma A.14.1.} \ \ {\it Given a positive function ${\cal K} \in C^{1, \,\, \alpha} \,(S^n)$\,,\, let $K (y)= {\cal K} \circ {\dot{\cal P}}^{-1} (y)$ be defined for\,} $y \in \R^n$.\, {\it Suppose\,} $v \in {\cal D}^{1, \,\, 2} \,\cap \, C^2 (\R^n)$    {\it is a positive solution of}
$$
\Delta_o\, v \,+ \,({\tilde c}_n \,K)\, v^{{n + 2}\over {n - 2}} \ =\  0 \ \ \ \ \ \ {\it{in}} \ \ \ \R^n.
$$
{\it Set }
$$
u \,(x)\ =\ v\,(y) \cdot \left(\, {{1 + |y|^{\,2}}\over 2} \right)^{{n - 2}\over 2} \ \ \ \ \ \ {\it for} \ \  x \,= \, {\dot{\cal P}}^{-1} (y) \in S^n \setminus \{ {\bf N} \}\,.\leqno (A.14.2)
$$
{\it Then $u$ has a removable singularity at north pole $\,{\bf N}\,,\,$  and satisfies the equation }
$$
\Delta_1\, u - \,{\tilde c}_n\,n \,(n - 1) + ({\tilde c}_n \,{\cal K}) \,u^{{n + 2}\over {n - 2}} \ =\  0 \ \ \ \ {\it in} \ \ S^n. \leqno (A.14.3)
$$

\vspace*{0.15in}

{\bf Proof.} \ \ As $\,v \in {\cal D}^{1\, \, 2}$\,,\, via the Sobolev inequality (\,see for instance (2.9) in Part I \cite{I}\,)\,,\, there exists a positive number
$C$ such that
$$\int_{|y| \ge 1} | \btd v (y) |^{\,2} \,\, + \ \int_{|y| \ge 1} v^{{2n}\over {n - 2}}
\,(y) \,\, \le C\,. \leqno (A.14.4)$$
By using the  Kelvin transform
$$
\tilde y \ = \ {y\over {|y|^{\,2}}} \ \ \ \ \ \ {\mbox{for}} \ \ \ |\,y| \,\ge\, 1\,,$$  we let
$$\tilde v \,(\tilde y) \,\,=\,\,  v (y) \cdot  {1\over {|\tilde y|^{n - 2} }}  \ \ \ \ \ \ {\mbox{with}} \ \ \ \tilde y \ = \ {y\over
{|y|^{\,2}}}\,. \leqno (A.14.5)$$ The equation
$$\Delta_o  \,v (y) \,+\, ({\tilde c}_n\,K) \,v^{{n + 2}\over {n - 2}} (y) \,\,= \,\,0  \ \ \ \
{\mbox{for}} \ \ |\,\tilde y|
\ge 1$$ is transformed into the equation
$$\Delta_o \, \tilde v \,(\tilde y) +  \left[ {\tilde c}_n\cdot K\,\left(\, {{\tilde y}\over {|{\tilde y}|^{\,2}}} \right)\right]\, {\tilde v}^{{n + 2}\over {n -
2}} (\tilde y)\ =\ 0
\
\
\
\ {\mbox{for}}
\
\ 0 \,< \,|\,y|
\,\le \,1\,.$$ Furthermore, from (A.14.4), there exists a positive number $C_1$ such
that
$$
\int_{|\tilde y| \le 1} | \btd \tilde v \,(\tilde y) |^{\,2} \,d\tilde y \ + \ \int_{|\tilde y|\, \le\, 1} {\tilde v}^{{2n}\over {n - 2}}
\,(\tilde y) \,d\tilde y \ \le \  C_1\,.
$$
By a result of Brezis and Kato \cite{Brezis-Kato}, we
have
$\tilde v
\,\in\, L^\infty (B_o ({1\over 2}))$.
It follows that
$$
v \,(y) \,\,\le \,\,{{C_2}\over {|\,y|^{n - 2}}} \ \ \ \ \mfor  \ \  |\,y|\,\gg\, 1\,.
$$
That is, $\,u\,$ is bounded in a deleted neighborhood of $\,{\bf N}$\,.\, Standard elliptic theory can be applied to show  that $u$ has a removable singularity at $\,{\bf N}\,$,\, and via the standard conformal transform (A.14.2), we show that equation (A.14.3) is fulfilled. \qed
%
%

%
%



\begin{center} {\LARGE   {\bf Construction of Blow-up Sequences for  the}}
\medskip \medskip \smallskip \smallskip \\ {\LARGE {\bf  Prescribed Scalar Curvature Equation}} \medskip \medskip \smallskip \smallskip\\ {\LARGE {\bf  on $S^n$.\,   I. \,Uniform Cancelation}}

\vspace{0.43in}

{\Large {Man Chun  {\LARGE L}}}{EUNG}\\[0.165in]

{\large {National University of Singapore}}${\,}^1$\\[0.1in]
{\tt matlmc@nus.edu.sg}\\
\end{center}
\vspace{0.33in}
\begin{abstract}
\vspace{0.25in}
\noindent For $\,n \ge 6\,,\,$ using the Lyapunov\,-\,Schmidt reduction method, we  describe how to construct (scalar curvature) functions on $S^n\,$,\, so that each of them enables the conformal scalar curvature equation to have an  infinite number of positive solutions, which form a blow-up sequence.   The prescribed  scalar curvature function  is shown  to have  $\,C^{\,n - 1\,, \, \, \beta}$\, smoothness. We present the argument in two parts. In this first part, we discuss the uniform cancelation property in the Lyapunov\,-\,Schmidt reduction method for the scalar curvature equation. We also explore relation between the Kazdan\,-Warner condition and the first order derivatives of the reduced functional, and symmetry in the second order derivatives of the reduced functional.

\end{abstract}

\vspace*{1.8in}

{\bf Key Words}\,:      Scalar Curvature Equation\,;\, Blow-up\,;\, Critical Points\,;\, Sobolev Spaces. \\[0.1in]
{\bf 2000 AMS MS Classification}\,: Primary 35J60\,;
\ \ Secondary 53C21. \\[0.01in]
%
%
%
%

$\overline{\ \ \ \ \ \ \ \ \ \ \ \ \ \ \ \ \ \ \ \ \  \ \ \ \ \ \ \ \ \ \ \ \ \ \ \ \ \ \ \ \ \ \ \ }$\\[0.05in]
\noindent{$\!\!{\!}^* {\small{\mbox{\, Department of Mathematics, National University of Singapore, 10, Lower Kent Ridge Rd.,}}}$}\\[0.03in]
 {\small{\mbox{\ \  \,Singapore 119076, Republic of Singapore\,.}}}

\newpage

{\bf \large {\bf  \S\,1. \ \ Introduction.}}\\[0.2in]
We apply the   Lyapunov\,-\,Schmidt reduction method to find non-constant function  $\,{\cal K}\,$ so that the equation
$$
\Delta_1 \, u \,-\, {\tilde c}_n \,n \,(n - 1) \,u \,+ \,({\tilde c}_n \,{\cal K})\, u^{{n + 2}\over {n - 2}} \ =\  0 \ \ \ \ {\mbox{in}} \ \ \,S^n \leqno (1.1)
$$
has infinite number of positive solutions $\,\{ u_i \,\}_{i = 1}^\infty$\,,\, which compose  a blow-up sequence of  solutions. We refer to \S\,1\,c for the (rather standard) notations we use.  \bk
 The blow\,-\,up phenomenon is endogenous to equation (1.1), because of the critical Sobolev exponent linked to $\,(n + 2)/(n - 2)\,$ (cf. \S\,2\,b.1)\,.\, It is systematically studied (see for examples  \cite{Chen-Lin-2} \cite{Chen-Lin-3}  \cite{Xu-Invent} \cite{Compact} \cite{Leung-blow-up} \cite{Leung-Supported} \cite{Li-1} \cite{Li-2}  \cite{Schoen-Notes})\,.\, When $\,{\cal K}\,$ equals to a positive constant, say, $\,(\,{\tilde c}_n\,{\cal K}) \,\equiv\, n\,(n - 2)$\,,\, equation (1.1) has a family of positive solutions (via the stereographic projection $\,\dot{\cal P}\,$ which sends the north pole $\,{\bf N}\,$ to $\,\infty\,$\,)\,:\\[0.01in]
 $$
 \left\{  \,  U_{\lambda, \  \xi} \,(x) \ := \   \left( {\lambda\over {\lambda^2 + |\,y - \xi|^{\,2}}} \right)^{{n-2}\over 2} \!\cdot \, \left( {{1 + |\,y|^2}\over 2} \right)^{{n - 2}\over 2}\ \ \ \ \  \mbox{with} \ \ (\,\lambda, \ \xi)\in \R^+\times \R^n \,\right\}.
$$

\vspace*{0.05in}

Here $\,y = \dot{\cal P}\, (x)$\, for $\,x \in S^n \setminus {\bf N}$\,,\,\, and $\,U_{\lambda, \  \xi} \,({\bf N})\, := \, [\,\lambda/2]^{\,{{n - 2}\over 2}}\,.\,\,$
 In particular, any sequence $\,\{U_{\lambda_i, \  \xi_i}\}_{i = 1}^\infty\,$ with $\,\lambda_i \to 0^+\,$ \,and\,\, $\xi_i \to \vec{\,0}\,$ is a blow\,-\,up sequence. The simple idea on perturbing $\,\{U_{\lambda, \  \xi}\}\,$ to form a blow\,-\,up sequence encounters a remarkable obstacle on stabilizing and fixing the scalar curvature function  ${\cal K}\,$ (non-\,constant this time)\,.\, \bk
  A significant development comes about when Ambrosetti, Berti and Malchiodi, using the Lyapunov\,-\,Schmidt reduction method,  construct many $C^k$\,-\,metrics $g$  on $S^n$ such that the Yamabe equation
$$
\Delta_g \, u \,-\, ({\tilde c}_n R_g) \,u \,+\, {\tilde c}_n \,n\,(n - 1) \,u^{{n + 2}\over {n - 2}} \ =\ 0 \ \ \ \ {\mbox{in}} \ \ S^n \leqno (1.2)
$$
has infinite number of positive solutions [\,here $g$ is non\,-\,conformal  to $g_1$ -- the standard metric on $S^n$\,;\, $R_g$ is the scalar curvature of $(S^n, \, g)$\,]\,.\,   See the renowned monograph \cite{Progress-Book}\,.\, This project stems from a study on these wonderful results and seek to apply them to the scalar curvature equation (1.1) so as to construct blow\,-\,up sequences.\bk
As we work through the Lyapunov\,-\,Schmidt method, we encounter two of the main features of equation (1.1)\,,\, namely, the balance and flexibility  (refer to \S\,6\,a\, and\, \S\,7\,)\,.\, The key object is the collection of solutions
$$
 (1.3)\ \ \ {\bf Z} \ :=\   \left\{ V_{\lambda, \  \xi} \in {\cal D}^{1\,,\, 2} \ \,\bigg\vert   \ \  V_{\lambda,  \  \xi} \,(y) \, :=\, \left( {\lambda\over {\lambda^2 + |\,y - \xi|^{\,2}}} \right)^{\!\!{{n-2}\over 2}} \, \  \mbox{with} \ \,(\,\lambda, \ \xi)\in \R^+\times \R^n \right\} \ \
 $$
of the equation
$$
\Delta_o \,v_o \ + \ n\, (n - 2) \,v_o^{{n + 2}\over {n - 2}}\  = \ 0 \ \ \ \ {\mbox{in}} \ \ \R^n.\leqno (1.4)
$$
`Sitting'\, inside the Hilbert space ${\cal D}^{1\,,\, 2}\,$ [\,see (2.7) for the definition]\,,\,  the pulled-back metric on ${\bf Z}\,$ is isometric (up to a rescaling) to the   hyperbolic metric on the upper half space $\R^+\times \R^n$ (this is shown in \S\,3\,b).\, We find it  curious that ${\bf Z}\,$  naturally carries such a   structure. This fact, although not directly linked to the main aim of the paper, provides valuable intuition on the vast amount of space for the solutions to be separated (see \S\,5\,a), and the correct order to be expected in the derivatives (see \S\,2   in \cite{II}\,,\, in particular,\, Lemma 2.65 in \cite{II}\,)\,.\,\bk
Via the stereographic projection $\dot{\cal P}$\, (see \S\,2\,a for more detail), equation (1.1) is transformed into

\vspace*{-0.3in}

$$
\Delta_o \,v + ({\tilde c}_n \,K)v^{{n + 2}\over {n - 2}}  \ = \ 0 \ \ \ \ \ \ \ {\mbox{in}} \ \ \ \ \R^n\,. \leqno (1.5)
$$
See (2.4) and (2.6) for the relations between $u$ and $v$, and between ${\cal K}$ and $K$\,.\, Let us write

\vspace*{-0.35in}

$$
K \ = \ 4n\, (n - 1) \ +\  \varepsilon H. \leqno (1.6)
$$
Note that $\,{\tilde c}_n \cdot  [\,4n \,(n - 1)\,] \ = \ n \,(n - 2)\,.\,$ Consider the `\,interaction' \,of \,$H$ with ${\bf Z}\,:$
$$
G_{|_{\bf Z}} \,({\bf z}) \ = \ {\bar c}_{-1} \int_{\R^n} H (y) \left[ {\lambda\over {\lambda^2 + |\,y - \xi|^{\,2}}} \,\right]^n \, \ \ \ \  \mfor {\bf z} \ =\  V_{\lambda,\,\ \xi} \in {\bf Z}\,.\leqno (1.7)
$$
Here\, ${\bar c}_{-1}$ is a {\it negative}\, constant defined in (2.15).
For $\,\varepsilon \in \R\,$ small, the question on finding a positive solution of equation (1.5) is reduced   to finding a {\it stable\,} critical point of $G_{|_{\bf Z}}\,.\ $ See \cite{Progress-Book}\,.\,
 This is discussed in \S\,2\,f.\bk
The balance condition is expressed by the   existence of a critical point for $\,G_{|_{\bf Z}}\,,\,$ and is linked to the Kazdan\,-Warner condition (\S\,6\,c).
To study   the flexibility part, we show that the Hessian matrix of $\,G_{|_{\bf Z}}\,$ at a critical point is always trace\,-\,free, thus can never be positive or negative definite. Instead of looking for strict local    maximum or minimum (as in the Yamabe equation), we seek {\it saddle points\,.\,} The contrast here is that one cannot `catch' a saddle point by looking at the values along the boundary and comparing them to the value of an inside point   (cf. the proofs of Theorem 2.16 in \cite{Progress-Book} and Proposition 24 in \cite{Brendle-1})\,.\bk
In \cite{Ambrosett-Scalar}\,,\, using the Lyapunov-\,Schimdt method on the scalar curvature equation (1.1), (finite number of) stable critical points of the $\,G_{|_{\bf Z}} $ are found by using a degree counting method.  The research is carried on in \cite{Ambrosetti-Malchiodi},  where the authors explore  symmetries  (see  also  \S\,7\,c). Each of these methods does not readily juxtapose   to produce a blow\,-\,up sequence of solutions. In this article and its sequel \cite{II}, we introduce annular domains (see   \S\,1\,a),\, and determine precisely the critical points of the {\it reduced functional\,} $\,G_{|_{\bf Z}}$\,,\, and show that the Hessian matrix at the critical point is non\,-\,degenerate (hence a stable critical point).
 By superimposing  concentric annular domains, and carefully estimating the gradient interference, we are able to find infinite number of stable critical points via degree theory for maps.

\vspace*{0.15in}


{\bf \S\,1\,a. \  } {\it Description of the main result.}\ \
Consider the sequence of  annular domains
$$
B_o \,\left( {{1+ \eta}\over {\bf a}}\right)\,\bigg\backslash \  \overline{B_o \,\left( {{1- \eta}\over {\bf a}}\right)\,} \,, \ \  \cdot \ \cdot \ \cdot\ , \ \ \ B_o \,\left( {{1+ \eta}\over {{\bf a}^k}}\right)\,\bigg\backslash \   \overline{B_o \,\left( {{1- \eta}\over {{\bf a}^k}}\right)\,}\ , \  \ \cdot\  \cdot\  \cdot\ .
$$
We describe the   numbers $\,{\bf a} \,> \,1\,$ and $\,\eta \,\in\, (0, \ 1)\,$  in (1.9)\,.\,
Choose a small positive number $\,\sigma\,$  and fix a number
  $\,\tau\, \in\, (n- 1\,, \ n)\,.\,$ Let $\,H^S\,$ be given by

  \vspace*{-0.075in}

$$
(1.8) \ \ \ \ \ \  H^S\, (y) \ = \ \cases{ {\displaystyle{1\over {{\bf a}^{\tau k} }} }   \ \ & {{\mbox{if}}} \ \ \   $\displaystyle{ y \in B_o \,\left( {{1+  \eta  }\over {{\bf a}^k}}\right)\,\bigg\backslash \  B_o \,\left( {{1-  \, \eta  }\over {{\bf a}^k}}\right)\,, \ \ \ \ \ k = 1, \ 2, \cdot \cdot \cdot\,,}\ \ \ \ \ \ \ \ \ $\cr
\ & \ \cr
0 & {\mbox if } \ \ $\displaystyle{ y \,\not\in \  {{\bigcup}_{k = 1}^\infty}   \left\{ B_o \,\left( {{1+ (\eta  + \sigma)}\over {{\bf a}^k}}\right)\,\bigg\backslash \, B_o \,\left( {{1- (\eta  + \sigma)}\over {{\bf a}^k}}\right)\right\}}$\,.}
$$

\vspace*{0.125in}

In between
$\ \,
B_o \,\left( {{1\, + \, (\, \eta  \, +\,  \sigma)}\over {{\bf a}^k}}\right) \big\backslash \,B_o \,\left( {{1\, +\,   \eta  }\over {{\bf a}^k}}\right)\, $\  \ \ \, and \ $\ B_o \,\left( {{1\,-\,  \eta  }\over {{\bf a}^k}}\right) \big\backslash \,B_o \,\left( {{1\, - \,  (\, \eta \,  +\,  \sigma)}\over {{\bf a}^k}}\right)\,,\
$
we properly  \\[0.1in]
smooth out $\,H^S\,$ so that $\ H^S \, \in\, C^{\, n-1, \ \beta} (\R^n)\,$\ \  (not necessarily in a rotationally symmetric manner\,,\, see \S\,3\,g in Part II \cite{II}\,)\,.\,
Together with Part II \cite{II}, We show the following result.

 \vspace*{0.15in}

 {\bf Main Theorem.} \ \   {\it For $\,n \ge 6\,$\,,\, let\,}  $H^S \in C^{\,n - 1\,, \ \beta} \,(\R^n)\,$ {\it be as described in }\,(1.8)\,,\, {\it  with the parameter $\,\tau\,$  and $\,\eta\,$ satisfying\,}

 \vspace*{-0.1in}

$$
\tau \,\in\, (n - 1\,, \ n)\,, \ \ \ \  1 - A^2 \ >\ \eta \ >\  B^2 \ >\  0  \ \ \ \ {\it{and}} \ \ \ \
  {{1 + \eta}\over {1 - \eta}} \ \le \   {5\over 2}\ .\leqno (1.9)
  $$
{\it Here $\,A\,$ and $\,B\,$ are positive numbers\,.\, There exist positive constants \,$C$\,,\, $C_1$\,,\, $c$ \,and\, $\varepsilon_o$ so that if the parameters \,${\bf a}$\, and $\,\sigma$ in\,} (1.8) {\it satisfy}
$$
{\bf a}\  > \ C^2 \ \ \ \ {\it{and}} \ \ \ \ 0 \ < \ \sigma\  <\  c^2\,,
$$

\vspace*{-0.25in}

{\it then the equation }
$$
\Delta_o \,v + {\tilde c}_n  \left[\, 4n\,(n - 1) \,+\, \varepsilon \, H^S\,\right]  v^{{n + 2}\over {n - 2}} \ \,= \  0 \ \ \ \  \   {\it{in}} \ \ \  \R^n$$
{\it has an infinite number of positive solutions\,} $\,\{ v_m \}_{m = 1}^\infty\, \subset\  C^{\,2\,,\ \bar\beta} \,(\R^n)$ {\it \,whenever\,} $|\,\varepsilon|\le \varepsilon_o\,.\,$ {\it Moreover,}

\vspace*{-0.3in}

 $$
 \Vert \,v_m - V_{\lambda_m\,, \ \xi_m} \Vert_\btd\ \,\le \  C_1\cdot\varepsilon \ \ \ \ \ \ {\it{for}} \ \ \ m = 1, \ 2, \cdot \cdot \cdot\,. \leqno (1.10)
 $$
{\it Here $\,\lambda_m \to 0\,$ and $\,|\,\xi_m| \to 0\,$ as $\,m \to \infty\,.\,$ As a result, \,$0\,$ is a blow\,-\,up point for $\,\{ v_m \}_{m = 1}^\infty$\,.}\\[0.2in]
\hspace*{0.5in}In (1.10), $\,\Vert \ \,\Vert_\btd\,$ represents the  $\,L^2$-\,norm on gradient for the Hilbert space ${\cal D}^{1\,, \ 2}$.\, See (2.7)\,.\, We also demonstrate in \cite{II}  how to transfer these solutions $\,\{ v_i \}_{i = 1}^\infty\,$ back to $S^n$ as solutions of equation (1.1).  [\,In (1.9), the number $5/2$\,    appears naturally in the calculations. It is, however, not claimed to be  sharp.\,]

\newpage

 \hspace*{0.5in}In this first part, we discuss the uniform cancelation property in the Lyapunov\,-\,Schmidt reduction method for   equation (1.1)\,.\, This uniformity enables us to find an infinite number of solutions. Already used in the Yamabe problem, and mentioned in \cite{Progress-Book}, this property is known by  experts working on the area. Nevertheless, it is an advantage to have the details recorded down. This is done in \S\,4\,.\, Given the available references, we select  those  (relatively) routine arguments and put them in the e\,-Appendix, which can be downloaded at

 \vspace*{-0.2in}

$$
{\mbox{www.math.nus.edu.sg/}}\, \tilde{\ }\, {\mbox{matlmc/e\,-Appendix.pdf}}
$$

\vspace*{0.15in}


{\bf \S\,1\,b. \  } {\it Additional highlights of the article.}\\[0.1in]
 $*_1$ \ \ In \S\,6\,,\, we present the first derivatives of the reduced functional $G_{|_{\bf Z}}$\,,\, and provide geometric interpretations of the formulas. The discussion culminates in showing that the existence of a critical point for $\,G_{|_{\bf Z}}$\,  implies the fulfillment of  the Kazdan\,-\,Warner condition  (refer to  Theorem 6.10 for full detail\,).\\[0.1in]
 $*_2$ \ \  In \S\,7\,,\, we   apply the formulas  on second derivatives  of \,$G_{|_{\bf Z}}$\,  to obtain existence result [\,for a single solution of equation (1.1)\,]\,  under symmetry conditions. See Theorem 7.18.

\vspace*{0.2in}

{\bf \S\,1\,c. \ } {\it Conventions.}\ \
Throughout this work,\, we assume  that the dimension $n \,\ge\, 3\,,\,$ except when otherwise is specifically mentioned, and let \, ${\tilde c}_n = (n - 2)/[\,4\,(n - 1)\,]\,.\,$
We observe the practice of  using $\,C,\,$ possibly with sub-indices, to denote various positive constants,
which may be rendered {\it differently\,} from line to line according to the contents. {\it Whilst we use $\,{\bar c}\,$ and  $\,{\bar C},\,$ possibly with sub-index, to denote a\,}   fixed\,
{\it positive constant which always keeps the same value as it is first defined\,}.\, [\,The {\it negative\,} constant $\,{\bar c}_{-1}\,$ is defined in (2.15).\,]\\[0.1in]
$\bullet_1$ \ \, Denote by $B_y \,(r)$ the open ball in $\R^n$ with center at $y$ and radius $r > 0\,,\,$  and $\,\Vert \,  S^n \, \Vert\,$ the measure of  $\,S^n\,$ in $\,\R^{n + 1}\,$ with respect to the standard metric.\\[0.1in]
$\bullet_2$ \ \, $\Delta_g$ is the Laplace-Beltrami operator associated with the metric $g$ on $S^n$\,.\, Likewise, $\,\Delta_o$ is the Laplace\,-\,Beltrami operator associated with the Euclidean metric $g_o$ on $\R^n,\,$ and   $\,\Delta_1$ is the
Laplace\,-\,Beltrami operator associated with the standard metric $g_1$ on $S^n$.\\[0.1in]
$\bullet_3$ \ \, Whenever there is no risk of misunderstanding, we suppress \,$dy$\, from the integral expressions.


  \newpage

{\bf \large \S\,2. \ \  The Lyapunov\,-\,Schmidt reduction\,: the case of  }\\[0.05in]
\hspace*{0.565in}{\bf \large  the conformal scalar curvature equation on $\R^n$.}\\[0.2in]
{\bf \S\,2\,a. \ } {\it Stereographic projection onto\,}  $\R^n$. \ \
Consider the stereographic projection
\begin{eqnarray*}
(2.1) \ \ \  \ \ \     \dot{\cal P} \,: \ S^n \setminus \{ {\bf N} \}\, &\longrightarrow &  {\R}^n\,\\
x & \mapsto & y \, =\, \dot{\cal P}(x)\,,\ \ \
{\mbox{where}} \ \    y_i \ = \ {{x_i}\over{1 - x_{n + 1}}}\,,
\ \ \ 1 \,\le\, i \,\le\, n\,.\ \ \  \ \ \ \ \ \  \ \ \
\end{eqnarray*}

Here $x \, =\, (x_1\,,\, \cdot \cdot \cdot, \ x_{n + 1})
\ \in \ S^n \setminus \{ {\bf N} \} \subset \R^{n + 1},\,\,$ and ${\bf N}\,=\, (0, \ \cdot \cdot \cdot\,, \ 0, \ 1)\,.\,$
Conversely,\,

$$\displaystyle{ x_i \ = \ {{2\,y_i}\over{1 + r^2}}\ , \ \ \ \ \  1 \,\le\, i \,\le\, n\,,\ \ \ \ \ {\mbox{and}} \ \ \
x_{n + 1} \ = \ {{r^2 -1 }\over{r^2 + 1}}\,,} \ \ \ \
{\mbox{where}} \ \  r \, = \, |\,y|\,.\, \leqno (2.2)$$
It is known that $\,\dot{\cal P}\,$ is a conformal map between $(S^n \setminus \{ {\bf N} \}\,, \ g_1)$ and $\,(\R^n,\ g_o)\,.\,$ The conformal factor is given by
$$
\ \ g_1 \,(x) \ = \ \left[ \,{{4}\over {(1 + r^2)^2}} \right]\,g_o \,(y) \ \ \ \ \ \ \ \ \ \ \ \ \ \ \ \ \ \ \mfor y \ = \ \dot{\cal P}(x)\,.\leqno (2.3)
$$
Set

\vspace*{-0.35in}

$$
v \, (y) \ :=\    u  \, ({\dot{\cal P}}^{-1} (y))\, \cdot \left( {2\over {1+ |y|^{\,2}}} \right)^{{n - 2}\over 2}  \ \  \ \ \mfor y \in \R^n. \leqno (2.4)
$$
Suppose $u$ is a positive solution of equation (1.1). Then  we have
$$
\Delta_o  \,v + ({\tilde c}_n \,K)\,v^{{n + 2}\over {n - 2}} \ = \  0 \ \ \ \ \ \ {\mbox{in}} \ \ \R^n\,. \leqno (2.5)
$$
We denote
$$
{\cal K} \ = \ 4n\,(n - 1) + \varepsilon\, {\cal H}\,.
$$
In (2.5), for $y \in \R^n\,,\,$ we set
$$
 K
(y) \  := \ {\cal K}  \,({\dot{\cal P}}^{-1} (y))\ = \ 4n\,(n - 1) \,+\, \varepsilon\, H (y) \ \  \ \  {\mbox{where}} \ \    H\,
(y)\ := \ {\cal H}  \,({\dot{\cal P}}^{-1} (y)) \,.
 \leqno (2.6)
$$

\vspace*{0.1in}

{\bf \S\,2\,b. \ } {\it  The Hilbert space and the}\, ({\it non-linear}) {\it  functional.}\ \
As there is an excellent book \cite{Progress-Book} on the topic, we  follow   closely \cite{Progress-Book}   and set up the basic notations for later study.
Throughout this article we work on the Hilbert space
$$
{\cal D}^{1, \  2}   \ =\ {\cal D}^{1, \ 2} (\R^n)\, :=\ \left\{ f\in L^{{2n}\over {n - 2}} (\R^n) \:\bigcap \: W^{1, \, 2}_{\small\mbox{loc}}\, (\R^n) \ \,  \bigg\vert \      \int_{\R^n} \langle \btd  \,f\,,\, \btd  \,f\rangle < \infty\right\}. \leqno (2.7)
$$
The inner product is defined by
$$
\langle \,f\,, \ \psi \,\rangle_{\,\btd}\ := \ \int_{\R^n} \langle \btd  \,f\,,\, \btd  \,\psi\rangle  \mfor \ f\,,\ \psi \in {\cal D}^{1, \  2}\,, \leqno (2.8)
$$
which induces the norm denoted by $\Vert \ \Vert_{\btd}\,.\,$
We recall the following properties.\\[0.2in]
{\bf \S\,2\,b.1.} \ \ {\it Sobolev's embedding with critical exponent\,.} \ \
(See, for instance, \cite{Brendle-1}.) There exists a dimensional constant ${\bar {\cal S}}_n\,$ (see \cite{Aubin}) such that
$$
\left( \,\int_{\R^n} |\,f|^{{2n}\over {n - 2}}  \right)^{\!\!{{n - 2}\over {2n}}} \ \le \ {\bar {\cal S}}_n \left(\, \int_{\R^n} \Vert \btd f \Vert^{\,2}  \right)^{\!{1\over 2}} \ \ \ \ \ \ {\mbox{for \ \ all}} \ \ \,f \in {\cal D}^{1, \, 2}. \leqno (2.9)
$$

\vspace*{0.1in}

 {\bf \S\,2.b.2.} \ \ {\it Integration by parts formula.} \ \
The following condition is sufficient for our use here.\,   Let $h \in {\cal D}^{1, \, 2} \,\cap\, C^2 \,(\R^n)$  with
$$
\int_{\R^n}   | \,\Delta_o \,h | ^{{2n}\over {n + 2}} \ <\  \infty\,.\leqno (2.10)
$$
Then we have
$$
\int_{\R^n} \langle \btd  \,f\,,\, \btd  \,h\rangle  \ = \ - \int_{\R^n}   f\cdot (\Delta_o \,h )  \ \ \mfor  \   {\mbox{all}} \ \ \,f \,\in\, {\cal D}^{1, \, 2}.\leqno (2.11)
$$
This can be shown by using  Sobolev inequality (2.9), and the fact that   $\,{\cal D}^{1,\,2}$ coincides with
the completion of $\,C^\infty_o\,(\R^n)\,$ with respect to
the $L^2$\,-\,norm of the gradient \cite{Ambrosett-Scalar}. See \S\,A.1 in the e\,-\,Appendix.

\vspace*{0.2in}

{\bf \S\,2.b.3.} \ \  {\it  The non\,-\,linear functional.}\ \
Corresponding to equation (1.5) with (1.6), consider the functional $\,I_\varepsilon\!:   {\cal D}^{1,\, 2} \to \R\,$ given by\\[0.05in]
(2.12)

\vspace*{-0.25in}

$$
 I_\varepsilon (f)  \, :=\,   \int_{\R^n} \left[\, {1\over 2}\, \langle \btd  \,f\,,\, \btd \,f \,\rangle \,-\,  {{n - 2}\over {2n}} \cdot  n\,(n - 2) \,f_+^{{2n}\over {n - 2}} \right]  \ -\  \varepsilon \cdot {{n-2}\over {2n}} \int_{\R^n}\  (\,{\tilde c}_n\,H\,) \,f_+^{{2n}\over {n - 2}}\,\,.
$$
Here $\,f_+\,$ denotes the positive part of $f\,.\,$
Separately, we define
$$
I_o\,(f) ={1\over 2} \int_{\R^n} \left[\,  \langle \btd_o \,f\,,\, \btd_o\,f \rangle -   (n - 2)^2   \,f_+^{{2n}\over {n - 2}} \,\right] \ \ \mfor f \in {\cal D}^{1, \, 2}\,, \leqno (2.13)
$$


$$
G \,(H)\,(f)\ = \  {\bar c}_{-1}\int_{\R^n} H \,f_+^{{2n}\over {n - 2}} \ \ \mfor f \in {\cal D}^{1, \, 2} \ \ \Longrightarrow \ \  I_\varepsilon  \, = \, I_o \,  + G (H)\,.\leqno (2.14)
 $$
Here

\vspace*{-0.3in}

$$\displaystyle{ \,{\bar c}_{-1}\ := \ -\,   {{ n\,-\,2 }\over {2\,n}}\cdot {\tilde c}_n \ = \ - \,{{(n\, - \,2)^2}\over {8n\,(n - 1)}}\ .\, }\leqno (2.15)$$
(The subindex of a negative sign is to remind ourselves that the constant is negative.)
In most cases there is no risk of confusion, we simply write
$$
G \,(H)\,(f)\ = \ G \,(f)\,.
$$
Given ${\cal K} \in C^\beta \,(S^n)$\, for $\beta \in (0, \ 1)\,,\,$
 (2.6) implies that $H$ is bounded and $\,H\in C^\beta \, (\R^n)$\,.\, It can be shown that any critical point $v$ of the functional $\,I_\varepsilon\,$ satisfies equation \,(2.1)\,.\, Moreover, standard elliptic theory and maximum principle imply that $v \in C^2 \,(\R^n)$  and $v > 0\,.\,$ See \cite{Gilbarg-Trudinger} and  \cite{Trudinger}\,.

\newpage

{\bf \S\,2\,c. \ } {\it  First order information --  $I_o'$\, and its kernel.}\ \
For   $\,f \,\in\, {\cal D}^{1, \, 2},\,$ a calculation using (2.13)  shows that the Fr\'echet derivative of $\,I_o$ is given by
$$
I_o' (f)\, [\,h] \ = \ \int_{\R^n} \!\left[ \,  \langle \btd \,f\,,\, \btd \, h \rangle -      n(n - 2) \,f^{{n + 2}\over {n - 2}}_+ \cdot h\,\right] \ \ \ \  \mfor h \in
{\cal D}^{1, \  2}. \leqno (2.16)
$$
We know that   (see \cite{Progress-Book}), \,for\, $v_o \in {\cal D}^{1, \  2}\,,\,$
$$
I_o' \,(v_o) \ \equiv \  0 \ \ \ \Longrightarrow \ \ v_o > 0 \ \ \, {\mbox{is \ \ smooth \ \ and \ \ }} \Delta \, v_o  + n\,(n - 2) \,v^{{n + 2}\over {n + 2}}_o\ = \ 0 \ \ \ \ {\mbox{in}} \ \ \R^n.
$$
The classification theorem obtained by Caffarelli-Gidas-Spruck/Gidas-Ni-Nirenberg \cite{Caffarelli-Gidas-Spruck} \cite{Gidas-Ni-Nirenberg.1} \cite{Gidas-Ni-Nirenberg.2}  shows that if $\,v_o\,$ is as above, then
$$
v_o \,(y)\ = \ V_{\lambda\,, \ \xi} \,(y)\ := \ \left( {\lambda\over {\lambda^2 + |\,y - \xi|^{\,2}}} \right)^{{n - 2}\over 2} \ \ \ \mfor \ y \in \R^n. \leqno (2.17)
$$
Here $\xi \in \R^n$ is fixed and $\lambda > 0\,.\,$ As set out in (1.3),
$${\bf Z}\ = \   {\mbox{Ker}}  \ I_o'\,.\, \leqno (2.18)$$
With the parametrization

\vspace*{-0.3in}

\begin{eqnarray*}
(2.19) \ \ \ \ \ \ \ \ \  \ \ \ \ \ \ \ \ \ \ \ \ \ \ \ \ \ \ \   \Phi\,: \R^+ \times \R^n & \to & {\bf Z}\\
 (\,\lambda\,, \, \xi) & \mapsto &  V_{\lambda\,, \  \xi}\,,\ \ \ \ \ \ \ \ \ \ \ \ \ \ \ \ \ \ \ \ \ \ \ \ \ \ \ \ \ \ \ \ \ \ \ \ \ \ \ \ \ \ \ \ \ \ \ \ \ \ \ \ \ \ \
\end{eqnarray*}
$\,{\bf Z}$\, is a $(n + 1)$\,-\,dimension  submanifold in  ${\cal D}^{1,\, 2}.$

\vspace*{0.2in}

{\bf \S\,2\,d. \ } {\it Second order information --  $I_o''$\, and its kernel.}\ \
Here we treat the second Fr\'echet derivative of $I_o$ at $f \in {\cal D}^{1, \  2}$ as a bilinear map on ${\cal D}^{1, \  2}$\,.\, See pp. 23 \,in \cite{Primer}.  The expression is given by
$$
(I_o'' (f)\, [\,  \phi] \  h)\ = \ \int_{\R^n} \!\left[\,  \langle \btd \,\phi\,,\, \btd \, h \rangle - n\,(n + 2) \,f_+^{{4}\over {n - 2}}\, \phi\cdot h\, \,\right] \mfor \phi\,, \ h \in {\cal D}^{1, \  2}\,. \leqno (2.20)
$$
In particular,   for $\,{\bf z} \in {\bf Z}\,,$
$$
(I_o'' ({\bf z})\, [\,\phi] \  h)\ = \ \int_{\R^n} \!\left[\,  \langle \btd \,\phi\,,\, \btd \, h \rangle - n\,(n + 2) \,{\bf z}^{{4}\over {n - 2}}\, \phi\cdot h\, \,\right]\mfor \phi\,, \ h \in {\cal D}^{1, \  2}\,. \leqno (2.21)
$$
Thus the condition
$$
(I_o'' ({\bf z})\, [\,\phi] \  h) \ \equiv \ 0 \ \ \  \ \ \  {\mbox{for \ \, all}} \, \ \ h \in {\cal D}^{1,\, 2}
$$
implies that
 $$
 \Delta_o \,\phi \ + \ n \,(n + 2) \,{\bf z}^{4\over {n - 2}} \,\phi\ = \ 0 \ \ \ \ {\mbox{in}} \ \ \R^n. \leqno (2.22)
$$
Equation (2.22) is also studied in \cite{Progress-Book}  \cite{Compact}\,.\bk
Using Riesz Representation Theorem (see, for example,  \S\,A.2 in the e\,-Appendix), we identify
\begin{eqnarray*}
(2.23) \ \ \ \      I_o'' ({\bf z}) \,: {\cal D}^{1, \  2} &\to& {\cal D}^{1, \  2}\\[0.15in]
f &\mapsto & {{ (I_o'' ({\bf z})\,[f] \,\,h) }\over {\Vert h \Vert_{\,\btd}^2}}\cdot h \ \ \ \ \ \ \  \  \{\,{\mbox{here}} \ \ h\not\equiv 0\,; \ \ \ h \,\perp {\mbox{Ker}}\,  I_o'' ({\bf z})\,[f]\,\}\,. \ \ \ \ \ \ \ \ \ \ \ \ \ \ \
\end{eqnarray*}
For the reduction method to work, we need the following   matching on the kernel of $ I_o'' ({\bf z})\,$ (refer to \cite{Progress-Book}\,)\,.\\[0.2in]
{\bf Lemma 2.24.} \ \ {\it For any ${\bf z} \in {\bf Z}$\,,\,   $I''_o\,({\bf z})$, as rendered in\,} (2.23), {\it is an index $\,0$ Fredholm operator\,.\, Moreover,  \,}
$$
{\mbox{Ker}}  \, I''_o\,({\bf z})\ = \ T_{\bf z} \,{\bf Z} \ \subset \ {\cal D}^{1, \  2}\,.
$$

\vspace*{-0.15in}

{\it Let }
$$
{\bf W}_{\bf z}\ := \ [\,{\mbox{Ker}}\,\, I''_o({\bf z})\,]^\perp \ = \ [\,T_{\bf z} \,{\bf Z}]^\perp \ \subset \ {\cal D}^{1, \  2}\,. \leqno (2.25)
$$
{\it Then }

\vspace*{-0.25in}

$$
 I_o'' \,({\bf z}) \,: \  {\bf W}_{\bf z} \to  {\bf W}_{\bf z}
$$
 {\it is an isomorphism.}\\[0.2in]
 %
%
%
 %
 \hspace*{0.5in}See the proof of Lemma 2.10 in pp. 21 (\cite{Progress-Book}), and Lemma 5.2 in  pp. 61 (\cite{Progress-Book}). See also  \S\,A.2 in the e\,-Appendix. Cf. \cite{Chen-Lin-3} and   Lemma 4.2 in \cite{Compact}\,.\, We first  observe that, in (2.23) and with the notations used there,
\begin{eqnarray*}
 (2.26) \ \ \ \ \ \ \  \ \ {\mbox{Ker}}\, I''_o \,({\bf z})\ = \ T_{\bf z} \,{\bf Z} \ \ & \Longrightarrow & \ \ T_{\bf z}\, {\bf Z} \,\subset\, {\mbox{Ker}}\, I''_o \,({\bf z}) \,[f] \\[0.1in]
 &\ \Longrightarrow & \ \  h \,\in\, {\bf {\bf W_z}}\  \ \   \   \ \ \ \   \ \ \ \ \ \ \ (\,{\mbox{as}} \ \ h \perp {\mbox{Ker}}\, I''_o \,({\bf z})\, [f]\,)\,. \ \ \ \ \ \ \  \ \ \ \ \ \ \
\end{eqnarray*}

\vspace*{0.1in}

{\bf \S\,2\,e. \ } {\it   The flow chart.}\ \
The simplicity of the Lyapunov\,-\,Schmidt reduction method is exposed  in the book \cite{Progress-Book}.
We keep the notations close to those used in \cite{Progress-Book}. We use the following chart to present the set\,-\,up of the Lyapunov\,-\,Schmidt reduction method, and  encourage  readers to refer to the excellent book \cite{Progress-Book}.

\newpage

\begin{eqnarray*}
I_\varepsilon \,(f)\ = \ I_o\,(f) &+& \varepsilon\, G \,(f)\ \ \  \ \ {\it for} \ \ f \in {\cal D}^{1, \  2}\,.\\
\downarrow  \ \ &\,&\\
{\bf Z}\ = \ {\mbox{Ker}} &\,&  \!\!\!\!\!\!\!\!\!\!\!\!\!I_o'\ \  \ \ \ \ \ \  \ \ \ \ \ \ \ {\it Refer \ \ to \ \ } (2.18)\,.\\
\downarrow  \ \ &\,&\\
T_{\bf z}\,{\bf Z} \ \ &\,& \ \ \longleftarrow \ \ {\mbox{Ker}}\,I''_o\,({\bf z})\ = \ T_{\bf z}\,{\bf Z} \ \ \ \ \ \ \ \  ({\mbox{Lemma \ 2.24}}\,.)\\
\downarrow  \ \ &\,&\\
{\bf W}_{\bf z}\ := \ (\,T_{\bf z}\,{\bf Z})^\perp\!\!\!\!\!\!\!\!\!\! &\,&\ \ \ \,((2.25)\,; \ \ {\it{write}} \ \ \,{\cal D}^{1, \  2}\ = \ T_{\bf z} \, {\bf Z} \oplus W_{\bf z}\,. ) \\
\downarrow  \ \ &\,&\\
P_{\bf z}: {\cal D}^{1, \  2}   \to {\bf W}_{\bf z}  & \ & \ \ \ \,({\it{Projection \ unto \ the \ normal.}})\\
\downarrow  \ \ &\,&\\
P_{\bf z}\!\circ I_\varepsilon' ({\bf z} + w_\varepsilon ({\bf z}))\ = \ 0  \!\!\!\!\!\!\!\!\!\! &\,& \ \ \ \ \leftarrow \ \ {\it The \ \ auxiliary \   \ equation\,. \ \ Solution}\,: w_\varepsilon ({\bf z}) \in W_{\bf z}\,.\\  |\ \,\,\,&\,&  \ \ \ \ ({\it Cancelation \ \ along \ \ the \ \ normal \ \ directions.)} \\
\downarrow  \ \ &\,&\\
 \Phi_\varepsilon\, ({\bf z})\ := \ I_\varepsilon\, ({\bf z} + w_\varepsilon\, ({\bf z})) \!\!\!\!\!\!\!\!\!\!\!\!\!\!\!\!&\,& \ \ \  \  ({\it The  \ \ finite \ \ dimensional \ \ reduction\,; \ \ {\bf z} \in {\bf Z}\,.})\\
\downarrow  \ \ &\,&\\
\Phi_\varepsilon' \,({\bf z}_\varepsilon + w_\varepsilon \,({\bf z}_\varepsilon))\ = \ 0  \!\!\!\!\!\!\!\!\!\! &\,& \ \ \ \ ({\it Critical \ \ point \ \ {\bf z}_\varepsilon  \ \ of \ \ the \ \ finite \ \ dimensional \ \ functional}.)  \\
\downarrow  \ \ &\,&\\
I_\varepsilon' \,({\bf z}_\varepsilon + w_\varepsilon \,({\bf z}_\varepsilon))\ = \ 0  \!\!\!\!\!\!\!\!\!\! &\,& \ \ \ \ (\,{\it If \ \  such  \ \  } {\bf z}_\varepsilon \ \ {\it exists \ \ for}  \ \ |\,\varepsilon| \ \ {\it small \ \ enough\,.\,}) \\
\downarrow  \ \ &\,&\\
 z_\varepsilon + w_\varepsilon\, ({\bf z}_\varepsilon)   \!\!\!\!\!\!\!\!\!\! &\,& \ \  \ \, {\it is  \ \ a   \ \ solution \ \ of \ \ equation} \ \ {\mbox{(1.5)}} \ \ {\it with}  \ \ {\mbox{(1.6)}}\,. \\
 & \,&  \ \ \ \,({\it Refer \ \ to} \ \ {\mbox{Theorem\  2.27}}\,.)
\end{eqnarray*}

\vspace*{0.15in}

{\bf Theorem 2.27. } \ \ {\it Keeping the notations in the chart, \, given a compact subset $\,{\bf Z_{\,c}} \subset {\bf Z}$\,,\, assume that $\,\Phi_\varepsilon\,$ has\,, for $\,|\,\varepsilon|\,$ sufficiently small, a critical point ${\bf z}_\varepsilon \in {\bf Z}_{\,c}\,.\,$  Then $\,{\bf z}_\varepsilon\, +\, w_\varepsilon ({\bf z}_\varepsilon)\,$ is a critical point of $\,I_\varepsilon$\,. }\\[0.2in]
\hspace*{0.5in}For the proof of Theorem 2.27, we refer to \cite{Progress-Book} (\,in particular, Theorem 2.12, Lemma 5.2 [\,pp. 61] and the discussion in pp. 59\,\,--\ 60\,,\, loc. cit.\,)\,.

\vspace*{0.2in}

{\bf \S\,2\,f.  } {\it Restriction of \,\,$G$\, to \,${\bf Z}$.}\ \
With the structure
$$
\Phi_\varepsilon ({\bf z})\ = \ I_o \,({\bf z} + w_\varepsilon ({\bf z})) + \varepsilon \,G \, ({\bf z} + w_\varepsilon ({\bf z}))\,,
$$
the finite dimension reduction   can be simplified further. Consider the restriction of $\,G\,$ to $\,{\bf Z}\,$\,:

\newpage

\begin{eqnarray*}
(2.28) \ \ \ \ \  \ \ \ \ \  \ \ G_{|_{\bf Z}} & : & \ \ {\bf Z}  \ \to \  \R\\
G_{|_{\bf Z}}\, ({\bf z}) &\ = \ & G_{|_{\bf Z}} \,(V_{\lambda\,, \ \xi}) \ = \ {\bar c}_{-1} \int_{R^n} H (y) \left[\, V_{\lambda\,, \, \xi} \right]^{{2n}\over {n - 2}} \, \ \ \ \ \ \ \ \ \ \ \ \ \   \  \ \     ({\mbox{see}} \ \ (2.14))\\
&\ = \ & {\bar c}_{-1} \int_{R^n} H (y) \left[ {\lambda\over {\lambda^2 \,+\, |\,y - \xi|^{\,2} }} \right]^n \, \mfor V_{\lambda\,, \ \xi} \,\in \,{\bf Z}\ \ \ \    ({\mbox{see}} \ \ (2.17)). \
\end{eqnarray*}
It is shown in \cite{Progress-Book} that $\,G_{|_{\bf Z}}\,$ and $\,G_{|_{\bf Z}}'\,$ are good approximations of $\,\Phi_\varepsilon\,$ and $\,\Phi_\varepsilon'\,$ on $\,\bf Z$\,, respectively. Precisely,\, for $\,{\bf z} \,\in\, {\bf Z}$\,,
\begin{eqnarray*}
\Phi_\varepsilon ({\bf z}) &= &I_o \,({\bf z}) \ + \ \varepsilon\, G_{|_{\bf}} \,({\bf z})\  + \ O \,(\varepsilon) \ \ \ \ \ \ \ [\,{\mbox{refer \   to \   Lemma\  2.15 \  in \   \cite{Progress-Book}\,, \ pp. 24}}\,]\\[0.13in]
\Phi_\varepsilon' ({\bf z}) &= &  \varepsilon\, G_{|_{\bf Z}}' ({\bf z}) \ + \ O \,(\varepsilon)\,, \ \ \ \ \ \ \ \  (\,{\mbox{the \ \ derivatives \ \ act \ \ on \ \ }} T_{\bf z} \,{\bf Z}\,; \ \ {\mbox{refer \ to \ the }}\\
 & \ & \ \ \ \  \ \ \ \ \ \ \ \ \ \ \ \ \ \ \ \ \ \ {\mbox{ \ \ \ \ \  \ \ \ \ \ \ \ \ \  \ \ \ \ \ \ \ proof \ of \ Theorem \  2.17 \  in \   \cite{Progress-Book}\,, \ pp. 26}})\,.
\end{eqnarray*}

\vspace*{0.1in}

{\bf Theorem 2.29.} \ \ {\it Let $\,G_{|_{\bf Z}}'\,(\bar {\bf z})\ = \ 0\,.\,$   Suppose one of the following conditions\,  holds}\,. \\[0.1in]
{\bf (a)} \ \ \,\,{\it $\bar {\bf z}\,$ is a strict local maximum or minimum of} $\,G_{|_{\bf Z}}\,.$ \\[0.1in]
{\bf (b)} \ \ {\it There is an open and bounded set $\,{\cal N} \subset {\bf Z}\,$ such that $\,\bar{\bf z} \,\in\, {\cal N}\,$ and} $${\mbox{deg}} \,\, (\,G_{|_{\bf Z}}'\,, \ {\cal N}\,, \ {\bf 0}) \  \not=\  0\,.$$
{\it Then for $|\,\varepsilon|$ small enough, $\,\Phi_\varepsilon$ has a critical point.\,} (\,{\it It follows from Theorem  }2.27 {\it that $I_\varepsilon$ also has a critical point.}\,)\\[0.2in]
\hspace*{0.5in}See Theorem 2.16 and Theorem 2.17 in \cite{Progress-Book}. Via the parametrization (2.19), we identify $G_{|_{\bf Z}}$ as a map from $\R^+ \times \R^n$ to $\,\R$.\, In particular, \, $\bar{\bf z}\ = \ V_{\bar\lambda\,, \,\,\bar \xi}$ \ is a critical point of  $\,G_{|_{\bf Z}}\,$ if and only if
 $$
 {{\partial G_{|_{\bf Z}} }\over {\partial \lambda}} \, (\bar\lambda\,, \,\bar \xi)\ = \ 0\,, \ \ \ \  {\mbox{and}} \ \ \ \  {{\partial G_{|_{\bf Z}} }\over {\partial\,\xi_j}} \, (\bar\lambda\,, \,\bar \xi)\ = \ 0 \mfor j\ = \ 1, \ 2,\ \cdot \cdot \cdot, \ n\,. \leqno (2.30)
  $$
Refer to \S\,A.4.2 in the e\,-appendix. With this, we can bring the discussion to $\R^+ \times \R^n$\,.\, In  particular, we have the following.

\vspace*{0.2in}

 {\bf Definition 2.31.} \ \ {\it A critical point\,}  $\,{\bf p}\, := \ (\bar\lambda\,, \ \bar\xi)\, \in \,\R^+ \times \R^n$  \,{\it of $\ G_{|_{\bf Z}}\,$ is called\,} stable {\it if there exists a ball}\, $B_{\,{\bf p}} \,(\rho) \,\subset\,  \R^+ \times \R^n$ {\it such \,that}\, ${\mbox{deg}} \, (\,\btd_{(\,\lambda\,,\,\, \xi)} \, G_{|_{\bf Z}}\,, \ B_{\,{\bf p}}\, (\rho)\,, \ \,\vec{\,0}\,\,) \,\not= \,0\,.$ \\[0.2in]
 \hspace*{0.5in}Using the parametrization (2.19), one infers that a stable critical point of $\,G_{|_{\bf Z}}\,$ fulfills condition {\bf (b)}\, in  Theorem  2.29. We refer to \S\,A.3 in the e\,-\,Appendix for selected fundamental properties on the degree of a map. See also \S\,A.4.2 in the e\,-appendix concerning the Hessian and the Jacobian.

\newpage

{\bf \large \S\,3.     Tangent space, hyperbolic structure, and the normal space.}\\[0.2in]
{\bf \S\,3\,a. \ } {\it The tangent space. }\ \
For a point ${\bf z}\ = \ V_{\lambda\,, \, \xi} \,\in\, {\bf Z}$\,,\, via (2.19), the  tangent space $\,T_{\bf z}\,{\bf Z}\,$ is spanned by the $(n + 1)$ functions\,:
\begin{eqnarray*}
(3.1) \  \ \  \varphi_o:= {{\partial V_{\lambda\,, \, \xi}}\over {\partial \lambda}}\ = \ {{\partial}\over {\partial \lambda}} \left[\left({\lambda\over {\lambda^2 + |\, y - \xi|^{\,2}}} \right)^{\!\!{{n - 2}\over 2}} \,\right] &\! = \! &\!-{{n - 2}\over 2} \cdot \lambda^{{n - 4}\over 2} \cdot {{ (\,\lambda^2 - |\, y - \xi|^{\,2})}\over {(\,\lambda^2 + |\, y - \xi|^{\,2})^{{n }\over 2} }}\,,\\
\varphi_j:= {{\partial V_{\lambda\,, \, \xi}}\over {\partial\,\xi_j}}\ = \ {{\partial}\over {\partial\,\xi_j}} \left[\left({\lambda\over {\lambda^2 + |\, y - \xi|^{\,2}}} \right)^{\!\!{{n - 2}\over 2}} \,\right] &\! = \! & -{{n - 2}\over 2} \cdot \lambda^{{n - 2}\over 2} \cdot {{ 2\,(\xi_j - y_j)}\over {(\,\lambda^2 + |\, y - \xi|^{\,2})^{{n }\over 2} }} \ \ \ \
\end{eqnarray*}
for $\,j\ = \ 1, \ 2, \cdot \cdot \cdot, \ n\,.\,$
 A direct calculation  shows that the collection is orthogonal. That is
$$
\langle \,\varphi_o\,, \ \varphi_j \rangle_{\,\btd}\ = \ 0\ = \ \langle\, \varphi_i\,, \ \varphi_j \rangle_{\,\btd} \ \ \ \ \   \mfor 1 \, \le \,i\, \not= j \,\le \,n\,.$$
Moreover,

\vspace*{-0.3in}

\begin{eqnarray*}
(3.2) \ \ \  \ \ \  \ \  \Vert \,\varphi_o \Vert_{\,\btd}^2 &\ = \ & \int_{\R^n} \bigg\langle \!\btd \, {{\partial V_{\lambda, \, \xi} }\over {\partial \lambda}} \ , \ \ \btd  \, {{\partial V_{\lambda, \, \xi} }\over {\partial \lambda}}\,\bigg\rangle\ = \  {{ {\bar c}_1^{\,2} }\over {\lambda^2}}\,,\\[0.1in]
(3.3) \ \ \  \ \ \   \ \   \Vert\, \varphi_j \Vert_{\,\btd}^2 &\ = \ & \int_{\R^n} \bigg\langle \!\btd \, {{\partial V_{\lambda, \, \xi} }\over {\partial\,\xi_j}} \ , \ \ \btd  \, {{\partial V_{\lambda, \, \xi} }\over {\partial\,\xi_j}}\,\bigg\rangle\ = \  {{ {\bar c}_2^{\,2} }\over {\lambda^2}} \ \ \ \ \mfor \ \  1 \,\le \,j \,\le \,n\,. \ \ \  \ \ \  \ \ \  \ \ \  \ \ \  \ \ \  \ \ \  \ \ \  \ \ \  \ \ \
\end{eqnarray*}
Here the positive constants $\,{\bar c}_1 \,$ and $\,{\bar c}_2\,$ depend on the dimension $n$ only.

\vspace*{0.2in}

 {\bf \S\,3\,b.    } {\it Geometry of $\ {\bf Z}\,$ and hyperbolic space.}\ \ 
 On $\R^+ \times \R^n$ with coordinates
\begin{eqnarray*}
& \ & (\,\lambda\,, \ \xi_1\,, \cdot \cdot \cdot\,, \ \xi_n)\,,\\
& \ & \!\!\!\!\!\!\!\!\!\!\!\!\!\!\!\!\!\!\!\!\!\!\!\!\!\!\!\!\!\!\!\!\!\!{\mbox{indexed zero}}   \ \,\uparrow \ \
\end{eqnarray*}
 using (2.19) to pull back the metric defined in (3.2) and (3.3), we obtain an induced Riemannian metric  on $\R^+ \times \R^n\,:$
 $$
 {{ c_\ell  }\over {\lambda^2}} \cdot \delta_{\ell\, k} \ \ \ \ {\mbox{with}} \ \ \ \ c_o\ = \ {\bar c}_o\,, \ \ \ c_1\ = \ c_2\ = \ \cdot  \cdot \cdot\ = \ c_n\ = \ {\bar c}_1 \ \ \ \ \ \ (\,\lambda > 0)\,, \leqno (3.4)
 $$
 where $\,0 \le \ell\,, \ k \le n\,.\,$
 The curvature tensor   calculation  shows that the induced metric has {\it constant\,} sectional  curvature equal to $\displaystyle{\,-\, {1\over {{\bar c}_o }}}\,.\ $ See \S\,A.4 in the e\,-Appendix.\, The standard hyperbolic metric on $\R^+ \times \R^n\,$ is
 $$
 {{1 }\over {t^2}} \,\, (dt^2 \otimes \,dy^2) \ \ \ \ \ {\mbox{at}} \ \ \, (t\,,\, y) \in \R^+ \times \R^n.
 $$
Thus we have the isometry\,:
 \begin{eqnarray*}
 \left( \R^+ \times \R^n\,, \ \ \ {{ c_\ell  }\over {\lambda^2}} \cdot  \delta_{\ell\, k}\right)    &\longrightarrow &  \left( \R^+ \times \R^n\,, \  \ \ {\bar c}_o \cdot \left[ \, {{1 }\over {t^2}} \,\, (\,dt^2 \times dy^2\,) \right]\,\right)\\
 (\,\lambda\,, \ \xi) & \mapsto &  \ (t\,, \ y)\ = \ \left(\,\lambda\,,\, \ \,{{ {\bar c}_1}\over { {\bar c}_o}} \cdot \xi\right)\,.
 \end{eqnarray*}
 The (rescaled) hyperbolic structure provides a  complete (non\,-\,compact) metric on \,{\bf Z}\,.

\newpage

{\bf  \S\,3 c.\ } {\it   Orthonormal Basis.}\ \
Following (3.1), let
\begin{eqnarray*}
(3.5) \ \ \ \ \ \ \ \ \ \ \ \  q_o &:=& {{\partial V_{\lambda, \ \xi} }\over {\partial \lambda}}\bigg/\, \bigg\Vert\, {{\partial V_{\lambda, \ \xi} }\over {\partial \lambda}} \,\bigg\Vert_{\,\btd}\ = \ {\bar C}_o\, \lambda \cdot \lambda^{{n - 4}\over 2} \cdot {{ (\,\lambda^2 - |\, y - \xi|^{\,2}) }\over {(\,\lambda^2 + |\, y - \xi|^{\,2})^{n\over 2} }}\,,\\[0.1in]
q_j &:= & {{\partial V_{\lambda, \  \xi} }\over {\partial\,\xi_j}}\bigg/\, \bigg\Vert \,{{\partial V_{\lambda, \  \xi} }\over {\partial\,\xi_j}} \,\bigg\Vert_{\,\btd}\ = \ {\bar C}_1\, \lambda \cdot \lambda^{{n - 2}\over 2} \cdot {{ 2 \,(\,\xi_j - y_j) }\over {(\,\lambda^2 + |\, y - \xi|^{\,2})^{n\over 2} }}\ \ \ \ \ \ \ \ \ \ \ \  \ \ \ \ \ \  \ \ \ \ \ \
\end{eqnarray*}
for $\,j\ = \ 1, \ 2, \cdot \cdot \cdot, \ n\,.\,$
Here $\,{\bar C}_o\,$ and $\,{\bar C}_1\,$ are positive constants depending on $n$ only.
Thus $\,\{ \,q_\ell \,\}_{\,0 \,\le \,\ell \,\le \,n}\ $ forms an orthonormal basis for $\,T_{\bf z}{\bf Z}\,$ with $\,{\bf z} = V_{\lambda\,,\ \xi}\,.\ $  Using (2.11), we find
$$
\bigg\Vert \,{{\partial q_\ell}\over {\partial \lambda}} \bigg\Vert^2_{\,\btd}\ = \ \bigg\langle\! \btd \, {{\partial q_\ell}\over {\partial \lambda}}  \,, \ \btd  \, {{\partial q_\ell}\over {\partial \lambda}} \bigg\rangle_{\,\btd} \ = \ \bigg\langle \!\Delta \, {{\partial q_\ell}\over {\partial \lambda}}  \,, \   {{\partial q_\ell}\over {\partial \lambda}} \bigg\rangle_{\,\btd}\ = \ \bigg\langle \!   {{\partial (\Delta \,q_\ell)}\over {\partial \lambda}}  \,, \   {{\partial q_\ell}\over {\partial \lambda}} \bigg\rangle_{\,\btd} \ ,
$$
leading to

\vspace*{-0.3in}

$$
\bigg\Vert \,D_o\, q_\ell  \bigg\Vert^2_{\,\btd}\ = \ \bigg\Vert \, {{\partial q_\ell}\over {\partial \lambda}} \bigg\Vert^2\ = \ {{{\bar C}_2}\over {\lambda^2}}\ . \leqno (3.6)
$$
\vspace*{-0.2in}

Likewise,

\vspace*{-0.2in}

$$
\bigg\Vert \,D_j\, q_\ell  \bigg\Vert^2_{\,\btd}\ = \ \bigg\Vert\, {{\partial q_\ell}\over {\partial\,\xi_j}} \bigg\Vert^2\ = \ {{{\bar C}_3}\over {\lambda^2}} \ \ \ \ \ \ \mfor j\ = \ 1, \ 2, \cdot \cdot \cdot, \ n\,. \leqno (3.7)
$$
In (3.6) and (3.7),  $\,{\bar C}_2\,$ and $\,{\bar C}_3\,$ are positive constants depending on $n$ only.
Refer to \S\,A.4.3 in the e\,-\,Appendix for  more details.

\vspace*{0.2in}


 {\bf \S\,3\,d.   \ } {\it     The normal space ${\bf W}_{\bf z}\ = \ (T_{\bf z} \,{\bf Z})^\perp.$}\ \
By definition,
\begin{eqnarray*}
{\bf W}_{\bf z}  \ = \   (T_{\bf z}\, {\bf Z})^\perp &=& \left\{ w \in {\cal D}^{1, \  2} \ \ \bigg\vert \ \ \int_{\R^n} \langle  \btd \,w\,, \btd \,\varphi_\ell\,  \rangle \,\ = \ 0 \mfor \ell\ = \ 0, \ 1, \cdot \cdot \cdot, \ n \,\right\}\\
&\ = \ &  \left\{ w \in {\cal D}^{1, \  2} \ \ \bigg\vert \ \int_{\R^n} w \, \Delta \,\varphi_i    \,\ = \ 0 \ \ \ \ \ \ \ \ \, \mfor \ell \ = \ 0, \ 1, \cdot \cdot \cdot, \ n \,\right\}\,.
\end{eqnarray*}
Here we apply (2.11) and use (3.1) to check  condition (2.10)\,.\,
Define
\begin{eqnarray*}
(3.8) \ \ \ \ \ \  \ \    \psi_o&\ = \ &\lambda^{{{n }\over 2} + 1}\cdot
 {{ (\,\lambda^2 - |\,y - \xi|^{\,2})}\over {(\,\lambda^2 + |\,y - \xi|^{\,2})^{{n \over 2}+2} }}\ ,\ \ \ \ \ \
 \psi_j  =  \lambda^{{{n }\over 2} + 2}\cdot
 {{ 2\,(y_i - \xi_i)}\over {(\,\lambda^2 + |\,y - \xi|^{\,2})^{{n \over 2}+2} }} \ \ \ \ \ \ \ \ \ \ \   \end{eqnarray*}
 for $\,j\ = \ 1, \ \cdot \cdot \cdot, \ n\,.\, $
These are constants times $\Delta\, \varphi_\ell$\,.\,  See \S\,A.4.3 in the e\,-Appendix.\, The power on $\lambda$ are chosen so that their $L^{{2n}\over {n - 2}}\, (\R^n)$\,-\,norms are constant (independent on $\lambda$ and $\xi$)\,.\, We can express
$$
{\bf W}_{\bf z}\ = \ \left\{ w \in {\cal D}^{1, \  2} \ \bigg\vert \ \int_{\R^n} w \, \psi_\ell    \,\ = \ 0 \ \ {\mbox{for}} \ \  \ell \ = \ 0, \ 1, \cdot \cdot \cdot, \ n \right\} \ \ {\mbox{for}} \ \  {\bf z}\ = \ V_{\lambda, \,\, \xi} \in  {\bf Z}\,. \leqno (3.9)
$$
Let us continue with the observation\,:

\newpage

\begin{eqnarray*}
(3.10) \ \ \int_{\R^n} V_{\lambda,\ \xi}^{{2n}\over {n - 2}}\ = \ {\mbox{Const.}} & \ &\!\!\!\!\!\!
\Longrightarrow \ \ D_\ell \, \left[ \int_{\R^n} V_{\lambda,\ \xi}^{{2n}\over {n - 2}} \right]\ = \ 0 \ \ \ \ \mfor \ell\ = \ 0\,, \ 1\,, \ \cdot \cdot \cdot, \ n\\
& \ &\!\!\!\!\!\!
\Longrightarrow \ \ \int_{\R^n} V_{\lambda,\ \xi}^{{n + 2}\over {n - 2}} \,[ D_\ell \,V_{\lambda,\ \xi} ]\ = \ 0 \mfor \ell\ = \ 0\,, \ 1\,, \ \cdot \cdot \cdot, \ n\\
& \ &\!\!\!\!\!\! \Longrightarrow \ \ \int_{\R^n} [\,\Delta_o \,V_{\lambda,\ \xi}] \,[ D_\ell \,V_{\lambda,\ \xi} ]  \ = \ 0 \ \ \ \ \ \ \ \  [\,\mbox{equation} \ \ (1.4)\,]\\[0.1in]
& \ &\!\!\!\!\!\! \Longrightarrow \ \ \langle \,V_{\lambda,\ \xi}\,, \ \phi_\ell \,\rangle_{\,\btd}\ = \ 0 \ \ \ \ \ \ \ \ \ \ \ \ \ \ \ \   [\,\mbox{via} \ \ (2.8) \ \& \ (2.11)\,] \ \ \ \ \\[0.1in]
& \ &\!\!\!\!\!\!\Longrightarrow \ \ {\bf z}\ = \ V_{\lambda\,,\  \xi} \in {\bf W}_{\bf z}\,.
\end{eqnarray*}
This reminds us the sphere in $\R^n$, where any point on the sphere is consider as a vector which is perpendicular to the tangent space of the sphere at that point. \bk
 The following property concerning ${\bf W}_{\bf z}$ is a first step toward {\it uniform cancelation}\, in the Lyapunov\,-\,Schmidt reduction method for equation (1.1).  Refer to \cite{Progress-Book}, \cite{Brendle-1} and the references therein. \\[0.15in]
%
%
%
%
%
%
{\bf Proposition 3.11.} \ \ {\it There exists a positive number $\,{\bar c}_2\,$ such that  the inequality }
$$
\int_{\R^n} \!\!\left(\, |\btd\, w|^{\,2} - n\,(n + 2)\, V_{\lambda\,, \, \xi}^{4\over {n - 2}} \, w^2\right)\ \ge  \ 2\, {\bar c}_2\, \Vert w \Vert_{\,\btd}^2 \ -\  {{16n^2}\over {{\bar c}_2}} \left(\int_{\R^n} V_{\lambda\,, \, \xi}^{{n + 2}\over {n - 2}} \, w\right)^2 \leqno (3.12)$$
({\it for all\, } $w \in {\bf W}_{\bf z}\,$) {\it holds uniformly for all\,} ${\bf z}\, = \, V_{\lambda\,, \,\, \xi} \in {\bf Z}\,$ ({\it that is, \ ${\bar c}_2\,$ does not depend on\ } ${\bf z} \in {\bf Z}\,,\, \  {\it{or}} \ \, w \in {\bf W}_{\bf z}\,)\,.$

\vspace*{0.2in}


{\bf Corollary 3.13. } \ \ {\it For\,} $w \perp ( \,{\mbox{span}} \, \{{\bf z} \} \oplus T_{\bf z}\, {\bf Z}\,)$\,  {\it we have }
$$
\int_{\R^n}   \left[\, |\!\btd  w|^{\,2} \,- \,n\,(n + 2)\, V_{\lambda\,, \, \xi}^{4\over {n - 2}} \, w^2\right]\  \ge \ 2 \,{\bar c}_2\, \Vert \,w \Vert_{\,\btd}^2\,,  \leqno (3.14)
$$
{\it uniformly for all\,} ${\bf z}\ = \ V_{\lambda\,,\,\, \xi} \in {\bf Z}\,.\, $
({\it Here \,${\bar c}_2\,$ is the same constant appeared in Proposition\,}  3.11.)

\vspace*{0.1in}

\hspace*{0.5in}As the proof of the corollary is well\,-\,known among experts, we place it in the e\,-Appendix (\,\S\,A.4.\,d) for the interested readers to view\,.\, Finally, consider  $w\ = \ {\bf z} =V_{\lambda\,,\,\, \xi}$\, [\,recall that (3.10) shows that ${\bf z} \in {\bf W_z}\,.]$\, Let us check what happens to the left hand side of (3.14)\,.\,  To see  precisely how  (3.14) fails to hold when $\,w\ = \ {\bf z}\ =\ V_{\lambda\,,\,\, \xi}\,,$\, we  observe
\begin{eqnarray*}
& \ &V_{\lambda\,,\  \xi} \,\,\Delta\,V_{\lambda\,, \  \xi} \ + \  n \,(n - 2)\, V_{\lambda\,, \, \xi}^{{2n}\over {n - 2}}\ = \ 0\\[0.1in]
  \Longrightarrow   & \ &  \int_{\R^n}  \!\left[\, |\btd\, V_{\lambda\,, \, \xi}|^{\,2} \ - \ n (n - 2) V_{\lambda\,, \, \xi}^{4\over {n - 2}} \cdot V_{\lambda\,, \  \xi}^2\right]\,\ = \ 0 \ \ \ \ \ ({\mbox{integration \ by \ parts}}) \\
& \ & \\
    \Longrightarrow   & \ &   \int_{\R^n} \!\!\left[\, |\btd\, V_{\lambda\,, \, \xi}|^{\,2} \ - \ n\,(n + 2) V_{\lambda\,, \, \xi}^{4\over {n - 2}} \cdot  V_{\lambda\,, \  \xi}^2\right]\ = \ - \,4n \int_{\R^n}V_{\lambda\,, \ \xi}^{{2n}\over {n - 2}} \\[0.1in]
     & \ &  \  \ \ \  \ \  \     = \ {{4n}\over {n(n - 2)}} \int_{\R^n} (\Delta\, V_{\lambda\,, \, \xi}) \,V_{\lambda\,, \  \xi}  \ = \  -\,{4\over {n - 2}} \int_{\R^n} |\!\btd\,{\bf  z}|^{\,2}\ = \ -\,{4\over {n - 2}}\, \Vert \, {\bf z} \Vert_{\,\btd}\,.\ \  \ \ \ \  \ \ \ \  \ \ \ \  \ \ \ \  \ \ \ \  \ \  \ \  \ \ \ \  \ \
\end{eqnarray*}

 \newpage


{\bf \large \S\,4. \ \ Uniform cancelation.}\\[0.2in]
{\bf \S\,4\,a. \ } {\it  Uniform solvability of the auxiliary equation.}\ \
The following result is from Lemma 2.21 in \cite{Progress-Book}, pp. 27. See also   (2.26) in \S\,2\,d\,.\\[0.2in]
{\bf  Lemma 4.1.} \ \  {\it Refer to}\, (2.23), Lemma 2.24 {\it and}\, (2.26).  {\it Assume  conditions\,} ({\bf i}) -- ({\bf iii}). \\[0.15in]
({\bf i}) \ \ $
 I_o'' ({\bf z}) : {\bf {\bf W_z}}   \to {\bf {\bf W_z}}
$
{\it is uniformly invertible in $\,{\bf z}$. That is,}
$$
\Vert \,( I_o''({\bf z}))^{-1} \Vert_{\,L \,({\bf {\bf W_z}} , \, \,{\bf {\bf W_z}} )}\ \le \ C \ \ \ \ {\it{for \ \ all}} \ \  {\bf z} \in {\bf Z}\,. \leqno (4.2)
$$

({\bf ii}) \ {\it \,The remainder in}
$$
  I_o'({\bf z}  + w) -I_o'' ({\bf z}) \,[w] \ \ \ \ {\it { is \ \  uniformly \ \ small \ \ in \ \ } } {\bf z}\,.\leqno (4.3)
$$

({\bf iii}) {\it  \,$\Vert \,P_{\bf z}\, G'\, ({\bf z}  + w) \Vert_\btd \ \le \ C$\, uniformly   in \,${\bf z}$} \ \ (\,{\it for all\,}  $\,\Vert w \Vert_{\,\btd} \,\le\, 1$)\,. \\[0.2in]
 {\it Then there exists $\,{\bar \varepsilon}_1 > 0$ such that for every $|\,\varepsilon|\le {\bar \varepsilon}_1\,,\,$ and for every $\,{\bf z} \in {\bf Z}$, the auxiliary equation}
$$
P_{\bf z}\circ I_\varepsilon' \,({\bf z} + w)\ = \ 0\leqno (4.4)
$$
{\it has a unique solution $w_\varepsilon ({\bf z})$ satisfying}
$$\Vert w_\varepsilon ({\bf z}) \Vert_{\,\btd} \to 0 \ \ \ \ {\it  uniformly \ \ as} \ \ |\,\varepsilon| \to 0\,.\leqno (4.5)
 $$

\hspace*{0.5in}In the next two sections (\S\,4\,b\, and \, \S\,4\,c)\,,\, we verify conditions {\bf (i)\,--\ (iii)}\, for the functional $I_\varepsilon$ defined in (2.12)\,.

\vspace*{0.15in}

{\bf \S\,4\,b. \ } {\it  Uniform invertibility of $\,I''_o \,({\bf z})$. } \ \
The following result is a direct consequence of Proposition 3.11 and Corollary 3.13. For the proof, see \cite{Progress-Book}, or \cite{Brendle-1}, or \,\S\,A.5.\,a\, in the e\,-Appendix.\\[0.2in]
{\bf Proposition 4.6.} \ \ {\it There exists a positive constant $\,{\bar c}_3$ such that for all $\,{\bf z} \in {\bf Z}$, }
$$
\Vert\, I''_o \,({\bf z})\, [f]\, \Vert_{\,\btd} \ \ge \ {\bar c}_3 \, \Vert f \Vert_{\,\btd} \ \ \ \ \ \ {\it{for \ \ all}} \ \ \ f \in {\bf W_z}\,. \leqno (4.7)
$$

\vspace*{0.2in}

%
%
%
%
%
%
%
{\bf \S\,4\,c. \ } {\it Uniform approximation and bounded projection.  }\ \
Likewise, the demonstrations of conditions ({\bf ii})  and ({\bf iii})
 in Lemma 4.1 are rather standard. We prefer to put the details  in \S\,A.5 and \S\,A.6 in the e\,-Appendix. Precisely, we show in  the e\,-Appendix that
$$
 R_{\,\bf z} (w) \ := \ I'_o\,({\bf z}  + w) \ - \ I''_o\,({\bf z} )\, [w] \ \ \ \ {\mbox{satisfies}} \ \ \  \Vert \,R_{\,\bf z}  (w) \Vert \ =\  o\,(\Vert  w \Vert_{ \btd} )\leqno (4.8) $$
 as $\Vert\,  w \Vert_{\,\btd} \to 0\,,\,$  uniformly  for $\,{\bf z}  \in {\bf Z}\,.$\,
Moreover,

 \vspace*{-0.175in}

$$
\Vert \,G' ({\bf z} + w) \Vert \ \le \ C \ \ \ \ \mfor {\mbox{all}} \ \ \ {\bf z}  \,\in\, {\bf Z} \ \ \ {\mbox{and}} \ \ w \ \ {\mbox{with}} \ \  \Vert\, w \Vert_{\,\btd}\ \le \ 1\,. \leqno (4.9)
$$

\vspace*{0.2in}

{\bf \S\,4\,d. \ } {\it  Uniform bounds in the expansions of $\,I_o$ and $\,G$\,.}\ \
We also have the following uniform bounds.
\begin{eqnarray*}
   (4.10)\   \ \ \ \ \ \   \ \ \ \ \ \   \ \ \     I_o\,({\bf z} + w_\varepsilon ({\bf z}))
        &\ = \ &     {\bar c}_4 \ + \   o\,(\Vert   w_\varepsilon ({\bf z})  \Vert\,)\,; \\[0.05in]
       \!\!\!(4.11) \ \  \ \ \ \ \   \ \ \ \ \ \  \  \ \ \, G ({\bf z} + w_\varepsilon ({\bf z}))     &= &G ({\bf z})
        + [ \,G' ({\bf z}) \,| \,w_\varepsilon ({\bf z})\,] \ + \ o\,(\Vert w_\varepsilon ({\bf z})  \Vert\,)\,; \\[0.05in]
          (4.12)\ \ \ \ \ \   \ \ \ \ \ \  \ \ \,  I_o' ({\bf z} \ + \ w_\varepsilon ({\bf z})) &\ = \ &    [\,I_o''({\bf z})\,| \,w_\varepsilon ({\bf z})\,]
           \ + \ o\,(\Vert  w_\varepsilon ({\bf z}) \Vert\,)\,;   \\[0.075in]
          (4.13) \ \ \ \ \ \   \ \ \ \ \ \  \ \  G' ({\bf z} \ + \  w_\varepsilon ({\bf z}))     &= &G' ({\bf z})
           + [ \,G'' ({\bf z}) \,| \,w_\varepsilon ({\bf z})\,]\  + \  o\,(\Vert w_\varepsilon
           ({\bf z})  \Vert\,)\,.\ \ \ \ \ \   \ \ \ \ \ \  \ \ \ \  \ \ \ \ \ \  \ \ \ \
         \end{eqnarray*}
Here $w_\varepsilon ({\bf z})$ is the solution to the auxiliary equation (4.4) described in Lemma 4.1, and
$$
{\bar c}_4\ = \ (n - 2) \int_{\R^n} \left({1\over {1 + |\,y|^2}} \right)^n. \leqno (4.14)
$$
We refer to \S\,A\,.5 and \S\,A\,.7 in the e\,-Appendix for more detail.

\vspace*{0.2in}

{\bf \S\,4\,e. \ } {\it  Uniform cancelation.}\\[0.15in]
 {\bf Proposition 4.15.} \ \  {\it There exists a positive constant $\,{\bar \varepsilon}_2$ such that for all $\,\varepsilon\,$ with $\,|\,\varepsilon|\le {\bar \varepsilon}_2\,,\,$ if $\Phi_\varepsilon$ has a sequence of critical points at $\,{\bf z}_{\,\varepsilon, \,i} \in {\bf Z}$\,,\, $i\ = \ 1, \ 2, \cdot\, \cdot\, \cdot\,\,,$\, } {\it then $\,{\bf z}_{\,\varepsilon, \,i}\, +\, w_\varepsilon ({\bf z}_{\,\varepsilon, \,i})$ are critical points of $\,I_\varepsilon\,$ for all $\,i\,.\,$  Here $\,w_\varepsilon ({\bf z}_{\,\varepsilon, \,i})\,$ is the solution of the auxiliary equation\,} (4.4) {\it as described in  Lemma\,} 4.1.\\[0.2in]
  \hspace*{0.5in}See \S\,A\,.\,8 in the e\,-\,Appendix for the argument, which  closely proceeds as in   the proof of Theorem 2.12 in\, \cite{Progress-Book}.
\vspace*{0.2in}

{\bf \S\,4\,f. \ } {\it Uniform levelness.}\ \
We show that the solutions of the auxiliary equation (4.4) as described in Lemma 4.1 stay  uniformly close to ${\bf Z}\,$,\, and `` is almost parallel\," to $\,{\bf Z}$\,.\, Precisely we have the following result.

\vspace*{0.2in}

{\bf Lemma 4.\,16.} \ \ {\it For $\,n \ge 6$\,,\, let $\,w_\varepsilon ({\bf z})\,$ be the unique solution of the auxiliary equation\,} (4.4) {\it as described in Lemma\,} 4.1\,.\, {\it There exists a positive number $\,{\bar \varepsilon}_3\,$ such that for every $\,varepsilon\,$ satisfying  $\,|\,\varepsilon|\le {\bar\varepsilon}_3$\,,\, we have    }

\vspace*{-0.3in}

\begin{eqnarray*}
(4.17) \ \ \Vert\, w_\varepsilon\, ({\bf z}) \Vert_{\,\btd}  \,\le \, {\bar C}_4 \,|\,\varepsilon| \ \ \  {\it{and}} \ \ \  \Vert\, D_\ell\, w_\varepsilon\, ({\bf z}) \Vert_{\,\btd} \,\le \,{{{\bar C}_4}\over \lambda}  \,|\,\varepsilon|^{4\over {n - 2}}    \ \ \    {\it for } \ \ \ell\ = \ 0\,, \ 1,\, \cdot \cdot \cdot\,, \ n\,.
\end{eqnarray*}
{\it Here $\,{\bar C}_4\,$ is a positive constant independent on   $\,{\bf z} \, = \, V_{\lambda\,, \ \xi} \, \in \,  {\bf Z}\,.$}
\vspace*{0.2in}

{\bf Proof.} \ \
We argue as in the proof of Lemma 2.11 in \cite{Progress-Book}, pp. 21. Let us begin with
$$  I'_\varepsilon ({\bf z} + w_\varepsilon ({\bf z})) \ \ \ \ \ \ \ \ \ \ \left[\,\ = \ I'_o \,({\bf z} + w_\varepsilon ({\bf z})) + \varepsilon\cdot G' \, ({\bf z} + w_\varepsilon ({\bf z})) \,\right] \leqno (4.18)
$$
\begin{eqnarray*}
&\ &  \!\!\!\!\!\!\!\!\!\!\!\!\!\!\!\!\!\!  = \   I'_o \,({\bf z})\  + \ I''_o \,({\bf z})\, [\,w_\varepsilon ({\bf z})\,]
+\ \varepsilon\,\{\,G' ({\bf z}) + \,G'' ({\bf z} ) [\,w_\varepsilon ({\bf z})\,]\,\} \ + \ o\, (\Vert \,w_\varepsilon ({\bf z}) \Vert_{\,\btd})\\[0.075in]
& \ & \ \ \ \ \ \ \ \ \ \ \ \ \ \ \ \ \ \ \ \ \ \ \ \ \ \ \ \ \ \ \ \ \ \ \ \ \ \ \ \   \ \ \  \ \ \ \ \ \ \ \ \ \ \ \ \ \ \ \     [\,{\mbox{via}} \ (4.12) \ \& \ (4.13)\,] \\[0.075in]
&\ &  \!\!\!\!\!\!\!\!\!\!\!\!\!\!\!\!\!\!    = \ \varepsilon\,G' ({\bf z}) \ + \ I''_o \,({\bf z})\, [w_\varepsilon ({\bf z})\,] \ +\ \varepsilon\,G'' ({\bf z} ) \,[w_\varepsilon ({\bf z})\,] \ + \ o \,(\Vert \,w_\varepsilon ({\bf z}) \Vert_{\,\btd}) \ \ \ {\mbox{uniformly \ \ in}} \ \ {\bf z}\,.
\end{eqnarray*}
From the auxiliary equation (4.4), we have the following.

\vspace*{-0.3in}

\begin{eqnarray*}
(4.19) \  & \ & P_{\bf z} \circ I'_\varepsilon ({\bf z} + w_\varepsilon ({\bf z}))\ = \ 0  \\[0.075in]
& \ &\!\!\!\!\!\!\!\!\!\!\!\!\!\!\!\!\!\!\!\!\!\!\!\!\!\! \Longrightarrow  \    \varepsilon\,P_{\bf z} \circ G' ({\bf z}) \,+\, P_{\bf z} \circ I''_o \,({\bf z})\, [\,w_\varepsilon ({\bf z})\,] +\varepsilon\, P_{\bf z} \circ G'' ({\bf z} ) \, [w_\varepsilon ({\bf z})\,] \,+\, o \,(\Vert\, w_\varepsilon ({\bf z}) \Vert_{\,\btd})\,=\,0\\[0.075in]
& \ &\!\!\!\!\!\!\!\!\!\!\!\!\!\!\!\!\!\!\!\!\!\!\!\!\!\! \Longrightarrow  \ P_{\bf z} \circ G' ({\bf z})  +\,  P_{\bf z} \circ I''_o \,({\bf z}) \, [\,\varepsilon^{-1}\, w_\varepsilon ({\bf z})\,]  +\, P_{\bf z} \circ G'' ({\bf z} ) [\,w_\varepsilon ({\bf z})\,]  + \, o\, ( \Vert \,\varepsilon^{-1} w_\varepsilon ({\bf z}) \Vert_{\,\btd} )\,=\, 0. \ \ \ \
\end{eqnarray*}
Via (4.5), we already know that $\,\Vert \,w_\varepsilon ({\bf z}) \Vert_{\,\btd}\,$ is uniformly small, hence
$$
P_{\bf z} \circ I''_o \,({\bf z})\, [\,\varepsilon^{-1} w_\varepsilon ({\bf z})\,]\ = \ P_{\bf z} \circ G' ({\bf z})  \ + \  o \,(1) \ + \ o\, ( \Vert \, \varepsilon^{-1} w_\varepsilon ({\bf z}) \Vert_{\,\btd} )\,. \leqno (4.20)
$$
Compare with (4.19), the small order terms in (4.20) are in ${\bf {\bf W_z}}$\,.\, The operator is identified via Riesz Representation Theorem with a vector in ${\cal D}^{1,\ 2}$\,.\, We check in \S\,A.9 in the e\,-Appendix that $G' ({\bf z})$ (as an operator) is uniformly bounded.
 Since
$$
P_{\bf z} \circ I''_o \,({\bf z}) \,: {\bf W_z} \to {\bf W_z}
$$
is uniformly invertible [\,refer to Proposition 4.6, (2.23)  and (2.26)\,]\,,\, it follows  that

\vspace*{-0.15in}

$$
\Vert \, \varepsilon^{-1} \cdot w_\varepsilon ({\bf z}) \Vert_{\,\btd}
$$

\vspace*{-0.15in}

is uniformly bounded\,.\bk
As for the second estimates in (4.17), Lemma 2.11 in \cite{Progress-Book} asserts that $w_\varepsilon ({\bf z})$ is $\,C^1\,$ in $\,{\bf z}$\,.\, Via the parametrization (2.19) it is differentiable ($C^1$) with respective to $\lambda$ and $\xi$\,.\, We make use of the inequality\,:

\vspace*{-0.2in}

$$
n \ \ge\  6 \ \ \Longleftrightarrow \ \ {4\over {n - 2}} \ \le \ 1 \ \Longrightarrow \ \bigg\vert\, |\,a + b\,|^{4\over {n - 2}} - a^{4\over {n - 2}}\bigg\vert \ \le \ |\,b|^{4\over {n - 2}} \leqno (4.21)
$$
for $\,a \, > \, 0$\, and $\,b \in \R\,.\,$
Refer to Lemma 6.18 in \cite{Progress-Book}, pp. 93 \,(\,for the definition of $2^*$ used in the book, see (6.9), pp. 75, loc. cit.)\,.\bk
Differentiating the auxiliary equation (4.4), as $P_{\bf z}$ is linear, we  obtain
$$ P_{\bf z} \circ I''_o \,({\bf z} + w_\varepsilon ({\bf z})) \,[\,D_\ell \,{\bf z} + D_\ell \,w_\varepsilon ({\bf z})\,] \ +\  \varepsilon\, P_{\bf z} \circ G'' ({\bf z} + w_\varepsilon ({\bf z})) \,[\,D_\ell\, \,{\bf z} + D_\ell \,w_\varepsilon ({\bf z})\,]\ = \ 0  \leqno (4.22)
$$
for $\,\ell\ = \ 0\,,\ 1, \ \cdot \cdot \cdot, \ n\,.\,$
  From (2.20) we recognize
\begin{eqnarray*}
(4.23)  \!\!\!\!\! &\ & (I_o'' ({\bf z} + w_\varepsilon ({\bf z}))\, [D_\ell \,{\bf z} + D_\ell \,w_\varepsilon ({\bf z})\,] \  h) \ \ \ \ \ \ \ \ \ \ \ \ \ \ \ \ \ \ \ \ \ \ \ \ \ \  \ \ \ \ \ \ \ \ \ \ \ \ \  \ \  (h \,\in\, {\cal D}^{1\,, \ 2}\,)\\
 & \ & \\
&  \,& \!\!\!\!\! \!\!\!\!\!\!\!\!\!\!= \int_{\R^n} \left[\,  \langle \btd \,[D_\ell \,{\bf z} + D_\ell \,w_\varepsilon ({\bf z})\,]\,,\, \btd \, h \rangle \,-\, n\,(n + 2) \,[\,{\bf z} + w_\varepsilon ({\bf z})\,]_+^{{4}\over {n - 2}}\, \,[\,D_\ell \,{\bf z} + D_\ell \,w_\varepsilon ({\bf z})\,]\cdot h\, \,\right]\\[0.15in]
& \ &  \ = \ \ \int_{\R^n}  \langle \btd \,[ D_\ell \,w_\varepsilon ({\bf z})\,]\,,\, \btd \, h \rangle - n\,(n + 2) \,\int_{\R^n}{\bf z}^{{4}\over {n - 2}}\, [ \,D_\ell \,w_\varepsilon ({\bf z})\,]\cdot h\, \,\\
& \ & \ \ \ \ \ \ \ \ \ \ +\, \int_{\R^n}  \langle \btd \,[ D_\ell \,{\bf z}]\,,\, \btd \, h \rangle - n\,(n + 2) \,\int_{\R^n}{\bf z}^{{4}\over {n - 2}}\, [ \,D_\ell \,{\bf z}]\cdot h\\
& \ & \ \ \ \ \ \ \ \ \ \ \ \ \ \ -\,n\,(n + 2) \int_{\R^n} \left[({\bf z}+ w)_+^{{4}\over {n - 2}} - {\bf z}^{{4}\over {n - 2}} \right]\, [ \,D_\ell \,{\bf z} + D_\ell \,w_\varepsilon ({\bf z})\,]\cdot h\,\\
& \ & \\
& \ &  \ = \ (I_o'' ({\bf z})\, [D_\ell \,w_\varepsilon ({\bf z})\,] \  h) \ \ + \ \ \ \   (I_o'' ({\bf z})\, [\,D_\ell \, z ] \  h) \\ & \ &
 \hspace*{2.4in}\ \  \ \  (\,\uparrow \ \ = \ 0 \ \ {\mbox{via \ \ Lemma \ 2.24}})\\[0.075in]
& \ & \ \ \ \ \ \  \ \ \  \ \ \ -\,n\,(n + 2) \int_{\R^n} \left[({\bf z}+ w_\varepsilon ({\bf z}))_+^{{4}\over {n - 2}} - {\bf z}^{{4}\over {n - 2}} \right]\, [ \,D_\ell \ {\bf z} + D_\ell \,w_\varepsilon ({\bf z})\,]\cdot h\,.
\end{eqnarray*}
Consider the last term in (4.23)\,:
\begin{eqnarray*}
(4.24) & \ & \bigg\vert\,\int_{\R^n} \left[({\bf z}+ w_\varepsilon ({\bf z}))_+^{{4}\over {n - 2}} \ - \ {\bf z}^{4\over {n - 2}} \right]\, [ \,D_\ell \,{\bf z} \ + \ D_\ell \,w_\varepsilon ({\bf z})\,]\cdot h \,\bigg\vert \\
& \le & \int_{\R^n} \bigg\vert\,({\bf z}+ w_\varepsilon ({\bf z}))_+^{{4}\over {n - 2}} \ - \ {\bf z}^{4\over {n - 2}} \bigg\vert\,  \,( \,|\,D_\ell \,{\bf z} | \ + \ |\,D_\ell \,w_\varepsilon ({\bf z})|\,) \cdot |\,h|\\
& \ &     \!\!\!\!\!\!  \bigg\{\, {\mbox{if \ \ at \ \ a \ \ point \ \ }}[\,{\bf z}+ w_\varepsilon ({\bf z})\,]_+\ = \ 0 \ \ \Longrightarrow \ \ |\,w_\varepsilon ({\bf z})| \ \ge \  {\bf z} =\  V_{\lambda\,, \ \xi} > 0\, \,; \\
& \ & \!\!\!\!\!\!\!\!\!\!\!\!    {\mbox{whereas\ \,  at \   a \   point \   }} \     [\,{\bf z}+ w_\varepsilon ({\bf z})\,]_+ \ > \ 0  \ \Longrightarrow \ \ ({\bf z}+ w_\varepsilon ({\bf z}))_+\ = \ |\,{\bf z} \ + \ w_\varepsilon ({\bf z})|\,\bigg\}\\
& \le & C_1 \int_{\R^n}  |\,w_\varepsilon ({\bf z})|^{{4}\over {n - 2}}  \,  [ \,|\,D_\ell \,{\bf z} | \ + \ |D_\ell \,w_\varepsilon ({\bf z})|\,]\cdot |\,h| \ \ \ \ \ \ \ \ \ \ \ \ \ \ \ \ \ \ \  \ \ [\,{\mbox{using}} \ \ (4.21)\,]\\
& \le & C_2\, \Vert \,h \Vert_{\,\btd} \cdot \left( \Vert D_\ell \,{\bf z} \Vert_{\,\btd}  \ + \  \Vert D_\ell \,w_\varepsilon ({\bf z}) \Vert_{\,\btd} \right)\left( \int_{\R^n}  |\,w_\varepsilon ({\bf z})|^{{2n}\over {n - 2}} \right)^{\!\!{\!2\over n}} \\
& \ & \ \ \ \ \ \ \       ({\mbox{apply \   H\"older's \   inequality \   twice\,; \    cf. \  \S\, A.5 \ in \ the \ e\,-\,Appendix}})\\[0.075in]
& \le & C_3\, \Vert \,h \Vert_{\,\btd} \cdot \left( \Vert D_\ell \,{\bf z} \Vert_{\,\btd}  +  \Vert D_\ell \,w_\varepsilon ({\bf z}) \Vert_{\,\btd} \right) \cdot \Vert   w_\varepsilon ({\bf z}) \Vert_{\,\btd}^{4\over {n - 2}}\ \ \      \left[\  {2\over n}\, = \left( {{n - 2}\over {2n}} \right)\!\cdot\! {4\over {n - 2}}\right]\\
& \le & C_4\,  \Vert\, h \Vert_{\,\btd} \cdot \left[ \,{1\over \lambda} + \Vert \,D_\ell \,w_\varepsilon ({\bf z}) \Vert_{\,\btd} \right]\cdot|\,\varepsilon|^{4\over {n - 2}} \ \ \ \ \ \ \ \ \ \ \ \ \ \ \ \  ({\mbox{as \ \ ususal}} \ \ {\bf z}\ = \ V_{\lambda\,, \ \xi}\,)\,.
\end{eqnarray*}
In the above, we also apply (3.1)\,--\ (3.3), and the first estimate in (4.17).\bk
For the second term in the right hand side of (4.22), from (2.14), we have
\begin{eqnarray*}
& \ & (\,G'' (f) \,[\,h] \, \phi) \ = \   - \ {\tilde c}_n  \cdot {{n + 2}\over {n - 2}} \cdot \!  \int_{\R^n} H \cdot f_+^{4\over {n - 2}}\, g \cdot h \ \   \mfor \ f\,, \ h\,, \ \phi \in {\cal D}^{1, \  2}\\[0.1in]
 \Longrightarrow      & \ & (\,G'' \,({\bf z}+w_\varepsilon ({\bf z})) \,[\,D_\ell \,{\bf z} + D_\ell \,w_\varepsilon ({\bf z})\,] \  h)
 \\[0.1in]
 & \ & \!\!\!\!\!\!\!=\ {{n + 2}\over {n - 2}} \ {\tilde c}_n\int_{\R^n} H \cdot [({\bf z}+w_\varepsilon ({\bf z})\,]^{4\over {n - 2}}_+\,\,(D_\ell \,{\bf z} + D_\ell \,w_\varepsilon ({\bf z})) \cdot h\,.
\end{eqnarray*}
Let us continue as in (4.24)\,:
\begin{eqnarray*}
(4.25) \ \  \ \  & \ & \bigg\vert \int_{\R^n} H \cdot [\,({\bf z}+w_\varepsilon ({\bf z})\,]^{4\over {n - 2}}_+\,\,(D_\ell \,{\bf z} + D_\ell \,w_\varepsilon ({\bf z})) \cdot h\,\bigg\vert\\[0.1in]
& \le & C_1 \int_{\R^n} \left[ \,{\bf z}^{4\over {n - 2}} + [\,w_\varepsilon ({\bf z})\,]^{4\over {n - 2}}\right]\cdot(| D_\ell \,{\bf z} |+ |D_\ell \,w_\varepsilon ({\bf z})|\,) \cdot| \,h|\\[0.1in]
& \le & C_2\, \left[ {1\over \lambda} + \Vert D_\ell \,w_\varepsilon ({\bf z}) \Vert_{\,\btd} + |\,\varepsilon|^{4\over {n - 2}} \left({1\over \lambda} + \Vert D_\ell \,w_\varepsilon ({\bf z}) \Vert_{\,\btd}  \right) \right]\cdot  \Vert \,h \Vert_{\,\btd}\,.\ \  \ \ \ \  \ \ \ \  \ \  \ \  \ \  \ \ \ \  \ \  \ \  \ \
\end{eqnarray*}
Combining the estimates (4.24) and (4.25), we obtain
\begin{eqnarray*}
& \ &  P_{\bf z} \circ I''_o \,({\bf z}) [D_\ell \,w_\varepsilon ({\bf z})\,]  \ = \ O \left( {{ |\,\varepsilon|^{4\over {n - 2}} }\over \lambda}\right) + O \left( |\,\varepsilon|^{4\over {n - 2}} \cdot \Vert D_\ell \, w_\varepsilon ({\bf z}) \Vert_{\,\btd}\right) \mfor n \ge 6\,.
\end{eqnarray*}
Write
$$
D_\ell \,w_\varepsilon ({\bf z})\ = \ [D_\ell \,w_\varepsilon ({\bf z})\,]_\perp + [D_\ell \,w_\varepsilon ({\bf z})\,]_{\,T_{\bf z} ({\bf Z})}\,.
$$
It follows from (2.21) and (2.23) that
$$
\Vert \,I''_o \,({\bf z}) \,([\,D_\ell \,w_\varepsilon ({\bf z})\,]_{\,T_{\bf z} ({\bf Z})}) \Vert_{\,\btd} \le C \Vert \,D_\ell \,w_\varepsilon ({\bf z})\,]_{\,T_{\bf z} ({\bf Z})}\Vert_{\,\btd}  \ \ \ \ \ \ \  \ \ \ \ \ \ \ \ \ \ \ \ \ \, ({\mbox{uniformly \ \ in}} \ \ {\bf z})\,.
$$
From (4.5) in Lemma 4.1,  we obtain [\,see also (A.8.1) in   the e\,-\,Appendix\,],
$$
\ \ \ \ \ \ \ \     \Vert   \,[\,D_\ell \,w_\varepsilon ({\bf z})\,]_{\,T_{\bf z} ({\bf Z})} \Vert_{\,\btd} \ \le \  C_1\cdot{{|\,\varepsilon|}\over {\lambda}}\ \ \ \ \ \ \ \ \ \ \ \ \ \ \ \ \ \ \ \ \    \ \ \ \ \ \ \ \ \ \ \ \ \ \  \  \ \ \    \ ({\mbox{uniformly \ \ in}} \ \ {\bf z})
$$
$$
\Longrightarrow \ \ \Vert\, I''_o \,({\bf z}) \,([\,D_\ell \,w_\varepsilon ({\bf z})\,]_{\,T_{\bf z} ({\bf Z})}) \Vert_{\,\btd} \ \le \ C_2\cdot {{|\,\varepsilon|}\over {\lambda}} \ \ \ \ \ \ \ \ \ \ \ \ \ \
\ \ \ \ \   \  \ \ \ \ \ \ \ \ \ \ \  \ \, ({\mbox{uniformly \ \ in}} \ \ {\bf z})\,.
$$
Hence
$$
P_{\bf z} \,\circ\, I''_o \,({\bf z})\, [\,D_\ell \,w_\varepsilon ({\bf z})\,]_\perp  \ = \ O \left( {{ |\,\varepsilon|^{4\over {n - 2}} }\over \lambda}\right) + O \left( |\,\varepsilon|^{4\over {n - 2}} \cdot \Vert\, [\,D_\ell \,w_\varepsilon ({\bf z})\,]_\perp \Vert_{\,\btd}\right)
$$
uniformly in $\,{\bf z}\,$.\,
[\,$I''_o \,({\bf z})$\,,\, considered as a linear map from ${\cal D}^{1, \  2}$\, to itself\,,\, is uniformly bounded in ${\bf z}$.]
 We obtain the desired estimate on $\Vert\, I''_o \,({\bf z}) \,(\,[\,D_\ell \,w_\varepsilon ({\bf z})\,]_\perp\,) \Vert_{\,\btd}\,$ by using  (2.26) and the uniform invertibility (4.7)\,.\qedwh
{\bf \S\,4\,g. \ } {\it Uniform existence for stable critical points.}\ \
The key merit of the following result is that we can allow  $\bar \lambda \to 0^+$.

\vspace*{0.2in}

 {\bf Theorem 4.26.}   \ \ {\it For $n \ge 6$\,,\, there exists positive constants $\,{\bar C}_5$ and $\,{\bar\varepsilon}_4$ such that the following holds. Suppose that $\,{\bar {\bf z}}\ = \ V_{{\bar \lambda}\,, \  {\bar \xi}} \in {\bf Z}\,$ is a stable critical point for $\,G_{|_{\bf Z}}\,,\,$ and there exists $\,\bar{\rho} > 0$ such that}
 $$
 B_{{\bar \lambda}\,, \,\, {\bar \xi}} \ (\bar \rho) \subset \R^+ \times \R^n  \ \ \ \ {\it{and}} \ \ \ \ \min_{(\,\lambda, \,\, \xi) \, \in \,  \partial B_{{\bar \lambda}\,, \  {\bar \xi}} \, (\bar \rho) } \,|\btd \,G_{|_{\bf Z}} (\,\lambda, \, \xi)| > {{{\bar C}_5}\over \lambda}\,. \leqno (4.27)
 $$
 [\,{\it  Here $\,\partial B_{{\bar \lambda}\,, \, {\bar \xi}} \,(\bar \rho)  $ is the boundary of the ball.\,}]\ \,{\it  Then for $|\,\varepsilon| < {\bar \varepsilon}_4$\,,\, the functional $\,{\bf I}_\varepsilon\,$ has a critical point at \,  $ \tilde{\bf z} + w_\varepsilon \,(\tilde{\bf z})\,$ with $\,\tilde{\bf z}\ = \ V_{\tilde\lambda\,,\,\, \tilde\xi}$\,,\, where $(\tilde\lambda\,,\,\, \tilde\xi) \in B_{{\bar \lambda}\,, \,\, {\bar \xi}} \ (\bar \rho)$\,.\,}  {\it Here $\,w_\varepsilon \, (\tilde{\bf z} )\,$ is the solution  of the auxiliary equation\,} (4.4) {\it as described in   Lemma\,} 4.1.\\[0.2in]
 {\bf Proof.} \ \
 We first verify the claim that, for all ${\bf z} \in {\bf Z}\,,\,$  we have
 $$
 D_\ell \, \Phi_\varepsilon  ({\bf z})  \ = \ \varepsilon\cdot [\, D_\ell \, G ({\bf z})\,]   \ + \ o \left({{\varepsilon}\over {\lambda}}\right) \ \ \ \ \ \ \ {\mbox{uniformly \ \ as}} \ \ \ \varepsilon \,\to 0\,.\leqno (4.28)
 $$
Here\, $\ell\ = \ 0, \ 1, \cdot \cdot \cdot, \ n\,.\,$
We start with  the definition in (2.12)\,:
\begin{eqnarray*}
& \ & \Phi_\varepsilon ({\bf z})  \ = \   I_\varepsilon \,({\bf z} + w_\varepsilon ({\bf z}))\\[0.075in]
&\ = \ & \int_{\R^n} \left\{\, {1\over 2}\, \big\langle \btd_o \,({\bf z} + w_\varepsilon ({\bf z}))\,,\, \btd_o\,({\bf z} + w_\varepsilon ({\bf z})) \,\big\rangle \ - \  {{n - 2}\over {2n}} \cdot  n\,(n - 2) \,[\,z + w_\varepsilon ({\bf z})\,]_+^{{2n}\over {n - 2}} \ \right\}\ \ \\
 & \ & \ \ \ \ \ \ \ \ \ \ \ \  \ \ \ \ \ \ \ \ \ \ \ \ + \, \varepsilon  \cdot  {\bar c}_{-1} \int_{\R^n} H \,[\,{\bf z} + w_\varepsilon ({\bf z})\,]_+^{{2n}\over {n - 2}}\end{eqnarray*}
 \begin{eqnarray*}
\Longrightarrow  & \ & D_\ell \,\Phi_\varepsilon ({\bf z})=   \int_{\R^n} \bigg\{\,   \big\langle \btd_o \,[\,{\bf z} \, + \, w_\varepsilon ({\bf z})\,]\,,\ \, \btd_o\,(D_\ell \,{\bf z} \ + \  D_\ell \,w_\varepsilon ({\bf z})) \,\big\rangle
 \\
 & \ & \ \ \ \ \ \ \ \ \ \ \ \ \ \ \ \ \ \ \ \  \ \ \ \ \ \ \ \ \left. -   \  n(n - 2) \,(D_\ell \,{\bf z} \, + \, D_\ell \,w_\varepsilon ({\bf z}))\cdot [\,{\bf z}\, + \,  w_\varepsilon ({\bf z})\,]_+^{{n+2}\over {n - 2}} \,\right\}\\
 & \ & \ \ \ \  \ \ \ \ \ \ \ \ \ \ \ \ \ \ \ \ \ \ \ \ \ \ \ \ \ \ \ \ \ \ \ \  - \, \varepsilon \cdot  {\tilde c}_n \int_{\R^n} H \cdot (D_\ell \,{\bf z} \, + \, D_\ell \,w_\varepsilon ({\bf z})) \cdot [\,{\bf z} \, + \, w_\varepsilon ({\bf z})\,]_+^{{n + 2}\over {n - 2}} \ \ \ \ \ \ \ \
\end{eqnarray*}

 \vspace*{-0.3in}

 \begin{eqnarray*}
\Longrightarrow  & \ & D_\ell \,\Phi_\varepsilon ({\bf z})=
 \int_{\R^n} \left\{\,   \langle \btd_o \,{\bf z}\,,\, \btd_o\, D_\ell \,{\bf z}   \rangle \,- \,n \,(n - 2) (D_\ell \,{\bf z})\,{\bf z}^{{n + 2}\over {n - 2}}\, \right\} \ \ \left(\ = \ 0 \right)\\
 & \ & \ \ \ \ \ \ \ \ \ \ \ \ \ \ \ \ \ \ \ \ \ \ \ \  \ \ \ \ \ \ \ \ \ \ \ \ \ \ \ \ \ \ \ \ \ \ \,
   [\,{\mbox{from \ \ (1.3) \ \,\&\, \ (1.4)\,}} \ \ \uparrow \ ]\\
 & \ & \ \  \ \ \ \ \ +\  \int_{\R^n}   \big\langle \btd_o \,{\bf z}\,,\, \btd_o\, D_\ell \,w_\varepsilon ({\bf z})   \big\rangle \ +\  \int_{\R^n}   \big\langle \btd_o \,w_\varepsilon ({\bf z})\,,\, \btd_o\, D_\ell \, {\bf z}    \big \rangle \\
 & \ & \ \ \ \ \ \ \ \ \ \ \ \ \ \ \ \ +\ \int_{\R^n}   \big\langle \btd_o \,w_\varepsilon ({\bf z})\,,\, \btd_o\, D_\ell \,w_\varepsilon ({\bf z})   \big\rangle\\
 & \ & \\
 & \ & \!\!\!\!\!\!\! \!\!\!\!\!\!\!\!\!\!\!-\, n\,(n - 2) \,\int_{\R^n} D_\ell \,{\bf z}  \,\left[ ({\bf z} \,+\, w_\varepsilon ({\bf z}))_+^{{n+2}\over {n - 2}} - {\bf z}^{{n + 2}\over {n - 2}} \right]- \,n\,(n - 2) \,\int_{\R^n} D_\ell \,w_\varepsilon ({\bf z})  \, [\,{\bf z} + w_\varepsilon ({\bf z})\,]_+^{{n+2}\over {n - 2}}  \\
 & \ & \ \ \ \  \   - \, \varepsilon \cdot {\tilde c}_n   \int_{\R^n} H \,[\, D_\ell \,{\bf z} ]\cdot {\bf z}^{{n + 2}\over {n - 2}} \,- \,\, \varepsilon \cdot {\tilde c}_n \int_{\R^n} H\,[ D_\ell \,{\bf z} ]\cdot \left[ \,({\bf z}+ w_\varepsilon ({\bf z}))_+^{{n + 2}\over {n - 2}} - {\bf z}^{{n + 2}\over {n - 2}}\right]\\
 & \ & \ \ \ \  \ \ \ \ \ \ \ \ \ \ \   - \, \varepsilon \cdot {\tilde c}_n   \int_{\R^n} H \cdot  D_\ell \,w_\varepsilon ({\bf z})  [\,{\bf z} + w_\varepsilon ({\bf z})\,]_+^{{n + 2}\over {n - 2}}
\end{eqnarray*}
for $\,\ell\ = \ 0, \ 1, \cdot \cdot \cdot, \ n\,.\,$
Together with  equation (1.4) and expressions (3.1)\,--\ (3.3), we have\\[0.1in]
(4.29)
\begin{eqnarray*}
   D_\ell \,\Phi_\varepsilon ({\bf z})&= & \varepsilon D_\ell \, G_{|_{\bf Z} } ({\bf z}) \, +\, \int_{\R^n}   \langle \btd_o \,w_\varepsilon ({\bf z})\,,\, \btd_o\, D_\ell \, {\bf z} \,   \rangle -n\,(n + 2) \cdot {\bf z}^{4\over {n - 2}} \,[\,D_\ell\, {\bf z}] \cdot w_\varepsilon ({\bf z}) \\
 & \ &  \ \ \ \ \ \ \ \ \ \ \ \ \ \ \ \ \ \ \ \ \ \ \ \ \ \ \ \ \ \ \ \ \
   [\,\uparrow \   \,= I''_o \,({\bf z}) \,[D_\ell\, {\bf z}\,] \,w_\varepsilon ({\bf z})= 0  \ ( {\mbox{Lemma 2.24}})\,]\,\\
    & \ & \\
  & \ & \ \ \ \ \   \ -\, n(n - 2) \,\int_{\R^n} D_\ell \,{\bf z}   \left[ ({\bf z} + w_\varepsilon  ({\bf z}))_+^{{n+2}\over {n - 2}} - {\bf z}^{{n + 2}\over {n - 2}} - {{n + 2}\over {n - 2}}\,\, {\bf z}^{4\over {n - 2}}   \,w_\varepsilon  ({\bf z})  \right]\\
  & \ & \   -\, n(n - 2) \,\int_{\R^n} D_\ell \,w_\varepsilon ({\bf z})  \, ({\bf z} + w_\varepsilon ({\bf z}))_+^{{n+2}\over {n - 2}}   +\, \int_{\R^n}   \langle \,\btd_o \,{\bf z}\,,\, \btd_o\, D_\ell \,w_\varepsilon ({\bf z}) \,  \rangle \\
   & \ & \ \ \ \ \ \ \ \ +\, \int_{\R^n}   \langle \btd_o \,w_\varepsilon ({\bf z})\,,\, \btd_o\, D_\ell \,w_\varepsilon ({\bf z})   \rangle\\
   & \ & \ \ \ \ \ \ \ \ \ \ \ \ \ \ - \,  \varepsilon \cdot {\tilde c}_n   \int_{\R^n} H \,[ D_\ell \,{\bf z} ]\cdot \left[ ({\bf z}+ w_\varepsilon ({\bf z}))_+^{{n + 2}\over {n - 2}} - {\bf z}^{{n + 2}\over {n - 2}}\right]\\& \ & \ \ \ \  \ \ \ \ \ \ \ \ \ \ \ \ \ \ \ - \, \varepsilon \cdot {\tilde c}_n \int_{\R^n} H \,  D_\ell \,w_\varepsilon ({\bf z})  [\,{\bf z} + w_\varepsilon ({\bf z})\,]_+^{{n + 2}\over {n - 2}}\\[0.1in]
   &\ = \ &   \varepsilon D_\ell \, G_{|_{\bf Z} } ({\bf z}) \ + \ {\bf I} \ + \ {\bf II} \ + \ {\bf III}  \ + \ {\bf IV}\,,  \end{eqnarray*}
   
   \vspace*{-0.1in}
   
   where
   \begin{eqnarray*}
    {\bf I} & =& -\,n\,(n - 2) \,\int_{\R^n} D_\ell \,\,{\bf z}  \,\left[ ({\bf z} + w _\varepsilon ({\bf z}))_+^{{n+2}\over {n - 2}} - {\bf z}^{{n + 2}\over {n - 2}} - {{n + 2}\over {n - 2}} \ {\bf z}^{4\over {n - 2}}\, w_\varepsilon ({\bf z})  \right]\,,\\[0.1in]
     {\bf II} & =&   -\, n\,(n - 2) \,\int_{\R^n} D_\ell \,w_\varepsilon ({\bf z})  \,\left[ ({\bf z} + w_\varepsilon ({\bf z}))_+^{{n+2}\over {n - 2}} - {\bf z}^{{n + 2}\over {n - 2}} - {{n + 2}\over {n - 2}} {\bf z}^{4\over {n - 2}} w_\varepsilon ({\bf z})  \right]\,, \\[0.1in]
     {\bf III} & =&   \, (I''_o \,({\bf z})\, [w_\varepsilon ({\bf z})\,]\, ( D_\ell\, w_\varepsilon ({\bf z})) - \ \varepsilon \cdot {\tilde c}_n   \int_{\R^n} H \,[ D_\ell \,{\bf z} ]\cdot \left[ ({\bf z}+ w_\varepsilon ({\bf z}))_+^{{n + 2}\over {n - 2}} - {\bf z}^{{n + 2}\over {n - 2}}\right]\,,\\[0.1in]
     {\bf IV} & =&   - \ \varepsilon \cdot {\tilde c}_n    \int_{\R^n} H \,  D_\ell \,w_\varepsilon ({\bf z})  [\,{\bf z} + w_\varepsilon ({\bf z})\,]_+^{{n + 2}\over {n - 2}}\ .
 \end{eqnarray*}
 Here we use (1.4) and (2.21)\,.\,  For the terms marked ({\bf I}) and ({\bf II}) in (4.29), see (A.5.5) in \S\,A\,.5 in the e\,-\,Appendix\,.\, The first term in ({\bf III}) can be estimated by using the uniform bound on $\Vert\,I_o'' ({\bf z})\Vert\,$ and  the bounds in Lemma 4.16\,.\, Whereas the second term in  ({\bf III}) is handled as in (4.24)\,,\, and ({\bf IV}) as in (A.6.2) in the e\,-Appendix.\, Together with (3.2) and (3.3),  we obtain
 
 \vspace*{-0.35in}
 
 \begin{eqnarray*}
 D_{\ell}\, \Phi_\varepsilon ({\bf z}) \ =\  \varepsilon D_\ell \, G_{|_{\bf Z} } ({\bf z}) \,+\, O\left(\! |\,\varepsilon|\cdot {{ |\,\varepsilon|}\over {\lambda}}\,\right) \,+\, O \left(\! |\,\varepsilon| \cdot {{ |\,\varepsilon|^{4\over {n - 2}} }\over {\lambda}}\, \right)    \mfor \ \ \ell\ = \ 0, \ 1, \cdot \cdot \cdot, \ n\,.
 \end{eqnarray*}
 Hence we establish the claim.
Thus we can find \,${\bar\varepsilon}_4\,$ small so that  for $|\,\varepsilon| < {\bar\varepsilon}_4$\,,\, we have
$$
\Vert\,\btd_{(\,\lambda\,, \, \xi)} \Phi_\varepsilon ({\bf z})\ -\ \varepsilon \btd_{(\,\lambda\,, \, \xi)} G_{|_{\bf Z}} ({\bf z}) \Vert \ \le \ c_1\cdot  {{ |\,\varepsilon|}\over {\lambda}}  \ < \ {\bar C}_5 \cdot {{ |\,\varepsilon|}\over {\lambda}} \   \mfor  {(\,\lambda, \, \,\xi) \, \in \,  \partial B_{\,{\bar \lambda}\,, \ {\bar \xi}}  \, ({\bar \rho}) } \,.
$$
Here, as usual, ${\bf z}\ = \ V_{\lambda\,, \ \xi}\,.\,$
Applying the degree theory \cite{Fonseca-Gangbo} (\,also described in \S\,A.3.3 and \S\,A.3.4 [\,in particular, (A.3.8)\,]),   we
obtain a critical point $\,\tilde{\bf z} \in B_{{\bar \lambda}\,, \ {\bar \xi}} \ (\bar \rho)$\, for \,$\Phi_\varepsilon ({\bf z})$. Applying Proposition 4.15\,,\, we complete the proof of the theorem. {{\hfill {$\rlap{$\sqcap$}\sqcup$}}}

\vspace*{0.3in}


{\bf \large \S\,5. \ \ Separation and blow\,-\,up.}\\[0.1in]
In this section, we provide criteria for a sequence of critical points $\,\{ {\bf z}_i \,+\, w_\varepsilon ({\bf z}_i ) \} \,$ for $I_\varepsilon$ to be distinct, and to form a blow\,-\,up sequence when $\,\lambda_i \to 0$\,.\, Take $\,{\bf z}_i\ = \ V_{\lambda_i\,,\  \xi_i} \,\in\, {\bf Z}$\,.\,
%
%


\vspace*{0.15in}

{\bf \S\,5\,a.} \ \ {\it Separation lemma.}\ \
Consider $\, V_{\lambda_o,\ \vec{\,0}}\ , \  \, V_{{\lambda_1},\,\, \xi_1} \,\in\, {\bf Z}$\,,\, where $\lambda_1 < \lambda_o\,.\,$ Via the parametrization (2.19), the pull-back metric (3.4) on $\R^+ \times \R^n$,\, and the identification with the (rescaled) hyperbolic metric in \S\,3\,b\,,\, we find that

\newpage

$$
\Vert \,V_{\lambda_o\,,\ \vec{\,\,0}} - V_{{\lambda_1}\,,\ \xi_1}  \Vert_{\,\btd} \ \ge \ \Vert \,V_{\lambda_o,\,\, 0} - V_{{\lambda_1},\  \vec{\,0}}  \Vert_{\,\btd} \\ = \ \ \int_{\lambda_1}^{\lambda_o} {{ {\bar c}_o}\over {\lambda}} \, d\lambda \\ = \ \ {\bar c}_o \, \ln \left( {{\lambda_o}\over {\lambda_1}}  \right)\,.\, \leqno (5.1)
$$
Recall that vertical lines  $\,(\,t\,, \  \bar \xi)\,$ [\,with $\,\bar \xi$\, fixed and $\,t \in \R^+$\,]\, are geodesics in the hyperbolic space. Thus the (curved) geodesic joining $\,(\,\lambda_o\,,\ \vec{\,0})\,$ and $\,(\,\lambda_1\,,\ \xi_1)\,$ in the hyperbolic space has longer length then that joining $\,(\,\lambda_o\,,\ \vec{\,0})\,$ and $\,(\,\lambda_1\,,\ \vec{\,0})\,$ whenever $\,\xi_1 \,\not=\, \ \vec{\,0}\,.$\\[0.2in]
{\bf   Lemma 5.2.} \ \ {\it Let }
$$
{\bf z}_o\ = \ V_{\lambda_o\,,\ \vec{\,0}} \ \ \ \ {\it and} \ \ \ \ {\bf z}_1\ = \ V_{\lambda_1\,,\,\,\xi_1}\,. \ \ \ \ {\it{Assume \ \ that }} \ \ \ {{\lambda_o}\over {\lambda_1}} \ \ge \ 2\,. \leqno (5.3)
$$
{\it For $\,{\bar \varepsilon}_5\,$ small enough, and for all  $\,|\, \varepsilon|\le {\bar \varepsilon}_5\,,\,$ we have }
$$
{\bf z}_o \,+\, w_\varepsilon ({\bf z}_o) \ \not=\  {\bf z}_1 \,+\,  w_\varepsilon ({\bf z}_1)\,.
$$
{\it Here $\,w_\varepsilon ({\bf z}_o)\,$ and $\,w_\varepsilon ({\bf z}_1) $ are the solutions of the auxiliary equation\,} (4.4) {\it as described in\,} Lemma 4.1.  \\[0.15in]
{\bf Proof.} \ \ Suppose the result is false. Then
\begin{eqnarray*}
& \ & \ \ \ \ \ \ \ \ \ \ \ \ \ \ \ \ \ \ \ \ \ {\bf z}_o \,+\, w_\varepsilon ({\bf z}_o)\ = \ {\bf z}_1 \,+\, w_\varepsilon ({\bf z}_1)\\
& \Longrightarrow &   {\bar c}_o \, \ln \,2 \ <  \ \Vert\, V_{\lambda, \,\, \vec{\,0}} - V_{\lambda_1, \,\, \xi_1} \Vert_{\,\btd}\  = \, \Vert w_\varepsilon ({\bf z}_1)  - w_\varepsilon ({\bf z}_o)  \Vert_{\,\btd} \le \Vert \,w_\varepsilon ({\bf z}_1)    \Vert_{\,\btd} + \Vert\, w_\varepsilon ({\bf z}_1)  \Vert_{\,\btd}\\
& \Longrightarrow &   {\bar c}_o \, \ln \,2 \  \le \  2 \, {\bar C}_4\,|\,\varepsilon|\,.
\end{eqnarray*}
Here we use (5.1) and  the uniform bound in Lemma 4.16.
Once we choose ${\bar \varepsilon}_5$ to be small enough,  we have a contradiction. \qedwh
 {\bf \S\,5\,b.} \ \  {\it Blow\,-\,up lemma.}\\[0.15in]
 {\bf Lemma 5.4.} \ \ {\it Let $\,\{ v_i\}_{i =1}^\infty\,$ be a sequence of positive $C^{2\,,\,\, \beta}$ solutions of equation\,}  (2.5)\,    {\it with $K \in C^{1, \ \tilde\beta}$\,.\, Suppose that}
 $$
 \Vert \,v_i- V_{\lambda_i,\, \,\xi_i} \Vert_{\,\btd} \le   \, C\,\tilde\varepsilon \ \ \ \ \ \ \ {\it{for}} \ \ \ i \gg 1\,, \leqno (5.5)
 $$
{\it where $\,\xi_i \to \xi_o\,$ and $\,\lambda_i \to 0^+$\,.\, There is a positive number \,${\bar\varepsilon}_6$\,} ({\it depending on $\,n\,$ only\,})  {\it such that   for all $\,\delta\,$ small enough}\, ({\it depending on $\,n\,$ and the sequence $\,\{ v_i\}_{i =1}^\infty\,$}),\, {\it the sequence $\{ v_i\}_{i =1}^\infty\,$ cannot stay bounded in $\,B_{\xi_o} \,(\delta)$\,,\, once}\, (5.5) {\it is satisfied for } $ \tilde\varepsilon \ \le\ {\bar\varepsilon}_6$\,.

\vspace*{0.2in}

{\bf Proof.} \ \ We compute
\begin{eqnarray*}
(5.6)&\,&\!\!\!\!\!\!\!\!\!\!\!\!\!\!\!\int_{B_{\xi_o}\,(\delta)} |\btd V_{\lambda_i\,,\,\, \xi_i} |^{\,2} \  \ge \  \int_{B_{\xi_i}\,\left({\delta\over 2}\right)}|\btd V_{\lambda_i\,,\,\, \xi_i} |^{\,2}  \!\! \mfor i \gg 1 \  \,  {\mbox{so \   that \   }} {B_{\xi_i} \left(\delta/2\right)} \subset B_{\xi_o}\,(\delta) \ \ \ \ \ \\[0.075in]
& \ge&  (n - 2)^2\int_{B_o\,\left({\delta\over 2}\right)} {{\lambda_i^{n - 2}\, r^2}\over {(\,\lambda_i^2 + r^2)^n }}\\[0.075in]
 & \ &  \!\!\!\!\!\!\!\!\!\!\!\! [\,{\mbox{via \ \ a \ \ translation \ \ }} (y - \xi_i) \to y\,,\, {\mbox{  find \ the \  gradient \  in \   polar \   coordinates}}\,] \\[0.075in]
& =  & \Vert S^{n - 1} \Vert \cdot (n - 2)^2\int_0^{\delta \over 2}  {{\lambda_i^{n - 2}\, r^{n + 1}}\over {(\,\lambda_i^2 + r^2)^n }}\, \,dr\\
&\ = \ &  (n - 2)^2\int_0^{\arctan  \left( {{\delta}\over{ 2\,\lambda_i}}\!\right) } \, {{\,[\,\tan \, \theta]^{n + 1} \,\,\sec^2 \theta \, d\theta\,}\over {[\,\sec  \theta\,]^{2n} }} \ \ \ \ \ \ \ \ \ \ \ \ \ \ \ \ \ \ \ \ [\,r\ = \ \lambda_i \tan \theta]\\
&\ = \ &  (n - 2)^2\int_0^{\arctan \left( {{\delta}\over{2\, \lambda_i}}\!\right) }\  [\,\sin \, \theta]^{n + 1}\,[\,\cos \, \theta]^{n -3}  \, d\theta \ge C_o\,, \ \     {\mbox{as \ long \ as }} \ {{\delta}\over {2\,\lambda_i }} \,\ge\, 1\,.
\end{eqnarray*}
Here $C_o$ depends on $n$ only.
Suppose that $\{ v_i \}_{i\ = \ 1}^\infty\,$ is uniformly bounded in $B_{\xi_o} \, (\delta)$. Via equation (2.5) and elliptic regularity theory on gradient estimates \cite{Gilbarg-Trudinger}, one infers that
$$
\Vert \,\btd \,v_i \Vert  \le C_1 \ \ \ {\mbox{in}} \ \ B_{\xi_o}\, (\delta) \ \ \ \ {\mbox{for \ \ all}} \ \  \ i \gg 1\,. \leqno (5.7)
$$
Observe that the constant $C_1$ may depend on $\{ v_i \}$\,.\,
(5.6) leads to

\vspace*{-0.25in}

\begin{eqnarray*}
(5.8) \ \ \ \ \    \ \ \ \ \  C_o & \le & \int_{B_{\xi_o} (\delta)} \Vert\,\btd \,V_{\lambda_i,\, \xi_i} \Vert^{\,2} \ = \   \int_{B_{\xi_o} (\delta)} \Vert\,\btd  \,v_i \ + \ (\, \btd_o\,V_{\lambda_i,\, \xi_i} \ - \ \btd  \,v_i)\,\Vert^2 \,\ \ \ \ \    \ \ \ \ \ \ \ \ \ \    \ \ \ \ \  \\
& \le &  2\int_{B_{\xi_o} (\delta)}\Vert\btd  \,v_i \,\Vert^{\,2} \,\, + \ 2\int_{B_{\xi_o} (\delta)}\Vert  \btd \,V_{\lambda_i,\, \xi_i}\ - \ \btd \,v_i\,\Vert^{\,2} \\
  &\le &   {C}_1^2 \cdot {{\Vert \,S^n\Vert}\over n} \cdot  \delta^n \ + \  2 \, {\tilde\varepsilon}^{\,2}\,.
\end{eqnarray*}
Combining (5.7) and (5.8), we obtain
$$
C_o \ \le\  { C}_1^2 \, {{\Vert \,S^n\Vert}\over n} \cdot  \delta^n \ + \  2 \, {\tilde\varepsilon}^{\,2}
$$
Thus once we choose ${\bar \varepsilon}_5$ and $\delta$ so that
$$
{{C_o }\over 2} \ > \ 2 \, {\bar\varepsilon}_{\,5} \ \ \ \ {\mbox{and}} \ \ \ \ {{C_o }\over 2}\  > \ {C}_1^2  \,{{\Vert \,S^n\Vert}\over n} \cdot  \delta^n\ ,
$$
we have a contradiction. Thus the conclusion of the lemma stands.\qed

\vspace*{0.05in}


{\large {\bf \S\,6. \ \ First order derivatives and the Kazdan\,-Warner condition.}}\\[0.2in]
From expressions (1.7) and (2.19), we obtain
\begin{eqnarray*}
(6.1) \ \ \ \ \    {{\partial G_{|_{\bf Z}} }\over {\partial \lambda }}\,
  (\,\lambda\,,\  \xi)
  &\ = \ &{{ n\,{\bar c}_{-1} }\over {\lambda}}
  \left\{ \, \int_{\R^n} H (y)
  \left[ {{\lambda^n\,  }\over {\left(\,\lambda^2 + |\,y - \xi|^{\,2}
  \right)^n}}\right]\,\,\right.\\
  & \ & \ \ \ \   \ \ \ \  \ \ \ \  \ \ \ \  \ \ \ \ \ \   - \,2 \, \left.
  \int_{\R^n} H (y) \left[
  {{ \lambda^{n+2}\, }\over {\left(\,\lambda^2 + |\,y - \xi|^{\,2}
  \right)^{n+1}}} \right] \right\}\,. \ \ \ \ \ \ \ \ \ \ \ \
\end{eqnarray*}

\vspace*{-0.2in}

This leads to

\vspace*{-0.2in}

$${{\partial G_{|_{\bf Z}} }\over {\partial \lambda }} \,(\,\lambda\,,\  \xi)\ = \  n\,{\bar c}_{-1}\,\lambda^{n-1} \int_{\R^n} H (y) \,  {{1}\over {\left(\,\lambda^2 + |\,y - \xi|^{\,2}\right)^n}}\left[  {{|\, y - \xi|^{\,2} - \lambda^2 }\over { |\, y - \xi|^{\,2} + \lambda^2 }}  \right] \,\,\,. \leqno (6.2)
 $$
Likewise,

\vspace*{-0.25in}

$$
{{\partial G_{|_{\bf Z}}}\over {\partial\,\xi_j}}\, (\,\lambda\,,\ \xi)\ = \ -\,2 n\,{\bar c}_{-1} \int_{\R^n} H (y) \,  {{\lambda^n\, (\xi_j - y_j)}\over {\left(\,\lambda^2 + |\,y - \xi|^{\,2}\right)^{n+1}}} \,\, \leqno (6.3) $$
for $\,j\ = \ 1, \ 2, \cdot \cdot \cdot, \ n\,.\,$
Note that the integrals in (6.1)
are independent on $\,\lambda\,$ and $\,\xi\,.$ This highlights the presence of $\lambda^{-1}\,$ in (6.1)\,.\,  Cf. (4.17) and (4.27).

\vspace*{0.2in}

{\bf Lemma 6.4.} \ \ {\it In\,} (6.3)\,,\, {\it suppose $\,\xi \, = \, 0\,,\,$ and $\,H \,(y)\, = \,H\, (|\,y|\,)\,$ is radially symmetric and bounded, then \,}

\vspace*{-0.1in}

$$
{{\partial  G_{|_{\bf Z}}}\over {\partial\,\xi_j }} \,(\,\lambda \,,\, 0)\ = \ 0   \ \ \ \ \ \ {\it for} \ \ \ \lambda \,> \,0 \ \ \ {\it and} \ \ \ 1\,\le\, j \,\le\, n\,.
$$

\vspace*{0.1in}


%
%
%
{\bf Proof.} \ \ This is because the integrals
$$
\int_{\R^n} H\, (|\,y|)  \cdot {{  y_j  }\over {\left(\,\lambda^2 + |\,y|^{\,2}\right)^{n+1}}} \ \ \ \ \  \ (j\ = \ 1, \ \,2, \cdot \cdot \cdot, \ n)
$$
vanish, due to antisymmetry. \qedwh
%
%
%
%
%
{\bf \S\,6\,a.} \ \ {\it Geometric interpretation.} \ \
Via the change of variables
$\,
\lambda \,\bar y\ = \  y - \xi\,,\,
$
we have
\begin{eqnarray*}
(6.5) \ \ \  \   {{\partial G_{|_{\bf Z}} }\over {\partial \lambda }} (\,\lambda\,,\, \xi) &=& {1\over {2^n}}\cdot{{ n\, {\bar c}_{-1} }\over {\lambda}} \left[ \,\int_{\R^n} H (\,\lambda \,\bar y \,+\, \xi) \, \left( 1 -      {{2  }\over { 1 + |\,\bar y|^{\,2} }}\,\right) \,{{2^n }\over {\left(1 + |\,\bar y |^{\,2}\right)^n}} \ d\bar y \,\right] \ \ \ \ \ \ \\
  & \ &\\
 \Longrightarrow  \ \  {{\partial G_{|_{\bf Z}} }\over {\partial \lambda }} \,(\,\lambda\,,\, \xi) &=&  {1\over {\lambda }}\cdot{{ n\, {\bar c}_{-1} }\over {2^n}} \int_{S^n}\! H (\,\lambda \,{\cal P}(x)\,+\, \xi)\, x_{n + 1} \ dS_x \ \ \ \  \ [\,{\mbox{via}} \ \ (2.2) \ \& \ \ (2.3)\,]\\
   = \ {1\over {\lambda }}\cdot{{ n\, {\bar c}_{-1} }\over {2^n}} \,\cdot & \ &\!\!\!\!\!\!\!\!\!\!\!\!\!\!\!  \left[ \ \int_{S^n_+} \! H (\,\lambda \,{\cal P}(x)\,+\, \xi)\, x_{n + 1} \ dS_x \ + \ \int_{S^n_-}\! H (\,\lambda \,{\cal P}(x)\,+\, \xi)\, x_{n + 1} \ dS_x \right]
\end{eqnarray*}
Intuitively, one can say  that a critical point is a balance between the upper \,($S^n_+$)\, and lower \,($S^n_-$)\, hemisphere. Likewise,
\begin{eqnarray*}
(6.6) \ \ \ \ \ \ \  \    {{\partial G_{|_{\bf Z}}}\over {\partial\,\xi_j}}\, (\,\lambda\,,\, \xi) &\ = \ &  \,2 n\,{\bar c}_{-1} \int_{\R^n} H (y) \,  {{\lambda^n\, (y_j - \xi_j)}\over {\left(\,\lambda^2 + |\,y - \xi|^{\,2}\right)^{n+1}}} \ \,\\
 &\ = \ &  \,{{ n\,{\bar c}_{-1}}\over {\lambda}} \int_{\R^n} H (\,\lambda \,\bar y + \xi) \, \left( {{2\,{\bar y}_j  }\over { 1+ |\,\bar y |^{\,2} }}  \right)\,  {{1 }\over {\left(1+ |\,\bar y |^{\,2}\right)^n}}\ \,\\
 &\ = \ &  \,{{ n\,{\bar c}_{-1}}\over {\lambda}} \int_{\R^n} H (\,\lambda \,\bar y + \xi) \, \left( {{2\,{\bar y}_j  }\over { 1+ |\,\bar y |^{\,2} }}  \right)\,  {{1 }\over {\left(1+ |\,\bar y |^{\,2}\right)^n}}\ \,\\[0.1in]
 \Longrightarrow \ \ {{\partial G_{|_{\bf Z}} }\over {\partial\,\xi_j }}\, (\,\lambda\,,\, \xi) &\ = \ &  {1\over {\lambda }}\cdot{{ n\, {\bar c}_{-1} }\over {2^n}} \int_{S^n} H (\,\lambda \,{\cal P}(x)+ \xi)\, x_j \ dS_x\ \ \ \  [\,{\mbox{via}} \ \ (2.2) \ \& \ \ (2.3)\,] \ \ \
\end{eqnarray*}
for $\,j\ = \ 1, \ 2, \cdot \cdot \cdot, \ n\,.\,$
It follows that if $\,{\cal H}\,$ is antipodal symmetric, i.e., $\,{\cal H} (-x)\ = \ {\cal H} (x)\,$,\, then $\,(1, \ {\vec{\,0}}\,)\,$  is a critical point of $\,G_{|_{\bf Z}}$\,.\, [\,Recall that \,${\cal H}\,$ and $\,H\,$ are related via  (2.6).\,] \bk
Next, we observe that  $\,x_1, \,\cdot \cdot \cdot, \ x_{n + 1}\,$ are first eigenfunctions of the standard Laplacian on $S^n$.  Using (6.5) and (6.5), we have the following.\\[0.2in]
%
%
%
{\bf Proposition 6.7.} \ \ {\it Let ${\cal H}$ and $H$ be related by } (2.6)\,.\, {\it Suppose that ${\cal H} (x)$ is orthogonal\,} ({\it with respect to the $L^2$-\,inner product})\,{\it to the first eigenspace of $\,\Delta_1\,.\,$ Then $(1, \ {\vec{\,0}}\,)$ is a critical point of $\,G_{|_{\bf Z}}$\,.\, }

\vspace*{0.2in}

{\bf \S\,6\,b.} \ \ {\it Homogeneous harmonic polynomials\,.} \ \  In order to apply Proposition 6.7 we consider in particular order $k$ homogeneous harmonic polynomials $\,P_\Delta^k\,$ on $\,{\R}^{n + 1}\,,$  where $\,k > 1\,$ is an  integer.  It is well known that
$$
\Delta_1\, P_\Delta^k + k (k + n - 1) P_\Delta^k\ = \ 0 \ \ \ \ {\mbox{in}} \ \ S^n\,.
$$
For $k > 1$, each $\,{\cal H}\, = \, P_\Delta^k\,$ is orthogonal to the first eigenspace, i.e., the condition in Proposition 6.7 is fulfilled.\\[0.2in]
{\bf \S\,6\,c.} \ \ {\it Kazdan\,-\,Warner condition.}\ \
A deep relation is revealed when we differentiate $\,{\cal K}\,$ with respect to a {\it conformal Killing vector field\,} {\bf  X} \, (generating a family of conformal transformations). In this way we obtain the renowned Kazdan\,-\,Warner formula \cite{BourguignonEzin}
$$ \ \
\int_{S^n} {\bf  X} \,({\cal K}) \, u^{{2n}\over {n - 2}} \  dS\ = \ 0  \ \ \ \  \ \ {\mbox{[\,refer \ to \ equation \ (1.1), \ here}} \ \ {\bf  X} \, ({\cal K})\ = \ \langle{\bf  X} \, \  \btd \,{\cal K}\rangle_{g_1}\,]\,.
$$
This is a necessary (but not sufficient) condition on $\,{\cal K}\,$ for equation (1.1)   has a positive solution (\,see \cite{Han-Li}).\\[0.2in]
%
%
%
%
%
{\bf Definition 6.8.} \ \ {\it A function\,} ${\cal K} \in C^1 (S^n)\,$ {\it is said to satisfy the K\,-W condition if there exists a\,} ({\it single}) {\it positive  function\,} $\psi \in C^o \,(S^n)$ {\it such that}
$$
\int_{S^n}  {\bf  X} \,({\cal K} ) \, \psi^{{2n}\over {n - 2}} \, dS\ = \ 0  \ \ \ \ \ \ {\it{for \ \ all \ \ conformal \ \ Killing \ \ vector \ \ fields \ \ }} \   {\bf  X} \,\,. \leqno (6.9)
$$

\vspace*{0.15in}

{\bf Theorem 6.10.} \ \ {\it For a given $\,{\cal H} \in C^1\,(S^n)\,$,\, suppose $\,G_{|_{\bf Z}}\,$ has a critical point at $\,(\,\lambda_c,\,\, \xi_c)\,$,\, then }
$$
{\cal K}\ = \ 4n(n -1) + \varepsilon \,{\cal H}
$$
{\it fulfills the $K\!-\!W$ conditions\,}\, (6.9)\,  {\it for any\,}  $\varepsilon \in \R$, {\it with}
$$
\psi \,(x)\ = \ \left({{\lambda_c}\over {\lambda_c^2 + |\, y - \xi_c|^{\,2}}}  \right)^{{n - 2}\over 2} \!\!\cdot \left( 1 + |\,y|^{\,2}\right)^{{n - 2}\over 2} \ \ \ \ \ \ \ \ {\it for} \ \  \ y\ = \ \dot{\cal P} \,(x)\,.
$$

\vspace*{0.15in}

\hspace*{0.5in}As the proof of Theorem 6.10 is well\,-known\,,\, we  direct the interested readers to  \S\,A.11 in the e\,-Appendix for reference. Nevertheless, Theorem 6.10 highlights the balance achieved by a critical point of $G_{|_{\bf Z}}$\,.\,

\vspace*{0.35in}

{\large {\bf \S\, 7. \ \  Second derivatives, symmetry and stability conditions.}}\\[0.15in]
We continue from (6.2)\,--\,(6.3),
$$
{{\partial^2 G_{|_{\bf Z}}}\over {\partial\,\xi_j\, \partial\,\xi_\ell\,}} \,(\,\lambda\,,\  \xi)\ = \  2^2 \,n\,(n + 1)\,{\bar c}_{-1} \int_{\R^n} H(y) \cdot  {{\lambda^n\, (\xi_j - y_j)(\xi_\ell - y_\ell)}\over {\left(\,\lambda^2 + |\,y - \xi|^{\,2}\right)^{n+2}}}\,\,, \leqno (7.1)
$$
\begin{eqnarray*}
(7.2) & \ &  {{\partial^2 G_{|_{\bf Z}}}\over {\partial \lambda\, \partial\,\xi_\ell\,}} \,(\,\lambda\,,\ \xi) \ = \  -\ 2 n^2\ {\bar c}_{-1} \int_{\R^n} H (y)\cdot  {{ \lambda^{n-1}\, (\xi_\ell - y_\ell)}\over {\left(\,\lambda^2 + |\,y - \xi|^{\,2}\right)^{n+1}}}\ + \\[0.075in]
& \  & \ \ \ \ \ \ \ \ \  + \ 2 n\,(n + 1)\,{\bar c}_{-1} \int_{\R^n} H(y)   \cdot {{2\,\lambda^{n+1}\, (\xi_\ell - y_\ell)}\over {\left(\,\lambda^2 + |\,y - \xi|^{\,2}\right)^{n+2}}}
  \mfor  1 \,\le\, j \,\not=\, \ell \,\le\, n\,.\ \ \ \ \ \
\end{eqnarray*}
Moreover,
\begin{eqnarray*}
(7.3) \ \ \ \ \ \ \ \ \ \ {{\partial^2 {G_{|_{\bf Z}} }}\over {\partial\,\xi_\ell^2}}\, (\,\lambda\,,\  \xi) &\ = \ & -\,2 n\,{\bar c}_{-1} \int_{\R^n} H (\,y)  \cdot {{\lambda^n\, }\over {\left(\,\lambda^2 + |\,y - \xi|^{\,2}\right)^{n+1}}}\,\\
& \  & \ \ \ \ \ +\ 4 n\,(n + 1)\,{\bar c}_{-1} \int_{\R^n} H (\,y)  \cdot {{ \lambda^n\, (\xi_\ell - y_\ell)^2}\over {\left(\,\lambda^2 + |\,y - \xi|^{\,2}\right)^{n+2}}}\,\\
\Longrightarrow \ \ {{\partial^2 {G_{|_{\bf Z}} }}\over {\partial\,\xi_\ell^2}} (\,\lambda\,,\, 0) &\ = \ & -\,2\, n\,{\bar c}_{-1} \int_{\R^n} H (\,y)  \cdot  {{\lambda^n\, }\over {\left(\,\lambda^2 \,+\, |\,y|^{\,2}\right)^{n+1}}}\, \,\\
& \  & \ \ \ \ \ +\ 4 \,n\,(n + 1)\,{\bar c}_{-1} \int_{\R^n} H (\,y) \cdot  {{ \lambda^n\, (y_\ell)^2}\over {\left(\,\lambda^2 + |\,y |^{\,2}\right)^{n+2}}}\,\,\,.\ \ \ \ \ \ \ \ \ \ \ \ \ \ \ \
\end{eqnarray*}

\vspace*{-0.25in}

It follows from (6.1) that
\begin{eqnarray*}
(7.4) \ \   {{\partial^2 {G_{|_{\bf Z}} }}\over {\partial \lambda^2 }} \, (\,\lambda\,,\  \xi) &= &  n\,(n - 1)\,{\bar c}_{-1} \int_{\R^n} H(\,y)   \cdot  {{\lambda^{n-2}\,  }\over {\left(\,\lambda^2 + |\,y - \xi|^{\,2}\right)^n}}\, \,\,\\
& \  & \ \ \ \   - \ 2n^2\, \,{\bar c}_{-1} \int_{\R^n} H(\,y)  \cdot  {{ \lambda^{n}\, }\over {\left(\,\lambda^2 + |\,y - \xi|^{\,2}\right)^{n+1}}}\, \,\,\\
& \  & \ \ \ \ \  \ \ \ \   - \ 2n\,(n + 1)\, \,{\bar c}_{-1} \int_{\R^n}H(\,y) \cdot   {{ \lambda^{n}\, }\over {\left(\,\lambda^2 + |\,y - \xi|^{\,2}\right)^{n+1}}}\, \,\,\\
& \  & \ \ \ \ \ \ \ \ \ \ \ \ \ \  + \ 4\, n\,(n + 1)\, \,{\bar c}_{-1} \int_{\R^n} H(\,y) \cdot   {{ \lambda^{n+2}\, }\over {\left(\,\lambda^2 + |\,y - \xi|^{\,2}\right)^{n+2}}}\, \,\,. \ \ \ \ \ \ \ \ \ \ \ \ \
\end{eqnarray*}

\vspace*{0.15in}

{\bf \S\,7\,a.} \ \  {\it Second derivative at a critical point.}\ \
Let $(\,\lambda_c\,,\,\ \xi_c)$ be a critical point for $G_{|_{\bf Z}}$\,.\, From (7.1), (7.2) and (6.3) we obtain
$$ {{\partial^2 G_{|_{\bf Z}}}\over {\partial\,\xi_j\, \partial\,\xi_\ell}} \,(\,\lambda_c\,,\, \xi_c) \ =\    4\, n\,(n + 1)\,{\bar c}_{-1} \int_{\R^n} H(y) \cdot  {{\lambda^n_c\, ({\xi_c}_j - y_j)({\xi_c}_\ell - y_\ell)}\over {\left(\,\lambda^2 + |\,y - \xi_c|^{\,2}\right)^{n+2}}}\,\,,\leqno (7.5)
$$
 $$
 {{\partial^2 G_{|_{\bf Z}}}\over {\partial \lambda\, \partial\,\xi_\ell}} \,(\,\lambda_c\,,\, \xi_c) \ =\   2\, n\,(n + 1)\,{\bar c}_{-1} \int_{\R^n} H(y) \cdot  {{ \lambda^{n+1}_c\, ({\xi_c}_\ell - y_\ell)}\over {\left(\,\lambda^2_c + |\,y - \xi_c|^{\,2}\right)^{n+2}}}\,\,.\leqno (7.6)
 $$
Here $1 \,\le\, j \,\not=\, \ell \,\le\, n\,.\,$
Before we continue, we register the following observation.

\vspace*{0.2in}

{\bf Lemma 7.7.} \ {\it In\,} (7.5) {\it and}\, (7.6)\,,\, {\it suppose $\,\xi_c\ = \ 0\,$ and $\,H (y)\ = \ H\, (|\,y|)\,$ is radially symmetric, then \,}
$$
{{\partial^2 G_{|_{\bf Z}}}\over {\partial\,\xi_j\, \partial \, \xi_\ell\,}} \,(\,\lambda_c\,,\, 0)\ = \ 0\ = \ {{\partial^2 G_{|_{\bf Z}}}\over {\partial \lambda\, \partial\, \xi_\ell\,}} \,(\,\lambda_c\,,\, 0) \ \ \ \ \ \ \ \ \ {\it{for}} \ \ \ 1 \ \le\  j \  \not=\  \ell\  \le\  n\,.
$$
{\it That is, the Hessian matrix is diagonal.}\\[0.2in]
{\bf Proof.} \ \ This is because the integrals
$$
\int_{\R^n} H\, (|\,y|)  \cdot  {{  y_j\,y_\ell }\over {\left(\,\lambda^2 + |\,y|^{\,2}\right)^{n+2}}}\, \ \ \ \ (j \,\not=\,\ell) \ \ \ \   {\mbox{and}} \ \ \ \int_{\R^n} H\, (|\,y|)  \cdot  {{  y_\ell }\over {\left(\,\lambda^2 + |\,y|^{\,2}\right)^{n+2}}}\,
$$
vanish by antisymmetry argument. \qed
From (7.4), together with (6.1), we obtain
\begin{eqnarray*}
(7.8) \ \ \   {{\partial^2 {G_{|_{\bf Z}} }}\over {\partial \lambda^2 }} (\,\lambda_c\,,\, \xi_c) &=&  2\,n\,(n - 1)\,{\bar c}_{-1} \int_{\R^n} H(\,y)   \cdot {{\lambda^{n}_c\,  }\over {\left(\,\lambda^2_c + |\,y - \xi_c|^{\,2}\right)^{n+1}}}\, \,\,\\
& \  & \ \ \ \ \ -\  2\,n^2\, \,{\bar c}_{-1} \int_{\R^n} H(\,y)  \cdot {{ \lambda^{n}_c\, }\over {\left(\,\lambda^2_c + |\,y - \xi_c|^{\,2}\right)^{n+1}}}\, \,\,\\
& \  & \ \ \ \ \  \ \ \ \ \ - \ 2\,n\,(n + 1)\, \,{\bar c}_{-1} \int_{\R^n}H(\,y)   \cdot{{ \lambda^{n}_c\, }\over {\left(\,\lambda^2_c + |\,y - \xi_c|^{\,2}\right)^{n+1}}}\, \,\,\\
& \  & \ \ \ \ \ \ \ \ \ \ \ \ \  + \ 4\, n\,(n + 1)\, \,{\bar c}_{-1} \int_{\R^n} H(\,y) \cdot   {{ \lambda^{n+2}_c\, }\over {\left(\,\lambda^2_c + |\,y - \xi_c|^{\,2}\right)^{n+2}}}\, \,\,\\[0.1in]
&\ = \ & - 2\,n \, (n + 2) \,{\bar c}_{-1} \int_{\R^n} H (\,y)\cdot   {{ \lambda^{n}_c\, }\over {\left(\,\lambda^2_c + |\,y  - \xi_c |^{\,2}\right)^{n+1}}}\, \,\,\\
 & \ & \ \ \ \ \ +  \ 4\,n\,(n + 1)  \,{\bar c}_{-1} \int_{\R^n} H (\,y) \cdot  {{ \lambda^{n+2}_c\, }\over {\left(\,\lambda^2_c + |\,y  - \xi_c |^{\,2}\right)^{n+2}}}\, \,\,\,.  \ \ \  \ \ \  \ \ \  \ \ \ \ \ \
\end{eqnarray*}
[\,Check\,: $2n\,(n - 1) - 2n^2 - 2n\,(n + 1)\ = \ -2n\,(n + 2)\,.$\,]\, As for the second derivatives in $\,\xi_\ell\,$,\,  from (7.3), we have
\begin{eqnarray*}
(7.9) \ \ \ \ \ \ \sum_{\ell\, = \, 1}^n \,{{\partial^2 {G_{|_{\bf Z}} }}\over {\partial\,\xi_\ell^2 }} \,(\,\lambda_c\,,\  \xi_c) &\ = \ & -\,2 n^2\,{\bar c}_{-1} \int_{\R^n} H (\,y)  \cdot {{\lambda^n_c\, }\over {\left(\,\lambda^2_c + |\,y - \xi_c|^{\,2}\right)^{n+1}}}\, \,\,\\
& \  & \ \ \   +\  4\, n\,(n + 1)\,{\bar c}_{-1} \int_{\R^n} H (\,y)  \cdot {{ \lambda^n_c\, \sum_{\ell\, = \, 1}^n    ({\xi_c}_\ell - y_\ell)^2}\over {\left(\,\lambda^2_c + |\,y - \xi_c|^{\,2}\right)^{n+2}}}\, \,\,\\[0.1in]
&\ = \ & -\ 2 \,n^2\,{\bar c}_{-1} \int_{\R^n} H ( y )\cdot    {{\lambda^n_c\, }\over {\left(\,\lambda^2_c + |\,y - \xi_c|^{\,2}\right)^{n+1}}}\, \,\\
& \  & \ \ \ \ \ +\ 4 \,n (n + 1) \,{\bar c}_{-1} \int_{\R^n} H (y) \cdot  {{ \lambda^n_c\, |\,y -\xi_c|^{\,2}}\over {\left(\,\lambda^2_c + |\,y - \xi_c|^{\,2}\right)^{n+2}}}\,\,\\
&\ = \ & -\ 2 \,n^2\,{\bar c}_{-1} \int_{\R^n} H (y)\cdot    {{\lambda^n_c\, }\over {\left(\,\lambda^2_c + |\,y -\xi_c|^{\,2}\right)^{n+1}}}\, \,\\
& \  & \ \ \ \ \ +\ 4 \,n (n + 1) \,{\bar c}_{-1} \int_{\R^n} H (y)\cdot    {{ \lambda^n_c\, (\,\lambda_c^2 + |\,y -\xi_c|^{\,2})}\over {\left(\,\lambda^2_c + |\,y -\xi_c|^{\,2}\right)^{n+2}}}\,\,\\
& \  & \ \ \ \ \     -\  4 \,n (n + 1) \,{\bar c}_{-1} \int_{\R^n} H (y) \cdot   {{ \lambda^n_c \times \lambda_c^2}\over {\left(\,\lambda^2_c + |\,y -\xi_c|^{\,2}\right)^{n+2}}}\,\,\\[0.1in]
&\ = \ & -\ 2 \,n^2\,{\bar c}_{-1} \int_{\R^n} H (y)\cdot    {{\lambda^n_c\, }\over {\left(\,\lambda^2_c + |\,y -\xi_c|^{\,2}\right)^{n+1}}}\, \,\\
& \  & \ \ \ \ \ +\,4 \ n (n + 1) \,{\bar c}_{-1} \int_{\R^n} H (y)\cdot   {{ \lambda^n_c\,  }\over {\left(\,\lambda^2_c + |\,y -\xi_c|^{\,2}\right)^{n+1}}}\,\,\\
& \  & \ \ \ \ \   -\  4 \, n\,(n + 1) \,{\bar c}_{-1} \int_{\R^n} H (y) \cdot  {{ \lambda^n_c \times \lambda_c^2}\over {\left(\,\lambda^2_c + |\,y -\xi_c|^{\,2}\right)^{n+2}}}\,\,\\[0.1in]
&\ = \ & 2\,n\,(n + 2)\,{\bar c}_{-1} \int_{\R^n} H (y)\cdot    {{\lambda^n_c\, }\over {\left(\,\lambda^2_c + |\,y -\xi_c|^{\,2}\right)^{n+1}}}\, \,\\
& \  & \ \ \    -\ 4 \, n\,(n + 1) \,{\bar c}_{-1} \int_{\R^n} H (y) \cdot   {{ \lambda^{n+2}_c}\over {\left(\,\lambda^2_c + |\,y -\xi_c|^{\,2}\right)^{n+2}}}\,\,\,. \ \ \ \ \ \ \ \ \ \ \
\end{eqnarray*}
[\,Check\,: $4n\,(n + 1) - 2n^2\ = \ 2n^2 + 4n\,.\,$]\, Combining (7.8) and (7.9), we observe the following result.

\vspace*{0.2in}

{\bf Proposition 7.10.} \ \ {\it At a critical point $\,(\,\lambda_c\,, \,\, \xi_c)\,$ for $\,G_{|_{\bf Z}}\,$,\, the Hessian matrix is trace-free. That is,  }
$$
{{\partial^2 {G_{|_{\bf Z}} }}\over {\partial \lambda^2 }} \, (\,\lambda_c\,,\, \,\xi_c)\ = \ -\,\left[\,\sum_{\ell\ = \ 1}^n \,{{\partial^2 {G_{|_{\bf Z}} }}\over {\partial \,\xi_\ell^2 }} \,(\,\lambda_c\,,\,\, \xi_c)\right]\,. \leqno (7.11)
$$

\vspace*{0.15in}

\hspace*{0.5in}It is well-\,known that at a critical point, the Hessian matrix  being positive  definite (resp.   negative definite)  is a sufficient condition for local minimum (resp. local maximum). In some sense, Proposition 7.10 guides us to look for saddle points.

\vspace*{0.2in}

{\bf \S\,7\,b.} \ \ {\it Geometric interpretation of the second derivatives.}\ \
If $\,(1, \ {\vec{\,0}}\,)\,$ is a critical point, then
\begin{eqnarray*}
(7.12) \ \ \ \ \ \ \ \  {{\partial^2 G_{|_{\bf Z}}}\over {\partial\,\xi_j\,
\partial\,\xi_\ell}} \,(1, \ {\vec{\,0}}\,)
&\ = \ &  4\, n\,(n + 1)\,{\bar c}_{-1} \int_{\R^n} H(y)\cdot
 {{   y_\ell\, y_j}\over {\left(1 + |\,y |^{\,2}\right)^{n+2}}}\,\\
 &\ = \ &
  {{n\,(n + 1)}\over {2^n}}\,{\bar c}_{-1}\,\int_{S^n} x_\ell \, x_j  \,H(x)
   \, dS_x \ \ \ \ \   (1 \,\le \,i \,\not= \,j \,\le\, n)\,,\ \ \ \ \ \ \ \ \ \ \ \ \ \ \ \ \\
& \ & \\
(7.13) \ \ \ \ \ \ \ \  \ {{\partial^2 G_{|_{\bf Z}}}\over {\partial \lambda\,
\partial\,\xi_\ell}} \,(1, \ {\vec{\,0}}\,)
&\ = \ &    2\, n\,(n + 1)\,{\bar c}_{-1} \int_{\R^n}
H(y)  \cdot  {{  ( - y_\ell)}\over {\left(1 + |\,y |^{\,2}\right)^{n+2}}}\,\\
&\ = \ & -\,{{n\,(n + 1)}\over {2^{n+1}}}\,{\bar c}_{-1}\,\int_{S^n} x_\ell \, (1 - x_{n + 1}) \,H(x) \, dS_x\,,
\end{eqnarray*}
$$ {{\partial^2 {G_{|_{\bf Z}} }}\over {\partial \lambda^2 }} \,(1, \ {\vec{\,0}}\,)
  \,= \, {{n\,(n + 1)}\over {2^n}}\cdot {\bar c}_{-1} \left[\,   \int_{S^n}
  x_{n + 1}^2 \,H(x) \, dS_x \ - \ {1\over {n + 1}} \int_{S^n}
   \,H(x) \, dS_x\right]\,, \leqno (7.14)
$$
$${{\partial^2 {G_{|_{\bf Z}} }}\over {\partial\,\xi_\ell^2}}\,
(1, \ {\vec{\,0}}\,)
 \ = \   -\,{{n}\over {2^n}}\cdot {\bar c}_{-1}\,\int_{S^n} H(x) \, dS_x \ + \  {{n\,(n + 1)}\over {2^n}}\cdot {\bar c}_{-1}\,\int_{S^n} (x_\ell )^2
\,H(x) \, dS_x\,.\leqno (7.15)
$$
We show in \S\,A.11 the calculations in obtaining (7.14) and (7.15).
Observe that (7.11) is fulfilled once we combine (7.14) and (7.15)  with
$$
x_1^2 + \cdot \cdot \cdot + x_n^2\ = \ 1 - x_{n + 1}^2 \ \ \ \ {\mbox{on}} \ \ S^n.
$$

\vspace*{0.2in}

{\bf \S\,7\,c. \ \ }{\it Symmetry\,.} \ \
Consider the symmetry described by
$$
  H \,(y_1, \cdot \cdot \cdot, \ -\,y_\ell\,, \ \cdot \cdot \cdot, \ y_n )\,\ = \ \,  H \,(y_1, \cdot \cdot \cdot, \  \,y_\ell\,, \ \cdot \cdot \cdot, \ y_n ) \mfor \ell\ = \ 1, \, 2, \cdot \cdot \cdot, \  n\,,\leqno (7.16)
$$
$$ H\,(r, \ \vartheta) \ = \ H \left({1\over {r }}\, , \ \vartheta\right) \ \ \ \mfor \ r \in (0, \ 1) \  \ \ \ {\mbox{and}} \  \ \ \vartheta\ \,:= \,\ {y\over {|\,y|}} \,\in\, S^{n - 1}. \leqno (7.17)
$$
Here $\,(r\,,\ \vartheta)\,$ is the polar coordinates of $\,\R^n\,.\,$
In this case $H$ is determined by its values in the `truncated quarter'
$$\,\{\, y\ = \ (y_1,\, \cdot \cdot \cdot, \ y_n) \in \R^n \ \,| \,\  \,y_1 \ge 0\,, \cdot \cdot \cdot, \ y_n\ge 0  \ \ \& \ \ |y| \ \,\le 1\, \}\,.$$

\vspace*{0.15in}

{\bf Theorem 7.18.} \ \ {\it Assume that\,} $H \in C^1 (\R^n)$ {\it satisfies the symmetric conditions\,} (7.16)\,,\,  (7.17),\,   {\it and the stability conditions}
\begin{eqnarray*}
(7.19) \ \ \ \ \int_{\R^n}  H (y)\cdot  \left( {1\over{ 1 + |\, y|^{\,2}}} \right)^{n+1} \,\, &\not=& 2\,(n + 1)\int_{\R^n}  H (y)\cdot { y}_\ell^2 \left( {1\over{ 1 + |\,y|^{\,2}}} \right)^{n+2} \ \ \ \ \ \ \ \ \ \  \\
&  \ & \\
& \ &   \ \ \ \ \ \ \ \ \ \ \ \ \ \ \ \ \ \ \ {\it{for}} \ \ \ell\ = \ 1, \, 2, \cdot \cdot \cdot, \ n\,, \ \ \  {\it{and}} \ \ \ \ \ \  \ \ \  \ \ \   \\
& \ & \\
(7.20) \ \ \ \   \int_{\R^n}  H (  y)\cdot \left( {1\over{ 1 + |\, y|^{\,2}}} \right)^{n+1} \,\,
& \not= &{{ 2\,(n + 1)}\over {n + 2}} \,  \int_{\R^n}  H ( y) \cdot  \left( {1\over{ 1+ |\,y|^{\,2}}} \right)^{n+2}\,\,\,.
\end{eqnarray*}
{\it Then equation\,} (1.5) {\it with relation\,} (1.6) {\it has a positive\, $C^2$ solution when   $|\,\varepsilon|$ is small. } %

\vspace*{0.15in}

{\bf Proof.} \ \
Refer to equation (6.2) with  $\lambda\ = \ 1$ and $\xi\ = \ 0$\,.\, Using a change of variables, we obtain
\begin{eqnarray*}
& \ &   n\,{\bar c}_{-1}    \int_{\R^n \,\setminus B_o (1)} H (y) \cdot   {{1}\over {\left(1 + r^2\right)^n}}\left[  {{r^2 - 1 }\over { r^2 + 1 }}  \right] \ \ \ \ \ \ \ \ \ \ \  \  \ \ \ \ \ \ \ \ \ \ \  \   \ \ \ \ \ \ \ \ \ \ \  \     (r\ = \ |\, y |\,)\\&\ = \ &   n\,{\bar c}_{-1} \,   \Vert \, S^{n - 1} \Vert\, \int_{S^{n - 1}} \int_1^\infty H (r, \  \theta) \cdot   {{1}\over {\left(1+ r^2\right)^n}}\left[  {{r^2 - 1 }\over { r^2 + 1}}  \right] r^{n - 1} \,dr\, dS_\theta \\
&\ = \ & -\, n\,{\bar c}_{-1}  \,\Vert \, S^{n - 1} \Vert \int_{S^{n - 1}} \int_0^1 H \left( \,{1\over {\bar r}}\ , \ \theta\right) \cdot \left[  {{   |\, \bar r |^{\,2} -1   }\over {  |\, \bar r |^{\,2}  + 1 }}  \right] \left( {{1}\over { 1 + |\,\bar r |^{\,2} }} \right)^n \,{\bar r}^{n - 1} \ d\bar r \,dS_\theta \\[0.15in]
& \  & \!\!\!\!\!\!\!\!\!\!\!\
\Longrightarrow  \ \ \ {{\partial G_{|_{\bf Z}} }\over {\partial \lambda }} \,(1, \ {\vec{\,0}}\,)\ \ \ \ \ \ \ \   \ \ \ \ \ \ \ \  \ \ \ \ \ \ \ \  \ \ \ \ \ \ \ \   \ \ \ \ \ \ \ \   \ \ \ \ \ \ \ \  \ \ \ \ \ \ \ \  \ \ \ \ \  (\,{\mbox{here}} \ \   \bar r\ = \ r^{-1}\,)\\ \\ &\ &\!\!\!\!\!\!\!\!\!\!\!\!\!\!\!\!\!\!\!\!=\ n\,{\bar c}_{-1}\,  \Vert \, S^{n - 1} \Vert \int_{S^{n - 1}} \int_0^1 \left[ H (\,s\,, \ \theta ) \!- \!H \left(\,{{1}\over s}\ , \ \theta \right) \right]\,\left[  {{ 1 - s^2   }\over {1 + s^2    }}  \right] \left( {{1}\over { 1 +  \,s^2 }} \right)^n \,s^{n - 1}ds\,dS_\theta\, = \ 0\\
 & \ & \ \ \ \ \ \ \ \ \ \  \ \ \ \ \ \ \ \ \ \  \ \ \ \ \ \ \ \ \ \  \ \ \ \ \ \ \ \ \ \  \ \  \ \ \ \ \ \ \ \ \ \  \ \ \ \ \ \ \ \ \  \ \ \  \ \ \  \ \ \ \ ({\mbox{via \ \ condition \ \ (7.17))}}\,.
\end{eqnarray*}
In view of (7.16), we also have
$$
{{\partial G_{|_{\bf Z}}}\over {\partial \,\xi_i}} \,(1, \ {\vec{\,0}}\,) \ = \ 2\, n\,{\bar c}_{-1}  \int_{\R^n}   \left[  {{1}\over {\left(1 + |\,y |^{\,2}\right)^{n+1}}}\, \right]\,[ H (y)\cdot  y_i]\,\,\ = \ 0\,.
$$
That is, $\,(1, \,\,{\vec{\,0}}\,)$ is a critical point of $G_{|_{\bf Z}}\,.\,$ Consider the Hessian matrix at $(1, \ {\vec{\,0}}\,)$\,.\, From (7.12), (7.13) and (7.16) we see that
\begin{eqnarray*}
{{\partial^2 G_{|_{\bf Z}}}\over {\partial\,\xi_j\, \partial \, \xi_\ell}} \,(1, \ {\vec{\,0}}\,) &\ = \ &  2^2 \,n\,(n + 1)\,{\bar c}_{-1} \int_{\R^n} {{   1  }\over {\left(1+ |\,y |^{\,2}\right)^{n+2}}}\, [H(y) \cdot y_\ell \, y_j]\,\,\ = \ 0\\
& \ & \\
{{\partial^2 G_{|_{\bf Z}}}\over {\partial \lambda\, \partial \, \xi_\ell}} \,(1, \ {\vec{\,0}}\,) &\ = \ &  \,2 n^2\,{\bar c}_{-1} \int_{\R^n}    {{  1  }\over {\left(1 + |\,y |^{\,2}\right)^{n+1}}}\,\, [ \,H (y)\cdot  y_\ell\,]\,\,\\
& \  & \ \ \ \ \ \ \ \   -\ 4\, n\,(n + 1)\,{\bar c}_{-1} \int_{\R^n}   {{ 1  }\over {\left(1 + |\,y|^{\,2}\right)^{n+2}}}\,\, [\, H (y)\cdot y_\ell\,]\,\ = \ 0\,.
\end{eqnarray*}
Here we apply symmetry argument for $\, 1 \,\le \,j\,\not= \,\ell \,\le\, n\,.\,$ From (7.19) and (7.20), together with the first equalities in (7.14) and (7.15), we obtain
 $$
 {{\partial^2 G_{|_{\bf Z}} }\over {\partial \lambda^2}} \,(1, \ {\vec{\,0}}\,) \not= 0 \ \ \ \ {\mbox{and}} \ \ \ \   {{\partial^2 G_{|_{\bf Z}} }\over {\partial\,\xi_\ell^2}}\, (1, \ {\vec{\,0}}\,) \not= 0 \mfor \ell\ = \ 1, \ 2, \cdot \cdot \cdot, \ n\,.
 $$
That is, $\,(1, \ {\vec{\,0}}\,)\,$ is a stable critical point of $\,G_{|_{\bf Z}}\,.\,$ Using Theorem 2.17 in \cite{Progress-Book}, pp. 25, together with elliptic regularity theory, we obtain a positive $C^2$ solution of equation (1.5) with relation (1.6) for $|\,\varepsilon|$ small enough.  \qedwh

\newpage

%
%
%
%
%
%
%
%
%
{\large \bf e\,-Appendix}\, \ \ can be found in
$
{\mbox{www.math.nus.edu.sg/}}\, \tilde{\ }\, {\mbox{matlmc/e\,-Appendix.pdf}}
$
